\input amstex\documentstyle{amsppt}  
\pagewidth{12.5cm}\pageheight{19cm}\magnification\magstep1
\topmatter
\title Non-unipotent character sheaves as a categorical centre\endtitle
\author G. Lusztig\endauthor
\address{Department of Mathematics, M.I.T., Cambridge, MA 02139}\endaddress
\endtopmatter   
\document
\define\tYY{\ti{\YY}}
\define\tvp{\ti\vp}

\define\cdo{\cdot}

\define\Irr{\text{\rm Irr}}

\define\Bpq{\Bumpeq}

\define\tco{\ti{\co}}
\define\btco{\ov{\ti{\co}}}

\define\du{\dot u}
\define\dw{\dot w}
\define\dz{\dot z}

\define\ds{\dot s}
\define\dy{\dot y}
\define\dZ{\dot Z}

\define\bco{\bar{\co}}

\define\mpb{\medpagebreak}

\define\bV{\bar V}

\define\hc{\hat c}

\define\hT{\hat T}
\define\da{\dagger}

\define\frl{\forall}

\define\si{\sim}

\define\sqc{\sqcup}

\define\qua{\quad}

\define\hL{\hat L}

\define\dx{\dot x}
\define\dM{\dot M}
\define\dL{\dot L}

\define\bG{\bar G}

\define\bZ{\bar Z}

\define\lb{\linebreak}

\define\bin{\binom}
\define\op{\oplus}
   
\redefine\spa{\spadesuit}
\define\part{\partial}
\define\emp{\emptyset}

\define\ra{\rangle}
\define\n{\notin}
\define\iy{\infty}
\define\m{\mapsto}
\define\do{\dots}
\define\la{\langle}
\define\bsl{\backslash}

\define\lra{\leftrightarrow}

\define\sub{\subset}    
\define\bxt{\boxtimes}
\define\T{\times}
\define\ti{\tilde}
\define\nl{\newline}
\redefine\i{^{-1}}

\define\un{\underline}
\define\ov{\overline}
\define\ot{\otimes}
\define\bbq{\bar{\QQ}_l}

\define\Hom{\text{\rm Hom}}
\define\End{\text{\rm End}}

\define\sg{\text{\rm sgn}}
\define\tr{\text{\rm tr}}
\define\rk{\text{\rm rk}}

\define\supp{\text{\rm supp}}
\define\card{\text{\rm card}}

\define\di{\diamond}

\define\a{\alpha}
\redefine\b{\beta}
\redefine\c{\chi}
\define\g{\gamma}
\redefine\d{\delta}
\define\e{\epsilon}

\define\io{\iota}
\redefine\o{\omega}
\define\p{\pi}
\define\ph{\phi}

\define\r{\rho}
\define\s{\sigma}
\redefine\t{\tau}
\define\th{\theta}
\define\k{\kappa}
\redefine\l{\lambda}
\define\z{\zeta}
\define\x{\xi}

\define\vp{\varpi}
\define\vt{\vartheta}

\redefine\G{\Gamma}
\redefine\D{\Delta}
\define\Om{\Omega}
\define\Si{\Sigma}

\redefine\L{\Lambda}
\define\Ph{\Phi}

\define\boc{\bold c}

\define\kk{\bold k}

\define\pp{\bold p}

\define\rr{\bold r}
\redefine\ss{\bold s}

\define\ww{\bold w}

\define\yy{\bold y}

\redefine\AA{\bold A}
\define\BB{\bold B}

\define\DD{\bold D}
\define\EE{\bold E}
\define\FF{\bold F}

\define\HH{\bold H}

\define\JJ{\bold J}

\define\LL{\bold L}

\define\NN{\bold N}

\define\QQ{\bold Q}

\define\TT{\bold T}
\define\UU{\bold U}

\define\ZZ{\bold Z}
\define\XX{\bold X}
\define\YY{\bold Y}

\define\ca{\Cal A}
\define\cb{\Cal B}
\define\cc{\Cal C}
\define\cd{\Cal D}
\define\ce{\Cal E}
\define\cf{\Cal F}
\define\cg{\Cal G}
\define\ch{\Cal H}
\define\ci{\Cal I}

\define\cl{\Cal L}
\define\cm{\Cal M}

\define\co{\Cal O}
\define\cp{\Cal P}

\define\car{\Cal R}
\define\cs{\Cal S}

\define\cu{\Cal U}
\define\cv{\Cal V}
\define\cw{\Cal W}
\define\cz{\Cal Z}
\define\cx{\Cal X}
\define\cy{\Cal Y}

\define\fb{\frak b}
\define\fc{\frak c}

\define\fh{\frak h}

\define\fk{\frak k}

\define\fo{\frak o}

\define\fs{\frak s}

\define\fA{\frak A}

\define\fD{\frak D}

\define\fK{\frak K}
\define\fL{\frak L}

\define\fR{\frak R}
\define\fS{\frak S}

\define\fZ{\frak Z}

\define\tb{\ti b}
\define\tc{\ti c}

\define\ty{\ti y}

\define\tE{\ti E}

\define\tK{\ti K}
\define\tL{\ti L}
\define\tM{\ti M}

\define\tP{\ti P}

\define\tV{\ti V}

\define\sha{\sharp}
\define\Mod{\text{\rm Mod}}

\define\bc{\bar c}
\define\bp{\bar p}

\define\bS{\bar S}

\define\bul{\bullet}

\define\che{\check}

\define\cir{\circ}

\define\tcb{\ti{\cb}}

\define\BBD{BBD}
\define\CONV{L17}
\define\CATEG{L18}
\define\CSI{L7}
\define\CSII{L8}
\define\CSIII{L9}
\define\CSIV{L10}
\define\HEC{L11}
\define\CDGVI{L12}
\define\CDGVII{L13}
\define\CDGIX{L14}
\define\CDGX{L15}
\define\CLEAN{L16}
\define\YO{Yo}
\define\KL{KL}
\define\ORA{L6}
\define\CLASS{L3}
\define\CLASSII{L4}
\define\GRE{L5}
\define\DL{DL}

\define\BFO{BFO}
\define\JS{JS}
\define\MA{Ma}
\define\ENO{ENO}
\define\REFL{L2}
\define\GP{GP}
\define\MADIS{L1}

\define\MU{M1}
\define\MUG{M2}
\define\SPALT{Spa}

\head Introduction\endhead
\subhead 0.1\endsubhead
Let $\kk$ be an algebraically closed field of characteristic $p\ge0$
and let $G$ be a reductive connected group over $\kk$.
We fix a prime number $l$ different from $p$. The theory of character sheaves developed in 
\cite{\CSI} and its sequels associates to $G$ a collection of simple perverse $\bbq$-sheaves on $G$
which in many respects mimic the irreducible representations of the finite Chevalley groups of the
same type as $G$. The classification of character sheaves was given in \cite{\CSIV}. A few years ago,
Bezrukavnikov, Finkelberg and Ostrik \cite{\BFO} gave a less computational (and more categorical)
approach to the classification of character sheaves assuming that the centre of $G$ is connected and 
that $p=0$. For applications to the study of finite Chevalley groups it was desirable to include the
case when $p>0$, but it was not clear how to do that by the method of \cite{\BFO} which relied on 
certain results on Harish-Chandra modules that are not available when $p>0$. In \cite{\CONV},
I found a way to obtain the classification of {\it unipotent} character sheaves in categorical terms
assuming that $p>0$, using a functor (truncated restriction) whose definition was different from that 
in \cite{\BFO}; moreover, in \cite{\CATEG}, 
I extended this to a classification of unipotent representations of a finite Chevalley group in 
categorical terms. In this paper I will extend the method of \cite{\CONV} to obtain the 
classification of not necessarily unipotent character sheaves of $G$ in categorical terms assuming 
that $p>0$. 

\subhead 0.2. Notation\endsubhead
In the rest of this paper $\kk$ is an algebraic closure of the finite field $\FF_q$ with $q$ elements.
All algebraic varieties are over $\kk$. We denote by $\pp$ the algebraic variety consisting of a 
single point. 
For an algebraic variety $X$ we write $\cd(X)$ for the bounded derived category of constructible
$\bbq$-sheaves on $X$. Let $\cm(X)$ be the subcategory of 
$\cd(X)$ consisting of perverse sheaves on $X$. 
For $K\in\cd(X)$ and $i\in\ZZ$ let $\ch^iK$ be the $i$-th cohomology sheaf of $K$ and let $K^i$ be
the $i$-th perverse cohomology sheaf of $K$; if $x\in X$, let $\ch^i_xK$ be the stalk of $\ch^iK$ 
at $x$. Let $\fD(K)$ be the Verdier dual of $K$.
If $X$ has a fixed $\FF_q$-structure $X_0$, we denote by 
$\cd_m(X)$ what in \cite{\BBD, 5.1.5} is denoted by $\cd_m^b(X_0,\bbq)$. 

Note that if $K\in\cd_m(X)$ then $K$ can be viewed as an object of $\cd(X)$ denoted again by $K$.
If $K\in\cd_m(X)$ is a perverse sheaf and $h\in\ZZ$, we denote by $gr_h(K)$ the subquotient of
pure weight $h$ of the weight filtration of $K$.
If $K\in\cd_m(X)$ and $i\in\ZZ$ we write $K\la i\ra=K[i](i/2)$ where $[i]$ is a shift and $(i/2)$
is a Tate twist; we write $K^{\{i\}}=gr_i(K^i)(i/2)$. 
If $K$ is a perverse sheaf on $X$ and $A$ is a simple perverse sheaf on $X$ we write $(A:K)$ for
the multiplicity of $A$ in a Jordan-H\"older series of $K$.

Assume that $C\in\cd_m(X)$ and that $\{C_i;i\in I\}$ is a family of objects of $\cd_m(X)$. We shall 
write $C\Bpq\{C_i;i\in I\}$ if the following condition is satisfied: there exist distinct elements
$i_1,i_2,\do,i_s$ in $I$, objects $C'_j\in\cd_m(X)$ ($j=0,1,\do,s$) and distinguished triangles 
$(C'_{j-1},C'_j,C_{i_j})$ for $j=1,2,\do,s$ such that $C'_0=0$, $C'_s=C$; moreover, $C_i=0$ unless 
$i=i_j$ for some $j\in[1,s]$. (See \cite{\CDGVII, 32.15}.)

Let $\ca=\ZZ[v,v\i]$ where $v$ is an indeterminate. Let $\bar{}:\ca@>>>\ca$ be the ring
homomorphism such that $\ov{v^m}=v^{-m}$ for any $m\in\ZZ$. If $f\in\QQ[v,v\i]$ and $j\in\ZZ$ we write
$(j;f)$ for the coefficient of $v^j$ in $f$.

Let $\cb$ be the variety of Borel subgroups of $G$. For any $B\in\cb$ let 
$U_B$ be the unipotent radical of $B$. In this paper we fix $\BB\in\cb$ and a 
maximal torus $\TT$ of $\BB$; let $\UU=U_\BB$. Let $\nu=\dim\UU=\dim\cb$, 
$\r=\dim\TT$, $\D=\dim G=2\nu+\r$. 

For any algebraic variety $X$ let $\fL=\fL_X=\a_!\bbq\in\cd(X)$ where $\a:X\T\TT@>>>X$ is the obvious 
projection. When $X$ is defined over $\FF_q$, $\fL$ is naturally an object of $\cd_m(X)$.

Unless otherwise specified, all vector spaces are over $\bbq$; in particular all representations of
a finite group $\G$ are assumed to be in (finite dimensional) $\bbq$-vector spaces. Let $\Mod\G$ be the 
category of representations of $\G$.

\subhead 0.3\endsubhead
We now discuss the content of various sections in some detail.
The main difference between \cite{\CONV} and the present paper is that the study of $G$-equivariant
sheaves on $\cb\T\cb$ is replaced by that of monodromic sheaves that is, certain $G$-equivariant
sheaves on $\tcb^2=G/\UU\T G/\UU$. The role that the Hecke algebra played in \cite{\CONV} is now 
played by 
a monodromic analogue $\HH$ of the Hecke algebra which was introduced (as an endomorphism algebra of 
the representation of a Chevalley group over $\FF_q$ induced by the unit representation of a Sylow 
$p$-subgroup) by Yokonuma \cite{\YO} in 1967. In Section 1 we recall from \cite{\CDGVII} various 
notions for $\HH$ that were known earlier for ordinary Hecke algebras: the canonical basis, the left 
cells, the two-sided cells,
the $a$-function, the asymptotic version. (Something close to the canonical basis of $\HH$ and its 
connection to intersection cohomology was already discussed in \cite{\ORA, Ch.1}.) A key role in our
discussion is the fact (see \cite{\CDGVII}) that $\HH$ is a matrix ring over an ordinary
extended Hecke algebra. In Section 2 we study the $G$-equivariant sheaves on 
$\tcb^2$ with monodromy of finite order dividing a fixed number $n$; we define truncated
convolution of such sheaves, see 2.24. This differs from the non-monodromic case since it now
involves direct images with compact support of non-proper maps, which makes the analysis more
complicated. In this section and in the subsequent ones we refer several times to two technical
lemmas \cite{\CONV, 1.12} and \cite{\CONV, 8.2} but we apply them in various cases which, although
not explicitly contained in those references, are proved just in the references. In Section 3 we 
define truncated convolution of $G$-equivariant sheaves on $Z=\TT\bsl\tcb^2$ with monodromy 
of order dividing $n$. Most of this section is concerned with the study of a functor $\fb$ (see 3.13) 
from sheaves on $Z$ to sheaves on $\tcb^2$ and its truncated version.
In Section 4 we discuss the unit object and rigidity of the truncated monoidal category 
$\cc^\boc\tcb^2$ of sheaves on $\tcb^2$ associated to a two-sided cell $\boc$ in $\HH$. 
In Section 5 we define truncated induction from a certain category
of sheaves $\cc^\boc Z$ on $Z$ associated to a two-sided cell $\boc$ of $\HH$ to a certain
category of sheaves $\cc^\boc G$ on $G$ associated to $\boc$ and we define truncated restriction 
going in 
the opposite direction. We also define truncated convolution in $\cc^\boc G$. In Section 6 we show 
(Theorem 6.13) that truncated restriction 
provides an equivalence of monoidal categories between $\cc^\boc G$
and the categorical centre of $\cc^\boc\tcb^2$. To do this we first 
prove a weak form of the adjunction
between truncated induction and truncated restriction. The adjunction is proved in full only as
a consequence of Theorem 6.13. Another consequence of Theorem 6.13 is that the character sheaves of
$G$ associated to $\boc$ are in bijection with the simple objects of the 
categorical centre of $\cc^\boc\tcb^2$.

\head Contents\endhead
1. Study of the algebra $\HH$.

2. Truncated convolution of sheaves on $\tcb^2$.

3. Sheaves on the variety $Z$.

4. The monoidal category $\cc^\boc\tcb^2$ and its centre.

5. Truncated induction, truncated restriction, truncated convolution on $G$.

6. The main results.

\head 1. Study of the algebra $\HH$\endhead
\subhead 1.1\endsubhead
Let $N\TT$ be the normalizer of $\TT$ in 
$G$, let $W=N\TT/\TT$ be the Weyl group and let $\k:N\TT@>>>W$ be the obvious 
homomorphism. For $w\in W$ we set $G_w=\BB\k\i(w)\BB$ so that $G=\sqc_wG_w$; 
let $\co_w=\{(x\BB x\i,y\BB y\i);x\in G,y\in G,x\i y\in G_w\}$ so that
$\cb\T\cb=\sqc_w\co_w$. For $w\in W$ let $\bG_w$ be the closure of $G_w$ in 
$G$; we have $\bG_w=\cup_{y\le w}G_y$ for a well defined partial order $\le$
on $W$. Let $\bco_w$ be the closure of $\co_w$ in $\cb^2$. Now $W$ is a 
(finite) Coxeter group with length function $w\m|w|=\dim\co_w-\nu$ and with 
set of generators $S=\{s\in W;|s|=1\}$. It acts on $\TT$ by 
$w:t\m w(t)=\o t\o\i$ where $\o\in\k\i(w)$.

\subhead 1.2\endsubhead
Let $s\in S$. Let $\UU_s$ be the unique root subgroup of $\UU$ with respect to
$\TT$ such that $\UU^-_s:=\s\UU_s\s\i\not\sub\UU$ for some/any $\s\in\k\i(s)$.
For any $\x\in\UU_s-\{1\}$ there is a unique $\eta\in\UU^-_s-\{1\}$ such that 
$\x\eta\x=\eta\x\eta\in\k\i(s)\sub N\TT$; we set $\s_\x=\x\eta\x=\eta\x\eta$. 
We have $\s_\x^4=1$. Note that $\x\m\eta$ is an isomorphism of algebraic 
varieties $\UU_s-\{1\}@>\si>>\UU^-_s-\{1\}$.

\subhead 1.3\endsubhead
Following Tits we define a cross-section $W@>>>N\TT$, $w\m\dw$ of
$\k:N\TT@>>>W$ as follows. For each $s\in S$ we choose $\x_s\in\UU_s-\{1\}$.
Let $w\in W$. We write $w=s_1s_2\do s_r$ where $s_i\in S$, $r=|w|$ and we set 
$\dw=\s_{\x_{s_1}}\s_{\x_{s_1}}\s_{\x_{s_2}}\do\s_{\x_{s_r}}\in\k\i(w)$. It is known
 that $\dw$ is independent of the choice of $s_1,s_2,\do,s_r$. 
Clearly, if $w,w'\in W$ satisfy $|ww'|=|w|+|w'|$, then $(ww')\dot{}=\dw\dw'$.

\subhead 1.4\endsubhead
In this paper we fix an integer $n\ge1$ such that $n\ne0$ in $\kk$. Let 
$\TT_n=\{t\in\TT;t^n=1\}$, $\fs=\Hom(\TT_n,\bbq^*)$. We have 
$\sha(\TT_n)=\sha(\fs)=n^\r$. Define $\io:\TT@>>>\TT$ by $t\m t^n$; clearly, 
$\io_!\bbq$ is a local system on $\TT$, equivariant for the $\TT$-action 
$t_1:t\m t_1^nt$ on $\TT$, hence $\TT_n$ acts naturally on each stalk of 
$\io_!\bbq$. We have $\io_!\bbq=\op_{\l\in\fs}L_\l$, where for any $\l\in\fs$,
$L_\l$ (a local system of rank $1$ on $\TT$) is such that $\TT_n$ acts on each
stalk of $L_\l$ through the character $\l$. 

The $W$-action on $\TT$ restricts to a $W$-action on $\TT_n$ hence induces a 
$W$-action on $\fs$. We shall write $W\fs$ instead of $W\T\fs$ (without group
structure); for $w\in W,\l\in\fs$ we shall write $w\cdo\l$ instead of 
$(w,\l)$. The following result can be deduced from \cite{\CDGVI, 28.2(a)}.

(a) {\it If $w\cdo\l\in W\fs$ and $w(\l)=\l$ then $L_\l$ is $T$-equivariant 
for the $T$-action $t:t'\m w(t)\i t't$ on $T$.}

\subhead 1.5 \endsubhead
Let $\t\in\TT$. We define $g_\t:\TT@>>>\TT$ by $t\m\t t$. We show that
for $\l\in\fs$,
the local systems $g_\t^*L_\l,L_\l$ are isomorphic. More precisely, we show
that any $\t'\in\TT$ such that $\t'{}^n=\t$ defines an isomorphism of local 
systems $g_\t^*L_\l@>\si>>L_\l$. The induced map 
$(g_\t^*L_\l)_t=(L_\l)_{\t t}@>c_{\t',t}>>(L_\l)_t$ on stalks at any 
$t\in\TT$ can be described as follows. We have 
$$(L_\l)_t=\{f:\io\i(t)@>>>\bbq;f(t_1t')=\l(t')f(t_1)\frl t_1\in\io\i(t),
t'\in\TT_n\},$$
$$(L_\l)_{\t t}=\{f':\io\i(\t t)@>>>\bbq;
f'(t_2t')=\l(t')f'(t_2)\frl t_2\in\io\i(\t t),t'\in\TT_n\}.$$
We have $c_{\t',t}(f)=f'$ where for any $t_2\in\io\i(\t t)$ we have 
$f'(t_2)=f(\t'{}\i t_2)$.

\subhead 1.6\endsubhead
For any root $\a:\TT@>>>\kk^*$ we denote by $\che\a:\kk^*@>>>\TT$ the 
corresponding coroot and by $s_\a$ the corresponding reflection in $W$.

Let $\l\in\fs$. Let $R_\l$ be the set of roots $\a:T@>>>\kk^*$ such that 
$\l(\che\a(z))=1$ for all $z\in\kk^*,z^n=1$. Let $W_\l$ be the subgroup of $W$
generated by $\{s_\a;\a\in R_\l\}$. We have $W_\l=W_{\l\i}$. Let 
$W'_\l=\{w\in W;w(\l)=\l\}$. Note that $W_\l\sub W'_\l$. There is a unique 
Coxeter group structure on $W_\l$ with length function $W_\l@>>>\NN$, 
$w\m|w|_\l$ such that, if $w\in W_\l$ and $w=s_1s_2\do s_r$ is any reduced 
expression of $w$ in $W$, then 
$$|w|_\l=\card\{i\in[1,r];s_r\do s_{i+1}s_is_{i+1}\do s_r\in W_\l\}.\tag a$$
See \cite{\CSI, 5.3}. 

\subhead 1.7\endsubhead
As in \cite{\CDGVI, 31.2}, let $\HH_n$ be the associative $\ca$-algebra with 
with generators $T_w (w\in W)$, $1_\l (\l\in\fs)$ and relations:

$1_\l1_{\l'}=\d_{\l,\l'}1_\l$ for $\l,\l'\in\fs$;

$T_wT_{w'}=T_{ww'}$ if $w,w'\in W$ and $|ww'|=|w|+|w'|$;

$T_w1_\l=1_{w(\l)}T_w$ for $w\in W,\l\in\fs$;

$T_s^2=v^2T_1+(v^2-1)\sum_{\l;s\in W_\l}T_s1_\l$ for $s\in W,|s|=1$;

$T_1=\sum_{\l\in\fs}1_\l$.
\nl
The algebra $\HH_n$ is closely related to the algebra introduced by Yokonuma \cite{\YO}.
(It specializes to it under $v=\sqrt{q},n=q-1$.) Since $n$ is fixed, we shall often write $\HH$ 
instead of $\HH_n$. Note that $T_1$ is the unit element of $\HH$ and that 
$\{T_w1_\l;w\cdo\l\in W\fs\}$ is an $\ca$-basis of $\HH$. The $\ca$-linear map
$\ti{}:\HH@>>>\HH$, $T_w1_\l\m T_w1_{\l\i}$ is an algebra automorphism. 
The 
$\ca$-linear map $\HH@>>>\HH$, $h\m h^\flat$, given by 
$T_w1_\l\m 1_\l T_{w\i}$ is an algebra antiautomorphism. (See 
\cite{\CDGVII, 32.19}.)

\subhead 1.8\endsubhead
For $w\in W$ we set $\hT_w=v^{-|w|}T_w\in\HH$. There is a unique ring 
homomorphism $\bar{}:\HH@>>>\HH$ such that 
$\ov{f\hT_w1_\l}=\bar f\hT_{w\i}\i1_\l$ for any $w\cdo\l\in W\fs,f\in\ca$; it 
has square $1$. As in \cite{\CDGVII, 34.4}, for any $w\cdo\l\in W\fs$ there is
a unique element $c_{w\cdo\l}\in\HH$ such that
$$c_{w\cdo\l}=\sum_{y\in W}p_{y\cdo\l,w\cdo\l}\hT_y1_\l$$
where $p_{y\cdo\l,w\cdo\l}\in v\i\ZZ[v\i]$ if $y\ne w$, 
$p_{w\cdo\l,w\cdo\l}=1$ and $\ov{c_{w\cdo\l}}=c_{w\cdo\l}$. Since 
$\bar{}:\HH@>>>\HH$, $\ti{}:\HH@>>>\HH$ commute, for any $w\cdo\l\in W\fs$, 
the element
$$\widetilde{c_{w\cdo\l\i}}=\sum_{y\in W}p_{y\cdo\l\i,w\cdo\l\i}\hT_y1_\l$$
satisfies the definition of $c_{w\cdo\l}$ hence 
$$\widetilde{c_{w\cdo\l\i}}=c_{w\cdo\l}.$$
In particular we have $p_{y\cdo\l\i,w\cdo\l\i}=p_{y\cdo\l\i,w\cdo\l\i}$ for 
any $y\cdo\l\in W\fs$.

For $y',w'$ in $W_\l$ let $P^\l_{y',w'}$ be the polynomial defined in 
\cite{\KL} in terms of the Coxeter group $W_\l$; let 
$$p^\l_{y',w'}=v^{-|w'|_\l+|y'|_\l}P^\l_{y',w'}(v^2)\in\ZZ[v\i].$$
Let $w\cdo\l\in W\fs$. From \cite{\ORA, 1.9(i)} we see that $wW_\l$ contains a
unique element $z$ such that $|z|$ is minimum; we write $z=\min(wW_\l)$; we 
have $w=zw'$ with $w'\in W_\l$. We show:
$$c_{zw'\cdo\l}=\sum_{y'\in W_\l}p^\l_{y',w'}\hT_{zy'}1_\l.\tag a$$
Since $p^\l_{y',w'}$ is $1$ if $y'=w'$ and is in $v\i\ZZ[v\i]$ if $y'\ne w'$, 
it is enough to show that 
$$\sum_{y'\in W_\l}p^\l_{y',w'}\hT_{zy'}1_\l\text{ is fixed by }
\bar{}:\HH@>>>\HH.\tag b$$
We can find a sequence $s_1,s_2,\do,s_k$ in $S$ such that
$$z(\l)=s_1s_2\do s_k\l\ne s_2\do s_k\l\ne\do\ne s_k\l\ne\l.$$ We argue by 
induction on $k$. If $k=0$ we have $z(\l)=\l$ and (b) follows from the proof 
of \cite{\CDGVII, 34.7}. Assume now that $k\ge1$. We have $z(\l)\ne(s_1z)(\l)$
hence $z\i s_1z\l\ne\l$. This implies that $s_1z=\min(s_1zW_\l)$. We have 
$s_1z\l=s_2\do s_k\l\ne\do\ne s_k\l\ne\l$; hence by the induction hypothesis 
applied to $s_1z$ instead of $z$ we see that
$$\sum_{y'\in W_\l}p^\l_{y',w'}(v^2)\hT_{s_1zy'}1_\l\text{ is fixed by }
\bar{}:\HH@>>>\HH.\tag c$$
For $y'\in W_\l$ we have $\hT_{zy'}1_\l=\hT_{s_1}1_{(s_1z)(\l)}\hT{s_1zy'}1_\l$ 
(we use again that $z(\l)\ne(s_1z)(\l)$) and $\hT_{s_1}1_{(s_1z)(\l)}$ is 
fixed by $\bar{}:\HH@>>>\HH$ (using that $z(\l)\ne(s_1z)(\l)$); we see that 
(b) follows from (c) by left multiplication with $\hT_{s_1}1_{(s_1z)(\l)}$.
This completes the proof of (a).

From (a) we see that
$$p_{y\cdo\l,zw'\cdo\l}=p^\l_{y',w'}(v^2)\text{ if }y=zy',y'\in W_\l,$$
$$p_{y\cdo\l,zw'\cdo\l}=0\text{ if }y\n zW_\l.$$
In particular we have $p_{y\cdo\l,w\cdo\l}\in\NN[v\i]$.
We show:
$$1_{w(\l)}c_{w\cdo\l}=c_{w\cdo\l}.\tag d$$
Using (a) it is enough to show that $1_{w(\l)}\hT_{zy'}1_\l=\hT_{zy'}1_\l$ for
any $y'\in W_\l$. We have $\hT_{zy'}1_\l=1_{(zy')(\l)}\hT_{zy'}$ and it is 
enough to show that $(zy')(\l)=(zw')(\l)$; this follows from $y'(\l)=\l$, 
$w'(\l)=\l$.

Let $w\cdo\l\in W\fs$. By (a) we have
$$c_{w\cdo\l}=\sum_{y\in wW_\l}p_{y\cdo\l,w\cdo\l}\hT_y1_\l.$$
Similarly, we have
$$c_{w\i\cdo w(\l)}=\sum_{y\in w\i W_{w(\l)}}p_{y\cdo w(\l),w\i\cdo w(\l)}
\hT_y1_{w(\l)}.$$
It follows that
$$\align&(c_{w\i\cdo w(\l)})^\flat=
\sum_{y\in w\i W_{w(\l)}}p_{y\cdo w(\l),w\i\cdo w(\l)}1_{w(\l)}\hT_{y\i}\\&
=\sum_{y\in w\i W_{w(\l)}}p_{y\cdo w(\l),w\i\cdo w(\l)}\hT_{y\i}1_{yw(\l)}.\endalign$$
For each $y$ in the last sum we have $y=w\i u$ with $uw(\l)=w(\l)$ hence 
$yw(\l)=w\i uw(\l)=w\i w(\l)=\l$. Thus, we have
$$\align&(c_{w\i\cdo w(\l)})^\flat
=\sum_{y\in w\i W_{w(\l)}}p_{y\cdo w(\l),w\i\cdo w(\l)}\hT_{y\i}1_\l\\&
=\sum_{y; y\i\in w\i W_{w(\l)}}p_{y\i\cdo w(\l),w\i\cdo w(\l)}\hT_y1_\l.\endalign$$
The condition that $y\i\in w\i W_{w(\l)}$ is equivalent to $y\in W_{w(\l)}w$ 
and also to $y\in wW_\l$. Hence
$$(c_{w\i\cdo w(\l)})^\flat=\sum_{y;y\in wW_\l}
p_{y\i\cdo w(\l),w\i\cdo w(\l)}\hT_y1_\l.$$
Note that $p_{y\i\cdo w(\l),w\i\cdo w(\l)}$ is $1$ if $y=w$ and is in 
$v\i\ZZ[v\i]$ if $y\ne w$. Also, since $\bar{}:\HH@>>>\HH$, 
${}^\flat:\HH@>>>\HH$ commute, $(c_{w\i\cdo w\i(\l)})^\flat$ is fixed by 
$\bar{}:\HH@>>>\HH$. It follows that $(c_{w\i\cdo w(\l)})^\flat$ satisfies the
defining property of $c_{w\cdo\l}$, hence
$$(c_{w\i\cdo w(\l)})^\flat=c_{w\cdo\l}.\tag e$$
We see also that 
$$p_{y\i\cdo w(\l),w\i\cdo w(\l)}=p_{y\cdo\l,w\cdo\l}$$
for any $y\in W$.

\subhead 1.9\endsubhead
Let $\fA$ be a {\it based $\ca$-algebra} that is, an associative $\ca$-algebra
with $1$ with a given finite basis $\{b_i;i\in I\}$ as an $\ca$-module, a 
given involution $i\m i^!$ of $I$ such that the $\ca$-linear map $x\m x^\flat$
defined by $b_i^\flat=b_{i^!}$ for all $i\in I$ is an algebra antiautomorphism
(necessarily preserving $1$) and a given subset $I_0$ of $\{i\in I;i^!=i\}$.
For $i,i'$ in $I$ we write $b_ib_{i'}=\sum_{j\in I}h_{i,i',j}b_j$ where 
$h_{i,i',j}\in\ca$. Let $j\underset\text{left}\to\preceq i$ (resp. $j\preceq i$) be the 
preorder on $I$ generated by the relation 
$h_{i',i,j}\ne0\text{ for some }i'\in I$, resp. by the relation 
$$h_{i,i',j}\ne0\text{ or }h_{i',i,j}\ne0\text{ for some }i'\in I.$$
We say that $i\underset\text{left}\to\si j$ (resp. $i\si j$) if $i\underset\text{left}\to\preceq j$ and
$j\underset\text{left}\to\preceq i$ (resp. $i\preceq j$ and $j\preceq i$). This is an 
equivalence relation on $I$; the equivalence classes are called left cells 
(resp. two-sided cells). Note that any two-sided cell is a union of left cells.

If $\boc$ is a two-sided cell and $i\in I$ we write $i\preceq\boc$ (resp. 
$\boc\preceq i$) if $i\preceq i'$ (resp. $i'\preceq i$) for some $i'\in\boc$; 
we write $i\prec\boc$ (resp. $\boc\prec i$) if $i\preceq\boc$ (resp. 
$\boc\preceq i$) and $i\n\boc$. If $\boc,\boc'$ are two-sided cells, we write 
$\boc\preceq\boc'$ (resp. $\boc\prec\boc'$) if $i\preceq i'$ (resp. 
$i\preceq i'$ and $i\not\si i'$) for some $i\in\boc,i'\in\boc'$.

Let $j\in I$. We can find an integer $m\ge0$ such that 
$h_{i,i',j}\in v^{-m}\ZZ[v]$ for all $i,i'$; let $a(j)$ be the smallest such 
$m$. For $i,i',j$ in $I$ there is a well defined integer $h^*_{i,i',j}$ such 
that
$$h_{i,i',j^!}=h^*_{i,i',j}v^{-a(j^!)}+\text{ higher powers of }v.$$ 
We say that the based algebra $\fA$ is {\it excellent} if properties Q1-Q11
below hold.

Q1. If $j\in I_0$ and $i,i'\in I$ satisfy $h^*_{i,i',j}\ne0$ then $i'=i^*$.

Q2. If $i\in I$, there exists a unique $j\in I_0$ such that 
$h^*_{i^!,i,j}\ne0$.

Q3. If $i'\preceq i$ then $a(i')\ge a(i)$. Hence if $i'\si i$ then 
$a(i')=a(i)$.

Q4. If $j\in I_0$, $i\in I$ and $h^*_{i^!,i,j}\ne0$ then $h^*_{i^!,i,j}=1$.

Q5. For any $i,j,k$ in $I$ we have $h^*_{i,j,k}=h^*_{j,k,i}$.

Q6. Let $i,j,k$ in $I$ be such that $h^*_{i,j,k}\ne0$. Then $i\underset\text{left}\to\si j^!$, 
$j\underset\text{left}\to\si k^!$, $k\underset\text{left}\to\si i^!$.

Q7. If $i'\underset\text{left}\to\preceq i$ and $a(i')=a(i)$ then $i'\underset\text{left}\to\si i$.

Q8. If $i'\preceq i$ and $a(i')=a(i)$ then $i'\si i$.

Q9. Any left cell $\G$ of $I$ contains a unique element of $j\in I_0$. We 
have $h^*_{i^!,i,j}=1$ for all $i\in\G$.

Q10. For any $i\in I$ we have $i\si i^!$.

Q11. Let $v'$ be a second indeterminate and let $h'_{i,j,k}\in\ZZ[v',v'{}\i]$ 
be obtained from $h_{i,j,k}$ by the substitution $v\m v'$, If $i,i',j,k\in I$ 
satisfy $a(j)=a(k)$ then
$$\sum_{j'\in I}h'_{k,i',j'}h_{i,j',j}=\sum_{j'\in I}h_{i,k,j'}h'_{j',i',j}.$$
In the remainder of this subsection we assume that $\fA$ is excellent. 
Consider the free abelian group $\fA^\iy$ with basis $\{t_i;i\in I\}$. We 
define a $\ZZ$-bilinear multiplication $\fA^\iy\T\fA^\iy@>>>\fA^\iy$ by 
$t_it_{i'}=\sum_{j\in I}h^*_{i,i',j^!}t_j$. As in \cite{\HEC, 18.3}, we see
using Q3,Q6 that this defines an associative ring structure on $\fA^\iy$ and 
we see using Q1,Q2,Q4,Q5 that $\sum_{i\in I_0}t_i$ is a unit element for this 
ring structure. Also from Q1,Q5 we see that $t_it_{i'}=\d_{i,i'}t_i$ for 
$i,i'\in I_0$.

From the definitions we have $h_{i^!,i'{}^!,j^!}=h_{i',i,j}$ for any $i,i',j$ 
in $I$. It follows that $a(j)=a(j^!)$ for any $j$ (this also follows from Q3, 
xQ10) and that $h^*_{i^!,i'{}^!,j}=h^*_{i',i,j^!}$ for any $i,i',j$ in $I$. 
Hence the $\ZZ$-linear map ${}^\flat:\fA^\iy@>>>\fA^\iy$ defined by 
$t_i^\flat=t_{i^!}$ for all $i\in I$ is a ring antiautomorphism.

We define an $\ca$-linear map $\psi:\fA@>>>\ca\ot\fA^\iy$ by 
$$\psi(b_i)=\sum_{i'\in I,j\in I_0;a(i')=a(j)}h_{i,j,i'}t_{i'}.$$
Using Q1,Q2,Q3,Q4,Q6,Q11 we see as in \cite{\HEC, 18.9} that $\psi$ is an 
$\ca$-algebra homomorphism preserving $1$.

We define a group homomorphism $\t:\fA^\iy@>>>\ZZ$ by $\t(t_i)=1$ if 
$i\in I_0$, $\t(t_i)=0$ if $i\in I-I_0$.  We show:

(a) {\it For $i,j\in I$ we have $\t(t_it_j)=1$ if $j=i^!$ and $\t(t_it_j)=0$ 
if $j\ne i^!$.}
\nl
An equivalent statement is that $\sum_{k\in I_0}h^*_{i,j,k^!}$ is $1$ if 
$j=i^!$ and is $0$ if $j\ne i^!$. This follows immediately from Q1,Q2,Q4.

\mpb

For any two-sided cell $\boc$ in $I$ let 
$\fA^\iy_\boc=\sum_{i\in\boc}\ZZ t_i\sub\fA^\iy$. From Q6 we see 
that if $\boc,\boc'$ are two-sided cells then $\fA^\iy_\boc\fA^\iy_{\boc'}$ is
zero if $\boc\ne\boc'$ and is contained in $\fA^\iy_\boc$ if $\boc=\boc'$. 
Hence $\fA^\iy_\boc$ is a ring with unit $\sum_{i\in I_0\cap\boc}t_i$ and 
$\fA^\iy=\op_\boc\fA^\iy_\boc$ as rings.

\subhead 1.10\endsubhead
Let $\fA$ be a based $\ca$-algebra with basis $\{b_i;i\in I\}$ and with 
$i\m i^!$, $I_0$ as in 1.9. We assume that $\fA$ is excellent. We use the 
notation in 1.9. We fix a two-sided cell $\boc$ of $I$ and we set $a=a(i)$ for
any $i\in\boc$. Let $r\ge1$, let $(i_1,i_2,\do,i_r)\in I^r$. We write
$$b_{i_1}b_{i_2}\do b_{i_r}
=\sum_{i\in I,k\in\ZZ}N(i,k)v^kb_i\text{ where }N(i,k)\in\ZZ.$$
We show:

(a) {\it Assume that $i\in\boc$. If $N(i,k)\ne0$ then $k\ge-(r-1)a$. If 
$N(i,-(r-1)a)\ne0$ then $i_u\in\boc$ for all $u\in[1,r]$.}
\nl
If $r=1$ the result is obvious. Now assume that $r\ge2$. We have 
$$\sum_{k\in\ZZ}N(i,k)v^k
=\sum h_{j_1,i_2,j_2}h_{j_2,i_3,j_3}\do h_{j_{r-1},i_r,j_r}\tag b$$
where the last sum is taken over all $j_1,j_2,\do,j_{r-1},j_r$ in $I$ such that
$$i=j_r\preceq j_{r-1}\preceq\do\preceq j_3\preceq j_2\preceq j_1=i_1.$$
Assume that $N(i,k)\ne0$. From (b) we see that
$$k=k_2+k_3+\do+k_r\text{ where }k_2\ge-a(j_2),\do,k_{r-1}\ge-a(j_{r-1}),
k_r\ge-a(j_r)=-a$$ 
for some $j_1,j_2,\do,j_{r-1},j_r$ as above. Using $Q3$ we see that
$$a=a(j_r)\ge a(j_{r-1})\ge\do\ge a(j_3)\ge a(j_2),$$
hence $k_2\ge-a$, $\do$, $k_r\ge-a$ and $k\ge-(r-1)a$, as required.

Assume now that $N(i,-(r-1)a)\ne0$. Then for some 
$j_1,j_2,j_3,\do,j_{r-1},j_r$ as above, the inequalities used above must be 
equalities
$$-k_2=\do=-k_{r-1}=a=a(j_r)=\do=a(j_3)=a(j_2)$$
and 
$$h^*_{i_1,i_2,j_2^!}\ne0,h^*_{j_2,i_3,j_3^!}\ne0,\do,
h^*_{j_{r-1},i_r,j_r^!}\ne0,$$
so that, by $Q6,Q8$ we have
$$i_1\si i_2\si j_2\si i_3\si j_3\si\do\si j_{r-1}\si i_r\si j_r.$$
Thus, $i_1\in\boc,\do,i_r\in\boc$. This proves (a).

We show:

(c) {\it Assume that $i\in\boc$ and $i_1\in\boc,\do,i_r\in\boc$. Then 
$$N(i,-(r-1)a)
=\sum h^*_{j_1,i_2,j_2^!}h^*_{j_2,i_3,j_3^!}\do h^*_{j_{r-1},i_r,j_r^!}$$
where the sum is taken over all $j_1,j_2,j_3,\do,j_{r-1},j_r$ in $\boc$ such 
that $j_1=i_1$, $j_r=i$.}
\nl
Note that $N(i,-(r-1)a)$ is the coefficient of $v^{-(r-1)a}$ in
$$\sum h_{j_1,i_2,j_2}h_{j_2,i_3,j_3}\do h_{j_{r-1},i_r,j_r}$$
where sum is taken over all $j_1,j_2,\do,j_{r-1},j_r$ in $I$ such that
$$i=j_r\preceq j_{r-1}\preceq\do\preceq j_3\preceq j_2\preceq j_1=i_1.$$
Such $j_u$ must satisfy $j_u\in\boc$ for all $u$ (since 
$j_1\in\boc,j_r\in\boc$). Hence the sum is equal to
$$\align&\sum(h^*_{j_1,i_2,j_2^!}v^{-a}+\text{ higher powers of $v$ })
(h^*_{j_2,i_3,j_3^!}v^{-a}+\text{ higher powers of $v$ })\do\\&
(h^*_{j_{r-1},i_r,j_r^!}v^{-a}+\text{ higher powers of $v$ })\\&=
\sum h^*_{j_1,i_2,j_2^!}h^*_{j_2,i_3,j_3^!}\do 
h^*_{j_{r-1},i_r,j_r^!}v^{-(r-1)a}+\text{ higher powers of $v$ }\endalign$$
where both sums are taken over all 
$j_1,j_2,\do,j_{r-1},j_r$ in $\boc$ such that $j_1=i_1$, $j_r=i$. Now (c) 
follows.

From (c) we deduce:

(d) {\it Assume that $i_1\in\boc,\do,i_r\in\boc$. Then 
$$t_{i_1}t_{i_2}\do t_{i_r}=\sum_{i\in\boc}N(i,-(r-1)a)t_i$$
(in $\fA^\iy$) where $N(i,-(r-1)a)$ is as in (c).}
\nl
We show:

(e) {\it Assume that $i_u\in\boc$ for some $u\in[1,r]$ and that $i\in I$, 
$k\in\ZZ$ are such that $N(i,k)\ne0$. Then either $i\in\boc$, $k\ge-(r-1)a$, 
or $i\prec\boc$.}
\nl
If $r=1$, the result is obvious. We now assume that $r\ge2$. We have
$$\sum_{k'\in\ZZ}N(i,k')v^{k'}
=\sum h_{j_1,i_2,j_2}h_{j_2,i_3,j_3}\do h_{j_{r-1},i_r,j_r}$$
where the last sum is taken over all $j_1,j_2,\do,j_r$ in $I$ such that 
$j_1=i_1,j_r=i$. Since the left hand side is $\ne0$, so is the right hand 
side. Thus there exist $j_1,j_2,\do,j_r$ as above such that 
$$h_{j_1,i_2,j_2}\ne0,h_{j_2,i_3,j_3}\ne0,\do,h_{j_{r-1},i_r,j_r}\ne0$$ 
hence $j_r\preceq j_{r-1}\preceq\do\preceq j_2\preceq j_1$ and 
$j_u\preceq i_u$. In particular we have $j_r\preceq i_u$ that is, 
$i\preceq\boc$. If $i\prec\boc$, there is nothing to prove. Thus we may assume 
that $i\in\boc$. In this case we have $k\ge-(r-1)a$ by (a). This proves (e).

\subhead 1.11\endsubhead
We now give some examples of excellent based $\ca$-algebras.

(i) Let $\l\in\fs$. Let $\HH_\l$ be the Hecke algebra of the Coxeter group
$W_\l$ with its basis $\{c_w; w\in I=W_\l\}$ defined as in 
\cite{\HEC, 5.3} with $i\m i^!$ given by $w\m w\i$ and $I_0$ being the set of 
distinguished involutions of $W_\l$ (defined as the set $\cd$ in 
\cite{\HEC, 14.1} with $W$ replaced by $W_\l$). Then $\HH_\l$ is excellent by 
results in \cite{\HEC, \S14,\S15}. 

\mpb

(ii) Let $\l\in\fs$. As in \cite{\CDGVII, 34.2}, $W'_\l$ is a semidirect 
product $W_\l\Om_\l$ where $\Om_\l$ is an abelian subgroup of $W'_\l$ such 
that any $x\in\Om_\l$ satisfies $xW_\l x\i=W_\l$ and $|xwx\i|_\l=|w|_\l$ for 
any $w\in W_\l$. Let $\HH'_\l$ be the $\ca$-module $\HH_\l\ot_\ca\ca[\Om_\l]$ 
with basis $\{c_w\ot x;w\in W_\l,x\in\Om_\l\}$. We regard 
$\HH'_\l$ as an $\ca$-algebra with multiplication
$(c_w\ot x)(c_{w'}\ot x')=(c_wc_{xw'x\i})\ot(xx')$
for $w,w'$ in $W_\l$ and $x,x'$ in $\Om_\l$.
We take $I=W_\l\T\Om_\l$, $i\m i^!$ given by $(w,x)\m(x\i(w\i),x\i)$ and
$I_0$ to be the set of all $(d,1)$ where $d$ is a distinguished involution of
$W_\l$. Then $\HH'_\l$ is excellent. (This follows easily from Case (i)).

\mpb

(iii) Let $\l\in\fs$. Let $\fA=1_\l\HH1_\l$ viewed as a subalgebra of $\HH$ 
with unit element $1_\l$ and with the basis $\{c_{w\cdo\l};w\in W'_\l\}$. In 
this case we take $I=W'_\l$. The involution $i\m i^!$ is given by $w\m w\i$ 
for $w\in W'_\l$. This is induced by the antiautomorphism of $\fA$ which is 
the restriction of the antiautomorphism $h\m h^\flat$ of $\HH$. We take $I_0$ 
to be the set of distinguished involutions of the Coxeter group $W_\l$.
In \cite{\HEC, 34.7} it is shown that $\fA$ is canonically isomorphic as a 
based $\ca$-algebra to $\HH'_\l$ in (ii). It follows that $\fA$ is excellent. 

\mpb

(iv) Let $\fo$ be a fixed $W$-orbit on $\fs$ and let $\l_0\in\fo$. Let 
$\EE$ be the set of all formal sums $x=\sum_{\l,\l'\in\fo}x_{\l,\l'}$ where
$x_{\l,\l'}\in1_{\l_0}\HH1_{\l_0}$ regarded naturally as an $\ca$-module and 
as an $\ca$-algebra where the product $xy$ of $x,y\in\EE$ is given by 
$$(xy)_{\l,\l'}=\sum_{\ti\l\in\fo}x_{\l,\ti\l}y_{\ti\l,\l'}$$ (we used the 
product in $1_{\l_0}\HH1_{\l_0}$). Let 
$I=\{(w,\l,\l')\in W'_{\l_0}\T\fo\T\fo\}$. We view $\EE$ as a based 
$\ca$-algebra with basis $\{b_{w,\l,\l'};(w,\l,\l')\in I\}$ where 
$b_{w,\l,\l'}$ has $(\l,\l')$-coordinate $c_{w\cdo\l_0}$ and all other
coordinates zero. The involution $i\m i^!$ of $I$ is given by 
$(w,\l,\l')\m(w\i,l',\l)$. This is induced by the antiautomorphism $h\m h^\di$
of the algebra $\EE$ such that $b_{w,\l,\l'}^\di=b_{w\i,\l',\l}$. The subset 
$I_0$ of $I$ is the set of all $(w,\l,\l')\in I$ such that $\l=\l'$ and 
$w\in\cd$, the set of distinguished involutions of $W_{\l_0}$. For $w,w'$ in 
$W'_{\l_0}$ we write 
$$c_{w\cdo\l_0}c_{w'\cdo\l_0}
=\sum_{w''\in W'_{\l_0}}h_{w,w',w''}c_{w''\cdo\l_0}$$ 
where $h_{w,w',w''}\in\ca$. Then the coefficients $h_{i_1,i_2,i_3}$ for 
$i_1,i_2,i_3$ in $I$ are given by
$$h_{w,\l_1,\l_2,w',\l'_1,\l'_2,w'',\l''_1,\l''_2}
=\d_{\l_2,\l'_1}\d_{\l_1,\l''_1}\d_{\l'_2,\l''_2}h_{w,w',w''}.$$
We see that the $a$-function on $I$ is given by $a(w,\l,\l')=a(w)$ where 
$a(w)$ is computed in $1_{\l_0}\HH1_{\l_0}$. Moreover,
$$h^*_{w,\l_1,\l_2,w',\l'_1,\l'_2,w'',\l''_1,\l''_2}
=\d_{\l_2,\l'_1}\d_{\l'_2,\l''_1}\d_{\l''_2,\l_1}h^*_{w,w',w''}.$$
We show that Q1-Q11 hold for $\fA=\EE$, using that we already know that they 
hold for $1_{\l_0}\HH1_{\l_0}$.

We prove Q1. Assume that 
$h^*_{w,\l_1,\l_2,w',\l'_1,\l'_2,w'',\l''_1,\l''_2}\ne0$ where $w''\in\cd$,
$\l''_1=\l''_2$. We must have $\l_2=\l'_1$, $\l_1=\l''_2=\l''_1=\l'_2$, 
$w'=w\i$ hence $(w,\l_1,\l_2)^!=(w',\l'_1,\l'_2)$. Thus Q1 holds.

We prove Q2. Assume that 
$h^*_{w,\l_1,\l_2,w\i,\l_2,\l_1,w'',\l''_1,\l''_1}\ne0$ where $w''\in\cd$. Then
$l_1=\l''_1$ and $w''$ is the unique element of $\cd$ such that 
$h^*_{w,w\i,w''}\ne0$; thus the uniqueness in Q2 holds. The same proof shows 
the existence in Q2.

We prove Q3. If $i=(w,\l_1,\l_2),i'=(w',\l'_1,\l'_2)$ then we have 
$i\preceq i'$ (resp. $i\si i'$) in $\EE$ if and only if $w\preceq w'$ (resp. $w\si w'$) in $1_{\l_0}\HH1_{\l_0}$. Hence Q3 for $\EE$ follows from Q3 for 
$1_{\l_0}\HH1_{\l_0}$. 

We prove Q4. Assume that 
$h^*_{w,\l_1,\l_2,w\i,\l_2,\l_1,w'',\l''_1,\l''_1}\ne0$ where $w''\in\cd$. Then
$h^*_{w,w\i,w''}\ne0$ hence by Q4 for $1_{\l_0}\HH1_{\l_0}$ we have 
$h^*_{w,w\i,w''}=1$. It follows that 
$h^*_{w,\l_1,\l_2,w\i,\l_2,\l_1,w'',\l''_1,\l''_1}=1$ as required.

We prove Q5. We must show that
$$\d_{\l_2,\l'_1}\d_{\l'_2,\l''_1}\d_{\l''_2,\l_1}h^*_{w,w',w''}=
\d_{\l'_2,\l''_1}\d_{\l''_2,\l_1}\d_{\l_2,\l'_1}h^*_{w',w'',w}.$$
This clearly follows from Q5 for $1_{\l_0}\HH1_{\l_0}$. 

We prove Q6. If $i=(w,\l_1,\l_2),i'=(w',\l'_1,\l'_2)$ then we have 
$i\underset\text{left}\to\preceq i'$ (resp. $i\underset\text{left}\to\si i'$) in $\EE$ if and only if 
$\l_2=\l'_2$ and $w\underset\text{left}\to\preceq w'$ (resp. $w\underset\text{left}\to\si w'$) in 
$1_{\l_0}\HH1_{\l_0}$. Assume that
$h^*_{w,\l_1,\l_2,w',\l'_1,\l'_2,w'',\l''_1,\l''_2}\ne0$. Then
$\d_{\l_2,\l'_1}\d_{\l'_2,\l''_1}\d_{\l''_2,\l_1}h^*_{w,w',w''}\ne0$ hence  
$\l_2=\l'_1,\l'_2=\l''_1,\l''_2=\l_1$ and (by Q6 for $1_{\l_0}\HH1_{\l_0}$) we
have $w\underset\text{left}\to\si w'{}\i$, $w'\underset\text{left}\to\si w''{}\i$, 
$w''\underset\text{left}\to\si w\i$. Thus Q6 
holds for $\EE$.

We prove Q8. Assume that $i=(w,\l_1,\l_2),i'=(w',\l'_1,\l'_2)$ and 
$i\preceq i'$, $a(i)=a(i')$. Then $w\preceq w'$ and $a(w)=a(w')$ in 
$1_{\l_0}\HH1_{\l_0}$ so that by Q8 for $1_{\l_0}\HH1_{\l_0}$ we have
$w\si w'$ and $i\si i'$. This proves Q8. The proof of Q7 is entirely similar 
to that of Q8.

We prove Q9. Now $\G$ is the set of all $(w,\l_1,\l_2)$ where $\l_2$ is 
fixed, $\l_1$ runs through $\fo$ and $w$ runs through a left cell $\G_0$ of
$W'_{\l_0}$. Let $w$ be the unique element in $\cd\cap\G_0$. Then 
$(w,\l_2,\l_2)$ is the unique element of $I_0\cap\G$. If 
$i=(w_1,\l_1,\l_2)\in\G$ then $w_1\in\G_0$ and 
$$h^*_{w_1\i,\l_2,\l_1,w,\l_1,\l_2,w,\l_2,\l_2}=h^*_{w_1\i,w_1,w}$$
and this is $1$ by Q9 for $1_{\l_0}\HH1_{\l_0}$.

We prove Q10. Let $(w,\l,\l')\in I$. It is enough to show that $w\si w\i$ in 
$1_{\l_0}\HH1_{\l_0}$; this follows from Q10 for $1_{\l_0}\HH1_{\l_0}$.

We prove Q11. We write $i=(w,\l_1,\l_2)$, $i'=(w',\l_3,\l_4)$, 
$j=(u,\l_5,\l_6)$, $k=(z,\l_7,\l_8)$. We have $a(z)=a(u)$. We must show 
$$\align&\d_{\l_8,\l_3}\d_{\l_7,\l_2}\d_{\l_1,\l_5}\d_{\l_4,\l_6}
\sum_{u'}h'_{z,w',u'}h_{w,u',u}\\&=
\d_{\l_8,\l_3}\d_{\l_7,\l_2}\d_{\l_1,\l_5}\d_{\l_4,\l_6}
\sum_{u'}h_{w,z,u'}h'_{u',w',u}.\endalign$$
This follows from Q11 for $1_{\l_0}\HH1_{\l_0}$.

We see that $\EE$ is excellent.

\mpb

(v) Let $\fo$ be a fixed $W$-orbit on $\fs$. Let $\HH_{\fo}$ be the 
$\ca$-subalgebra of $\HH$ with $\ca$-basis 
$\{c_{w\cdo\l};w\in W,\l\in\fo\}$. We view $\HH_{\fo}$ as a based 
$\ca$-algebra. We take $I=\{w\cdo\l\in W\fs;\l\in\fo\}$. The involution 
$i\m i^!$ of $I$ is $w\cdo\l\m w\i\cdo w(\l)$. This is induced by the 
antiautomorphism of $\HH_{\fo}$ which is the restriction of the 
antiautomorphism $h\m h^\flat$ of $\HH$. We take $I_0$ to be the set of all 
$w\cdo\l$ where $\l\in\fo$ and $w$ is a distinguished involution of 
$W_\l$. We show that $\HH_{\fo}$ is excellent. We fix $\l_0\in\fo$ and 
let $\EE$ be as in (iv) above. For any $\l\in\fo$ we choose a sequence 
$\ss_\l=(s_1,s_2,\do,s_r)$ in $S$ such that
$$\l_0\ne s_1\l_0\ne s_2s_1\l_0\ne\do\ne s_r\do s_2s_1(\l_0)=\l$$
and we set $\t_\l=\hT_{s_1}\hT_{s_2}\do\hT_{s_r}\in\HH$, 
$[\ss_\l]=s_1s_2\do s_r\in W$. Note that
$\t_\l^flat=\hT_{s_r}\do\hT_{s_2}\hT_{s_1}\in\HH$. We define an $\ca$-linear 
map $\Psi:\HH_{\fo}@>>>\EE$ by 
$$\Psi(h)_{\l_1,\l_2}
=\t_{\l_1}1_{\l_1}h1_{\l_2}\t_{\l_2}^\flat\in 1_{\l_0}\HH1_{\l_0}$$
for any $\l_1,\l_2$ in $\fo$.
In \cite{\HEC, 34.10} it is shown that $\Psi$ is an isomorphism of 
$\ca$-algebras and $\Psi\i$ carries the basis element $b_{w,\l_1,\l_2}$ of 
$\EE$ onto the basis element $c_{[\ss_{\l_1}]\i w[\ss_{\l_2}]\cdo\l_2}$ of 
$\HH_{\fo}$. We show that $\Psi(h^\flat)=(\Psi(h))^\di$ for all 
$h\in\HH_{\fo}$. Indeed for $\l_1,\l_2$ in $\fo$ we have
$$\align&((\Psi(h))^\di)_{\l_1,\l_2}=(\Psi(h)_{\l_2,\l_1})^\flat
=(\t_{\l_2}1_{\l_2}h1_{\l_1}\t_{\l_1}^\flat)^\flat\\&=
\t_{\l_1}1_{\l_1}h^\flat1_{\l_2}\t_{\l_2}^\flat=\Psi(h^\flat)_{\l_1,\l_2}.
\endalign$$
If $w$ is a distinguished involution of $W_{\l_0}$ and $\l\in\fo$, then 
$$\Psi\i(b_{w,\l,\l})=c_{[\ss_\l]\i w[\ss_\l]\cdo\l};$$
 note that conjugation by
$[\ss_\l]\i$ is a Coxeter group isomorphism $W_{\l_0}@>\si>>W_\l$ hence 
$[\ss_\l]\i w[\ss_\l]$ is a distinguished involution of $W_\l$. This argument 
shows that $\Psi\i$ induces a bijection from $I_0$ defined in terms of $\EE$ 
to $I_0$ defined in terms of $\HH_{\fo}$. Using the fact that $\EE$ is 
excellent we now deduce that $\HH_{\fo}$ is excellent.

\mpb

(vi) We consider the $\ca$-algebra $\HH$ with its $\ca$-basis 
$\{c_{w\cdo\l};w\cdo\l\in W\fs\}$. We view $\HH$ as a based $\ca$-algebra. We
take $I=W\fs$. The involution $i\m i^!$ of $I$ is $w\cdo\l\m w\i\cdo w(\l)$. 
This is induced by the antiautomorphism $h\m h^\flat$ of $\HH$. We take $I_0$ 
to be the set of all $w\cdo\l$ where $\l\in\fs$ and $w$ is a distinguished 
involution of $W_\l$. We have $\HH=\op_{\fo}\HH_{\fo}$ (as algebras) 
where $\fo$ runs over the set of $W$-orbits in $\fs$ and $\HH_{\fo}$ is 
as in (v). Using the fact that each $\HH_{\fo}$ is excellent, it follows 
immediately that $\HH$ is excellent.

In this case we shall write $\DD$ instead of $I_0$.

\mpb

In particular, in case (vi) the two-sided cells of $W\fs$ and the $a$-function
$a:W\fs@>>>\NN$ are well defined. Note that each two-sided cell of $W\fs$ (in 
case (vi)) is equal to a two-sided cell of $W\T\fo$ (in case (v)) for a 
unique $W$-orbit $\fo$ in $\fs$. Moreover for any two-sided cell $\fc$ of 
$W\T\fo$ (with $\l\in\fo$), the subset $\{w\in W'_\l;w\cdo\l\in\fc\}$ is
a two-sided cell of $W'_\l$ (in case (ii)) and this gives a bijection between 
the set of two-sided cells of $W\T\fo$ and the set of two-sided cells of 
$W'_\l$ in case (ii), which in turn is in bijection with the set of orbits of 
the conjugation action of $\Om_\l$ on the set of two-sided cells of $W_\l$ in 
case (i).

\mpb

The based algebras in (i)-(vi) have the additional properties that

(a) $h_{i,j,k}\in\NN[v,v\i]$, $\ov{h_{i,j,k}}=h_{i,j,k}$ for any $i,j,k\in I$; 

(b) $h^*_{i,j,k}\in\NN$ for any $i,j,k\in I$.
\nl
Indeed, (a) is well known in the cases (i),(ii); from this we deduce by the 
arguments in (iii)-(vi) that (a) holds in each case (iii)-(vi). Clearly, (b) 
follows from (a). From (a) we see that for the based algebras in (iii)-(vi) we
have
$$h_{i,i',j^!}=h^*_{i,i',j}v^{a(j^!)}+\text{ lower powers of }v.$$ 

\mpb

In (i) the ring $\HH_\l^\iy$ has $\ZZ$-basis $\{t_w;w\in W_\l\}$ in 
natural bijection with the $\ca$-basis $(c_w)$ of $\HH_\l$. 

In (ii) the ring $(\HH'_\l)^\iy$ has $\ZZ$-basis 
$\{t_w\ot x;w\in W_\l,x\in\Om_\l\}$ in natural bijection with the $\ca$-basis 
$(c_w\ot x)$ of $\HH'_\l$. The multiplication is given by
$$(t_w\ot x)(t_{w'}\ot x')=\sum_{z\in W_\l}c_zt_z\ot(xx')$$
where $t_wt_{xw'x\i}=\sum_{z\in W_\l}c_zt_z$, $c_z\in\ZZ$, is the product in 
$\HH_\l^\iy$.

In (iii) we have an identification $(1_\l\HH1_\l)^\iy=(\HH'_\l)^\iy$ (as rings)
for which the basis element $t_{(wx)\cdo\l}$ (with $w\in W_\l,x\in\Om_\l$) of
$(1_\l\HH1_\l)^\iy$ corresponds to the basis element $t_w\ot x$ of
$(\HH'_\l)^\iy$.

In (iv) the ring $\EE^\iy$ has $\ZZ$-basis $\{t_{w,\l,\l'};
w\in W'_{\l_0},\l\in\fo,\l'\in\fo\}$ in natural bijection with the
$\ca$-basis $(b_{w,\l,\l'})$ of $\EE$. This ring is canonically isomorphic to
a matrix ring with entries in the ring $(1_{\l_0}\HH1_{\l_0})^\iy$ with its
natural basis.

In (v) we have an identification $\HH_{\fo}^\iy=\EE^\iy$ (as rings) for 
which the basis element $t_{w,\l_1,\l_2}$ of $\EE^\iy$ corresponds to the 
basis element $t_{[\ss_{\l_1}]\i w[\ss_{\l_2}]\cdo\l_2}$ of $\HH_{\fo}^\iy$.

In (vi) we have an identification $\HH^\iy=\op_{\fo}\HH_{\fo}^\iy$ (as
rings) for which the basis elements $t_{w\cdo\l}$ in the two sides correspond 
to each other.

\subhead 1.12\endsubhead
For a based $\ca$-algebra $\fA$ as in 1.9 we set $\fA^v=\bbq(v)\ot_\ca\fA$;
we set $\fA^1=\bbq\ot_\ca\fA$ where $\bbq$ is viewed as an 
$\ca$-algebra via $v\m1$. 
 By definition, $\HH^1$ is the associative 
$\bbq$-algebra with generators $T_w (w\in W)$, $1_\l (\l\in\fs)$ and relations:

$1_\l1_{\l'}=\d_{\l,\l'}1_\l$ for $\l,\l'\in\fs$;

$T_wT_{w'}=T_{ww'}$ if $w,w'\in W$;

$T_w1_\l=1_{w(\l)}T_w$ for $w\in W,\l\in\fs$;

$T_1=\sum_{\l\in\fs}1_\l$.
\nl
The elements $\{T_w1_\l;w\cdo\l\in W\fs\}$ form an $\bbq$-basis of $\HH^1$ and 
$T_1$ is the unit element.
Consider the group algebra $\bbq[W\TT_n]$ where $W\TT_n$ is the semidirect 
product of $W$ and $\TT_n$ with $\TT_n$ normal and $W$ acting on $\TT_n$ by 
$w:t\m w(t)$. We define a $\bbq$-linear map $\bbq[W\TT_n]@>>>\HH^1$ by 
$wt\m\sum_{\l\in\fs}\l(t)T_w1_\l$. From the definitions we see that this is an
isomorphism of $\bbq$-algebras; we shall use it to identify
$\bbq[W\TT_n]=\HH^1$. For $\l\in\fs$ we set
$(1_\l\HH1_\l)^1=\bbq\ot_\ca(1_\l\HH1_\l)$; under the identification above we
have $(1_\l\HH1_\l)^1=\bbq[W'_\l]$.

\mpb

Recall that we have $\HH^\iy=\op_{\boc}\HH^\iy_\boc$ as rings. Here, for any 
two-sided cell $\boc$, $\HH^\iy_\boc$ has basis 
$\{t_{w\cdo\l};w\cdo\l\in\boc\}$; it is a 
ring with unit element $\sum_{w\cdo\l\in\DD_\boc}t_{w\cdo\l}$ where 
$\DD_\boc=\DD\cap\boc$. We set $\JJ=\bbq\ot\HH^\iy$. 
We have $\JJ=\op_\boc\JJ_\boc$ (as algebras) where for any two-sided cell 
$\boc$ of $W\fs$ we set $\JJ_\boc=\bbq\ot\HH^\iy_\boc$. 

Now $\psi:\HH@>>>\ca\ot\HH^\iy$ and 
$\psi:1_\l\HH1_\l@>>>\ca\ot(1_\l\HH1_\l)^\iy$ 
induce by extension of scalars isomorphisms of split semisimple 
$\bbq(v)$-algebras $\psi^v:\HH^v@>\si>>\bbq(v)\ot\HH^\iy$,
$\psi^v:(1_\l\HH1_\l)^v@>\si>>\bbq(v)\ot(1_\l\HH1_\l)^\iy$ and isomorphisms of 
semisimple $\bbq$-algebras $\psi^1:\HH^1@>\si>>\JJ$, 
$\psi^1:\bbq[W'_\l]@>\si>>\bbq\ot(1_\l\HH1_\l)^\iy$. (See
\cite{\CDGVII, 34.12(b),(c)}.)

Let $\Irr(W\TT_n)$ be a set of representatives for the isomorphism classes of 
simple $\bbq[W\TT_n]=\HH^1$-modules. For any $W\TT_n$-module $E$ let $E^\iy$ be 
the corresponding $\JJ$-module (via $\psi^1$) and let $E^v$ be the 
$\HH^v$-module corresponding to $\bbq(v)\ot E^\iy$ under $\psi^v$. 
For any $W$-orbit $\fo$ on $\fs$ let $\Irr_{\fo}(W\TT_n)$ be the set of 
all $E\in\Irr(W\TT_n)$ such that $1_{\l'}E=0$ for all $\l'\n\fo$. We have 
$\Irr(W\TT_n)=\sqc_{\fo}\Irr_{\fo}(W\TT_n)$. If 
$E\in\Irr_{\fo}(W\TT_n)$,
then for any $\l\in\fo$, $1_\l E$ is a simple $(1_\l\HH1_\l)^1$-module, 
that is, a simple $W'_\l$-module. Moreover $E\m1_\l E$ is a bijection between 
$\Irr_{\fo}(W\TT_n)$ and a set of representatives $\Irr(W'_\l)$ for the 
isomorphism classes of simple $\bbq[W'_\l]$-modules.

For any $E\in\Irr_{\fo}(W\TT_n)$ and $\l\in\fo$ let $(1_\l E)^\iy$ be 
the $\bbq\ot(1_\l\HH1_\l)^\iy$-module corresponding to $1_\l E$ via $\psi^1$ 
and let $(1_\l E)^v$ be the $(1_\l\HH1_\l)^v$-module corresponding to 
$\bbq(v)\ot(1_\l E)^\iy$ under $\psi^v$. Note that $(1_\l E)^v=1_\l(E^v)$.

If $E\in\Irr_{\fo}(W\TT_n)$, $\l\in\fo$ and $w\in W'_\l$ then we have
$$\tr(c_{w\cdo\l},E^v)=\tr(c_{w\cdo\l},(1_\l E)^v).$$
For any left cell $\L$ of $W\fs$ contained in $\fo$ we denote by $[\L]$ the
$\bbq[W\TT_n]$-module such that $[\L]^\iy$ is the $\bbq$-subspace of $\JJ$ 
spanned by $\{t_{w\cdo\l};w\cdo\l\in\L\}$ (a left ideal of $\JJ$). We show:

(a) {\it Let $z\cdo\l$ be the unique element of $\L\cap\DD$. Then for any 
$E\in\Irr_{\fo}W\TT_n$, $\tr(t_{z\cdo\l},E^\iy)$ is equal to the 
multiplicity of $E^\iy$ in the $\JJ$-module $[\L]^\iy$.}
\nl
An equivalent statement is that 
$$\dim(t_{z\cdo\l}E^\iy)=\dim\Hom_\JJ([\L]^\iy,E^\iy).$$ 
It is enough to show that the $\bbq$-linear map 
$\Hom_{\JJ}([\L]^\iy,E^\iy)@>>>t_{z\cdo\l}E^\iy$, $\x\m\x(t_{z\cdo\l})$ is an
isomorphism. The proof is along the same lines as that of \cite{\HEC, 21.3}.

\subhead 1.13\endsubhead
Let $\fA$ be one of the based $\ca$-algebra $\HH_\l,\HH'_\l$ (with $\l\in\fs$)
or $\HH$. Note that in these cases $I$ is $W_\l,W'_\l,W\fs$ respectively and 
$\fA^1$ is the group algebra $\bbq[\cw]$ where $\cw$ is $W_\l,W'_\l,W\TT_n$ 
respectively. Note that $\psi:\fA@>>>\ca\ot\fA^\iy$ induces an isomorphism
$\bbq[\cw]=\fA^1@>\si>>\bbq\ot\fA^\iy$. Under this isomorphism an irreducible
$\cw$-module $E$ corresponds to a simple $\bbq\ot\fA^\iy$-module $E^\iy$.
We have $\bbq\ot\fA^\iy=\op_\boc(\bbq\ot\fA^\iy_\boc)$ where $\boc$ runs over
the two-sided cell of $I$. Hence if $E$ is an irreducible $\cw$-module then 
there is a unique two-sided cell $\boc_E$ of $I$ such that 
$\bbq\ot\fA^\iy_\boc$ acts as zero on $E^\iy$ for any $\boc\ne\boc_E$; thus 
$E^\iy$ can be viewed as a simple $\bbq\ot\fA^\iy_{\boc_E}$-module.
For an irreducible $\cw$-module $E$ let $a_E\in\NN$ be the constant value of 
the restriction of $a:I@>>>\NN$ to $\boc_E$.

\subhead 1.14 \endsubhead
Since $\ti{}:\HH@>>>\HH$ permutes the elements in the basis $\{c_{w\cdo\l}\}$ 
according to the involution $w\cdo\l\m\widetilde{w\cdo\l}=w\cdo\l\i$ of 
$W\fs$, we see that the image of a two-sided cell $\boc$ of $W\fs$ under this 
involution is again a two-sided cell $\ti\boc$ of $W\fs$ and the value of the 
$a$-function on $\ti\boc$ is equal to the value of the $a$-function on $\boc$.

\subhead 1.15\endsubhead
Applying ${}^\flat$ to the equation 
$$c_{x\cdo\l}c_{y\cdo\l'}=
\sum_{z\cdo\l''\in W\fs}h_{x\cdo\l,y\cdo\l',z\cdo\l''}c_{z\cdo\l''},$$
we get
$$\align&c_{y\i\cdo y(\l')}c_{x\i\cdo x(\l)}=
\sum_{z\cdo\l''\in W\fs}h_{x\cdo\l,y\cdo\l',z\cdo\l''}c_{z\i\cdo z(\l'')}\\&=
\sum_{z\cdo\l''\in W\fs}h_{x\cdo\l,y\cdo\l',z\i\cdo z(\l'')}c_{z\cdo\l''}\endalign$$
hence 
$$h_{x\cdo\l,y\cdo\l',z\i\cdo z(\l'')}
=h_{y\i\cdo y(\l'),x\i\cdo x(\l),z\cdo\l''}.$$
This shows that the involution $z\cdo\l\m z\i\cdo z(\l)$ of $W\fs$ preserves 
the preorder $\preceq$ hence it maps any two-sided cell onto a two-sided cell.
(In fact, it maps each two-sided cell $\boc$ onto itself. Indeed, it is enough
to show that some element $z\cdo\l$ of $\boc$ satisfies
$z\cdo\l=z\i\cdo z(\l)$; we can take $z\cdo\l$ to be any element of the 
nonempty subset $\DD_\boc$ of $\boc$.) We also see that the $a$-function on 
$W\fs$ is constant on the orbits of our involution and that the group 
isomorphism ${}^\flat:\HH^\iy@>>>\HH^\iy$ given by 
$t_{z\cdo\l}\m t_{z\i\cdo z(\l)}$ 
is a ring antiautomorphism. Note that our involution restricts to the identity
permutation of the subset $\DD$ of $W\fs$ and that the algebra homomorphism 
$\psi:\HH@>>>\ca\ot\HH^\iy$ is compatible with the antiautomorphisms 
${}^\flat$ of $\HH$ and of $\ca\ot\HH^\iy$.

\subhead 1.16\endsubhead
Let $()^\spa:\bbq@>>>\bbq$ be a field automorphism which maps any root of $1$ 
in $\bbq$ to its inverse. The field automorphism $\bbq(v)@>>>\bbq(v)$ which 
maps $v$ to $v$ and $x\in\bbq$ to $x^\spa$ is denoted again by ${}^\spa$.

\subhead 1.17\endsubhead
We can view naturally $W$ as a subgroup of $GL(V)$ where
$$V:=\bbq\ot\Hom(\TT,\kk^*)(-1).$$
For any $i\ge0$ let $S^iV$ be the $i$-th symmetric power of the vector space 
$V$. Then $SV=\op_{i\ge0}S^iV$ is naturally a commutative algebra. Now $W$ acts
naturally on $S^iV$ for any $i$. 

Let $\l\in\fs$. If $E,\tE\in\Mod W_\l$
(resp. $E',\tE'\in\Mod W'_\l$) with $E$ (resp. $E'$) irreducible, we set 
$\tE^E=\Hom_{W_\l}(E,\tE)$ (resp. $\tE'{}^{E'}=\Hom_{W'_\l}(E',\tE')$). For 
$E$ (resp. $E'$) as above there exists $i\ge0$ such that $(S^iV)^E\ne0$ (resp.
$(S^iV)^{E'}\ne0$); let $b_E$ (resp. $b_{E'}$) be the smallest such $i$. If 
$E$ (resp. $E'$) are as above we say that $E$ (resp. $E'$) is univalent if
$\dim((S^{b_E}V)^E)=1$ (resp. $\dim((S^{b_{E'}}V)^{E'})=1$). We show:

(a) {\it Let $E\in\Mod W_\l$ be irreducible and univalent. There exists
$E'\in\Mod W'_\l$ irreducible such that $E$ appears in $E'|_{W_\l}$ and 
$b_{E'}=b_E$; moreover, $E'$ is uniquely determined up to isomorphism by these
properties and is univalent.}
\nl
Let $\un E$ be the unique $W_\l$-submodule of $S^{b_E}V$ that is isomorphic to
$E$. Let $E'=\sum_{x\in\Om_\l}x\un E\sub S^{b_E}V$ (notation of 1.11(ii)) 
where we have used the $W$-action on 
$S^{b_E}V$. Then $E'$ is a $W'_\l$-submodule of $S^{b_E}V$; moreover, for each 
$x\in\Om_\l$, $x\un E$ is an irreducible $W_\l$-submodule of $S^{b_E}V$. If 
$\ce$ is an irreducible $W_\l$-submodule of $E'$ then $\ce$ is isomorphic to 
$x\un E$ as an irreducible $W_\l$-submodule (for some $x\in\Om_\l$). But
$x\un E$ is a univalent $W_\l$-submodule hence we have necessarily 
$\ce=x\un E$.
Since any irreducible $W_\l$-submodule of $E'$ is equal to $x\un E$ for some
$x\in\Om_\l$, we see that any nonzero $W'_\l$-submodule of $E'$ contains
$x\un E$ for some $x\in\Om_\l$; being stable under the action of $\Om_\l$, it
is equal to $E'$. Thus $E'$ is an irreducible $W'_\l$-submodule of $S^{b_E}V$. 
Clearly, $E'$ appears with multiplicity $1$ (resp. $0$) in the $W'_\l$-module 
$S^{b_E}V$ (resp. $S^iV$ with $0\le i<b_E$). Thus $b_{E'}=b_E$ and $E'$ is 
univalent. Thus the existence of $E'$ in (a) is proved. Now let $E'_1$ be an
irreducible $W'_\l$-module such that $E$ appears in $E'_1|_{W_\l}$ and 
$b{E'_1}=b_E$. We can find a $W'_\l$-submodule $\un E'_1$ of $S^{b_E}V$ that
is isomorphic to $E'_1$. By assumption we can find a $W_\l$-submodule of
$\un E'_1$ which is isomorphic to $E$; this is necessarily equal to $\un E$.
For any $x\in\Om_\l$ we must have $x\un E\sub\un E'_1$ so that $E'\sub E'_1$.
Since $E'_1$ is irreducible as a $W'_\l$-module we have $E'_1=E'$. This proves
(a).

\subhead 1.18\endsubhead
If $E,\tE\in\Mod W\TT_n$ with $E$ irreducible, we set $\tE^E=\Hom_{W\TT_n}(E,\tE)$.
For any $i\ge0$, $W\TT_n$ acts on $S^iV\ot\bbq[\TT_n]$ ($V$ as in 1.17) by
$wx_1:v\ot x\m w(v)\ot w(x_1x)$ (with $w\in W, v\in V,x_1\in\TT_n,x\in\TT_n$).
If $E\in\Mod W\TT_n$ is irreducible, there exists $i\ge0$ such that 
$(S^iV\ot\bbq[\TT_n])^E\ne0$. Indeed, there exists $\l\in\fs$ and 
$E_1\in\Mod W'_\l$ irreducible such that $E$ is induced by the representation 
$E_1\ot\l$ of the subgroup $W'_\l\TT_n$. Then $E_1$ appears in the 
$W'_\l$-module $S^{b_{E_1}}V$ hence $E$ appears in the $W\TT_n$-module 
$S^{b_{E_1}}V\ot\bbq[\TT_n]$. Thus we can take $i=b_{E_1}$. Let $b_E$ be the 
smallest $i\ge0$ such that $(S^iV\ot\bbq[\TT_n])^E\ne0$. Note that 
$b_E\le b_{E_1}$. Conversely, assume that $(S^iV\ot\bbq[\TT_n])^E\ne0$. Since 
the $W\TT_n$-module $S^iV\ot\bbq[\TT_n]$ is induced by the $W$-module $S^iV$ 
we see (using Frobenius reciprocity) that
$\dim(S^iV\ot\bbq[\TT_n])^E=\dim\Hom_W(E,S^iV)$. Since the $W$-module $E$ is
induced by the $W'_\l$-module $E_1$, the last dimension is equal to
$\dim(S^iV)^{E_1}$. Thus $\dim(S^iV)^{E_1}\ne0$ so that $b_{E_1}\le i$. We 
see that $b_E=b_{E_1}$. This argument shows also that
$(S^{b_E}V\ot\bbq[\TT_n])^E=(S^{b_{E_1}}V)^{E_1}$.
We say that $E$ is univalent if $\dim(S^{b_E}V\ot\bbq[\TT_n])^E=1$ or
equivalently if $E_1$ is univalent.

\subhead 1.19\endsubhead
Let $\fA$, $I$, $\cw$ be as in 1.13. Thus $\cw$ is $W_\l,W'_\l$ (with 
$\l\in\fs$) or $W\TT_n$. Let $\sg$ be the (one dimensional) sign representation
of $W$. The composition of $\sg$ with the obvious homomorphism $\cw@>>>W$ (the 
inclusion if $\cw$ is $W_\l$ or $W'_\l$, the projection if $\cw=W\TT_n$) is 
denoted again by $\sg$. If $E$ is an irreducible $\cw$-module then $E\ot\sg$ is
again an irreducible $\cw$-module. 
An irreducible representation $E$ of $W$ is said to be {\it special} if
$a_E=b_{E\ot\sg}$. We show:

(a) {\it If $E$ is an irreducible $\cw$-module then $a_E\le b_{E\ot\sg}$.}

(b) {\it For any two-sided cell $\boc$ of $I$ there exists a unique (up to 
isomorphism) irreducible special representation $E$ of $\cw$ such that 
$\boc_E=\boc$. Moreover, $E\ot\sg$ is univalent.}
\nl
In the case where $\fA=\HH_\l$, (a),(b) are known from \cite{\CLASS}. 

We prove (a) for $\fA=\HH'_\l$. Let $E$ be an irreducible $W'_\l$-module and 
let $\boc=\boc_E$ (a two-sided cell of $W'_\l$). Let $E_0$ be an irreducible 
$W_\l$-module appearing in $E|_{W_\l}$ and let $\boc_0=\boc_{E_0}$ (a 
two-sided cell of $W_\l$). We have $\boc_0\sub\boc$ and the $a$-function of
$W_\l$ takes the same value on $\boc_0$ as the $a$-function of $W'_\l$ on
$\boc$. Hence $a_E=a_{E_0}$. Now $E_0\ot\sg$ appears in $(E\ot\sg)|_{W_\l}$
hence $b_{E_0\ot\sg}\le b_{E\ot\sg}$. Since $a_{E_0}\le b_{E_0\ot\sg}$ is
already known we see that $a_E\le b_{E\ot\sg}$. Thus (a) holds for 
$\fA=\HH'_\l$. 

We prove (b) for $\fA=\HH'_\l$.
Let $\boc$ be a two-sided cell of $W'_\l$. We can find a two-sided cell
$\boc_0$ of $W_\l$ such that $\boc_0\sub\boc$. We can find an irreducible
$W_\l$-module $E_0$ such that $\boc_{E_0}=\boc_0$ and $a_{E_0}=b_{E_0\ot\sg}$.
By 1.17(a) we can find an irreducible $W'_\l$-module $E'$ such that
$E_0\ot\sg$ appears in $(E'\ot\sg)|_{W_\l}$ and $b_{E'\ot\sg}=b_{E_0\ot\sg}$
Then $E_0$ appears in $E'|_{W_\l}$ hence $\boc_{E_0}\sub\boc_{E'}$. Thus 
$\boc_0\sub\boc_{E'}$ so that $\boc_{E'},\boc$ have nonempty intersection and 
$\boc_{E'}=\boc$; we see also that $a_{E_0}=a_{E'}$ so that
$a_{E'}=b_{E'\ot\sg}$. Thus the existence part of (b) is proved. Assume now 
that $E''$ is an irreducible $W'_\l$-module such that $\boc_{E''}=\boc$ and 
$a_{E''}=b_{E''\ot\sg}$. Let $E_1$ be an irreducible $W_\l$-module which
appears in $E''|_{W_\l}$ and let $\boc_1=\boc_{E_1}$ so that $\boc_1\sub\boc$
and $a_{E_1}=a_{E''}$. Replacing $\boc_1$ by $x\boc_1x\i$ for some 
$x\in\Om_\l$, we can assume that $\boc_1=\boc_0$. Now $E_1\ot\sg$ appears in
$(E''\ot\sg)|_{W_\l}$. Hence $b_{E_1\ot\sg}\le b_{E''\ot\sg}=a_{E''}=a_{E_1}$.
Since $a_{E_1}\le b_{E_1\ot\sg}$ by (a), it follows that
$a_{E_1}=b_{E_1\ot\sg}=b_{E''\ot\sg}$. Similarly we have 
$a_{E_0}=b_{E_0\ot\sg}$. By the uniqueness in (b) for $W_\l$ we see that 
$E_1\cong E_0$ as $W_\l$-modules; moreover $E_0\ot\sg$ is univalent. Now 
$E_0\ot\sg$ appears in $(E'\ot\sg)|_{W_\l}$ and 
$b_{E'\ot\sg}=b_{E_0\ot\sg}$; moreover, $E_0\ot\sg$ appears in 
$(E''\ot\sg)|_{W_\l}$ and $b_{E''\ot\sg}=b_{E_0\ot\sg}$. By the uniqueness in 
1.17(a) we see that $E''\ot\sg\cong E'\ot\sg$ so that $E''\cong E'$; from 
1.13(a) we see also that $E''\ot\sg$ is univalent. Thus (b) holds for 
$\fA=\HH'_\l$. 

We prove (a) for $\fA=\HH$. Let $E$ be an irreducible $W\TT_n$-module and 
let $\boc=\boc_E$ (a two-sided cell of $W\fs$). We can find $\l\in\fs$ such 
that $1_\l E\ne0$. Then $1_\l E$ is an irreducible $(1_\l\HH1_\l)^1$-module 
hence an irreducible $W'_\l$-module. Let $\boc_1=\boc_{1_\l E}$, a two-sided
cell of $W'_\l$. Then $\{w\cdo\l;w\in\boc_1\}\sub\boc$ and $a_{\boc_1}=a_\boc$
hence $a_{1\l E}=a_E$. Now $1_\l(E\ot\sg)=(1_\l E)\ot\sg$ hence by an argument
in 1.18 we have $b_{E\ot\sg}=b_{(1_\l E)\ot\sg}$. Since 
$a_{1_\l E}\le b_{(1_\l E)\ot\sg}$ is already known we see that
$a_E\le b_{E\ot\sg}$. Thus (a) holds for $\fA=\HH$.

We prove (b) for $\fA=\HH$. Let $\boc$ be a two-sided cell of $W\fs$. Note 
that $\boc$ is also a two-sided cell of $W\T\fo$ for some $W$-orbit 
$\fo$ in $W\fs$. We can find $\l\in\fo$ and a two-sided cell $\boc_0$ of
$W'_\l$ such that $\{w\cdo\l;w\in\boc_0\}\sub\boc$ and $a_{\boc_0}=a_\boc$. We
can find an irreducible $W'_\l$-module $E_0$ such that $\boc_{E_0}=\boc_0$ and 
$a_{E_0}=b_{E_0\ot\sg}$. Let $E'$ be the $W\TT_n$-module induced by the 
$W'_\l\TT_n$-module $E_0\ot\l$; note that $E'$ is irreducible. We have 
$a_{E'}=a_{E_0}$. Moreover, $E'\ot\sg$ is the $W\TT_n$-module induced by the 
$W'_\l\TT_n$-module $(E_0\ot\sg)\ot\l$; hence by an argument in 1.18 we have
$b_{E'\ot\sg}=b_{E_0\ot\sg}$. Thus we have $a_{E'}=b_{E'\ot\sg}$. Thus, the
existence part of (b) is proved.

Assume now that $E''$ is an irreducible $W\TT_n$-module such that 
$\boc_{E''}=\boc$ and $a_{E''}=b_{E''\ot\sg}$. We can find $\l'\in\fs$ such 
that $1_{\l'}E''\ne0$. Since $\boc_{E''}=\boc$ we must have $\l'\in\fo$.
Replacing $\l'$ by $w(\l')$ for some $w\in W$, we can assume that $\l'=\l$ so
that $1_\l E''\ne0$. Then $\{w\in W;w\cdo\l\in\boc\}$ is a two-sided cell of
$W'_\l$, necessarily equal to $\boc_0$; moreover, $\boc_{1_\l E''}=\boc_0$
hence $a_{1\l E''}=a_{E''}$. Now $E''$ is the $W\TT_n$-module induced by the
$W'_\l\TT_n$-module $(1_\l E'')\ot\l$ hence $E''\ot\sg$ is the $W\TT_n$-module
induced by the $W'_\l\TT_n$-module $((1_\l E'')\ot\sg)\ot\l$ hence by the
argument in 1.18 we have $b_{E''\ot\sg}=b_{(1_\l E'')\ot\sg}$. 
It follows that $a_{1_l E''}=b_{(1_\l E'')\ot\sg}$. Using this and
$a_{E_0}=b_{E_0\ot\sg}$ and also the uniqueness part in (b) for $W'_\l$ we see
that $E_0\cong 1_\l E''$ as $W'_\l$-modules. Since $E'$ (resp. $E''$) is
induced by the $W'_\l\TT_n$-module $E_0\ot\l$ (resp. $(1_\l E'')\ot\l$) we
deduce that $E'\cong E''$ as $W\TT_n$-modules. From 1.18 we see also that 
$E''$ is univalent. Thus (b) holds for $\fA=\HH$. This completes the proof of 
(a),(b) in all cases.

\mpb

The special representation of $\cw$ associated to $\boc$ in (b) 
is denoted by $E_\boc$. It is well defined up to isomorphism. By (b), the special 
representations of $\cw$ (up to isomorphism) are in natural bijection with the 
two-sided cells of $I$.

\subhead 1.20\endsubhead
Let $V,SV=\op_{i\ge0}S^iV$ be as in 1.17. For any $i\ge0$ we set
$$\ci^i=\sum_{i'>0,i''\ge0;i'+i''=i}(S^{i'}V)^WS^{i''}V\sub S^iV,$$
$\bS^iV=S^iV/\ci^i$, where $(S^{i'}V)^W$ is the space of $W$-invariants in
$S^{i'}V$ (we have used the algebra structure of $SV$). Let 
$\ci=\op_{i\ge0}\ci^i$, $\bS V=\op_{i\ge0}\bS V^i=SV/\ci$. Note that $\ci$ is an
ideal in $SV$ hence $\bS V$ is a (graded) algebra. Note also that the $W$-action
on $SV$ preserves $\ci$ hence it induces a $W$-action on $\bS V$ which is 
compatible with the grading and with the algebra structure. 

The following
property is well known:

(a) {\it $\bS^iV=0$ for $i>\nu$; $\bS^\nu$ is isomorphic to $\sg$ as a 
$W$-module (in particular it is $1$-dimensional). For $i\in[0,\nu]$, the 
bilinear pairing $\bS^iV\T\bS^{\nu-i}@>>>\bS^\nu$ given by multiplication in 
$\bS V$ is perfect.}
\nl
From (a) we deduce that for $i\in[0,\nu]$ we have
$\bS^iV\cong\bS^{\nu-i}\ot\sg$ as $W$-modules. (We use that any $W$-module is
isomorphic to its dual.) Hence if $\l\in\fs$ we have

(b) {\it $\bS^iV\cong\bS^{\nu-i}\ot\sg$ as $W'_\l$-modules. ($W'_\l$ acts by 
restriction of the $W$-action.) In particular. for any irreducible 
representation $E$ of $W'_\l$ we have 
$(\bS^iV)^E\cong(\bS^{\nu-i})^{E\ot\sg}$.}
\nl
Clearly, if $E$ is an irreducible $W'_\l$-module and $0\le i\le b_E$ then
$(\ci^i)^E=0$ hence $(\bS^iV)^E=(S^iV)^E$. In particular we have

(c) $(\bS^iV)^E=0$ for $0\le i<b_E$ and $(\bS^{b_E}V)^E=(S^{b_E}V)^E\ne0$.
\nl
Using (b),(c) we see that:

(c${}'$) {\it If $E$ is an irreducible $W'_\l$-module then $(\bS^iV)^E=0$ for
$i>\nu-b_{E\ot\sg}$. Moreover, $\dim(\bS^{\nu-b_{E\ot\sg}}V)^E$ is $1$ if $E$ 
is special and $\ge1$ if $E$ is not special.}
\nl
Since $a_E\le b_{E\ot\sg}$ (see 1.19(a)) with equality if and only if $E$ is 
special, we deduce:

(d) {\it If $E$ is an irreducible $W'_\l$-module then $(\bS^iV)^E=0$ for
$i>\nu-a_E$. Moreover, $\dim(\bS^{\nu-a_E}V)^E$ is $1$ if $E$ is special and
$0$ if $E$ is not special.}

\subhead 1.21\endsubhead
The $W\TT_n$-action on $S^iV\ot\bbq[\TT_n]$ (see 1.18) leaves
$\ci^i\ot\bbq[\TT_n]$ stable hence it induces a $W\TT_n$-action on 
$\bS^iV\ot\bbq[\TT_n]$. We show:

(a) {\it Let $E$ be an irreducible $W\TT_n$-module. We have 
$(\bS^iV\ot\bbq[\TT_n])^E=0$ for $i>\nu-a_E$. Moreover, 
$\dim(\bS^{\nu-a_E}V\ot\bbq[\TT_n])^E$ is $1$ if $E$ is special and $0$ if 
$E$ is not special.}
\nl
We can find $\l\in\fs$ such that $1_\l E\ne0$. Then $1_\l E$ is an irreducible
$W'_\l$-module and $E$ is induced by the representation $(1_\l E)\ot\l$ of the
subgroup $W'_\l\TT_n$. Since for $i\ge0$ the $W\TT_n$-module 
$\bS^iV\ot\bbq[\TT_n]$ is induced by the $W$-module $\bS^iV$ we see (using 
Frobenius reciprocity) that
$\dim(\bS^iV\ot\bbq[\TT_n])^E=\dim\Hom_W(E,\bS^iV)$. Since the $W$-module $E$ 
is induced by the $W'_\l$-module $1_\l E$, the last dimension is equal to
$\dim(\bS^iV)^{1_\l E}$. Now (a) follows from 1.20(d) applied to $1_\l E$
instead of $E$, using the equality $a_E=a_{1_\l E}$ and the fact that $E$ is
special if and only if $1_\l E$ is special.

\head 2. Truncated convolution of sheaves on $\tcb^2$\endhead
\subhead 2.1\endsubhead
For $w\in W$ and $\o\in\k\i(w)$ we define $G_w@>>>\TT$, $g\m g_\o$, by
$g\in\UU\o g_\o\UU$, $g_\o\in\TT$.
Let $\tcb=G/\UU$. Now $G\T\TT^2$ acts on $\tcb^2$ by 
$$(g,t_1,t_2):(x\UU,y\UU)\m(gxt_1^n\UU,gyt_2^n\UU).$$
The orbits of this action are indexed by $W$: to $w\in W$ corresponds the orbit 
$\tco_w=\{(x\UU,y\UU)\in\tcb^2;x\i y\in G_w\}$. The closure of $\tco_w$ in 
$\tcb^2$ is $\btco_w=\cup_{y\le w}\tco_w$.

Let $w\in W$, $\o\in\k\i(w)$. We define $j_\o:\tco_w@>>>\TT$ by 
$j_\o(x\UU,y\UU)=(x\i y)_\o$.
 Let $\l\in\fs$. We set $L_\l^\o=j_\o^*L_\l$. Now 
$L_\l$ is equivariant for the $G\T\TT^2$-action 
$(g,t_1,t_2):t\m w\i(t_1)^{-n}tt_2^n$ on $\TT$ and this action is compatible 
under $j_\o$ with the $G\T\TT^2$-action on $\tco_w$ (as above); hence 
$L_\l^\o$ is a $G\T\TT^2$-equivariant local system of rank $1$ on $\tco_w$
such that the induced action of $\TT_n^2$ on any stalk is via the character 
$(t_1,t_2)\m w(\l)\i(t_1)\l(t_2)$. (Note that $\TT_n^2$ acts trivially on 
$\tcb^2$.) Now let 
$$\G_w=\{(g,t_0,t_1)\in G\T\TT^2;gt_0^n\UU=\UU,g\o t_1^n\UU=\o\UU\}$$
be the stabilizer in $G\T\TT^2$ of $(\UU,\o\UU)\in\tco_w$. Setting $g=tu$ 
where $t\in\TT$, $u\in U_1:=\UU\cap\o\UU\o\i$, we can identify
$$\G_w=\{(u,t,t_0,t_1)\in U_1\T\TT^3;tt_0^n=1,w\i(t)t_1^n=1\}.$$
The subgroup 
$$\{(u,t,t_0,t_1)\in U_1\T\TT^3;t_0=w(t_1),tt_0^n=1\}$$
of $\G_w$ is clearly connected and has the same dimension as $\G$ (namely 
$\r+\dim U_1$) hence it is the identity component $\G_w^0$ of $\G_w$. We can 
view $\TT_n^2$ as the kernel of the surjective homomorphism 
$\G_w@>>>\TT\T U_1$, $(u,t,t_0,t_1)\m(u,t_0^n)$, whose restriction to $\G_w^0$
must also be surjective. It follows that $\G_w=T_n^2\G_w^0$ hence
$$\G_w/\G_w^0=\TT^2/(\TT^2\cap\G_w^0)=\TT^2/\{(t_0,t_1)\in\TT^2;t_0=w(t_1)\}.$$
Now the $G\T\TT^2$-equivariant local systems on $\tco_w$ correspond to 
representations of $\G_w/\G_w^0$ hence to representations of $\TT_n^2$ which 
are trivial on $\{(t_1,t_2)\in\TT_n^2;t_1=w(t_2)\}$.
We see that the local systems $L_\l^{\dw},\l\in\fs$ form a set of representatives for the
isomorphism classes of irreducible $G\T\TT^2$-equivariant local systems on 
$\tco_w$. 

\mpb

We define $\ti\fh:\tcb^2@>>>\tcb^2$ by $(xU,yU)\m(yU,xU)$. 
Let $w\in W$, $\o\in\k\i(w)$, $\l\in\fs$. Define $\x:\TT@>>>\TT$ by $\x(t)=w(t\i)$.
From the definitions we have $j_\o\ti\fh=\x j_{\o\i}:\tco_{w\i}@>>>\TT$.
Hence $\ti\fh^*L_\l^\o=L_{w(\l\i)}^{\o\i}$. (We use that $\x^*L_\l=L_{w(\l\i)}$.)

\subhead 2.2\endsubhead
Let $w\in W$, $\o\in\k\i(w)$. For $l\in\fs$ we shall view $L_\l^\o$ as a 
constructible sheaf on $\tcb^2$ which is $0$ on $\tcb^2-\tco_w$. Let 
$L_\l^{\o\sha}$ be its extension to an intersection cohomology complex on 
$\btco_w$ viewed as a complex on $\tcb^2$, equal to $0$ on $\tcb^2-\btco_w$. 
Let $\LL_\l^\o=L_\l^{\o\sha}\la|w|+\nu+2\r\ra$, a simple perverse sheaf on 
$\tcb^2$. Note that $L_\l^\o$ (resp. $\LL_\l^\o$) is (noncanonically) isomorphic to $L_\l^{\dw}$
(resp. $\LL_\l^{\dw}$). (We use 1.5). We have 
$$\ti\fh^*\LL_\l^\o=\LL_{w(\l\i)}^{\o\i}.\tag a$$

\subhead 2.3 \endsubhead
For $i<j$ in $[0,2]$ we define $p_{ij}:\tcb^3@>>>\tcb^2$ by 
$(x_0\UU,x_1\UU,x_2\UU)\m(x_i\UU,x_j\UU)$. For $L,L'$ in $\cd(\tcb^2)$ we set 
$L\cir L'=p_{02!}(p^*_{01}L\bxt p^*_{12}L')\in\cd(\tcb^2)$. This operation is
associative. Hence for ${}^1L,{}^2L,\do,{}^mL$ in $\cd(\tcb^2)$, 
${}^1L\cir{}^2L\cir\do\cir{}^mL\in\cd(\tcb^2)$ is defined.

We have $p_{02}=\bp_{02}f$ where $f:\tcb^3@>>>\tcb\T\cb\T\tcb$ is 
$(x_0\UU,x_1\UU,x_2\UU)\m(x_0\UU,x_1\BB x_1\i,x_2\UU)$ and 
$\bp_{02}:\tcb\T\cb\T\tcb@>>>\tcb^2$ is $(x_0\UU,B,x_2\UU)\m(x_0\UU,x_2\UU)$. 
We show:

(a) {\it Let $w,w'\in W$, $\o\in\k\i(w),\o'\in\k\i(w')$, $\l,\l'\in\fs$. If
$w'(\l')\ne\l$ then $L_\l^\o\cir L_{\l'}^{\o'}=0$.}
\nl
It is enough to show that 
$f_!(p^*_{01}L_\l^\o\bxt p^*_{12}L_{\l'}^{\o'})=0$. Hence it is enough to 
show that for any $(x_0,x_1,x_2)\in G^3$ and any $i$ we have
$$H^i_c(f\i(x_0\UU,x_1\BB x_1\i,x_2\UU),
p^*_{01}L_\l^\o\bxt p^*_{12}L_{\l'}^{\o'})=0.$$
We have 
$$f\i(x_0\UU,x_1\BB x_1\i,x_2\UU)=\{(x_0\UU,x_1\t\UU,x_2\UU);\t\in\TT\}$$ 
hence this fibre of $f$ is empty unless $\x_0\i x_1\in\UU\o t_0\UU$, 
$x_1\i x_2\in\UU\o't'_0\UU$ for some $t_0,t'_0$ in $\TT$ (which we now 
assume) so that the fibre can be identified with $\TT$. The restriction of 
$p_{01}$ (resp. $p_{12}$ to this fibre can be identified with $\t\m t_0\t$ 
(resp. $\t\m w'{}\i(\t\i)t'_0$). Then 
$p^*_{01}L_\l^\o\bxt p^*_{12}L_{\l'}^{\o'}$ becomes the local system 
$L_\l\ot L_{w'(\l')\i}=L_{\l w'(\l')\i}$ on $\TT$. It remains to use that 
$H^i_c(\TT,L_{\l_1})=0$ if $\l_1\in\fs-\{1\}$.

\subhead 2.4\endsubhead
Let $w,w'\in W$ be such that $|ww'|=|w|+|w'|$, let 
$\o\in\k\i(w),\o'\in\k\i(w')$ and let $\l,\l'\in\fs$. We show: 

(a) {\it If $w'(\l')=\l$, then we have canonically 
$L_\l^\o\cir L_{\l'}^{\o'}=L_{\l'}^{\o\o'}\ot\fL$.}
\nl
Let 
$$Y=\{(x\UU,y\UU,t,t')\in \tcb\T\tcb\T\TT\T\TT;x\i y\in\UU\o\o'w'{}\i(t)t'\UU\}.$$
Define $h:\TT^2@>>>\TT$ by $h(t,t')=w'{}\i(t)t'$. Define $j:Y@>>>\tco_{ww'}$ by 
$$j(x\UU,y\UU,t,t')=(x\UU,y\UU).$$
 Define $j_1:Y@>>>\TT$ by 
$j_1(x\UU,y\UU,t,t')=(t,t')$. Let $j'=j_{\o\o'}:\tco_{ww'}@>>>\TT$ be as in 2.1.
From the definitions we have 
$$L_\l^\o\cir L_{\l'}^{\o'}=j_!(j_1^*(L_\l\bxt L_{\l'}))=
j'{}^*(h_!(L_\l\bxt L_{\l'})).$$
To prove (a) it remains to show that $h_!(L_\l\bxt L_{\l'})=L_{\l'}\ot\fL$. 
Replacing $\l$ by $w'{}\i(\l)$ and $h$ by $h':\TT^2@>>>\TT$, $h'(t,t')=tt'$ 
we see that it is enough to show that $h'_!(L_\l\bxt L_\l)=L_\l\ot\fL$. We 
have $h_!(L_\l\bxt L_\l)=h'_!h'{}^*L_\l$ and it remains to use the equality 
$h'_!\bbq=\fL$.

\subhead 2.5\endsubhead
Let $s\in S$, $\l'\in\fs$. Let $L'$ be the local system of rank $1$ on 
$\tcb^2\T(\UU_s-\{1\})$ whose restriction to $\tcb^2\T\{\x\}$ is 
$L_{\l'}^{\s_\x}$ for any $\x\in\UU_s-\{1\}$ (see 1.2). Let
$\hL_{\l'}^s=c_!L'\in\cd(\tcb^2)$ where $c:\tcb^2\T(\UU_s-\{1\})@>>>\tcb^2$ is
the obvious projection. Clearly, we have 
$$\hL_{\l'}^s\Bpq\{\ch^2\hL_{\l'}^s[-2],\ch^1\hL_{\l'}^s[-1]).$$
Moreover, if $s\n W_{\l'}$ then $\ch^2\hL_{\l'}^s[-2]=0$, 
$\ch^1\hL_{\l'}^s[-1]=0$ hence $\hL_{\l'}^s=0$. If $s\in W_{\l'}$ then 
$\ch^2\hL_{\l'}^s[-2]=L_{\l'}^{\ds}[-2](-1)$, 
$\ch^1\hL_{\l'}^s[-2]=L_{\l'}^{\ds}[-1]$. 

\subhead 2.6\endsubhead
Let $s\in S$ and let $\l,\l'\in\fs$ be such that $s(\l')=\l$. From the 
definitions we see that: 
$$L_\l^{\ds}\cir L_{\l'}^{\ds\i}\Bpq\{L_{\l'}^1[-2](-1)\ot\fL,\hL_{\l'}^s\}.$$
Using the results in 2.5 we deduce:

(a) {\it If $s\n W_\l$, then 
$L_\l^{\ds}\cir L_{\l'}^{\ds\i}=L_{\l'}^1[-2](-1)\ot\fL$.}

(b) {\it If $s\in W_\l$, then}
$$L_\l^{\ds}\cir L_{\l'}^{\ds\i}\Bpq\{L_{\l'}^1[-2](-1)\ot\fL,
L_{\l'}^{\ds}[-2](-1)\ot\fL,L_{\l'}^{\ds}[-1]\ot\fL\}.$$
\nl
(Note that the conditions $s\in W_\l$ and $s\in W_{\l'}$ are equivalent.)

\subhead 2.7\endsubhead
Let $s\in S,w\in W$ be such that $|sw|<|w|$ and let $\o\in\k\i(w)$, 
$\l,\l'\in\fs$ be such that $w(\l')=\l$. We show:

(a) {\it If $s\n W_\l$ then 
$L_\l^{\ds}\cir L_{\l'}^\o\ot\fL=L_{\l'}^{\ds\o}[-2](-1)\ot\fL\ot\fL$.}

(b) {\it If $s\in W_\l,$ then}
$$L_\l^{\ds}\cir L_{\l'}^\o\ot\fL\Bpq\{L_{\l'}^{\ds\o}[-2](-1)\ot\fL^{\ot2},
L_{\l'}^\o[-2](-1)\ot\fL^{\ot2},L_{\l'}^\o[-1]\ot\fL^{\ot2}\}.$$
Using 2.4(a) we have $L_{\l'}^\o\ot\fL=L_{\l'' }^{\ds\i}\cir L_{\l'}^{\ds\o}$ 
where $\l''=(sw)(\l')$. Hence 
$$L_\l^{\ds}\cir L_{\l'}^\o\ot\fL=L_\l^{\ds}\cir L_{\l'' }^{\ds\i}\cir L_{\l'}^{\ds\o}.$$
We now apply the results in 2.6 to describe $L_\l^{\ds}\cir L_{\l'' }^{\ds\i}$. In case (a), we obtain
$$L_\l^{\ds}\cir L_{\l'}^\o\ot\fL=L_{\l''}^1\cir L_{\l'}^{\ds\o}\ot\fL[-2](-1).$$
By 2.4(a) this equals $L_{\l'}^{\ds\o}\ot\fL^{\ot2}[-2](-1)$, proving (a). In 
case (b) we obtain
$$\align&L_\l^{\ds}\cir L_{\l'}^\o\ot\fL\Bpq\\&\{L_{\l''}^1\cir L_{\l'}^{\ds\o}[-2](-1)\ot\fL,
L_{\l''}^{\ds\i}\cir L_{\l'}^{\ds\o}[-2](-1)\ot\fL,
L_{\l''}^{\ds\i}\cir L_{\l'}^{\ds\o}[-1]\ot\fL\}.\endalign$$
Here we substitute $L_{\l''}^1\cir L_{\l'}^{\ds\o}=L_{\l'}^{\ds\o}\ot\fL$,
$L_{\l''}^{\ds\i}\cir L_{\l'}^{\ds\o}=L_{\l'}^\o\ot\fL$ (see 2.4(a)) and (b) follows.

\subhead 2.8\endsubhead
We choose an $\FF_q$-rational structure on $G$. We shall assume that $\BB$ 
(hence $\UU$) is defined over $\FF_q$, that $\TT$ is defined and split over 
$\FF_q$ and that the integer $n$ in 1.4 divides $q-1$. Then for each $s\in S$,
the subgroup $\UU_s$ is defined over $\FF_q$; we shall also assume that in 1.3
we have $\x_s\in\UU_s(\FF_q)-\{1\}$. We have induced $\FF_q$-structures on $\cb,\tcb$. For any 
$w\in W$, $\co_w,\bco_w,\tco_w,\btco_w$ inherit natural 
$\FF_q$-structures. For any $w\in W$ we write $\k\i_q(w)$ instead of
$\k\i(w)\cap G(\FF_q)$; note that $\dw\in\k\i_q(w),(\dw)\i\in\k\i_q(w\i)$.
Now the local system $\io_!\bbq$ in 1.4 is naturally pure 
of weight zero (since $\bbq$ is so) and each of its direct summands $L_\l$ is 
itself naturally pure of weight zero (since $n$ divides $q-1$). If 
$\o\in\k_q\i(w)$, it follows that the local system $L_\l^\o$ on $\tco_w$ 
is naturally pure of weight zero. Hence $L_\l^{\o\sha}$, $\LL_\l^\o$ are 
naturally pure of weight zero. In particular, $L_\l^{\dw}$, $L_\l^{\dw\sha}$, 
$\LL_\l^{\dw}$ are naturally pure of weight zero.

Let $\cd^\spa\tcb^2$ be the subcategory of $\cd(\tcb^2)$ consisting of objects
which are restrictions of objects in the $G\T\TT^2$-equivariant derived 
category. Let $\cd^\spa_m\tcb^2$ be the subcategory of $\cd_m(\tcb^2)$ 
consisting of objects which are restrictions of objects in the mixed 
$G\T\TT^2$-equivariant derived category. Let $\cm^\spa\tcb^2$ (resp. 
$\cm^\spa_m\tcb^2$) be the subcategory of $\cd^\spa\tcb^2$ (resp. 
$\cd^\spa_m\tcb^2$) consisting of objects which are perverse sheaves.

If $w\in W$, $\o\in\k\i_q(w)$ then $L_\l^\o$ (resp. $L_\l^{\o\sha}$, 
$\LL_\l^\o$) is (noncanonically) isomorphic to $L_\l^{\dw}$ (resp. 
$L_\l^{\dw\sha}$, $\LL_\l^{\dw}$) as objects of $\cd_m(\tcb^2)$.

\subhead 2.9\endsubhead
Let $L\in\cd^\spa_m(\tcb^2)$. For any $w\in W$, $i\in\ZZ$, $\ch^iL|_{\tco_w}$ 
is a $G\T\TT^2$-equivariant local system with an induced mixed structure. We
can write it as 
$$\ch^iL|_{\tco_w}=\op_{\l\in\fs}V_{L,i,w,\l}\ot L^{\dw}_\l$$
where $V_{L,i,w,\l}$ are mixed $\bbq$-vector spaces. For $j\in\ZZ$ let 
$V_{L,i,w,\l,j}$ be the subquotient of $V_{L,i,w,\l,j}$ of pure weight $j$. We
set
$$\g(L)=\sum_{w\in W,\l\in\fs}\sum_{i,j\in\ZZ}(-1)^i(-v)^j\dim V_{L,i,w,\l,j}
T_w1_\l\in\HH.$$
For example, if $w\in W$ and $\o\in\k\i_q(w)$, $\l\in\fs$ then
$$\g(L_\l^\o)=\g(L_\l^{\dw})=T_w1_\l.$$
Note that 

(a) {\it if $(L,L',L'')$ is a distinguished triangle in $\cd^\spa_m(\tcb^2)$,
then $\g(L')=\g(L)+\g(L'')$.}

\subhead 2.10\endsubhead
Let $w,w'\in W$, $\o\in\k\i_q(w)$, $\o'\in\k\i_q(w')$, $\l,\l'\in\fs$. We show:
$$\g(L_\l^\o\cir L_{\l'}^{\o'})=(v^2-1)^\r\g(L_\l^\o)\g(L_{\l'}^{\o'}).\tag a$$
The right hand side of (a) is $(v^2-1)^\r T_w1_\l T_{w'}1_{\l'}$. We prove (a)
by induction on $|w|$. If $|w|=0$ then by 2.4(a), 2.3(a), the left hand side 
of (a) is $(v^2-1)^\r T_{w'}1_{\l'}$ (if $w'(\l')=\l$) and $0$ otherwise; this
is clearly equal to the right hand side of (a). Now assume that $|w|\ge1$. We 
can find $s\in S$ such that $|w|=|ws|+1$. The right hand side of (a) is 
$$(v^2-1)^\r T_w1_\l T_{w'}1_{\l'}=(v^2-1)^\r T_{ws}T_s1_\l T_{w'}1_{\l'}.
\tag b$$
If $w'(\l')\ne\l$ then (b) is $0$. If $w'(\l')=\l$ and $|sw'|=|w'|+1$ then 
(b) is $(v^2-1)^\r T_{ws}T_{sw'}1_{\l'}$. If $w'(\l')=\l$ and $|sw'|=|w'|-1$,
$s\n W_\l$ then (b) is $v^2(v^2-1)^\r T_{ws}T_{sw'}1_{\l'}$. If $w'(\l')=\l$ 
and $|sw'|=|w'|-1$, $s\in W_\l$ then (b) is 
$$(v^2-1)^\r T_{ws}(v^2T_{sw'}+(v^2-1)T_{w'})1_{\l'}.$$ 
Let $\l''=s(\l)$. By 2.4(a) we have 
$L_\l^\o\ot\fL=L_{\l''}^{\o\ds\i}\cir L_\l^{\ds}$, hence
$$\fK:=L_\l^\o\cir L_{\l'}^{\o'}\ot\fL\ot\fL
=L_{\l''}^{\o\ds\i}\cir L_\l^{\ds}\cir L_{\l'}^{\o'}\ot\fL.\tag c$$
If $w'(\l')\ne\l$, then $\fK=0$ by 2.3(a); hence in this case (a) holds. Thus 
we can assume that $w'(\l')=\l$. If $|sw'|=|w'|+1$ we have (using 2.4(a))
$$\fK=L_{\l''}^{\o\ds\i}\cir L_{\l'}^{\ds\o'}\ot\fL\ot\fL$$
hence by the induction hypothesis
$$(v^2-1)^{2\r}\g(L_\l^\o)\cir L_{\l'}^{\o'})=(v^2-1)^{3\r}T_{ws}1_{\l''}T_{sw'}1_{\l'};$$
hence in this case (a) holds. We now assume that $w'(\l')=\l$, $|sw'|=|w'|-1$.
Using 2.7(a),(b) to describe $L_\l^{\ds}\cir L_{\l'}^{\o'}\ot\fL$ we deduce 
that
$$\fK=L_{\l''}^{\o\ds\i}\cir L_{\l'}^{\ds\o'}\ot\fL^{\ot2}[-2](-1)
\text{ if }s\n W_\l,$$
$$\align&\fK\Bpq\{L_{\l''}^{\o\ds\i}\cir L_{\l'}^{\ds\o'}\ot\fL^{\ot2}[-2](-1),
\\&L_{\l''}^{\o\ds\i}\cir L_{\l'}^{\o'}\ot\fL^{\ot2}[-2](-1),
L_{\l''}^{\o\ds\i}\cir L_{\l'}^{\ds\o'}\ot\fL^{\ot2}[-1]\}\text{ if }
s\in W_\l.\endalign$$
It follows that 
$$\g(\fK)=v^2(v^2-1)^{2\r}\g(L_{\l''}^{\o\ds\i}\cir L_{\l'}^{\ds\o'})
\text{ if }s\n W_\l,$$
$$\g(\fK)=(v^2-1)^{2\r}(v^2\g(L_{\l''}^{\o\ds\i}\cir L_{\l'}^{\ds\o'})
+(v^2-1)\g(L_{\l''}^{\o\ds\i}\cir L_{\l'}^{\o'}))\text{ if }s\in W_\l.$$
Using the induction hypothesis we see that 
$$\g(\fK)=v^2(v^2-1)^{3\r}T_{ws}1_{\l''}T_{sw'}1_{\l'}\text{ if }s\n W_\l,$$
$$\g(\fK)=(v^2-1)^{3\r}(v^2T_{ws}1_{\l''}T_{sw'}1_{\l'}
+(v^2-1)T_{ws}1_{\l''}T_{w'}1_{\l'})\text{ if }s\in W_\l.$$
Thus, (a) holds.

\subhead 2.11\endsubhead
Let $r\ge1$ and let ${}^1L,{}^2L,\do,{}^rL$ be objects of $\cd^\spa_m\tcb^2$.
We show:
$$\g({}^1L\cir{}^2L\cir\do\cir{}^rL)
=(v^2-1)^{(r-1)\r}\g({}^1L)\g({}^2L)\do\g({}^rL).\tag a$$
When $r=1$, (a) is obvious. For $r\ge2$, (a) follows easily by induction from 
the case when $r=2$. Thus we may assume that $r=2$. For $j=1,2$ we have
$${}^jL\Bpq\{\ch^i({}^jL)[-i];i\in\ZZ\}$$
hence (using 2.9(a)), $\g({}^jL)=\sum_i(-1)^i\g(\ch^i({}^jL))$. Moreover,
$${}^1L\cir{}^2L\Bpq\{\ch^i({}^1L)\cir\ch^{i'}({}^2L)[-i-i'];i,i'\in\ZZ\}$$
hence (using 2.9(a)),
$$\g({}^1L\cir{}^2L)
=\sum_{i,i'\in\ZZ}(-1)^{i+i'}\g(\ch^i({}^1L)\cir\ch^{i'}({}^2L)).$$ 
Thus we can assume that ${}^1L=L_\l^{\dw}$, ${}^2L=L_{\l'}^{\dw'}$ where 
$w,w'\in W$, $\l,\l'\in\fs$. In this case, (a) follows from 2.10(a).

\subhead 2.12\endsubhead
Let $\l\in\fs$. We choose for each $\eta\in W_\l$ an element 
$\ddot\eta\in\k\i(\eta)$ as follows. Assume first that $|\eta|_\l=1$. We write
$\eta=s_1s_2\do s_rs_{r+1}s_r\do s_1$ with $s_1,s_2\do,s_{r+1}$ in $S$; we set
$$\ddot\eta=\ds_1\ds_2\do\ds_r\ds_{r+1}\ds_r\i\do\ds_1\i.$$
Assume next that $|\eta|_\l=m$. We write $\eta=\eta_1\eta_2\do\eta_m$ with 
$\eta_i\in W_\l$ such that $|\eta_1|_\l=\do=|\eta_m|_\l=1$, $|\eta|_\l=m$ 
and we set $\ddot\eta=\ddot\eta_1\ddot\eta_2\do\ddot\eta_m$. (In particular, 
$\ddot1=1$.)  

We now define for each $w\in W$ an element $\ddot w\in\k\i_q(w)$ as 
follows. There is a unique $z\in W$ such that $z=\min(wW_\l)$. We have 
$w=z\eta$ for a unique $\eta\in W_\l$. We set $\ddot w=\dz\ddot\eta$.

Let $w,y\in W$. Let $z=\min(wW_\l)$. We write $w=z\eta$ with $\eta\in W_\l$. 
Let $i\in\ZZ$. The statements (a),(b) below can be deduced from 
\cite{\ORA, 1.24} in the same way as \cite{\CSIII, 12.4} was deduced from \cite{\ORA, 1.24}.

(a) {\it We have $\ch^iL_\l^{\ddot w\sha}|_{\tco_y}=0$ unless $i$ is even and 
$y\in wW_\l$.}

(b) {\it Assume that $i$ is even and $y\in wW_\l$. We write $y=z\eta'$ with $\eta'\in W_\l$. We have 
$$\ch^iL_\l^{\ddot w\sha}|_{\tco_y}\Bpq\{(L_\l^{\ddot y})_h(-i/2);h\in[1,n_{\l,\eta',\eta,i}]\}$$
where $(L_\l^{\ddot y})_h$ are copies of $L_\l^{\ddot y}$ and 
$n_{\l,\eta',\eta,i}$ is the coefficient of $X^{i/2}$ in 
$$X^{(1/2)(|w|-|y|-|\eta|_\l|+|\eta'|_\l)}P^\l_{\eta',\eta}(X),$$
see 1.8.}
\nl
From (a),(b) we deduce:
$$L_\l^{\ddot w\sha}\Bpq\{(L_\l^{(z\eta')\ddot{}})_h\la-i\ra;
\eta'\in W_\l,h\in[1,n_{\l,\eta',\eta,i}]\}.\tag c$$
This is compatible with the natural mixed structures. Using 2.9(a), we deduce
$$\g(L_\l^{\ddot w\sha})=\sum_{\eta'\in W_\l;i\in2\ZZ}n_{\l,\eta',\eta,i}
\g(L_\l^{(z\eta')\ddot{}})v^i,$$
that is
$$\g(L_\l^{\ddot w\sha})=\sum_{\eta'\in W_\l}
v^{|w|-|z\eta'|-|\eta|_\l+|\eta'|_\l}P^\l_{\eta',\eta}(v^2)T_{z\eta'}1_\l$$
hence, using 1.8(a),
$$\g(L_\l^{\o\sha})=v^{|w|}c_{w,\l},\tag d$$
for any $\o\in\k\i_q(w)$.

\subhead 2.13\endsubhead
Let $w,w'\in W$, $\o\in\k\i(w)$, $\o'\in\k\i(w')$ and $\l,\l'\in\fs$.  
We show:

(a) {\it If $w'(\l')\ne\l$ then $L_\l^{\o\sha}\cir L_{\l'}^{\o'\sha}=0$.}

(b) {\it If $w'(\l')\ne\l$ then $L_\l^\o\cir L_{\l'}^{\o'\sha}=0$, 
$L_\l^{\o\sha}\cir L_{\l'}^{\o'}=0$.} 
\nl
We prove (a). We write $w=zw_1$ (resp. $w'=z'w'_1$) where $z=\min(zW_\l)$ 
(resp. $z'=\min(z'W_{\l'})$) and $w_1\in W_\l$ (resp. $w'_1\in W_{\l'}$). 
Using 2.12(c) it is enough to show that $L_\l^{\dy_1}\cir L_{\l'}^{\dy'_1}=0$ 
for any $y_1\in W_\l$, $y'_1\in W_{\l'}$. Using 2.3(a) it is enough to show 
that for $y'_1\in W_{\l'}$ we have $z'y'_1(\l')\ne\l$. We have 
$y'_1(\l')=\l'$, $w'_1(\l')=\l$ hence 
$z'y'_1(\l')=z'(\l')=z'w'_1(\l')=w'(\l')$. It remains to use our assumption 
that $w'(\l')\ne\l$.

We prove (b). For the first (resp. second) equality in (b) we repeat the proof
of (a) but take $y_1=w_1$ (resp. $y'_1=w'_1$).

It is not difficult to prove the following strengthening of (a).

(c) {\it Assume that for some $j\in\ZZ$, $\LL_\eta^{\dz}$ is a composition 
factor of $(\LL_\l^{\dw}\cir\LL_{\l'}^{\dw'})^j$. Then $\eta=\l'=w'{}\i(\l)$.}

\subhead 2.14\endsubhead
{\it In the remainder of this paper we fix a two-sided cell $\boc$ of $W\fs$ and we set 
$a=a(w\cdo\l)$ for some/any $w\cdo\l\in\boc$. Let $\fo$ be the unique 
$W$-orbit on $\fs$ such that $w\cdo\l\in\boc\implies\l\in\fo$.}
\nl
Let $Y=\tcb^2$. Let $\cm^{\preceq}Y$ (resp. $\cm^{\prec}Y$) be the subcategory
of $\cd^\spa Y$ whose objects are perverse sheaves $L$ such that any 
composition factor of $L$ is of the form $\LL_\l^{\dw}$ for some 
$w\cdo\l\preceq\boc$ (resp. $w\cdo\l\prec\boc$). Let $\cm^{\preceq}_mY$ (resp. 
$\cm^{\prec}_mY$) be the subcategory of $\cd^\spa_mY$ whose objects are in 
$\cm^{\preceq}Y$ (resp. $\cm^{\prec}Y$). Let $\cd^{\preceq}Y$ (resp. 
$\cd^{\prec}Y$) be the subcategory of $\cd^\spa Y$ whose objects are complexes
$L$ such that $L^j$ is in $\cm^{\preceq}Y$ (resp. $\cm^{\prec}Y$) for any $j$.
Let $\cd^{\preceq}_mY$ (resp. $\cd^{\prec}_mY$) be the subcategory of 
$\cd^\spa_mY$ whose objects are also in $\cd^{\preceq}Y$ (resp. 
$\cd^{\prec}Y$). Let $\cc^\spa Y$ be the subcategory of $\cm^\spa Y$ 
consisting of semisimple objects. Let $\cc_0^\spa Y$ be the subcategory of 
$\cm^\spa_mY$ consisting of those $L$ such that $L$ is pure of weight zero. 
Let $\cc^\boc Y$ be the subcategory of $\cm^\spa Y$ consisting of objects 
which are direct sums of objects of the form $\LL_\l^{\dw}$ with 
$w\cdo\l\in\boc$. Let $\cc_0^\boc Y$ be the subcategory of $\cc^\spa_0Y$ 
consisting of those $L\in\cc^\spa_0Y$ such that, as an object of $\cc^\spa Y$,
$L$ belongs to $\cc^\boc Y$. For $L\in\cc_0^\spa Y$ let $\un{L}$ be the 
largest subobject of $L$ such that as an object of $\cc^\spa Y$, we have 
$\un{L}\in\cc^\boc Y$.

\subhead 2.15\endsubhead
Let $r\ge1$. We define an action of 
$\cg=G\T(\UU^*)^r\T\TT^{2r+1}\T\UU^{r+1}$ on $G^{r+1}$ by
$$\align&(g,u_1,u_2,\do,u_r,t_1,\do,t_r,t'_0,t'_1,\do,t'_r,u'_0,u'_1,\do,u'_r):\\&
(g_0,g_1,\do,g_r)\m(gg_0t'_0{}^nu'_0{}\i,u_1t_1^{-n}g_1t'_1{}^nu'_1{}\i,\do,
u_rt_r^{-n}g_rt'_r{}^nu'_r{}\i).\endalign$$
The orbits of this action are indexed by $W^r$; to $\ww=(w_1,\do,w_r)\in W^r$ 
corresponds the orbit $G^{r+1}_\ww=G\T G_{w_1}\T G_{w_2}\T\do\T G_{w_r}$. The 
restriction of the $\cg$-action to the subgroup
$$\align&\cg':=\{(g,u_1,u_2,\do,u_r,t_1,\do,t_r,t'_0,t'_1,\do,t'_r,u'_0,u'_1,\do,u'_r)\in\cg;\\&
g=1,t_1=\do=t_r=t'_0=\do=t'_r=1,u'_0=u_1,u'_1=u_2,\do,u'_{r-1}=u_r\}\endalign$$
(isomorphic to $\UU^{r+1}$) is free and the map $\th:G^{r+1}@>>>\tcb^{r+1}$ 
given by 
$$(g_0,g_1,\do,g_r)\m(g_0\UU,g_0g_1\UU,\do,g_0g_1\do g_r\UU)$$
identifies $\tcb^{r+1}$ with $\cg'\bsl G^{r+1}$. For 
$\ww=(w_1,\do,w_r)\in W^r$ and $J\sub[1,r]$ we define
$$G^{r+1,J}_\ww=\{(g_0,g_1,\do,g_r)\in G^{r+1};g_i\in\bG_{w_i} \frl i\in J,
g_i\in G_{w_i}\frl i\in[1,r]-J\},$$ 
$$\align&\tco^J_\ww=\{(x_0\UU,x_1\UU,\do,x_r\UU)\in\tcb^{r+1};\\&x_{i-1}\i x_i\UU\in
\bG_{w_i}\frl i\in J,x_{i-1}\i x_i\in G_{w_i}\frl i\in[1,r]-J\}.\endalign$$
Now $\th$ identifies $\tco^J_\ww$ with $\cg'\bsl G^{r+1,J}_\ww$ and 
$$\tco^J_\ww=\sqc_{\yy=(y_1,y_2,\do,y_r)\in W^r;y_i\le w_i\frl i\in J,y_i=w_i
\frl i\in[1,r]-J}\tco_\yy^\emp.$$
Note that $\tco_\ww^\emp$ is irreducible of dimension $\nu+(r+1)\r+|\ww|$ where
$$|\ww|=|w_1|+|w_2|+\do+|w_r|.$$
Until the end of 2.22 we fix $\ww=(w_1,\do,w_r)\in W^r$, 
$\pmb\o=(\o_1,\o_2,\do,\o_r)$ such that $\o_i\in\k\i_q(w_i)$ for 
$i=1,\do,r$ and $\pmb\l=(\l_1,\l_2,\do,\l_r)\in\fs^r$.

\mpb

Define $c:\tco_\ww^\emp@>>>\TT^r$ and $\tc:G^{r+1}_\ww@>>>\TT^r$ by
$$c(x_0\UU,x_1\UU,\do,x_r\UU)=((x_0\i x_1)_{\o_1},(x_1\i x_2)_{\o_2},\do,
(x_{r-1}\i x_r)_{\o_r}),$$
$$\tc(g_0,g_1,\do,g_r)=((g_1)_{\o_1},(g_2)_{\o_2},\do,(g_r)_{\o_r}),$$
so that $\tc=c\th$.
Let $M^{\pmb\o}_{\pmb\l}\in\cd_m(\tcb^{r+1})$ be the local system 
$c^*(L_{\l_1}\bxt\do\bxt L_{\l_r})$ on $\tco_\ww^\emp$ extended by $0$ on 
$\tcb^{r+1}-\tco_\ww^\emp$. 
Let $\tM^{\pmb\o}_{\pmb\l}\in\cd_m(G^{r+1})$ be the local system 
$\tc^*(L_{\l_1}\bxt\do\bxt L_{\l_r})$ on $G^{r+1}_\ww$ extended by $0$ on 
$G^{r+1}-G^{r+1}_\ww$. Note that 
$$\tM^{\pmb\o}_{\pmb\l}=\th^*M^{\pmb\o}_{\pmb\l}.$$
From the definitions we have
$$M^{\pmb\o}_{\pmb\l}=p_{01}^*L_{\l_1}^{\o_1}\ot p_{12}^*L_{\l_2}^{\o_2}
\ot\do\ot p_{r-1,r}^*L_{\l_r}^{\o_r}.$$
(Here $p_{ij}:\tcb^{r+1}@>>>\tcb^2$ are the obvious projections.)
Note that $\tM^{\pmb\o}_{\pmb\l}\in\cd_m(G^{r+1})$ is $\cg$-equivariant.
Indeed, $\cg$ acts on $\TT^r$ by
$$\align&(g,u_1,u_2,\do,u_r,t_1,\do,t_r,t'_0,t'_1,\do,t'_r,u'_0,u'_1,\do,u'_r):
\\&(t''_1,t''_2,\do,t''_r)\m(w_1\i(t_1^{-n})t''_1t'_1{}^n,w_2\i(t_2^{-n})t''_2t'_2{}^n,\do,
w_r\i(t_r^{-n})t''_rt'_r{}^n),\endalign$$
$\th$ is compatible with the $\cg$-actions and 
$L_{\l_1}\bxt\do\bxt L_{\l_r}$ is a $\cg$-equivariant local system. Let $J\sub[1,r]$. We set 
$$M^{\pmb\o,J}_{\pmb\l}=p_{01}^*{}^1L\ot p_{12}^*{}^2L\ot\do
\ot p_{r-1,r}^*{}^rL\in\cd_m(\tcb^{r+1}),$$
$$L^{\pmb\o,J}_{\pmb\l}=p_{0r!}M^{\pmb\o,J}_{\pmb\l}\la|\ww|\ra
={}^1L\cir{}^2L\cir\do\cir{}^rL\la|\ww|\ra\in\cd_m(\tcb^2),$$
where ${}^iL$ is $L_{\l_i}^{\o_i\sha}$ for $i\in J$ and $L_{\l_i}^{\o_i}$ for 
$i\n J$. Note that $M^{\pmb\o,\emp}_{\pmb\l}=M^{\pmb\o}_{\pmb\l}$. Moreover,

{\it $M^{\pmb\o,J}_{\pmb\l}$ is the intersection cohomology complex of 
$\tco_\ww^J$ with coefficients in $M^{\pmb\o}_{\pmb\l}$.}
\nl
To prove this, it is enough to show that $\th^*M^{\pmb\o,J}_{\pmb\l}$ is the 
intersection cohomology complex of $G^{r+1,J}_\ww$ with coefficients in 
$\tM^{\pmb\o}_{\pmb\l}$; this is immediate.

\mpb

Consider the free $\TT^{r-1}$-action on $\tcb^{r+1}$ given by
$$\align&(\t_1,\t_2,\do,\t_{r-1}):(x_0\UU,x_1\UU,\do,x_{r-1}\UU,x_r\UU)\m\\&
(x_0\UU,x_1\t_1\UU,\do,x_{r-1}\t_{r-1}\UU,x_r\UU).\endalign$$
Note that $\tco^J_\ww$ is stable under this $\TT^{r-1}$-action. We also have
a free $\TT^{r-1}$-action on $\TT^r$ given by
$$\align&(\t_1,\t_2,\do,\t_{r-1}):(t_1,t_2,\do,t_r)\m\\&
(t_1\t_1,w_2\i(\t_1\i)t_2\t_2,w_3\i(\t_2\i)t_3\t_3,\do,
w_{r-1}\i(\t_{r-2}\i)t_{r-1}\t_{r-1},w_r\i(\t_{r-1}\i)t_r).\endalign$$
Let ${}'\tcb^{r+1}=\TT^{r-1}\bsl\tcb^{r+1}$. Let 
${}'\tco^J_\ww=\TT^{r-1}\bsl\tco^J_\ww$ (a locally closed subvariety of 
${}'\tcb^{r+1}$). Let ${}'\TT^r=\TT^{r-1}\bsl\TT^r$. Note that 
${}'\tco^\emp_\ww=\TT^{r-1}\bsl\tco^\emp_\ww$ is an open dense smooth 
irreducible subvariety of ${}'\tco^J_\ww$. Now $c:\tco_\ww^\emp@>>>\TT^r$ is
compatible with the $\TT^{r-1}$-actions on $\tco_\ww^\emp,\TT^r$ hence it 
induces a map ${}'c:{}'\tco_\ww^\emp@>>>{}'\TT^r$. The homomorphism 
$c':\TT^r@>>>\TT$ given by
$$(t_1,t_2,\do,t_r)\m t_1w_2(t_2)w_2w_3(t_3)\do w_2w_3\do w_r(t_r)$$
is constant on each orbit of the $\TT^{r-1}$-action on $\TT^r$ hence it 
induces a morphism ${}'\TT^r@>>>\TT$ whose composition with ${}'c$ is 
denoted by $\bc:{}'\tco_\ww^\emp@>>>\TT$. Let ${}'M^{\pmb\o,\emp}_{\pmb\l}$ be
the local system $\bc^*L_{\l_1}$ on ${}'\tco_\ww^\emp$ extended by $0$ on 
${}'\tcb^{r+1}-{}'\tco_\ww^\emp$. Let 
${}'M^{\pmb\o,J}_{\pmb\l}\in\cd_m({}'\tcb^{r+1})$ be the intersection 
cohomology complex of ${}'\tco_\ww^J$ with coefficients in 
${}'M^{\pmb\o,\emp}_{\pmb\l}$ extended by $0$ on ${}'\tcb^{r+1}-{}'\tco_\ww^J$.
Let $\bp_{0r}:{}'\tco_\ww^J@>>>\tcb^2$ be the map induced by $p_{0r}:\tco_\ww^J@>>>\tcb^2$. We define 
${}'L^{\pmb\o,J}_{\pmb\l}\in\cd_m^\spa\tcb^2$ as follows: if
$$\l_k=w_{k+1}(\l_{k+1})\text{ for }k=1,2,\do,r-1$$
(in which case we say that $\pmb\l$ is $\ww$-adapted) we set 
$${}'L^{\pmb\o,J}_{\pmb\l}=\bp_{0r!}{}'M^{\pmb\o,J}_{\pmb\l}\la|\ww|\ra;$$
if $\pmb\l$ is not $\ww$-adapted, we set ${}'L^{\pmb\o,J}_{\pmb\l}=0$.
Let $\hc^J:\tco_\ww^J@>>>{}'\tco_\ww^J$ be the obvious (orbit) map. We show:

(a) {\it If $\pmb\l$ is $\ww$-adapted then 
$M^{\pmb\o,J}_{\pmb\l}=(\hc^J)^*{}'M^{\pmb\o,J}_{\pmb\l}$.}
\nl
Since $\hc^J$ is a $\TT^{r-1}$-bundle, it is enough to show that 
$$M^{\pmb\o,\emp}_{\pmb\l}=(\hc^\emp)^*{}'M^{\pmb\o,\emp}_{\pmb\l}$$
or that 
$$c^*(L_{\l_1}\bxt\do\bxt L_{\l_r})=(\hc^\emp)^*\bc^*L_{\l_1}.$$
We have a commutative diagram
$$\CD
\tco_\ww^\emp      @>c>>  \TT^r   \\
@V\hc^\emp VV                   @Vc'VV  \\
{}'\tco_\ww^\emp@>\bc>>\TT
\endCD$$
hence $(\hc^\emp)^*\bc^*=c^*c'{}^*$ and it is enough to show that
$$L_{\l_1}\bxt\do\bxt L_{\l_r}=c'{}^*L_{\l_1}.$$
This follows from the equality
$$\l_1(c'(t_1,t_2,\do,t_r))=
\l_1(t_1)\l_2(t_2)\do\l_r(t_r)\text{ for all }(t_1,t_2,\do,t_r)\in\TT_n^r$$
which is a consequence of $\pmb\l$ being $\ww$-adapted.

We now show:

(b) {\it We have $L^{\pmb\o,J}_{\pmb\l}=\fL^{\ot(r-1)}\ot{}'L^{\pmb\o,J}_{\pmb\l}$.}
\nl
If $\pmb\l$ is not $\ww$-adapted then from 2.3(a), 2.13(a),(b), we see that 
$L^{\pmb\o,J}_{\pmb\l}=0$ hence (b) holds. We now assume that $\pmb\l$ is 
$\ww$-adapted. Using (a) we have
$$\align&L^{\pmb\o,J}_{\pmb\l}=p_{0r!}M^{\pmb\o,J}_{\pmb\l}\la|\ww|\ra=
p_{0r!}(\hc^J)^*{}'M^{\pmb\o,J}_{\pmb\l}\la|\ww|\ra\\&=
\bp_{0r!}(\hc^J)_!(\hc^J)^*{}'M^{\pmb\o,J}_{\pmb\l}\la|\ww|\ra=
\bp_{0r!}((\hc^J)_!\bbq)\ot(\hc^J)^*{}'M^{\pmb\o,J}_{\pmb\l}\la|\ww|\ra)
\endalign$$
and it remains to use that $(\hc^J)_!\bbq=\fL^{\ot(r-1)}$.

\mpb

We prove the following result.

(c) {\it There is a natural bijection between $\fs^r$ and the set of 
isomorphism classes of irreducible $\cg$-equivariant local systems on the
$\cg$-orbit $G^{r+1}_\ww$: to \lb $\pmb\l'=(\l'_1,\l'_2,\do,\l'_r)\in\fs^r$ 
corresponds the local system
$\tM^{\pmb\o'}_{\pmb\l'}|_{G_{\ww'}}=(\th^*M^{\pmb\o'}_{\pmb\l'})_{G_{\ww'}}$
where $\pmb\o'=(\dw_1,\dw_2,\do,\dw_r)$.}
\nl
Let $\G$ be the stabilizer of $(1,\dw_1,\do,\dw_r)\in G^{r+1}_\ww$ in $\cg$.
We have
$$\align&\G=\{(g,u_1,u_2,\do,u_r,t_1,\do,t_r,t'_0,t'_1,\do,t'_r,u'_0,u'_1,\do,u'_r)
\in\cg;\\&g=u'_0t'_0{}^{-n},\dw_1\i u_1\dw_1=u'_1,\dw_2\i u_2\dw_2=u'_2,\do,
\dw_r\i u_r\dw_r=u'_r,\\&t_1^n=w_1(t'_1{}^n),\do,t_r^n=w_r(t'_r{}^n).\}\endalign$$
The closed subgroup 
$$\align&\{(g,u_1,u_2,\do,u_r,t_1,\do,t_r,t'_0,t'_1,\do,t'_r,u'_0,u'_1,\do,u'_r)\in
\cg;\\&g=u'_0t'_0{}^{-n},\dw_1\i u_1\dw_1=u'_1,\dw_2\i u_2\dw_2=u'_2,\do, 
\dw_r\i u_r\dw_r=u'_r,\\&t_1=w_1(t'_1),\do,t_r=w_r(t'_r)\}\endalign$$
of $\G$ is clearly connected of the same dimension as $\G$ (namely
$(r+1)\nu+(r+1)\r$) hence it is equal to the identity component $\G^0$ of $\G$.
We can view $\TT_n^{2r}$ as the kernel of the surjective homomorphism 
$\G@>>>G\T\UU^r\T\TT^{r+1}\T\UU^{r+1}$,
$$\align&(g,u_1,u_2,\do,u_r,t_1,\do,t_r,t'_0,t'_1,\do,t'_r,u'_0,u'_1,\do,u'_r)\m\\&
(u_1,u_2,\do,u_r,t_1^n,\do,t_r^n,t'_0,u'_0,u'_1,\do,u'_r)\endalign$$
whose restriction to $\G^0$ must also be surjective. It follows that 
$\G=\TT_n^{2r}\G^0$ hence 
$$\align&\G/\G^0=\TT_n^{2r}/(\TT_n^{2r}\cap\G^0)\\&=
\TT_n^{2r}/\{(t_1,\do,t_r,t'_1,\do,t'_r)\in\TT_n^{2r};
t_1=w_1(t'_1),\do,t_r=w_r(t'_r)\}.\endalign$$
Note that the irreducible $\cg$-equivariant local systems on $G^{r+1}_\ww$
correspond to irreducible representations of $\G/\G^0$ hence to 
representations of $\TT_n^{2r}$ which are trivial on 
$$\{(t_1,\do,t_r,t'_1,\do,t'_r)\in T_n^{2r};t_1=w_1(t'_1),\do,t_r=w_r(t'_r)\}.$$
Such representations are uniquely determined by their restriction to
$$\{(t_1,\do,t_r,t'_1,\do,t'_r)\in\TT_n^{2r};t_1=t_2=\do=t_r=1\}$$
hence they are in natural bijection with $\fs^r$. This proves (c).

Using (c) and the fact that the $\cg$-orbits on $G^{r+1}$ are indexed by
$W^r$, we deduce:

(d) {\it There is a natural bijection between $W^r\T\fs^r$ and the set of 
isomorphism classes of simple $\cg$-equivariant perverse sheaves on 
$G^{r+1}$: to \lb $\ww'=(w'_1,w'_2,\do,w'_r)\in W^r$ and
$\pmb\l'=(\l'_1,\l'_2,\do,\l'_r)\in\fs^r$ corresponds the simple perverse
sheaf $\th^*M^{\pmb\o',[1,r]}_{\pmb\l'}\la\dim G^{r+1}_{\ww'}\ra$ where 
$\pmb\o'=(\dw'_1,\dw'_2,\do,\dw'_r)$.}

\subhead 2.16\endsubhead
We preserve the setup of 2.15. We assume that $J=[1,r]$. In this case, 
$\bp_{0r}:{}'\tco_\ww^{[1,r]}@>>>\tcb^2$ is clearly a proper morphism. Hence, 
by Deligne's theorem,

(a) {\it ${}'L^{\pmb\o,[1,r]}_{\pmb\l}$ is pure of weight zero.}
\nl
We set $L=L^{\pmb\o,[1,r]}_{\pmb\l}$, ${}'L={}'L^{\pmb\o,[1,r]}_{\pmb\l}$. 
From (a) it follows that for $j\in\ZZ$, ${}'L^j$ is pure of weight $j$ hence
$${}'L^j=\op_{w\cdo\l\in W\fs}\bV_{w\cdo\l,j}\LL_\l^{\dw}\tag b$$
where $\bV_{w\cdo\l,j}$ are mixed $\bbq$-vector spaces of pure weight $j$. For
any $(w,\l)\in W\fs$ and any $j\in\ZZ$ we show:

(c) {\it We have $$\dim\bV_{w\cdo\l,j}=N(w\cdo\l,-j+\nu+2\r)$$ where 
$N(w,\l,k)=N(w,\l,-k)\in\NN$ are given by the equality (in $\HH$):}
$$c_{w_1\cdo\l_1}c_{w_2\cdo\l_2}\do c_{w_r\cdo\l_r}
=\sum_{w\cdo\l\in W\fs,k\in\ZZ}N(w\cdo\l,k)v^kc_{w\cdo\l}.$$
From 2.11(a) and 2.12(d) we have (setting $\d=(r-1)\r$):
$$\align&\g(L)=(v^2-1)^\d\g(L_{\l_1}^{\o_1\sha})\g(L_{\l_2}^{\o_2\sha})\do 
\g(L_{\l_r}^{\o_r\sha})v^{-|\ww|}\\&
=(v^2-1)^\d v^{|w_1|}c_{w_1\cdo\l_1}v^{|w_2|}c_{w_2\cdo\l_2}\do v^{|w_r|}
c_{w_r\cdo\l_r}v^{-|\ww|}\\&
=(v^2-1)^\d c_{w_1\cdo\l_1}c_{w_2\cdo\l_2}\do c_{w_r\cdo\l_r}
=(v^2-1)^\d\sum_{w\cdo\l\in W\fs,k\in\ZZ}N(w\cdo\l,k)v^kc_{w\cdo\l}.\endalign$$
From the definitions we have (using (b)):
$$\align&\g({}'L)=\sum_j(-1)^j\g({}'L^j)
=\sum_j(-1)^j\sum_{w\cdo\l\in W\fs}\dim\bV_{w\cdo\l,j}(-v)^j
\g(L_\l^{\dw\sha}\la|w|+\nu+2\r\ra)\\&=\sum_j\sum_{w\cdo\l\in W\fs}\dim\bV_{w\cdo\l,j}v^{|w|}
c_{w\cdo\l}v^{j-|w|-\nu-\r}.\endalign$$
From 2.15(b) we have $\g(L)=(v^2-1)^\d\g({}'L)$ hence
$$\align&\g(L)=(v^2-1)^\d\sum_j\sum_{(w,\l)\in W\fs}\dim\bV_{w\cdo\l,j}c_{w\cdo\l}
v^{j-\nu-2\r}\\&=(v^2-1)^\d\sum_{w\cdo\l\in W\fs,k\in\ZZ}N(w\cdo\l,k)v^kc_{w\cdo\l}.\endalign$$
Since $c_{w\cdo\l}$ are linearly independent in $\HH$, it follows that for any
$w\cdo\l$ we have
$$\sum_j\dim\bV_{w\cdo\l,j}v^{j-\nu-2\r}=\sum_{k\in\ZZ}N(w\cdo\l,k)v^k
=\sum_{k\in\ZZ}N(w\cdo\l,-k)v^k$$
hence for any $j$ we have $\dim\bV_{w\cdo\l,j}=N(w\cdo\l,-j+\nu+2\r)$, as 
required.

\subhead 2.17\endsubhead
We preserve the setup of 2.15; let $J\sub[1,r]$. We set 
$L^J=L^{\pmb\o,J}_{\pmb\l}$, ${}'L^J={}'L^{\pmb\o,J}_{\pmb\l}$. 
As in 2.16, we set $\d=(r-1)\r$.

We now analyze the complex $\fL^{\ot(r-1)}\in\cd_m(\text{point})$.
We can find free abelian groups $\cx_{2\d-i}$ of rank $\bin{\d}{i}$, 
($i\in\ZZ$) such that $\cx_{2\d}=\ZZ$, complexes 
$R_{\le2\d-i}\in\cd_m(\text{point})$ ($i\in[0,\d+1]$) and distinguished 
triangles 
$$(R_{\le2\d-i-1},R_{\le2\d-i},\cx_{2\d-i}\ot\bbq(i-\d)[i-2\d]),\qua (i\in[0,\d])$$ 
in $\cd_m(\text{point})$ such that $R_{\le2\d}=\fL^{\ot(r-1)}$, 
$R_{\le\d-1}=0$. It follows that for $i\in[0,\d]$ we have distinguished 
triangles in $\cd_m(\tcb^2)$:
$$(R_{\le2\d-i-1}\ot{}'L^J,R_{\le2\d-i}\ot{}'L^J,\cx_{2\d-i}(i-\d)\ot{}'L^J[i-2\d])$$
hence we have exact sequences
$$\align&\do@>>>\cx_{2\d-i}(i-\d)\ot({}'L^J)^{-2\d+i+j-1}@>>>
(R_{\le2\d-i-1}\ot{}'L^J)^j\\&@>>>(R_{\le2\d-i}\ot{}'L^J)^j@>>>
\cx_{2\d-i}(i-\d)\ot({}'L^J)^{-2\d+i+j}@>>>\do.\endalign$$
Thus, setting 
$$\car_{i,j}=(R_{\le2\d-i}\ot{}'L^J)^j\text{ for }i\in[0,\d+1],$$ 
$$\cp_{i,j}=\cx_{2\d-i}(i-\d)\ot({}'L^J)^{-2\d+i+j}\text{ for }i\in[0,\d],$$
we have $\car_{\d+1,j}=0$ for all $j$ and, for any $i\in[0,\d]$, we have an 
exact sequence in $\cm_m(\tcb^2)$:
$$\do@>>>\cp_{i,j-1}@>>>\car_{i+1,j}@>>>\car_{i,j}@>>>\cp_{i,j}@>>>
\car_{i+1,j+1}@>>>\car_{i,j+1}@>>>\do.\tag a$$
Note that for any $j$ we have
$$\car_{0,j}=(L^J)^j,\tag b$$
$$\cp_{0,j}=({}'L^J)^{j-2\d}(-\d).\tag c$$
Indeed, (c) is obvious; (b) follows from 2.15(b):
$$\car_{0,j}=(R_{\le2\d}\ot{}'L^J)^j=(\fL^{\ot(r-1)}\ot{}'L^J)^j=(L^J)^j.$$

\subhead 2.18\endsubhead
We preserve the setup of 2.15; there is no assumption on $J$. The restriction 
of $M:=M^{\pmb\o,[1,r]}_{\pmb\l}$ to $\tco^J_\ww$ (an open dense subset of 
$\tco^{[1,r]}_\ww$) is the same as the restriction of
$M^J:=M^{\pmb\o,J}_{\pmb\l}$ to $\tco^J_\ww$; the restriction of $M$ to 
$\tco^{[1,r]}_\ww-\tco^J_\ww$ (a closed subset of $\tco^{[1,r]}_\ww$),
extended by $0$ on $\tcb^{r+1}-(\tco^{[1,r]}_\ww-\tco^J_\ww)$, is denoted by 
$\dM^J$. We have a distinguished triangle 
$$(M^J,M,\dM^J)\tag a$$
in $\cd_m(\tcb^{r+1})$. We have the following result.

(b) {\it Let $h\in\ZZ$. Let $K$ be either $\dM^J$ or $M^J$. Any composition 
factor of $K^h\in\cm(\tcb^{r+1})$ is of the form 
$M^{\pmb\o',[1,r]}_{\pmb\l'}\la|\ww'|+\nu+(r+1)\r\ra$ for some 
$\ww'=(w'_1,w'_2,\do,w'_r)\in W^r$, $\pmb\l'=(\l'_1,\l'_2,\do,\l'_r)\in\fs^r$
such that $w_i=w'_i$, $\l_i=\l'_i$ for all $i\in J$; here 
$\pmb\o'=(\dw'_1,\dw'_2,\do,\dw'_r)$.}
\nl
It is enough to show that for any $h$, any composition factor of
$(\th^*K)^h$ ($\th$ as in 2.15) is of the form
$\th^*M^{\pmb\o',[1,r]}_{\pmb\l'}\la|\ww'|+\nu+(r+1)\r\ra$ for some 
$\ww',\pmb\l',\pmb\o'$ as in (b). To see this we use the fact that
$(\th^*K)^h$ is a $\cg$-equivariant perverse sheaf on $G^{r+1}$ (it is 
obtained from the $\cg$-equivariant object $\tM^{\pmb\o}_{\pmb\l}$ by 
operations which preserve $\cg$-equivariance: passage to an intersection
cohomology complex, restriction to a $\cg$-invariant subvariety, taking a
perverse cohomology sheaf) and that all simple $\cg$-equivariant perverse 
sheaves on $G^{r+1}$ are of the form 
$\th^*M^{\pmb\o',[1,r]}_{\pmb\l'}\la|\ww'|+\nu+(r+1)\r\ra$ with $\ww'$, 
$\pmb\l'$, $\pmb\o'$ as in (b) (but with $\ww'$ unrestricted), see 2.15(d).
The fact that the $\ww',\pmb\l'$ which appear are restricted as in (b) is immediate.

We show:

(c) {\it $(\dM^J\la|\ww|+\nu+(r+1)\r-1\ra)^j=0$ for any $j>0$.}
\nl
It is enough to show that $\dim\supp\ch^h(\dM^J[|\ww|+\nu+(r+1)\r-1])\le-h$ 
for any $h\in\ZZ$. Assume first that $h\le-|\ww|-\nu-(r+1)\r$. Since $M$ is an
intersection cohomology complex with support of dimension $|\ww|+\nu+(r+1)\r$,
we have $\dim\supp\ch^{h-1}(M[|\ww|+\nu+(r+1)\r])<-h+1$ hence
$$\dim\supp\ch^{h-1}(\dM^J[|\ww|+\nu+(r+1)\r])<-h+1$$
hence $\dim\supp\ch^{h-1}(\dM^J[|\ww|+\nu+(r+1)\r])\le-h$, hence
$$\dim\supp\ch^h(\dM^J[|\ww|+\nu+(r+1)\r-1])\le-h.$$
Next we assume that $h=-|\ww|-\nu-(r+1)\r+1$. Then
$$\align&\dim\supp\ch^{h-1}(\dM^J[|\ww|+\nu+(r+1)\r])\le
\dim(\tco^{[1,r]}_\ww-\tco^J_\ww)\le\\&|\ww|+\nu+(r+1)\r-1=-h,\endalign$$
hence
$$\dim\supp\ch^h(\dM^J[|\ww|+\nu+(r+1)\r-1])\le-h.$$
Finally, assume that $h\ge-|\ww|-\nu-(r+1)\r+2$. Then 
$\ch^{h-1}(M[|\ww|+\nu+(r+1)\r])=0,$ hence
$\ch^{h-1}(\dM^J[|\ww|+\nu+(r+1)\r])=0,$ hence 
$\ch^h(\dM^J[|\ww|+\nu+(r+1)\r-1])=0$. This proves (c).

\subhead 2.19\endsubhead
We preserve the setup of 2.15; there is no assumption on $J$. We shall need a 
variant of the results in 2.18.
The restriction of ${}'M:={}'M^{\pmb\o,[1,r]}_{\pmb\l}$ to ${}'\tco^J_\ww$ (an
open dense subset of ${}'\tco^{[1,r]}_\ww$) is the same as the restriction of 
${}'M^J:={}'M^{\pmb\o,J}_{\pmb\l}$ to ${}'\tco^J_\ww$; the restriction of 
${}'M$ to ${}'\tco^{[1,r]}_\ww-{}'\tco^J_\ww$ (a closed subset of 
${}'\tco^{[1,r]}_\ww$), extended by $0$ on 
${}'\tcb^{r+1}-({}'\tco^{[1,r]}_\ww-{}'\tco^J_\ww)$, is denoted by ${}'\dM^J$.
We have a distinguished triangle 
$$({}'M^J,{}'M,{}'\dM^J)\tag a$$
in $\cd({}'\tcb^{r+1})$. The following result can be deduced from 2.18(b).

(b) {\it Let $h\in\ZZ$. Let $'K$ be either ${}'\dM^J$ or ${}'M^J$. Any 
composition factor of $({}'K)^h\in\cm({}'\tcb^{r+1})$ is of the form 
${}'M^{\pmb\o',[1,r]}_{\pmb\l'}\la|\ww'|+\nu+2\r\ra$ for some 
$\ww'=(w'_1,w'_2,\do,w'_r)\in W^r$, $\pmb\l'=(\l'_1,\l'_2,\do,\l'_r)\in\fs^r$ 
such that $w_i=w'_i$, $\l_i=\l'_i$ for all $i\in J$, and $\pmb\l'$ is 
$\ww$-adapted; here $\pmb\o'=(\dw'_1,dw'_2,\do,\dw'_r)$.}
\nl
We note:

(c) {\it $({}'\dM^J\la|\ww|+\nu+2\r-1\ra)^j=0$ for any $j>0$.}
\nl
The proof is entirely similar to that of 2.18(c); alternatively it can be 
deduced from 2.18(c).

\subhead 2.20\endsubhead
We preserve the setup of 2.15. Assume that $w_u\cdo\l_u\in\boc$ for some 
$u\in J$. We set $L^J=L^{\pmb\o,J}_{\pmb\l}$, 
${}'L^J={}'L^{\pmb\o,J}_{\pmb\l}$. Let $\dM^J$ be as in 2.18; let ${}'\dM^J$ 
be as in 2.19. Let $\dL^J=p_{0r!}\dM^J\la|\ww|\ra\in\cd(\tcb^2)$, 
${}'\dL^J=\bp_{0r!}{}'\dM^J\la|\ww|\ra\in\cd(\tcb^2)$. Let $j\in\ZZ$. We have 
the following results, in which $\car_{i,j}$, $\cp_{i,j}$ are as in 2.17 with 
$J=[1,r]$ and $\d=(r-1)\r$. 

(a) {\it We have $(L^J)^j\in\cm^{\preceq}\tcb^2$. If $j>2\d+\nu+2\r+(r-1)a$ 
then $(L^J)^j\in\cm^{\prec}\tcb^2$.}

(b) {\it We have $(\dL^J)^j\in\cm^{\preceq}\tcb^2$. If 
$j\ge2\d+\nu+2\r+(r-1)a$ then $(\dL^J)^j\in\cm^{\prec}\tcb^2$.}

(c) {\it We have $({}'L^J)^j\in\cm^{\preceq}\tcb^2$. If $j>\nu+2\r+(r-1)a$ 
then $({}'L^J)^j\in\cm^{\prec}\tcb^2$.}

(d) {\it We have $({}'\dL^J)^j\in\cm^{\preceq}\tcb^2$. If $j\ge\nu+2\r+(r-1)a$
then $({}'\dL^J)^j\in\cm^{\prec}\tcb^2$.}

(e) {\it If $i\in[0,\d+1]$, $J=[1,r]$, then 
$\car_{i,j}\in\cm^{\preceq}\tcb^2$.}

(f) {\it If $i\in[0,\d+1]$, $j>2\d-i+\nu+2\r+(r-1)a$, $J=[1,r]$, then  
$\car_{i,j}\in\cm^{\prec}\tcb^2$.}

(g) {\it If $i\in[0,\d]$, $J=[1,r]$, then $\cp_{i,j}\in\cm^{\preceq}\tcb^2$. 
If $i\in[0,\d]$, $j>2\d-i+\nu+2\r+(r-1)a$, $J=[1,r]$, then 
$\cp_{i,j}\in\cm^{\prec}\tcb^2$.}
\nl
We prove (e) by descending induction on $i$. If $i=\d+1$ then, since 
$\car_{\d+1,j}=0$, there is nothing to prove. Now assume that $i\in[0,\d]$. 
Assume that $\LL_\l^{\dw}$ is a composition factor of $\car_{i,j}$ (without 
the mixed structure). We must show that $w\cdo\l\preceq\boc$. By the induction
hypothesis we can assume that $\LL_\l^{\dw}$ is not a composition factor of 
$\car_{i+1,j}$; hence by 2.17(a) it is a composition factor of $\cp_{i,j}$. 
Hence $\LL_\l^{\dw}$ is a composition factor of $({}'L^{[1,r]})^{-2\d+i+j}$.
Hence $\bV_{w\cdo\l,k}$ in 2.16 is $\ne0$ for some $k$. Using 2.16(c) we see 
that $N(w\cdo\l,k)\ne0$ for some $k$. Using the definition of $N(w\cdo\l,k)$ 
we see that $w\cdo\l\preceq\boc$ (recall that $w_u\cdo\l_u\in\boc$ for some $u$) and 
(e) is proved.

We prove (f) by descending induction on $i$. If $i=\d+1$ then, since 
$\car_{\d+1,j}=0$, there is nothing to prove. Now assume that $i\in[0,\d]$. 
Assume that $\LL_\l^{\dw}$ is a composition factor of $\car_{i,j}$ (without 
the mixed structure). We must show that $w\cdo\l\prec\boc$. By the induction 
hypothesis we can assume that $\LL_\l^{\dw}$ is not a composition factor of 
$\car_{i+1,j}$ (we have $j>2\d-i-1+\nu+2\r+(r-1)a$); hence by 2.17(a), 
$\LL_\l^{\dw}$ is a composition factor of $\cp_{i,j}$. Hence $\LL_\l^{\dw}$ is
a composition factor of $({}'L^{[1,r]})^{-2\d+i+j}$. Hence 
$\bV_{w\cdo\l,-2\d+i+j}$ in 2.16 is $\ne0$. Using 2.16(c) we see that 
$N(w\cdo\l,2\d-i-j+\nu+2\r)\ne0$. We have $2\d-i-j+\nu+2\r<-(r-1)a$. Using 
1.10(a) we deduce that $w\cdo\l\prec\boc$ and (f) is proved.

We prove (g). This follows from the exact sequence 2.17(a) (with $J=[1,r]$) 
using (e),(f).

We prove (a) assuming that $J=[1,r]$. From (e),(f) we have 
$\car_{0,j}\in\cm^{\preceq}\tcb^2$ and $\car_{0,j}\in\cm^{\prec}\tcb^2$ if 
$j>2\d+\nu+2\r+(r-1)a$. Using 2.17(b) we deduce that (a) holds (when 
$J=[1,r])$.

We prove (c) assuming that $J=[1,r]$. From (g) we have 
$\cp_{0,j}\in\cm^{\preceq}\tcb^2$ and $\cp_{0,j}\in\cm^{\prec}\tcb^2$ if 
$j>2\d+\nu+2\r+(r-1)a$. Using 2.17(c) we deduce that $({}'L^J)^{j-2\d}(-\d)$ 
is in $\cm^{\preceq}\tcb^2$ and is in $\cm^{\prec}\tcb^2$ if 
$j-2\d>\nu+2\r+(r-1)a$. We deduce that (c) holds (when $J=[1,r]$).

We prove (b). Assume that $j\in\ZZ$ and $w\cdo\l\in\fs$ is such that 
$\LL_\l^{\dw}$ is a composition factor of $(\dL^J)^j$ (without mixed 
structure). Then there exists $h$ such that $\LL_\l^{\dw}$ is a composition 
factor of $(p_{0r!}(\dM^J[|\ww|])^h[-h])^j$. We have $(\dM^J[|\ww|)^h\ne0$  
hence $(\dM^J[|\ww|+\nu+(r+1)\r-1])^{h-\nu-(r+1)\r+1}\ne0$, hence by 2.18(c), 
$h-\nu-(r+1)\r+1\le0$. From 2.18(b) we see that there exist 
$\ww'=(w'_1,w'_2,\do,w'_r)\in W^r$, $\pmb\l'=(\l'_1,\l'_2,\do,\l'_r)\in\fs^r$ 
such that $w_i=w'_i,\l_i=\l'_i$ for all $i\in J$ and  $\LL_\l^{\dw}$ is a 
composition factor of 
$$(p_{0r!}(M^{\pmb\o',[1,r]}_{\pmb\l}[|\ww'|+\nu+(r+1)\r][-h])^j
=L^{\pmb\o',[1,r]}_{\pmb\l'}[|\ww'|])^{j+\nu+(r+1)\r-h};$$
here $\pmb\o'=(\dw_1',\do,\dw'_r)$. From the part of (a) that is already 
proved (for $\ww',\pmb\l'$ instead of $\ww,\pmb\l$) we see that 
$w\cdo\l\preceq\boc$ and that if $j+\nu+(r+1)\r-h>2\d+\nu+2\r+(r-1)a$ that is, 
if $j>2\d+(r-1)\r+(r-1)a+h$, then $w\cdo\l\prec\boc$. If 
$j\ge2\d+\nu+2\r+(r-1)a$ then, using that $0>h-\nu-(r+1)\r$, we see that we 
have indeed $j>2\d+(r-1)\r+(r-1)a+h$. This proves (b).

The proof of (d) is entirely similar to that of (b); it uses the already 
proved part of (c) and it uses the results of 2.19 instead of those in 2.18.

We prove (a) without assumption on $J$. Applying $p_{0r}!$ to 2.18(a) we get a
distinguished triangle $(L^J,L^{[1,r]},\dL^J)$. This gives rise to an exact 
sequence
$$(\dL^J)^{j-1}@>>>(L^J)^j@>>>(L^{[1,r]})^j@>>>(\dL^J)^j.$$
Using this together with (b) and the already proved part of (a) we see that 
(a) holds in general.

We prove (c) without assumption on $J$. Applying $\bp_{0r}!$ to 2.19(a) we get
a distinguished triangle $({}'L^J,{}'L^{[1,r]},{}'\dL^J)$. This gives rise to
an exact sequence
$$({}'\dL^J)^{j-1}@>>>({}'L^J)^j@>>>({}'L^{[1,r]})^j@>>>({}'\dL^J)^j.$$
Using this together with (d) and the already proved part of (c) we see that 
(c) holds in general.

\subhead 2.21\endsubhead
Let $j\in\ZZ$. In (a),(b) below, $\car_{i,j}$ is as in 2.17 with arbitrary $J$.

(a) {\it If $i\in[0,\d+1]$ then $\car_{i,j}\in\cm^{\preceq}\tcb^2$.}

(b) {\it If $i\in[0,\d+1]$, $j>2\d-i+\nu+2\r+(r-1)a$ then 
$\car_{i,j}\in\cm^{\prec}\tcb^2$.}
\nl
Note that (a),(b) are generalizations of 2.20(e), 2.20(f) (which correspond to
the case $J=[1,r]$).

We prove (a),(b) by descending induction on $i$. If $i=\d+1$ then, since 
$\car_{\d+1,j}=0$, there is nothing to prove. Now assume that $i\in[0,\d]$. 
Assume that $\LL_\l^{\dw}$ is a composition factor of $\car_{i,j}$ (without 
the mixed structure). We must show that $w\cdo\l\preceq\boc$ and that, if 
$j>2\d-i+\nu+2\r+(r-1)a$, then $w\cdo\l\prec\boc$. Using 2.17(a), we see that
$\LL_\l^{\dw}$ is a composition factor of $\car_{i+1,j}$ or of $\cp_{i,j}$. In
the first case, using the induction hypothesis we see that 
$w\cdo\l\preceq\boc$ and that, if $j>2\d-i+\nu+2\r+(r-1)a$ (so that 
$j>2\d-i-1+\nu+2\r+(r-1)a$), then $w\cdo\l\prec\boc$. In the second case,
$\LL_\l^{\dw}$ is a composition factor of $({}'L^J)^{-2\d+i+j}$. Using 2.20(c)
we see that $w\cdo\l\preceq\boc$ and that, if $j>2\d-i+\nu+2\r+(r-1)a$ (so 
that $-2\d+i+j>\nu+2\r+(r-1)a$), then $w\cdo\l\prec\boc$. This proves (a),(b).

\subhead 2.22\endsubhead
We preserve the setup of 2.20. In (a),(b) below we take $j=2\d+\nu+2\r+(r-1)a$.
We have the following results.

(a) {\it We have canonically $\un{gr_j((L^J)^j)}@>\si>>\un{gr_j({}'L^J){}^{j-2\d}(-\d))}$.}  

(b) {\it We have canonically $\un{gr_j((L^J)^j)}@>\si>>\un{gr_j((L^{[1,r]})^j)}$.}

(c) {\it If $(w_1,\l_1)\in\boc,\do,(w_r,\l_r)\in\boc$, $(w,\l)\in\boc$ and 
$j=2\d+\nu+\r+(r-1)a$ then the multiplicity of $\LL_\l^{\dw}$ in 
$gr_j((L^J)^j)$ is 
$$\sum h^*_{z_1\cdo\l'_1,w_2\cdo\l_2,z_2\cdo\l'_2}
h^*_{z_2\cdo\l'_2,w_3\cdo\l_3,z_3\cdo\l'_3}\do  
h^*_{z_{r-1}\cdo\l'_{r-1},w_r\cdo\l_r,z_r\cdo\l'_r}$$
where the sum is taken over all $z_1\cdo\l'_1,z_2\cdo\l'_2,\do,z_r\cdo\l'_r$
in $\boc$ such that $z_1\cdo\l'_1=w_1\cdo\l_1$, $z_r\cdo\l'_r=w\cdo\l$.}

(d) {\it Assume that $i\in[0,\d+1]$. Then $\car_{i,j}$ (notation of 2.17) is 
mixed of weight $\le j-i$.}
\nl 
We prove (d) by descending induction on $i$. If $i=\d+1$ there is nothing to 
prove. Assume now that $i\le\d$. By Deligne's theorem, ${}'L^J$ is mixed of 
weight $\le0$; hence $({}'L^J)^{-2\d+i+j}$ is mixed of weight $\le-2\d+i+j$ 
and $\cx_{2\d-i}(i-\d)\ot({}'L^J)^{-2\d+i+j}$ is mixed of weight 
$\le-2\d+i+j-2(i-\d)=j-i$. In other words, $\cp_{i,j}$ (notation of 2.17) is 
mixed of weight $\le j-i$. Thus in the exact sequence 
$\car_{i+1,j}@>>>\car_{i,j}@>>>\cp_{i,j}$ coming from 2.17(a) in which
$\car_{i+1,j}$ is mixed of weight $\le j-i-1<j-i$ (by the induction 
hypothesis) and $\cp_{i,j}$ is mixed of weight $\le j-i$ we must have that 
$\car_{i,j}$ is mixed of weight $\le j-i$. This proves (d).

We prove (a). From 2.17(a) we deduce an exact sequence 
$$gr_j(\car_{1,j})@>>>gr_j(\car_{0,j})@>>>gr_j(\cp_{0,j})@>>>
gr_j(\car_{1,j+1}).$$
By (d) we have $gr_j(\car_{1,j})=0$. We have $gr_j(\car_{0,j})=gr_j((L^J)^j)$, 
$gr_j(\cp_{0,j})=gr_j(({}'L^J)^{-2\d+j}(-\d))$. Moreover, by 2.21(b) we have 
$\car_{1,j+1}\in\cd^{\prec}\tcb^2$ since $j+1>2\d-1+\nu+2\r+(r-1)a$. 
It follows that $gr_j(\car_{1,j+1})\in\cd^{\prec}\tcb^2$.
Thus the exact sequence above induces an isomorphism as in (a).

We prove (b). As in 2.20 we have an exact sequence
$$(\dL^J)^{j-1}@>>>(L^J)^j@>>>(L^{[1,r]})^j@>>>(\dL^J)^j.$$
This gives rise to an exact sequence
$$gr_j((\dL^J)^{j-1})@>>>gr_j((L^J)^j)@>>>gr_j((L^{[1,r]})^j)@>>>
gr_j((\dL^J)^j).$$
From 2.20(b) we have $gr_j((\dL^J)^j)\in\cm^{\prec}(\tcb^2)$. It is then 
enough to show that $gr_j((\dL^J)^{j-1})=0$. Since $\dM^J\la|\ww|\ra$ is the
restriction of a pure complex of weight $0$ to a subspace, it is mixed of 
weight $\le0$ (see \cite{\BBD, 5.1.14}). Hence 
$\dL^J=(p_{0r!}\dM^J\la|\ww|\ra$ is mixed of weight $\le0$ and $(\dL^J)^{j-1}$
is mixed of weight $le j-1$ (see \cite{\BBD, 5.4.1}). Thus, 
$gr_j((\dL^J)^{j-1})=0$ as required. This proves (b).

We prove (c). By (a) and (b), the multiplicity in (c) is equal to the 
multiplicity of $\LL_\l^{\dw}$ in
$$gr_j({}'L^{[1,r]}){}^{j-2\d}(-\d))=({}'L^{[1,r]})^{j-2\d}(-\d)$$
(we use the fact that ${}'L^{[1,r]}$ is pure of weight zero); thus it is equal
to $\bV_{w\cdo\l,j-2\d}=\bV_{w\cdo\l,\nu+2\r+(r-1)a}$ hence also to 
$N(w\cdo\l,-(r-1)a)$ (see 2.16(c)). It remains to use the equality 1.10(c).

\subhead 2.23\endsubhead
Let $L,L'\in\cd^\spa\tcb^2$. We show:

(a) {\it Assume that $L,L'\in\cm^\spa\tcb^2$ and that either $L$ or $L'$ is in
$\cd^{\preceq}\tcb^2$. If $j>a-\nu$ then $(L\cir L')^j\in\cm^{\prec}\tcb^2$.}
\nl
We can assume that $L=\LL_\l^{\dw},L'=\LL_{\l'}^{\dw'}$ and that either 
$w\cdo\l\in\boc$ or $w'\cdo\l'\in\boc$. According to 2.20(b) we have 
$$(L_\l^{\dw\sha}\la|w|\ra\cir L_{\l'}^{\dw'\sha}\la|w'|\ra)^{j'}\in
\cm^{\prec}\tcb^2$$ 
if $j'>4\r+\nu+a$. 
Hence 
$$(L_\l^{\dw\sha}\la|w|+\nu+2\r\ra\cir L_{\l'}^{\dw'\sha}\la|w'|+\nu+2\r\ra)^j\in\cm^{\prec}\tcb^2$$ 
if $j+2\nu+4\r>4\r+\nu+a$ that is, if $j>a-\nu$. This proves (a).

(b) {\it If $L\in\cd^{\preceq}\tcb^2$ or $L'\in\cd^{\preceq}\tcb^2$ then 
$L\cir L'\in\cd^{\preceq}\tcb^2$. If $L\in\cd^{\prec}\tcb^2$ or 
$L'\in\cd^{\prec}\tcb^2$ then $L\cir L'\in\cd^{\prec}\tcb^2$.}
\nl
The first assertion of (b) is shown in the same way as (a). The second 
assertion of (b) can be reduced to the first assertion.

\subhead 2.24\endsubhead  
For $L,L'\in\cc^\boc_0\tcb^2$ we set
$$L\un\cir L'=\un{(L\cir L')^{\{a-\nu\}}}\in\cc^\boc_0\tcb^2.$$
(For the notation ${}^{\{i\}}$ see 0.2.)
Now let $L,L',L''\in\cc^\boc_0\tcb^2$. By 2.23(b), we have 
$L'\cir L''\in\cd_m^{\preceq}\tcb^2$ and the functor 
$\Ph:\cd_m^{\preceq}\tcb^2@>>>\cd_m^{\preceq}\tcb^2$, $\Ph(K)=L\cir K$, is 
well defined. We note that, by 2.23(a), (i),(ii) below hold.

(i) {\it If $X_0\in\cm_m^{\preceq}\tcb^2$ then $(\Ph(X_0))^h\in\cm_m^{\prec}\tcb^2$
for any $h>a-\nu$.}

(ii) {\it $L'\cir L''$ is mixed of weight $\le0$ and $(L'\cir L'')^h\in\cm_m^{\prec}\tcb^2$ 
for any $h>a-\nu$.}
\nl
Similarly we have $L\cir L'\in\cd_m^{\preceq}\tcb^2$ and the functor 
$\Ph':\cd_m^{\preceq}\tcb^2@>>>\cd_m^{\preceq}\tcb^2$, $\Ph'(K)=K\cir L'$, is 
well defined. Moreover, (iii),(iv) below hold.

(iii) {\it If $X_0\in\cm_m^{\preceq}\tcb^2$ then $(\Ph'(X_0))^h\in\cm_m^{\prec}\tcb^2$ for any 
$h>a-\nu$.}

(iv) {\it $L\cir L'$ is mixed of weight $\le0$ and $(L\cir L')^h\in\cm_m^{\prec}\tcb^2$ for any 
$h>a-\nu$.}
\nl
We now apply \cite{\CONV, 1.12} with $Y_1,Y_2$ replaced by $\tcb^2$ we see that
$$\un{\Ph(\un{(L'\circ L'')^{\{a+\r-\nu\}}})^{\{a+\r-\nu\}}}=\un{(\Ph(L'\cir L''))^{\{2a-2\nu\}}},$$
$$\un{\Ph'(\un{(L\circ L')^{\{a+\r-\nu\}}})^{\{a+\r-\nu\}}}=\un{(\Ph'(L\cir L'))^{\{2a-2\nu\}}}.$$
Thus, we have
$$L\un{\cir}(L'\un{\cir}L'')=\un{(L\cir L'\cir L'')^{\{2a-2\nu\}}},$$
$$(L\un{\cir}L')\un{\cir}L''=\un{(L\cir L'\cir L'')^{\{2a-2\nu\}}}.$$
Hence
$$L\un{\cir}(L'\un{\cir}L'')=(L\un{\cir}L')\un{\cir}L''.$$
We see that $L,L'\m L\un{\cir}L'$ defines a monoidal structure on 
$\cc^\boc_0\tcb^2$. Hence if ${}^1L,{}^2,\do,{}^rL$ are in $\cc^\boc_0\tcb^2$ 
then ${}^1L\un{\cir}{}^2L\un{\cir}\do\un{\cir}{}^rL\in\cc^\boc_0\tcb^2$
is well defined; using \cite{\CONV, 1.12} repeatedly, we have
$${}^1L\un{\cir}{}^2L\un{\cir}\do\un{\cir}{}^rL=
\un{({}^1L\cir{}^2L\cir\do\cir{}^rL)^{\{(r-1)(a-\nu)\}}}.\tag a$$

\subhead 2.25\endsubhead
Let $L_0,L_1\in\cc^\boc_0\tcb^2$. We show that we have canonically
$$\fD(L_0\un{\cir}L_1)=\fD(L_0)\un{\cir}\fD(L_1).\tag a$$
(Note that in the right hand side, $\un\cir$ is relative to $\ti{\boc}$
instead of $\boc$, see below.)
We can assume that $L_i=\LL_{\l_i}^{\dw_i}$ ($i=0,1$) where 
$w_i\cdo\l_i\in\boc$ ($i=0,1$). Let $\tL_i=\fD(L_i)=\LL_{\l_i\i}^{\dw_i}$, 
$i=0,1$. Note that $w_i\cdo\l_i\i\in\ti{\boc}$ (see 1.14) and
$\tL_0,\tL_1\in\cc^{\ti\boc}_0\tcb^2$. It is enough to show that 
$$\fD(L_0\un\cir L_1)=\tL_0\un\cir\tL_1.$$
If $w_1(\l_1)\ne\l_0$ (that is, $w_1(\l_1\i)\ne\l_0\i$), we have 
$L_0\cir L_1=0$, $\tL_0\cir\tL_1=0$ hence both sides of (a) are zero. Now 
assume that $w_1(\l_1)=\l_0$. Let $L=L^{\pmb\o,J}_{\pmb\l}$, 
${}'L={}'L^{\pmb\o,J}_{\pmb\l}$, ${}'M={}'M^{\pmb\o,J}_{\pmb\l}$ where 
$\pmb\o=(\dw_0,\dw_1)$, $\pmb\l=(\l_0,\l_1)$, $J=\{1\}$. Let 
$\tL=L^{\pmb\o,J}_{\pmb\l'}$, ${}'\tL={}'L^{\pmb\o,J}_{\pmb\l'}$     
${}'\tM={}'\tM^{\pmb\o,J}_{\pmb\l'}$ where $\pmb\l'=(\l_0\i,\l_1\i)$. By 
definition we have 
$$L=(L_0\la-|w_0|-\nu-2\r\ra\cir L_1\la-|w_1|-\nu-2\r\ra)\la|w_0+|w_1|\ra$$
hence $L=(L_0\cir L_1)\la-2\nu-4\r\ra$. Thus,
$$L^{4\r+\nu+a}=(L_0\cir L_1)^{a-\nu}(-\nu-2\r)$$ 
and
$$gr_{4\r+\nu+a}(L^{4\r+nu+a})=gr_{a-\nu}(L_0\cir L_1)^{a-\nu}(-\nu-2\r).$$
By 2.22(a) we have 
$$\un{gr_{4\r+\nu+a}(L^{4\r+nu+a})}=\un{{}'L^{\nu+2\r+a}}(-\r)$$
hence
$$\un{gr_{a-\nu}(L_0\cir L_1)^{a-\nu}}((a-\nu)/2)=
\un{{}'L^{\nu+2\r+a}}((a+2\r+\nu)/2).$$
Similarly,
$$\un{gr_{a-\nu}(\tL_0\cir\tL_1)^{a-\nu}}((a-\nu)/2)
=\un{{}'\tL^{\nu+2\r+a}}((a+2\r+\nu)/2).$$
It is then enough to show that 
$$\fD(\un{{}'L^{\nu+2\r+a}}((a+2\r+\nu)/2))
=\un{{}'\tL^{\nu+2\r+a}}((a+2\r+\nu)/2).$$
This would follow from the stronger result that
$$\fD({}'L^{\nu+2\r+a})={}'\tL^{\nu+2\r+a}(a+2\r+\nu).$$
Recall that ${}'L=\bp_{0r!}({}'M)\la|w_0|+|w_1|\ra$ and similarly
${}'\tL=\bp_{0r!}({}'\tM)\la|w_0|+|w_1|\ra$ where 
${}'M\la|w_0|+|w_1|+\nu+2\r\ra$, ${}'\tM\la|w_0|+|w_1|+\nu+2\r\ra$ are 
perverse sheaves, each being $\fD$ of the other. Since $\bp_{01}$ is
proper, $\bp_{01!}$ commutes with $\fD$. It follows that 
$$\fD({}'L\la\nu+2\r\ra)={}'\tL\la\nu+2\r\ra,$$
hence 
$$\fD(({}'L\la\nu+2\r\ra)^j)=({}'\tL\la\nu+2\r\ra)^{-j},$$
that is
$$\fD({}'L^{\nu+2\r+j})={}'\tL^{\nu+2\r-j}(\nu+2\r)$$
for any $j$; in particular,
$$\fD({}'L^{a+2\nu+\r})={}'\tL^{-a+\nu+2\r}(\nu+2\r).$$
Thus it is enough to prove
$${}'\tL^{-a+\nu+2\r}(\nu+2\r)={}'\tL^{\nu+2\r+a}(a+2\r+\nu),$$
that is
$${}'\tL^{-a+\nu+2\r}={}'\tL^{\nu+2\r+a}(a).$$
From the hard Lefschetz theorem applied to the projective morphism $\bp_{01}$ 
and to ${}'\tM\la|w_0|+|w_1|+\nu+2\r\ra$ (a perverse sheaf of pure weight $0$)
we have canonically for any $i$:
$${}'\tL^{\nu+2\r-i}={}'\tL^{\nu+2\r+i}(i).$$
Taking $i=a$ we obtain the desired result. This proves (a).

\subhead 2.26\endsubhead
Let $r\ge1$ and let $\ww=(w_1,w_2,\do,w_r)\in W^r$,
$\pmb\l=(\l_1,\l_2,\do,\l_r)\in\fs^r$ and let $\pmb\o=(\o_1,\o_2,\do,\o_r)$ be
such that $\o_i\in\k\i_q(w_i)$ for $i=1,\do,r$. 

(a) {\it Assume that $w_i\cdo\l_i\in\boc$ for some $i\in[1,r]$ and that
$$\un{(\LL_{\l_1}^{\o_1}\cir\LL_{\l_2}^{\o_2}\cir\do\cir
\LL_{\l_r}^{\o_r})^{\{(r-1)(a-\nu)\}}}\ne0.$$
Then $w_i\cdo\l_i\in\boc$ for all $i\in[1,r]$.}
\nl
Let $j=\nu+2r\r+(r-1)a$. By assumption we have 
$\un{gr_j((L^{\pmb\o,[1,r]}_{\pmb\l})^j)}\ne0$. Hence by 2.22(a) we have 
$$\un{gr_j(({}'L^{\pmb\o,[1,r]}_{\pmb\l})^{j-2(r-1)\r}(-(r-1)\r))}\ne0.$$
Thus there exists $w\cdo\l\in\boc$ such that $\LL_\l^{\dw}$ has nonzero
multiplicity in $$({}'L^{\pmb\o,[1,r]}_{\pmb\l})^{j-2(r-1)\r},$$ that is,
$\bV_{w\cdo\l,j-2(r-1)\r}\ne0$ (notation of 2.16). Using 2.16(c) we see that \lb
$N(w\cdo\l,-j+2\r+\nu)\ne0$ that is $N(w\cdo\l,-(r-1)a)\ne0$. Using now 
1.10(a) we see that $w_i\cdo\l_i\in\boc$ for all $i\in[1,r]$.

\head 3. Sheaves on the variety $Z$\endhead   
\subhead 3.1\endsubhead
Let 
$$Z=\{(B,B',gU_B);(B,B')\in\cb^2,g\in G,gBg\i=B'\}.$$
We define $\e:\tcb^2@>>>Z$ by $(x\UU,y\UU)\m(x\BB x\i,y\BB y\i,y\UU x\i)$. Now 
$\e$ identifies $Z$ with $\TT\bsl\tcb^2$ where $\TT$ acts on $\tcb^2$ by 
$t:(x\UU,y\UU)\m(xt\UU,yt\UU)$. Note that $Z$ inherits an $\FF_q$-structure 
from $\cb\T\cb\T G$.

\subhead 3.2\endsubhead
The $G\T\TT^2$-action on $\tcb^2$ (see 2.1) induces a $G\T\TT^2$-action on
$\TT\bsl\tcb^2$ (see 3.1) hence a $G\T\TT^2$-action on $Z$ in which the 
subgroup $\{(1,t_1,t_2)\in G\T\TT^2;t_1=t_2\}$ acts trivially. For $w\in W$ 
let $Z_w=\{(B,B',gU_B);(B,B')\in\co_w,g\in G,gBg\i=B'\}$; this is a single 
$G\T\TT^2$-orbit on $Z$ with closure 
$$\bZ_w=\{(B,B',gU_B);(B,B')\in\bco_w,g\in G,gBg\i=B'\}$$
and we have $Z=\sqc_{w\in W}Z_w$. Note that $Z_w=\e(\tco_w)$, 
$\tco_w=\e\i(Z_w)$, $\bZ_w=\e(\btco_w)$, $\btco_w=\e\i(\bZ_w)$.

Let $\o\in\k\i_q(w)$. We have a diagram $\TT@<j_\o<<\tco_w@>\e_w>>Z_w$ 
where $\e_w$ is the restriction of $\e$ and $j_\o$ is as in 2.1. 
 Let $\l\in\fs$ be such that 
$w(\l)=\l$. The $\TT$-action on $\tcb^2$ in 3.1 is compatible under $j_\o$ 
with the $\TT$-action $t:t'\m w\i(t\i)t't$ on $\TT$ and $L_\l$ is equivariant 
for this action (by 1.4(a) with $w$ replaced by $w\i$) hence $j_\o^*L_\l$ is 
$\TT$-equivariant so that there is a well defined local system $\cl_\l^\o$ of
rank $1$ on $Z_w$ such that $\e_w^*\cl_\l^\o=j_\o^*L_\l$. Note that the 
induced action of $\TT_n^2$ (which acts trivially on $Z_w$) on any stalk of 
$\cl_\l^\o$ is via the character $(t_1,t_2)\m\l\i(t_1)\l(t_2)$. Moreover,
$\cl_\l^\o$ is naturally pure of weight zero. We have 
$\e^*_w\cl_\l^\o=L_\l^\o$. 

We show the converse:

(a) {\it Let $\cl$ be an irreducible $G\T\TT^2$-equivariant local system on 
$Z_w$. Then $\cl$ is isomorphic to $\cl_\l^\o$ for a unique $\l\in\fs$ such 
that $w(\l)=\l$.}
\nl
The local system $\e^*_w\cl$ on $\tcb^2$ is irreducible and 
$G\T\TT^2$-equivariant hence, by 2.1, is isomorphic to $L_\l^\o$ for a well 
defined $\l\in\fs$. Now the restriction of $\e^*_w\cl$ to any fibre of $\e_w$
is the constant sheaf. On the other hand the restriction of $L_\l^\o$ to any
fibre of $\e_w$ is (under an identification with $\TT$) of the form 
$L_{w(\l\i)\l}$ which is trivial if and only if $w(l)=\l$. We see that we must
have $w(\l)=\l$. We have $\e^*_w\cl\cong\e^*_w\cl_\l^\o$ (both are $L_\l^\o$)
hence $\cl\cong\cl_\l^\o$. This proves (a).

\mpb

We define $\fh:Z@>>>Z$ by $(B,B',gU_B)\m(B',B,g\i U_{B'})$.  Note that 
$\fh\e=\e\ti\fh:\tcb^2@>>>Z$. For $L\in\cd_mZ$ we set $L^\da=\fh^*L$.

\subhead 3.3\endsubhead
Let $\overset\smile\to{W\fs}=\{w\cdo\l\in W\fs;w(\l)=\l\}$,
$\overset\smile\to\boc=\overset\smile\to{W\fs}\cap\boc$.
For $w\cdo\l\in\overset\smile\to{W\fs}$ and 
$\o\in\k\i_q(w)$ we shall view $\cl_\l^\o$ as a constructible sheaf on 
$Z$ which is $0$ on $Z-Z_w$. Let $\cl_\l^{\o\sha}$ be its extension to an 
intersection cohomology complex of $\bZ_w$, viewed as a complex on $Z$, equal 
to $0$ on $Z-\bZ_w$. Let
$$\Bbb L_\l^\o=\cl_\l^{\o\sha}\la|w|+\nu+\r\ra,$$
a simple perverse sheaf on $Z$. 
Note that $\cl_\l^\o$ (resp. $\Bbb L_\l^\o$) is noncanonically isomorphic to 
$\cl_\l^{\dw}$ (resp. $\Bbb L_\l^{\dw}$.)

We define $\ti\e:\cd(Z)@>>>\cd(\tcb^2)$ and $\ti\e:\cd_m(Z)@>>>\cd_m(\tcb^2)$ by
$$\ti\e(L)=\e^*(L)\la\r\ra.$$
From the definition we have 
$$\e^*\cl_\l^{\o\sha}=L_\l^{\o\sha},\qua\ti\e\Bbb L_\l^\o=\LL_\l^\o.$$
Note that $\cl_\l^{\o\sha},\Bbb L_\l^\o$ are naturally pure of weight zero.

Let $\cd^\spa Z$ be the subcategory of $\cd(Z)$ consisting of objects which 
are restrictions of objects in the $G\T\TT^2$-equivariant derived category. 
Let $\cd^\spa_mZ$ be the subcategory of $\cd_m(Z)$ consisting of objects which
are restrictions of objects in the mixed $G\T\TT^2$-equivariant derived 
category. Let $\cm^\spa Z$ (resp. $\cm^\spa_mZ$) be the subcategory of 
$\cd^\spa Z$ (resp. $\cd^\spa_mZ$) consisting of objects which are perverse 
sheaves. We define $\cd^{\preceq}Z,\cd^{\prec}Z$, $\cm^{\preceq}Z$,
$\cm^\prec Z$, $\cd_m^{\preceq}Z$,$\cd_m^{\prec}Z$, $\cm_m^{\preceq}Z$,
$\cm_m^\prec Z$, $\cc^\spa Z$, $\cc_0^\boc Z$ as in 2.14, by replacing (in 2.14)
$Y$ by $Z$ and $\LL_\l^{\dw}$ by $\Bbb L_\l^{\dw}$ (with $w\cdo\l$ required to
be in $\overset\smile\to{W\fs}$). For $M\in\cc_0^\spa Z$ let $\un{M}$ be the largest 
subobject of $M$ such that as an object of $\cc^\spa Z$, we have 
$\un{M}\in\cc^\boc Z$.

\mpb

From 3.2(a) we see that, if $M\in\cm^\spa Z$, then any composition factor of 
$M$ is of the form $\Bbb L_\l^{\dw}$ for some $w\cdo\l\in\overset\smile\to{W\fs}$. From the 
definitions we see that $M\m\ti\e M$ is a functor $\cd^\spa Z@>>>\cd^\spa\tcb^2$ 
and also $\cd_m^\spa Z@>>>\cd_m^\spa\tcb^2$; moreover, it is a 
functor $\cm^\spa Z@>>>\cm^\spa\tcb^2$ and also 
$\cm_m^\spa Z@>>>\cm_m^\spa\tcb^2$. From the definitions we see that for 
$M\in\cm^\spa Z$ we have

(a) {\it $M\in\cm^{\preceq}Z$ if and only if $\ti\e M\in\cm^\preceq\tcb^2$; 
we have $M\in\cm^{\prec}Z$ if and only if $\ti\e M\in\cm^\prec\tcb^2$.}
\nl
Note that if $X\in\cd(Z)$ and $j\in\ZZ$ then 
$$(\e^*X)^{j+\r}=\e^*(X^j)[\r].\tag b$$
Moreover, if $Y\in\cm_m(Z)$ and $j'\in\ZZ$ then 
$$gr_{j'}(\ti\e Y)=\ti\e(gr_{j'}Y).\tag c$$

\mpb

Let $\l\in\fs$, $w\in W'_\l$, $\o\in\k\i(w)$. From 2.2(a) we deduce
$$(\Bbb L_\l^\o)^\da=\Bbb L_{\l\i}^{\o\i}.\tag d$$

\subhead 3.4\endsubhead
Let $r,f$ be integers such that $0\le f\le r-3$. Let
$$\cy=\{((x_0\UU,x_1\UU,\do,x_r\UU),g)\in\tcb^{r+1}\T G;g\in x_{f+3}\UU x_f\i,
g\in x_{f+2}\BB x_{f+1}\i\}.$$
Define $\vt:\cy@>>>\tcb^{r+1}$ by
$$((x_0\UU,x_1\UU,\do,x_r\UU),g)\m(x_0\UU,x_1\UU,\do,x_r\UU).$$
For $y,y'\in W$ let
$$\tcb^{r+1}_{[y,y']}=\{(x_0\UU,x_1\UU,\do,x_r\UU)\in\tcb^{r+1};
x_f\i x_{f+1}\in G_y,x_{f+2}\i x_{f+3}\in G_{y'}\}.$$
We show:

(a) {\it Let $\x\in\tcb^{r+1}_{[y,y']}$. If $yy'\ne1$ then $\vt\i(\x)=\emp$. If
$yy'=1$ then $\vt\i(\x)\cong\kk^{\nu-|y|}$.}
\nl
We set $\x=(x_0\UU,x_1\UU,\do,x_r\UU)$. If $\vt\i(\x)\ne\emp$ then 
$x_f\i x_{f+1}\in G_y$, $x_{f+2}\i x_{f+3}\in G_{y'}$ and 
$(x_{f+3}\UU x_f\i)\cap(x_{f+2}\BB x_{f+1}\i)\ne\emp$. Hence for some 
$u\in\UU$, $b\in\BB$ we have 
$$ux_f\i x_{f+1}=x_{f+3}\i x_{f+2}b\in G_y\cap G_{y'{}\i}$$ 
so that $yy'=1$. If we assume that $yy'=1$, then $\vt\i(\x)$ can be identified with
$$\{g\in G;g\in x_{f+3}\UU x_f\i,g\in x_{f+2}\BB x_{f+1}\i\}$$
hence with 
$$\{(u,b)\in\UU\T\BB;ux_f\i x_{f+1}=x_{f+3}\i x_{f+2}b\}.$$
We substitute $x_{f+3}\i x_{f+2}=u_0\dy t_0u'_0$, 
$x_f\i x_{f+1}=u_1\dy t_1u'_1$ where $u_0,u'_0,u_1,u'_1\in\UU$, $t_0\in\TT$. 
Then $\vt\i(\x)$ is identified with $\{(u,b)\in\UU\T\BB;uu_1\dy t_1u'_1=u_0\dy t_0u'_0b\}$.
The map $(u,b)\m u_0\i uu_1$ identifies this variety with
$\UU\cap\dy\BB\dy\i\cong\kk^{\nu-|y|}$. This proves (a).

\mpb

Now $\TT^2$ acts freely on $\cy$ by
$$\align&(t_1,t_2):((x_0\UU,x_1\UU,\do,x_r\UU),g)\m\\&((x_0\UU,x_1\UU,\do,x_f\UU,
x_{f+1}t_1\UU,x_{f+2}t_2\UU,x_{f+3}\UU,\do,x_r\UU),g).\endalign$$
Let
$${}^!\cy=\TT\bsl\{((x_0\UU,x_1\UU,\do,x_r\UU),g)\in\tcb^{r+1}\T G;
g\in x_{f+3}\UU x_f\i,g\in x_{f+2}\UU x_{f+1}\i\}$$
where $\TT$ acts freely by
$$\align&t:((x_0\UU,x_1\UU,\do,x_r\UU),g)\m\\&((x_0\UU,x_1\UU,\do,x_f\UU,x_{f+1}t\UU,
x_{f+2}t\UU,x_{f+3}\UU,\do,x_r\UU),g).\endalign$$
Clearly, the obvious map $\b:{}^!\cy@>>>\TT^2\bsl\cy$ is an isomorphism. We 
define ${}^!\eta:{}^!\cy@>>>Z$ by
$$((x_0\UU,x_1\UU,\do,x_r\UU),g)\m\e(x_{f+1}\UU,x_{f+2}\UU).$$
We define $\t:\cy@>>>{}^!\cy$ as the composition of the obvious map
$\cy@>>>\TT^2\bsl\cy$ with $\b\i$. Let $\eta={}^!\eta\t:\cy@>>>Z$. We have
$$\eta((x_0\UU,x_1\UU,\do,x_r\UU),g)=\e(x_{f+1}t\i\UU,x_{f+2}t'{}\i\UU)$$
where $t,t'$ in $\TT$ are such that $g\in x_{f+2}t'{}\i\UU tx_{f+1}\i$. 

\subhead 3.5\endsubhead
Let $z\cdo\l\in\overset\smile\to{W\fs}$. Let $P=\eta^*\cl_\l^{\dz\sha}$. Let 
$p_{ij}:\tcb^{r+1}@>>>\tcb^2$ be the projection to the $ij$ coordinates. We 
have the following result:
$$\vt_!P\Bpq\{p_{f,f+1}^*L_\l^{\dy}\ot p_{f+1,f+2}^*L_\l^{\dz\sha}\ot
p_{f+2,f+3}^*L_{y(\l)}^{\dy\i}\la2|y|-2\nu\ra;y\in W\}.\tag a$$
Define $e:\tcb^{r+1}@>>>\tcb^4$ by
$$(x_0\UU,x_1\UU,\do,x_r\UU)\m(x_f\UU,x_{f+1}\UU,x_{f+2}\UU,x_{f+3}\UU).$$
Then (a) is obtained by applying $e^*$ to the statement similar to (a) in
which $\{0,1,\do,r\}$ is replaced by $\{f,f+1,f+2,f+3\}$. Thus it is enough to
prove (a) in the special case where $r=3,f=0$. In the remainder of the proof
we assume that $r=3,f=0$.

For any $y,y'$ in $W$ let $\vt_{y,y'}:\vt\i(\tcb^4_{[y,y']})@>>>\tcb^4$ be the
restriction of $\vt$. Let $P^{y,y'}$ be the restriction of $P$ to
$\vt\i(\tcb^4)_{[y,y']}$. Clearly, we have
$$\vt_!P\Bpq\{\vt_{y,y'!}P^{y,y'};(y,y')\in W^2\}.$$
Since $\vt\i(\tcb^{r+1}_{[y,y']})=\emp$ when $yy'\ne1$, see 3.4(a), we deduce 
that
$$\vt_!P\Bpq\{\vt_{y,y\i!}P^{y,y\i};y\in W\}.$$
Hence to prove (a) it is enough to show for any $y\in W$ that
$$\vt_{y!}P_y=p_{01}^*L_\l^{\dy}\ot p_{12}^*L_\l^{\dz\sha}\ot 
p_{23}^*L_{y(\l)}^{\dy\i}\la2|y|-2\nu\ra\tag b$$
where we write $\vt_y,P_y$ instead of $\vt_{y,y\i},P^{y,y\i}$. We have a 
cartesian diagram
$$\CD
\tV_y@>\tb>>\ti\cv_y\\
@VVV     @VVV\\
V_y@>b>>\cv_y
\endCD$$
where 
$$V_y=\{(x_0\UU,x_1\UU,x_2\UU,x_3\UU)\in\tcb^4;x_0\i x_1\in G_y,
x_1\i x_2\in G_z,x_2\i x_3\in G_{y\i}\},$$
$$\align&\cv_y=\TT\bsl\{(x_0\UU,x_1\UU,x_2\UU,x_3\UU)\in\tcb^4;x_0\i x_1\in G_y,
x_1\i x_2\in G_z,x_2\i x_3\in G_{y\i},\\& (x_0\i x_1)_{\dy}=(x_3\i x_2)_{\dy}\}\endalign$$
with $\TT$ acting freely by simultaneous right multiplication on $x_1$ and 
$x_2$, $\tV_y=\vt\i(V_y)$ and
$$\align&\ti\cv_y=\TT\bsl\{((x_0\UU,x_1\UU,x_2\UU,x_3\UU),g)\in\tcb^4\T G;x_0\i x_1\in G_y,
x_1\i x_2\in G_z, \\& x_2\i x_3\in G_{y\i},g\in x_3\UU x_0\i,g\in x_2\UU x_1\i\};\endalign$$
we have 
$$b(x_0\UU,x_1\UU,x_2\UU,x_3\UU)=\TT-\text{orbit of }
(x_0\UU,x_1t\UU,x_2t'\UU,x_3\UU)$$
where $t,t'$ in $\TT$ are such that $(x_0\i x_1t)_{\dy}=(x_3\i x_2t')_{\dy}$,
$$\tb((x_0\UU,x_1\UU,x_2\UU,x_3\UU),g)=\TT-\text{orbit of }
((x_0\UU,x_1t\UU,x_2t'\UU,x_3\UU),g)$$
where $t,t'$ in $\TT$ are such that $g\in x_2t'\UU t\i x_1\i$; the vertical 
maps are the obvious ones. We also have a cartesian diagram
$$\CD
\tV'_y@>\tb'>>\ti\cv'_y\\
@VVV     @VVV\\
V'_y@>b'>>\cv'_y
\endCD$$
where $\tV'_y,\ti\cv'_y,V'_y,\cv'_y$ are defined in the same way as
$\tV_y,\ti\cv_y,V_y,\cv_y$ but the condition $x_1\i x_2\in G_z$ is replaced by
the condition $x_1\i x_2\in\bG_z$; the maps $\tb',b'$ are given by the same 
formulas as $\tb,b$; the vertical maps are the obvious ones.

Let $j:V'_y@>>>\tcb^4$ be the inclusion. It is enough to show that
$$j^*\vt_{y!}P_y=j^*(p_{01}^*L_\l^{\dy}\ot p_{12}^*L_\l^{\dz\sha}\ot 
p_{23}^*L_{y(\l)}^{\dy\i})\la2|y|-2\nu\ra.$$
By definition, $P|_{\tV'_y}$ is the inverse image of $\cl_\l^{\dz\sha}$ under 
the composition of $\tb'$ with $\ti\cv'_y@>>>\cv'_y@>{}^!\eta_y>>Z$ where the 
first map is the obvious one and
$${}^!\eta_y(x_0\UU,x_1\UU,x_2\UU,x_3\UU)=\e(x_1\UU,x_2\UU).$$
Hence $P|_{\tV'_y}$ is the inverse image of $\cl_\l^{\dz\sha}$ under the 
composition of $\eta_y:={}^!\eta_yb'$ with the obvious map 
$\vt'_y:\tV'_y@>>>V'_y$. Since $\vt_y$ is an affine space bundle with fibres 
of dimension $\nu-|y|$, it follows that 
$j^*\vt_{y!}P_y=\eta_y^*\cl_\l^{\dz\sha}\la2|y|-2\nu\ra$. Thus it is enough to show that 
$$\eta_y^*\cl_\l^{\dz\sha}=j^*(p_{01}^*L_\l^{\dy}\ot p_{12}^*L_\l^{\dz\sha}\ot 
p_{23}^*L_{y(\l)}^{\dy\i}).$$
Since $\eta_y$ is smooth as a map to $\bZ_z$, we see that
$\eta_y^*\cl_\l^{\dz\sha}$ is the intersection cohomology complex of $V'_y$ 
with coefficients in the local system $\eta_y^*\cl_\l^{\dz}$ on $V_y$.
Now, 
$$j^*(p_{01}^*L_\l^{\dy}\ot p_{12}^*L_\l^{\dz\sha}\ot 
p_{23}^*L_{y(\l)}^{\dy\i})$$
is the intersection cohomology complex of $V'_y$ with coefficients in the 
local system 
$$j^*(p_{01}^*L_\l^{\dy}\ot p_{12}^*L_\l^{\dz}\ot p_{23}^*L_{y(\l)}^{\dy\i})$$
on $V_y$. It is then enough to show that these two local systems on $V_y$ are 
the same.
One local system is $h^*L_\l$ with $h:V_y@>>>\TT$ given by
$$(x_0\UU,x_1\UU,x_2\UU,x_3\UU)\m(t_1x_1\i x_2t_2\i)_{\dz}$$
where 
$$t_1=(x_0\i x_1)_{\dy}\in\TT,t_2=(x_3\i x_2)_{\dy}\in\TT.$$
The other local system is $h'{}^*L_\l$ with $h':V_y@>>>\TT$ given by
$$(x_0\UU,x_1\UU,x_2\UU,x_3\UU)\m t'_1z(t'_2)(zy\i)(t'_3)$$
where 
$$t'_1=(x_0\i x_1)_{\dy}\in\TT, t'_2=(x_1\i x_2)_{\dz}\in\TT,
t'_3=(x_2\i x_3)_{\dy\i}\in\TT.$$
It is enough to show that $h^*L_\l=h'{}^*L_\l$. Since the map $u:\TT@>>>\TT$, 
$t\m z\i(t)$ satisfies $u^*L_\l=L_\l$ (recall that $z(\l)=\l$), we have 
$h'{}^*L_\l=h'{}^*u^*L_\l$ hence it is enough to show that 
$h(\x)=\z\i(h'(\x))$ for any $\x=(x_0\UU,x_1\UU,x_2\UU,x_3\UU)\in V$, 
or that, if $t_1,t_2,t'_1,t'_2,t'_3$ are associated to $\x$ as above, then
$$(t_1x_1\i x_2t_2\i)_{\dz}=z\i(t'_1z(t'_2)(zy\i)(t'_3)).$$
We have $t_1=t'_1$ and $x_3\i x_2\in\UU\dy t_2\UU$ hence 
$$x_2\i x_3\in\UU t_2\i\dy\i\UU=\UU\dy\i y(t_2\i)\UU$$ 
so that $t'_3=y(t_2\i)$ and $t_2\i=y\i(t'_3)$. We have
$$t_1x_1\i x_2t_2\i\in t_1\UU\dz t'_2\UU t_2\i=\UU\dz z\i(t_1)t'_2t_2\i\UU,$$
so that
$$(t_1x_1\i x_2t_2\i)_{\dz}=z\i(t_1)t'_2t_2\i=z\i(t'_1)t'_2y\i(t'_3),$$
as required. This completes the proof of (b) hence that of (a).

\subhead 3.6\endsubhead
Let 
$$(w_1,w_2,\do,w_f,w_{f+2},w_{f+4},\do,w_r)\in W^{r-2},$$
$$(\l_1,\l_2,\do,\l_f,\l_{f+2},\l_{f+4},\do,\l_r)\in\fs^{r-2}.$$ 
We set $z=w_{f+2},\l=\l_{f+2}$. We assume that $z(\l)=\l$.
Let $P$ be as in 3.5. Let 
$P'=\ot_{i\in[1,r]-\{f+1,f+2,f+3\}}p_{i-1,i}^*L_{\l_i}^{\dw_i\sha}\in
\cd_m(\tcb^{r+1})$, $\tP=P\ot\vt^*P'\in\cd_m(\cy)$. For any $y\in W$ we set
$$\ww_y=(w_1,w_2,\do,w_f,y,w_{f+2},y\i,w_{f+4},\do,w_r)\in W^r,$$
$$\pmb\o_y=(\dw_1,\dw_2,\do,\dw_f,\dy,\dw_{f+2},\dy\i,\dw_{f+4},\do,\dw_r),$$
$$\pmb\l_y=(\l_1,\l_2,\do,\l_f,\l_{f+2},\l_{f+2},y(\l_{f+2}),\l_{f+4},\do,
\l_r)\in\fs^r.$$
We set $\Xi=\vt_!\tP$. We have the following result.
$$\Xi\Bpq\{M^{\pmb\o_y,[1,r]-\{f+1,f+3\}}_{\pmb\l_y}\la2|y|-2\nu\ra;
y\in W\}\tag a$$
in $\cd_m(\tcb^{r+1})$. This follows immediately from 3.5(a) since 
$\Xi=P'\ot\vt_!(P)$.

\subhead 3.7\endsubhead
We preserve the setup of 3.6.
Let $\cs=\sqc_{\ww'}\tco_{\ww'}^\emp$ where the union is over all
$\ww'=(w'_1,\do,w'_r)\in W^r$ such that $w'_i=w_i$ for $i\n\{f+1,f+3\}$. This 
is a locally closed subvariety of $\tcb^{r+1}$. 
For $y\in W$ let $R_y$ be the restriction of 
$M^{\pmb\o_y,\emp}_{\pmb\l_y}$ to $\tco^\emp_{\ww_y}$ extended by $0$ on 
$\cs-\tco^\emp_{\ww_y}$ (a constructible sheaf on $\cs$). From the definitions 
we have
$$M^{\pmb\o_y,[1,r]-\{f+1,f+3\}}_{\pmb\l_y}|_\cs=R_y.$$
From 3.6(a) we deduce
$$\Xi|_{\cs}\Bpq\{R_y\la2|y|-2\nu\ra;y\in W\}.$$
We now restrict further to $\tco^\emp_{\ww_y}$ (for $y\in W$); we obtain
$$\Xi|_{\tco^\emp_{\ww_y}}
\Bpq\{R_{y'}\la2|y'|-2\nu\ra|_{\tco^\emp_{\ww_y}};y'\in W\}.$$
In the right hand side we have $R_{y'}\la2|y'|-2\nu\ra|_{\tco^\emp_{\ww_y}}=0$ 
if $y'\ne y$. It follows that $\Xi|_{\tco^\emp_{\ww_\l}}=R_y\la2|y|-2\nu\ra|_{\tco^\emp_{\ww_y}}$.
Since $R_y|_{\tco^\emp_{\ww_y}}$ is a local system we deduce for $y\in W$ the following result.

(a) {\it Let $h\in\ZZ$. If $h=2\nu-2|y|$ then
$\ch^h\Xi|_{\tco^\emp_{\ww_y}}=R_y|_{\tco^\emp_{\ww_y}}(|y|-\nu)$. If 
$h\ne2\nu-2|y|$ then $\ch^h\Xi|_{\tco^\emp_{\ww_y}}=0$.}  

\subhead 3.8\endsubhead
We preserve the setup of 3.6. We set
$$k=(\sum_{i\in[1,r]-\{f+1,f+3\}}|w_i|)+3\nu+(r+1)\r.\tag a$$ 
For $y\in W$ we set
$$K_y=M^{\pmb\o_y,[1,r]-\{f+1,f+3\}}_{\pmb\l_y}\la|\ww_y|+\nu+(r+1)\r\ra,$$
$$\tK_y=M^{\pmb\o_y,[1,r]}_{\pmb\l_y}\la|\ww_y|+\nu+(r+1)\r\ra.$$
From 3.6(a) we deduce:
$$\Xi\la k\ra\Bpq\{K_y;y\in W\}.\tag b$$
We show:

(c) {\it For any $j>0$ we have $(\Xi\la k\ra)^j=0$. Equivalently, $\Xi^j=0$ 
for any $j>k$.}
\nl
Using (b) we see that it is enough to show that for any $y\in W$ we have
$(K_y)^j=0$ for any $j>0$. Now $\tK_y$ is a (simple) perverse sheaf hence for 
any $i$ we have $\dim\supp\ch^i\tK_y\le-i$. Moreover $K_y$ is obtained by 
restricting $\tK_y$ to an open subset of its support and then extending the 
result (by zero) on the complement of this subset in $\tcb^{r+1}$. Hence 
$\supp\ch^iK_y\sub\supp\ch^i\tK_y$ so that $\dim\supp\ch^iK_y\le-i$. Since
this holds for any $i$ we see that $(K_y)^j=0$ for any $j>0$. 

\subhead 3.9\endsubhead
We preserve the notation of 3.6. We show:

(a) {\it Let $j\in\ZZ$ and let $X$ be a composition factor of $\Xi^j$. Then
$X\cong M^{\pmb\o',[1,r]}_{\pmb\l'}\la|\ww'|+\nu+(r+1)\r\ra$ for some 
$$\ww'=(w'_1,w'_2,\do,w'_r)\in W^r,\pmb\l'=(\l'_1,\l'_2,\do,\l'_r)\in\fs^r$$
such that $w'_i=w_i$, $\l'_i=\l_i$ for $i\in[1,r]-\{f+1,f+3\}$ and such that 
$$\l'_{f+1}=w'_{f+2}(\l'_{f+2}), \l'_{f+2}=w'_{f+3}(\l'_{f+3}).$$
Here $\pmb\o'=(\dw'_1,\dw'_2,\do,\dw'_r)$.}
\nl
From 3.6(a) we see that, for some $y\in W$, $X$ is a composition factor of 
$$(M^{\pmb\o_y,[1,r]-\{f+1,f+3\}}_{\pmb\l_y}\la2|y|-2\nu\ra)^j.$$
Using this and 2.18(b) we see that $X$ is as required except that the equalities above 
for $\l'_{f+1},\l'_{f+2}$ may not be 
satisfied. To see that they are in fact satisfied we note that 
$$(M^{\pmb\o_y,[1,r]-\{f+1,f+3\}}_{\pmb\l_y}\la2|y|-2\nu\ra)^j$$ 
is equivariant for the $\TT^2$-action
$$\align&(t_1,t_2):(x_0\UU,x_1\UU,\do,x_r\UU)\m  \\&
(x_0\UU,x_1\UU,\do,x_f\UU,x_{f+1}t_1\UU,x_{f+2}t_2\UU,x_{f+3}\UU,\do,x_r\UU)\endalign$$
hence so are its composition factors and this implies that the equalities
above for $\l'_{f+1},\l'_{f+2}$ do hold.

\subhead 3.10\endsubhead
From 3.8(c) we see that we have a distinguished triangle
$(\Xi',\Xi,\Xi^k[-k])$ where $\Xi'\in\cd_m(\tcb^{r+1})$ satisfies 
$(\Xi')^j=0$ for all $j\ge k$. We show:

(a) {\it Let $j\in\ZZ$ and let $K$ be one of $\Xi,\Xi^j,\Xi'$. For any 
$\ww'\in W^r$ and any $h\in\ZZ$, $\ch^hK|_{\tco^\emp_{\ww'}}$ is a local 
system.}
\nl
We prove (a) for $K=\Xi$. Using 3.6(a), we see that it is enough to show that 
$\ch^h(M^{\pmb\o_y,[1,r]-\{f+1,f+3\}}_{\pmb\l_y})|_{\tco^\emp_{\ww'}}$ is a 
local system for any $h$. This follows from the fact that
$\th^*M^{\pmb\o_y,[1,r]-\{f+1,f+3\}}_{\pmb\l_y}$ (see 2.15) is 
$\cg$-equivariant.

We prove (a) for $K=\Xi^j$. Using 3.6(a), we see that it is enough to show
that $\ch^h((M^{\pmb\o_y,[1,r]-\{f+1,f+3\}}_{\pmb\l_y})^j)|_{\tco^\emp_{\ww'}}$ 
is a local system for any $h$ and any $j$. This again follows from the
$\cg$-equivariance statement in the previous paragraph.

Now (a) for $K=\Xi'$ follows from (a) for $\Xi$ and $\Xi^k[-k]$ using the long
exact sequence for cohomology sheaves of $(\Xi',\Xi,\Xi^k[-k])$ restricted to
$\tco^\emp_{\ww'}$.

\mpb

(b) {\it Let $(y,y')\in W^2$, $i=2\nu-|y|-|y'|$. Let 
$$\ww_{y,y'}=(w_1,w_2,\do,w_f,y,w_{f+2},y',w_{f+3},\do,w_r)\in W^r.$$
The induced homomorphism $\ch^i\Xi|_{\tco^\emp_{\ww_{y,y'}}}@>>>
\ch^{i-k}(\Xi^k)|_{\tco^\emp_{\ww_{y,y'}}}$ is an isomorphism.}
\nl
We have an exact sequence of constructible sheaves
$$\ch^i\Xi'|_{\tco^\emp_{\ww_{y,y'}}}@>>>\ch^i\Xi|_{\tco^\emp_{\ww_{y,y'}}}@>>>
\ch^{i-k}(\Xi^k)|_{\tco^\emp_{\ww_{y,y'}}}@>>>\ch^{i+1}\Xi'|_{\tco^\emp_{\ww_{y,y'}}}.$$
Hence it is enough to show that 
$\ch^{i'}\Xi'|_{\tco^\emp_{\ww_{y,y'}}}=0$ if $i'\ge i$. Assume that 
$$\ch^{i'}\Xi'|_{\tco^\emp_{\ww_{y,y'}}}\ne0$$ for some $i'\ge i$.
Since $\ch^{i'}\Xi'|_{\tco^\emp_{\ww_{y,y'}}}$ is a local system (see (a)), 
we deduce that $\tco^\emp_{\ww_{y,y'}}$ is contained in $\supp(\ch^{i'}\Xi')$.
We have $(\Xi'[k-1])^j=0$ for any $j>0$ hence 
$\dim\supp(\ch^{i''}\Xi'[k-1])\le-i''$ for any $i''$. Taking $i''=i'-k+1$ we 
deduce that 
$$\dim\tco^\emp_{\ww_{y,y'}}\le\dim\supp(\ch^{i'}\Xi')\le-i'+k-1\le-i+k-1$$
hence
$$|\ww_{y,y'}|+\nu+(r+1)\r\le-i+k-1.$$
We have $|\ww_{y,y'}|+\nu+(r+1)\r=-i+k$ hence $-i+k\le-i+k-1$, contradiction.
This proves (b).

\subhead 3.11\endsubhead
For $(y,y')\in W^2$ we set
$$\pmb\o_{y,y'}=(\dw_1,\dw_2,\do,
\dw_f,\dy,\dw_{f+2},((y'{}\i)\dot{})\i,\dw_{f+3},\do,\dw_r)\in W^r,$$ 
$$\pmb\l_{y,y'}=(\l_1,\l_2,\do,\l_f,\l_{f+2},\l_{f+2},y'{}\i(\l_{f+2}),
\l_{f+4},\do,\l_r)\in\fs^r,$$
$$K_{y,y'}=M^{\pmb\o_{y,y'},\emp}_{\pmb\l_{y,y'}}
\la|\ww_{y,y'}|+\nu+(r+1)\r\ra\in\cm_m(\tcb^{r+1}),$$
$$\tK_{y,y'}=M^{\pmb\o_{y,y'},[1,r]}_{\pmb\l_{y,y'}}
\la|\ww_{y,y'}|+\nu+(r+1)\r\ra\in\cm_m(\tcb^{r+1}).$$
Note that when $y=y'$, $\ww_{y,y'},\pmb\o_{y,y'},\pmb\l_{y,y'},K_{y,y'},
\tK_{y,y'}$ become $\ww_y,\pmb\o_y,\pmb\l_y$ (see 3.6) and $K_y,\tK_y$ (see 
3.8). We show that we have canonically
$$gr_0(\Xi^k(k/2))=\op_{y\in W}\tK_y.\tag a$$
Since $gr_0(\Xi^k(k/2))$ is a semisimple perverse sheaf of pure weight zero,
it is a direct sum of simple perverse sheaves, necessarily of the form 
described in 3.9(a). Thus we have canonically
$$gr_0(\Xi^k(k/2))=\op_{(y,y')\in W^2}V_{y,y'}\ot\tK_{y,y'}$$
where $V_{y,y'}$ are mixed $\bbq$-vector spaces of pure weight $0$. By
\cite{\BBD, 5.1.14}, $\Xi$ is mixed of weight $\le0$ hence $\Xi^k(k/2)$ is
mixed of weight $\le0$. Hence we have an exact sequence in $\cm_m(\tcb^{r+1})$:
$$0@>>>\cw\i(\Xi^k(k/2))@>>>\Xi^k(k/2)@>>>gr_0(\Xi^k(k/2))@>>>0\tag a$$
that is,
$$0@>>>\cw\i(\Xi^k(k/2))@>>>\Xi^k(k/2)@>>>
\op_{(y,y')\in W^2}V_{y,y'}\ot\tK_{y,y'}@>>>0.$$
(Here $\cw\i(?)$ denotes the part of weight $\le-1$ of a mixed perverse sheaf.)
Hence for any $(\ty,\ty')\in W^2$ we have an exact sequence of (mixed) cohomology
sheaves restricted to $\tco^\emp_{\ww_{\ty,\ty'}}$ (where 
$h=2\nu-|\ty|-|\ty'|-k$): 
$$\align&\ch^h(\cw\i(\Xi^k(k/2)))@>\a>>\ch^h(\Xi^k(k/2))@>>>
\op_{(y,y')\in W^2}V_{y,y'}\ot\ch^h(\tK_{y,y'})@>>>\\&
\ch^{h+1}(\cw\i(\Xi^k(k/2))).\tag b\endalign$$
Moreover, by 3.10(b), we have
an equality of local systems on $\tco^{\emp}_{\ww_{\ty,\ty'}}$:
$$\ch^h(\Xi^k(k/2))=\ch^{h+k}(\Xi(k/2))=\ch^{2\nu-|y|-|y'|}(\Xi(k/2))$$
and this is $R_{\ty}(k/2+|\ty|-\nu)$ if $\ty\ty'=1$ (see 3.7(a)) and is $0$
if $\ty\ty'\ne1$ (see 3.4(a)) hence is pure of weight $-k-|\ty|-|\ty'|+\nu=h$.
On the other hand, $\ch^h(\cw\i(\Xi^k(k/2)))$ is mixed of weight $\le h-1$; it
follows that $\a$ in (b) must be zero.

Assume that $\ch^h(\tK_{y,y'})$ is nonzero on $\tco^\emp_{\ww_{\ty,\ty'}}$. 
Then, by 3.10(a), $\tco^\emp_{\ww_{\ty,\ty'}}$ is contained in 
$\supp\ch^h(\tK_{y,y'})$ which has dimension $\le-h$ (resp. $<-h$ if 
$(y,y')\ne(\ty,\ty')$); hence $-h=\dim \tco^\emp_{\ww_{\ty,\ty'}}$ is 
$\le-h$ (resp. $<-h$); we see that we must have $(y,y')=(\ty,\ty')$ and we 
have $\ch^h(\tK_{y,y'})=\ch^h(K_{y,y'})$ on $\tco^\emp_{\ww_{\ty,\ty'}}$.

Assume that $\ch^{h+1}(\cw\i(\Xi^k(k/2)))$ is not identically $0$ on 
$\tco^\emp_{\ww_{\ty,\ty'}}$. Then, by 3.10(a), $\tco^\emp_{\ww_{\ty,\ty'}}$ 
is contained in $\supp\ch^{h+1}(\cw\i(\Xi^k(k/2)))$ which has dimension 
$\le-h-1$; hence $-h=\dim\tco^\emp_{\ww_{\ty,\ty'}}\le-h-1$, a contradiction. 
We see that (b) becomes an isomorphism of local systems on 
$\tco^\emp_{\ww_{\ty,\ty'}}$:
$$0=V_{\ty,\ty'}\ot K_{\ty,\ty'} \text{ if }\ty\ty'\ne1,$$
$$R_{\ty}(-h/2)@>\si>>V_{\ty,\ty'}\ot\ch^h(K_{\ty,\ty'})\text{ if }\ty\ty'=1.$$
When $\ty\ty'=1$ we have $\ch^h(K_{\ty,\ty'})=R_{\ty}(-h/2)$ as local systems
on $\tco^\emp_{\ww_{\ty,\ty'}}$. It follows that $V_{\ty,\ty'}$ is $\bbq$ if
$\ty\ty'=1$ and is $0$ if $\ty\ty'\ne1$. This proves (a).

\subhead 3.12\endsubhead
Let $h\in[1,r]$.
Let ${}_h\cd^{\preceq}\tcb^{r+1}$ (resp. ${}_h\cd^{\prec}\tcb^{r+1}$) be the 
subcategory of $\cd\tcb^{r+1}$ consisting of objects $K$ such that for any 
$j\in\ZZ$, any composition factor of $K^j$ is of the form 
$M^{\pmb\o,[1,r]}_{\pmb\l}\la|\ww|+\nu+(r+1)\r\ra$ for some 
$\ww=(w_1,\do,w_r)\in W^r$, $\pmb\l=(\l_1,\l_2,\do,\l_r)\in\fs^r$ such that 
$w_h\cdo\l_h\preceq\boc$ (resp. $w_h\cdo\l_h\prec\boc$). (Here
$\pmb\o=(\dw_1,\dw_2,\do,\dw_r)$.)
\nl
Let ${}_h\cm^{\preceq}\tcb^{r+1}$ (resp. ${}_h\cm^{\prec}\tcb^{r+1}$) be the 
subcategory of ${}_h\cd^{\preceq}\tcb^{r+1}$ (resp. 
${}_h\cd^{\prec}\tcb^{r+1}$) consisting of perverse sheaves.

If $K\in\cm_m(\tcb^{r+1})$ is pure of weight $0$ and is also in
${}_h\cd^{\preceq}\tcb^{r+1}$ we denote by $\un K$ the sum of all simple 
subobjects
of $K$ (without mixed structure) which are not in 
${}_h\cd^{\prec}\tcb^{r+1}$.

\subhead 3.13\endsubhead
Let $Z@<\eta<<\cy@>\vt>>\tcb^4$ be as in 3.4 with $r=3,f=0$. We define 
$\fb:\cd(Z)@>>>\cd(\tcb^2)$ and $\fb:\cd_m(Z)@>>>\cd_m(\tcb^2)$ by 
$$\fb(L)=p_{03!}\vt_!\eta^*L.$$
We show:

(a) {\it If $L\in\cd^{\preceq}(Z)$ then $\fb(L)\in\cd^{\preceq}\tcb^2$.}

(b) {\it If $L\in\cd^{\prec}(Z)$ then $\fb(L)\in\cd^{\prec}\tcb^2$.}

(c) {\it If $L\in\cm^{\preceq}(Z)$ and $h>5\r+2\nu+2a$ then 
$(\fb(L))^h\in\cm^{\prec}\tcb^2$.}
\nl
We can assume that $L=\Bbb L_\l^{\dz}$ where $z\cdo\l\in\overset\smile\to{W\fs}$, 
$z\cdo\l\preceq\boc$. Applying 3.5(a) with $P=\eta^*\cl_\l^{\dz\sha}$ we see 
that
$$\fb(\cl_\l^{\dz\sha})\Bpq\{L^{\dy,\dz,\dy\i,\{2\}}_{\l,\l,y(\l)}
\la-|z|-2\nu\ra;y\in W\},$$
hence
$$\fb(\Bbb L_\l^{\dz\sha})\Bpq\{L^{\dy,\dz,\dy\i,\{2\}}_{\l,\l,y(\l)}
\la-\nu+\r\ra;y\in W\}.$$
To prove (a) it is enough to show that for any $y\in W$ we have
$$L^{\dy,\dz,\dy\i,\{2\}}_{\l,\l,y(\l)}\in\cd^{\preceq}\tcb^2.$$
When $z\cdo\l\in\boc$ this follows from 2.10(a). When $z\cdo\l\prec\boc$ this 
again follows from 2.10(a), applied to the two-sided cell containing $z\cdo\l$
instead of $\boc$. The same argument proves (b).
To prove (c) we can assume that $z\cdo\l\in\boc$; it is 
enough to prove that for any $y\in W$ we have 
$$(L^{\dy,\dz,\dy\i,\{2\}}_{\l,\l,y(\l)}\la-\nu+\r\ra)^h\in\cm^{\prec}\tcb^2$$
if $h>5\r+2\nu+2a$ or that $(L^{\dy,\dz,\dy\i,\{2\}}_{\l,\l,y(\l)})^j\in\cm^{\prec}\tcb^2$
if $j>6\r+\nu+2a$. This follows from 2.20(a). This completes the proof of (a),
(b),(c).

We define $\un\fb:\cc^\boc_0(Z)@>>>\cc^\boc_0(\tcb^2)$ by
$$\un\fb(L)=\un{gr_{5\r+2\nu+2a}((\fb(L))^{5\r+2\nu+2a})}((5\r+2\nu+2a)/2).
$$ 
We show:

(d) {\it Let $z\cdo\l\in\overset\smile\to\boc$. We have
$\un\fb(\Bbb L_\l^{\dz})=\op_{y\in W;y\cdo\l\in\boc}
\LL_\l^{\dy}\un{\cir}\LL_\l^{\dz}\un{\cir}\LL_{y(\l)}^{\dy\i}.$}
\nl
We shall apply \cite{\CONV, 1.12} with 
$\Ph:\cd_m(Y_1)@>>>\cd_m(Y_2)$ replaced by 
$p_{03!}:\cd_m(\tcb^4)@>>>\cd_m(\tcb^2)$ and with $\cd^{\preceq}(Y_1)$, 
$\cd^{\preceq}(Y_2)$ replaced by ${}_2\cd^{\preceq}(\tcb^2)$, 
${}_2\cd^{\preceq}(\tcb^4)$, see 3.12. We shall take $\XX$ in {\it loc.cit.} 
equal to $\vt_!\eta^*\Bbb L_\l^{\dz}$.
The conditions of {\it loc.cit.} are satisfied: those concerning $\XX$ are
satisfied with $c'=2\nu+3\r$. (For $h>|z|+3\nu+4\r$ we have $\Xi^h=0$ that is
$(\XX[-|z|-\nu-\r])^h=0$, with $\Xi$ as in 3.8(c). Hence if $j>2\nu+3\r$ we 
have $\XX^j=0$.) The conditions concerning $p_{03!}$ are satisfied with
$c=2\r+2a$. (This follows from 2.20(a).) Since 
$\fb(\Bbb L_\l^{\dz})=p_{03!}\XX$ and $c+c'=5\r+2\nu+2a$, we see that
$$\un\fb(\Bbb L_\l^{\dz})=\un{gr_{2\r+2a}(p_{03!}((\un{gr_{2\nu+3\r}
((\vt_!\eta^*\Bbb L_\l^{\dz})^{2\nu+3\r})}((2\nu+3\r)/2)))^{2\r+2a})}(\r+a).$$
Using 3.11(a) we see that (with $\Xi$ as in 3.11(a) and $k=|z|+3\nu+4\r$) we 
have
$$\align&\un{gr_{2\nu+3\r}((\vt_!\eta^*\Bbb L_\l^{\dz})^{2\nu+3\r})}((2\nu+3\r)/2)\\&=
\un{gr_{2\nu+3\r}((\Xi\la|z|+\nu+\r\ra)^{2\nu+3\r})}((2\nu+3\r)/2)\\&=
\un{gr_0(\Xi^k(k/2)}
=\op_{y\in W}M^{\dy,\dz,\dy\i,[1,3]}_{\l,\l,y(\l)}\la2|y|+|z|+\nu+4\r\ra.\endalign$$
Hence
$$\align&\un\fb(\Bbb L_\l^{\dz})=\un{gr_{2\r+2a}(\op_{y\in W}(p_{03!}
M^{\dy,\dz,\dy\i,[1,3]}_{\l,\l,y(\l)}\la2|y|+|z|+\nu+4\r\ra)^{2\r+2a})}(\r+a)\\&
=\un{gr_{2\r+2a}(\op_{y\in W}
(L^{\dy,\dz,\dy\i,[1,3]}_{\l,\l,y(\l)})^{6\r+\nu+2a}((\nu+4\r)/2))}(\r+a).\endalign$$
Using 2.26(a) we see that in the last direct sum the contribution of $y\in W$
is $0$ unless $y\cdo\l\in\boc$. For the terms corresponding to $y$ such that
$y\cdo\l\in\boc$, we may apply 2.24(a). Now (d) follows.

\subhead 3.14\endsubhead
Let $Z@<{}^!\eta<<{}^!\cy$ be as in 3.4 with $r=3,f=0$. 
Let ${}^!\tcb^4$ be the space of 
orbits of the free $\TT^2$-action on $\tcb^4$ given by 
$$(t_1,t_2):(x_0\UU,x_1\UU,x_2\UU,x_3\UU)\m(x_0\UU,x_1t_1\UU,x_2t_2\UU,x_3\UU);
$$
let ${}^!\vt:{}^!\cy@>>>{}^!\tcb^4$ be the map induced by $\vt$. We define 
$\fb':\cd(Z)@>>>\cd(\tcb^2)$ and $\fb':\cd_m(Z)@>>>\cd_m(\tcb^2)$ by 
$$\fb'(L)=p_{03!}{}^!\vt_!{}^!\eta^*L.$$
(The map ${}^!\tcb^4@>>>\tcb^2$ induced by $p_{03}:\tcb^4@>>>\tcb^2$ is
denoted again by $p_{03}$.) Let $\t:\cy@>>>{}^!\cy$ be as in 3.4 (it is a 
principal $T^2$-bundle). We have the following results.

(a) {\it If $L\in\cd^{\preceq}(Z)$ then $\fb'(L)\in\cd^{\preceq}\tcb^2$.}

(b) {\it If $L\in\cd^{\prec}(Z)$ then $\fb'(L)\in\cd^{\prec}\tcb^2$.}

(c) {\it If $L\in\cm^{\preceq}(Z)$ and $h>\r+2\nu+2a$ then 
$(\fb'(L))^h\in\cm^{\prec}\tcb^2$.}
\nl
We can assume that $L=\Bbb L_\l^{\dz}$ where $z\cdo\l\in\overset\smile\to{W\fs}$, 
$z\cdo\l\preceq\boc$. A variant of the proof of 3.5(a) gives:
$$\fb'(\cl_\l^{\dz\sha})\Bpq\{{}'L^{\dy,\dz,\dy\i,\{2\}}_{\l,\l,y(\l)}
\la-|z|-2\nu\ra;y\in W\},$$
hence
$$\fb'(\Bbb L_\l^{\dz\sha})\Bpq\{{}'L^{\dy,\dz,\dy\i,\{2\}}_{\l,\l,y(\l)}
\la-\nu+\r\ra;y\in W\}.$$
To prove (a) it is enough to show that for any $y\in W$ we have
$${}'L^{\dy,\dz,\dy\i,\{2\}}_{\l,\l,y(\l)}\in\cd^{\preceq}\tcb^2.$$
When $z\cdo\l\in\boc$ this follows from 2.10(c). When $z\cdo\l\prec\boc$ this 
again follows from 2.10(c), applied to the two-sided cell containing $z\cdo\l$
instead of $\boc$. The same argument proves (b). 
To prove (c) we can assume that $z\cdo\l\in\boc$; it is 
enough to prove that for any $y\in W$ we have 
$$({}'L^{\dy,\dz,\dy\i,\{2\}}_{\l,\l,y(\l)}\la-\nu+\r\ra)^h\in
\cm^{\prec}\tcb^2$$
if $h>\r+2\nu+2a$ or that $({}'L^{\dy,\dz,\dy\i,\{2\}}_{\l,\l,y(\l)})^j\in\cm^{\prec}\tcb^2$
if $j>2\r+\nu+2a$. This follows from 2.20(c). This completes the proof of (a),(b),(c).

\mpb

We define $\un{\fb'}:\cc^\boc_0(Z)@>>>\cc^\boc_0(\tcb^2)$ by
$$\un{\fb'}(L)=\un{gr_{\r+2\nu+2a}((\fb'(L))^{\r+2\nu+2a})}((\r+2\nu+2a)/2).
$$ 
In the remainder of this subsection we fix $z\cdo\l\in\overset\smile\to\boc$ and we
set $L=\Bbb L_\l^{\dz}$. We show:

(d) {\it We have canonically $\un{\fb'}(L)=\un\fb(L)$.}
\nl
The method of proof is similar to that of 2.22(a). It is based on the fact that
$$\fb(L)=\fb'(L)\ot\fL^{\ot2}$$
which follows from the definitions. We define $\car_{i,j}$ for $i\in[0,2\r+1]$
and $\cp_{i,j}$ for $i\in[0,2\r]$ as in 2.17, but replacing $L^J,{}'L^J,r,\d$ 
by $\fb(L),\fb'(L),3,2\r$. In particular, we have
$$\cp_{i,j}=\cx_{4\r-i}(i-2\r)\ot(\fb'(L))^{-4\r+i+j}\text{ for }i\in[0,2\r]
$$
where $\cx_{4\r-i}$ is a free abelian group of rank $\bin{2\r}{i}$ and 
$\cx_{4\r}=\ZZ$. We have for any $j$ an exact sequence analogous to 2.17(a):
$$\do@>>>\cp_{i,j-1}@>>>\car_{i+1,j}@>>>\car_{i,j}@>>>\cp_{i,j}@>>>
\car_{i+1,j+1}@>>>\car_{i,j+1}@>>>\do,\tag e$$
and we have
$$\car_{0,j}=(\fb(L))^j,\qua \cp_{0,j}=(\fb'(L))^{j-4\r}(-2\r).$$
We show:

(f) {\it If $i\in[0,2\r+1]$ then $\car_{i,j}\in\cm^{\preceq}\tcb^2$.}

(g) {\it If $i\in[0,2\r+1]$, $j>6\r-i+\nu+2a$ then 
$\car_{i,j}\in\cm^{\prec}\tcb^2$.}
\nl
We prove (f),(g) by descending induction on $i$ as in 2.21. If $i=2\r+1$ then,
since $\car_{2\r+1,j}=0$, there is nothing to prove. Now assume that 
$i\in[0,2\r]$. 
Assume that $\LL_\l^{\dw}$ is a composition factor of $\car_{i,j}$ (without 
the mixed structure). We must show that $w\cdo\l\preceq\boc$ and that, if 
$j>6\r-i+\nu+2a$ then $w\cdo\l\prec\boc$. Using (e), we see that 
$\LL_\l^{\dw}$ is a composition factor of $\car_{i+1,j}$ or of $\cp_{i,j}$. In
the first case, using the induction hypothesis we see that 
$w\cdo\l\preceq\boc$ and that, if $j>6\r-i+\nu+2a$ (so that 
$j>6\r-i-1+\nu+2a$), then $w\cdo\l\prec\boc$. In the second case,
$\LL_\l^{\dw}$ is a composition factor of $(\fb'(L))^{-4\r+i+j}$. Using (a),
(c), we see that $w\cdo\l\preceq\boc$ and that, if $j>6\r-i+\nu+2a$ (so that 
$-4\r+i+j>\nu+2\r+2a$), then $w\cdo\l\prec\boc$. This proves (f),(g).

We show:

(h) {\it Assume that $i\in[0,2\r+1]$. Then $\car_{i,j}$ is mixed of weight 
$\le j-i$.}
\nl 
We argue as in 2.22 by descending induction on $i$. If $i=2\r+1$ there is 
nothing to prove. Assume now that $i\le2\r$. By Deligne's theorem, 
$\fb'(L)$ is mixed of weight $\le0$; hence $(\fb'(L))^{-4\r+i+j}$ is mixed 
of weight $\le-4\r+i+j$ and $\cx_{4\r-i}(i-2\r)\ot(\fb'(L))^{-4\r+i+j}$ is 
mixed of weight $\le-4\r+i+j-2(i-2\r)=j-i$. In other words, $\cp_{i,j}$ is
mixed of weight $\le j-i$. Thus in the exact sequence 
$\car_{i+1,j}@>>>\car_{i,j}@>>>\cp_{i,j}$ coming from (e) in which
$\car_{i+1,j}$ is mixed of weight $\le j-i-1<j-i$ (by the induction 
hypothesis) and $\cp_{i,j}$ is mixed of weight $\le j-i$ we must have that 
$\car_{i,j}$ is mixed of weight $\le j-i$. This proves (h).

\mpb

We now prove (d). From (e) we deduce an exact sequence 
$$gr_j(\car_{1,j})@>>>gr_j(\car_{0,j})@>>>gr_j(\cp_{0,j})@>>>
gr_j(\car_{1,j+1}).$$
By (h) we have $gr_j(\car_{1,j})=0$. We have 
$gr_j(\car_{0,j})=gr_j(\fb(L)^j)$,
$gr_j(\cp_{0,j})=gr_j((\fb'(L))^{-4\r+j}(-2\r))$. Moreover, by (g) we have 
$\car_{1,j+1}\in\cd^{\prec}\tcb^2$ since $j+1>6\r-1+\nu+2a$. It follows that 
$gr_j(\car_{1,j+1})\in\cd^{\prec}\tcb^2$.
Thus the exact sequence above induces an isomorphism as in (d).

\mpb

We show:

(i) {\it Let $L\in\cd(Z)$. Let $L'\in\cm(\tcb^2)$ be $G$-equivariant. 
We have canonically} $$\fb'(L)\cir L'=L'\cir\fb'(L).$$
\nl
Let $R=\TT\bsl\{((x_0\UU,x_1\UU,x_2\UU,x_3\UU),g)\in\tcb^4\T G;g\in x_2\UU x_1\i\}$ where
$\TT$ acts freely by
$$t:((x_0\UU,x_1\UU,x_2\UU,x_3\UU),g)\m((x_0\UU,x_1t\UU,x_2t\UU,x_3\UU),g).$$

Define $c_0:R@>>>Z$ by $((x_0\UU,x_1\UU,x_2\UU,x_3\UU),g)\m\e(x_1\UU,x_2\UU)$.

Define $c_1:R@>>>\tcb^2$ by $((x_0\UU,x_1\UU,x_2\UU,x_3\UU),g)\m(x_0\UU,gx_3\UU)$.

Define $c_2:R@>>>\tcb^2$ by $((x_0\UU,x_1\UU,x_2\UU,x_3\UU),g)\m(g\i x_0\UU,x_3\UU)$.

Define $c_3:R@>>>\tcb^2$ by $((x_0\UU,x_1\UU,x_2\UU,x_3\UU),g)\m(x_0\UU,x_3\UU)$. We have 
$$L'\cir\fb'(L)=c_{3!}(c_1^*L'\ot c_0^*L), \fb'(L)\cir L'=c_{3!}(c_2^*L'\ot c_0^*L).$$ 
It is enough to show that $c_1^*L'=c_2^*L'$. This follows from the $G$-equivariance of $L'$.

\mpb

(j) {\it If $L\in\cc^\boc_0Z$, $L'\in\cc^\boc\tcb^2$, then we have canonically 
$\un\fb(L)\un\cir L'=L'\un\cir\un\fb(L)$.}
\nl
By (d), it is enough to prove that $\un\fb'(L)\un\cir L'=L'\un\cir\un\fb'(L)$.
Using (i) together with (a),(b),(c) and results in 2.23, we see that both sides
are equal to 
$$\align&\un{gr_{\r+\nu+3a}(c_{3!}(c_1^*L'\ot c_0^*L))^{\r+\nu+3a}}((\r+\nu+3a)/2)\\&
=\un{gr_{\r+\nu+3a}(c_{3!}(c_2^*L'\ot c_0^*L))^{\r+\nu+3a}}((\r+\nu+3a)/2).\endalign$$

\subhead 3.15\endsubhead
Let 
$$\fZ=\{(z_0\UU,z_1\UU,z_2\UU,z_3\UU),g)\in\tcb^4\T G;g\in z_2\BB z_1\i\}.$$
Define $\ti\vt:\fZ@>>>\tcb^4$ by
$((z_0\UU,z_1\UU,z_2\UU,z_3\UU),g)\m(z_0\UU,z_1\UU,z_2\UU,z_3\UU)$. Let 
$${}'\cy=\{((x_0\UU,x_1\UU,x_2\UU,x_3\UU,x_4\UU),g)\in\tcb^5\T G;
g\in x_3\UU x_0\i, g\in x_2\BB x_1\i\},$$
$${}''\cy=\{((x_0\UU,x_1\UU,x_2\UU,x_3\UU,x_4\UU),g)\in\tcb^5\T G;
g\in x_4\UU x_1\i, g\in x_3\BB x_2\i\},$$
Define ${}'\vt:{}'\cy@>>>\tcb^5$, ${}''\vt:{}''\cy@>>>\tcb^5$ by
$$((x_0\UU,x_1\UU,x_2\UU,x_3\UU,x_4\UU),g)\m
(x_0\UU,x_1\UU,x_2\UU,x_3\UU,x_4\UU).$$
We have isomorphisms ${}'\fc:{}'\cy@>\si>>\fZ$, ${}''\fc:{}''\cy@>\si>>\fZ$ given by
$${}'\fc:((x_0\UU,x_1\UU,x_2\UU,x_3\UU,x_4\UU),g)\m((x_0\UU,x_1\UU,x_2\UU,x_4\UU),g),$$
$${}''\fc:((x_0\UU,x_1\UU,x_2\UU,x_3\UU,x_4\UU),g)\m((x_0\UU,x_2\UU,x_3\UU,x_4\UU),g).$$
Define ${}'d:\tcb^5@>>>\tcb^4$, ${}''d:\tcb^5@>>>\tcb^4$ by
$${}'d:(x_0\UU,x_1\UU,x_2\UU,x_3\UU,x_4\UU)\m(x_0\UU,x_1\UU,x_2\UU,x_4\UU),$$
$${}''d:(x_0\UU,x_1\UU,x_2\UU,x_3\UU,x_4\UU)\m(x_0\UU,x_2\UU,x_3\UU,x_4\UU).$$
We fix $w,u$ in $W$ and $\l,\l'$ in $\fs$ such that $w(\l)=\l$.

The smooth subvarieties
$${}'\cu=\{((x_0\UU,x_1\UU,x_2\UU,x_3\UU,x_4\UU),g)\in{}'\cy;
x_1\i x_2\in G_w,x_3\i x_4\in G_u\},$$
$$\cu=\{((x_0\UU,x_1\UU,x_2\UU,x_3\UU),g)\in\fZ;x_1\i x_2\in G_w,
x_0\i g\i x_3\in G_u\},$$   
$${}''\cu=\{((x_0\UU,x_1\UU,x_2\UU,x_3\UU,x_4\UU),g)\in{}''\cy;
x_2\i x_3\in G_w,x_0\i x_1\in G_u\},$$
of ${}'\cy,\fZ,{}''\cy$ correspond to each other under the isomorphisms
${}'\cy@>{}'\fc>>\fZ@<{}''\fc<<{}''\cy$. Moreover, the maps 
${}'\s:{}'\cu@>>>Z$, $\s:\cu@>>>Z$, ${}''\s:{}''\cu@>>>Z$ given by
$$((x_0\UU,x_1\UU,x_2\UU,x_3\UU,x_4\UU),g)\m\e((x_1\UU,x_2\UU),$$
$$((x_0\UU,x_1\UU,x_2\UU,x_3\UU),g)\m\e((x_1\UU,x_2\UU),$$
$$((x_0\UU,x_1\UU,x_2\UU,x_3\UU,x_4\UU),g)\m\e(x_2\UU,x_3\UU),$$
correspond to each other under the isomorphisms
${}'\cy@>{}'\fc>>\fZ@<{}''\fc<<{}''\cy$.

Also the maps ${}'\ti\s:{}'\cu@>>>\tcb^2$, $\ti\s:\cu@>>>\tcb^2$, given by
$$((x_0\UU,x_1\UU,x_2\UU,x_3\UU,x_4\UU),g)\m(x_3\UU,x_4\UU),$$
$$((x_0\UU,x_1\UU,x_2\UU,x_3\UU),g)\m(gx_0\UU,x_3\UU)$$
correspond to each other under the isomorphism ${}'\cy@>{}'\fc>>\fZ$ and the 
maps $\ti\s_1:\cu@>>>\tcb^2$, ${}''\ti\s:{}''\cu@>>>\tcb^2$ given by
$$((x_0\UU,x_1\UU,x_2\UU,x_3\UU),g)\m(x_0\UU,g\i x_3\UU),$$
$$((x_0\UU,x_1\UU,x_2\UU,x_3\UU,x_4\UU),g)\m(x_0\UU,x_1\UU),$$
correspond to each other under the isomorphism $\fZ@<{}''\fc<<{}''\cy$. It 
follows that the local systems ${}'\s^*\cl_\l^{\dw}$, $\s^*\cl_\l^{\dw}$, 
${}''\s^*\cl_\l^{\dw}$ correspond to each other under the isomorphisms 
${}'\cy@>{}'\fc>>\fZ@<{}''\fc<<{}''\cy$; the local systems 
${}'\ti\s^*L_{\l'}^{\du}$, $\ti\s^*L_{\l'}^{\du}$ correspond to each other 
under the isomorphism ${}'\cy@>{}'\fc>>\fZ$; the local systems 
$\ti\s_1^*L_{\l'}^{\du}$, ${}''\ti\s^*L_{\l'}^{\du}$ correspond to each other 
under the isomorphism $\fZ@<{}''\fc<<{}''\cy$. Moreover, we have 
$\ti\s^*L_{\l'}^{\du}=\ti\s_1^*L_{\l'}^{\du}$ by the $G$-equivariance of 
$L_{\l'}^{\du}$. Let ${}'K, K,{}''K$ be the intersection cohomology complex of
the closure of ${}'\cu,\cu,{}''\cu$ respectively with coefficients in the 
local system
$${}'\s^*\cl_\l^{\dw}\ot{}'\ti\s^*L_{\l'}^{\du},
\s^*\cl_\l^{\dw}\ot\ti\s^*L_{\l'}^{\du},{}''\s^*\cl_\l^{\dw}\ot{}''\ti\s^*L_{\l'}^{\du},$$
on ${}'\cu,\cu,{}''\cu$ (respectively) extended by $0$ on the complement of 
this closure in ${}'\cy,\fZ,{}''\cy$. We see that ${}'K,K,{}''K$ correspond to
each other under the isomorphisms ${}'\cy@>{}'\fc>>\fZ@<{}''\fc<<{}''\cy$. 
Hence we have ${}'\fc_!({}'K)=K={}''\fc_!({}''K)$. Using this and the 
commutative diagram
$$\CD
{}'\cy@>{}'\fc>>\fZ@<{}''\fc<<{}''\cy\\
@V{}'\vt VV   @V\ti\vt VV     @V{}''\vt VV   \\
\tcb^5@>{}'d>>\tcb^4@<{}''d<<\tcb^5
\endCD$$
we see that 
$${}'d_!{}'\vt_!({}'K)={}''d_!{}''\vt_!({}''K).\tag a$$
(Both sides are equal to $\ti\vt_!K$.)

\subhead 3.16\endsubhead
In this subsection we study the functor ${}'d_!:\cd_m(\tcb^5)@>>>\cd_m(\tcb^4)$. 
Let $\ww=(w_1,w_2,w_3,w_4)$, $\pmb\l=(\l_1,\l_2,\l_3,\l_4)$,
$\pmb\o=(\o_1,\o_2,\o_3,\o_4)$ (with $\o_i\in\k\i_q(w_i)$). Assume that
$w_4\cdo\l_4\preceq\boc$. Let $K=M_{\pmb\l}^{\pmb\o,[1,4]}\la|\ww|+5\r+\nu\ra\in\cd_m(\tcb^5)$.
We show:

(a) {\it If $h>a+\r$ then $({}'d_!K)^h\in{}'\cm^{\prec}(\tcb^4)$. Moreover,}
$$\align&\un{gr_{a+\r}(({}'dK)^{a+\r})}((a+\r)/2)=\op_{y'\in W;y'{}\i\cdo\l_4\in\boc}
\Hom_{\cc^\boc\tcb^2}(\LL_{\l_4}^{\dy'{}\i},\LL_{\l_3}^{\o_3}\un{\cir}
\LL_{\l_4}^{\o_4})\\&\ot M_{\l_1,\l_2,\l_4}^{\o_1,\o_2,\dy'{}\i,[1,3]}\la
|w_1|+|w_2|+|y'|+4\r+\nu\ra.\endalign$$
We shall apply \cite{\CONV, 1.12} with 
$\Ph:\cd_m(Y_1)@>>>\cd_m(Y_2)$ replaced by 
$\Ph_0:\cd_m(\tcb^2)@>>>\cd_m(\tcb^4)$,
$M\m p_{01}^*L_{\l_1}^{\o_1\sha}\la|w_1|\ra\ot
p_{12}^*L_{\l_2}^{\o_2\sha}\la|w_2|\ra\ot p_{23}^*M\la\r-\nu\ra$
and with $\cd^{\preceq}(Y_1)$, $\cd^{\preceq}(Y_2)$ replaced by
${}_4\cd^{\preceq}(\tcb^2)$, ${}_4\cd^{\preceq}(\tcb^4)$, see 3.15. We shall 
take $\XX$ in {\it loc.cit.} equal to 
$\Xi:=L_{\l_3}^{\o_3\sha}\cir L_{\l_4}^{\o_4\sha}\la|w_3|+|w_4|+4\r+2\nu\ra$.
The conditions of {\it loc.cit.} are satisfied: those concerning $\XX$ are
satisfied with $c'=a-\nu$ (see 2.23); those concerning $\Ph_0$ are satisfied 
with $c=\r+\nu$ (using the definitions). Since ${}'d_!K=\Ph_0(\Xi)$ and 
$c+c'=a+\r$ we see that the first sentence in (a) holds; moreover, we see 
that, setting
$K_1=\un{gr_{a+\r}(({}'d_!K)^{a+\r})}((a+\r)/2)$, we have
$$K_1=\un{gr_{\r+\nu}((p_{01}^*L_{\l_1}^{\o_1\sha}\la|w_1|\ra\ot 
p_{12}^*L_{\l_2}^{\o_2\sha}\la|w_2|\ra\ot p_{23}^*M\la\r-\nu\ra)^{\r+\nu})}
((\r+\nu)/2)$$    
where 
$$M=\LL_{\l_3}^{\o_3}\un{\cir}\LL_{\l_4}^{\o_4}=\op_{y'{}\i\cdo\eta\in\boc}
\Hom_{\cc^\boc\tcb^2}(\LL_\eta^{\dy'{}\i},
\LL_{\l_3}^{\o_3}\un{\cir}\LL_{\l_4}^{\o_4})\LL_\eta^{\dy'{}\i}.$$
From 2.13(c) we see that $\eta$ above must satisfy $\eta=\l_4$. Thus we have
$$\align&K_1=\op_{y'\in W;y'{}\i\cdo\l_4\in\boc}\Hom_{\cc^\boc\tcb^2}
(\LL_{\l_4}^{\dy'{}\i},
\LL_{\l_3}^{\o_3}\un{\cir}\LL_{\l_4}^{\o_4}) \\& \ot\un{gr_{\r+\nu}
(M_{\l_1,\l_2,\l_4}^{\o_1,\o_2,\dy'{}\i,[1,3]}
\la|w_1|+|w_2|+|y'|+3\r\ra)^{\r+\nu})}((\r+\nu)/2).\endalign$$
It remains to use that
$$\align&\un{gr_{\r+\nu}(M_{\l_1,\l_2,\l_4}^{\o_1,\o_2,\dy'{}\i,[1,3]}
\la|w_1|+|w_2|+|y'|+3\r\ra)^{\r+\nu})}((\r+\nu)/2)\\&=\un{gr_{\r+\nu}
(M_{\l_1,\l_2,\l_4}^{\o_1,\o_2,\dy'{}\i,[1,3]}\la|w_1|+|w_2|+|y'|+4\r+\nu\ra)^0
(-(r+\nu)/2))}((\r+\nu)/2)\\&=
M_{\l_1,\l_2,\l_4}^{\o_1,\o_2,\dy'{}\i,[1,3]}\la|w_1|+|w_2|+|y'|+4\r+\nu\ra.\endalign$$
We state the following properties of the functor
${}'d_!:\cd_m(\tcb^5)@>>>\cd_m(\tcb^4)$.

(b) {\it If $K\in{}_4\cd^{\preceq}(\tcb^5)$ then 
${}'d_!(K)\in{}_4\cd^{\preceq}(\tcb^4)$.}

(c) {\it If $K\in{}_4\cd^{\prec}(\tcb^5)$ then 
${}'d_!(K)\in{}_4\cd^{\prec}(\tcb^4)$.}

(d) {\it If $K\in{}_4\cm^{\preceq}(\tcb^5)$ and $h>a+\r$ then 
$({}'d_!(K))^h\in{}_4\cm^{\prec}(\tcb^4)$.}
\nl
We prove (b). We can assume that $K$ is as in the first paragraph of this subsection. 
It is enough to show that
for $j\in\ZZ$ we have $(\Ph_0(\Xi))^j\in{}_4\cm^{\preceq}(\tcb^4)$ (with 
$\Ph_0,\Xi$ as above). It is enough to show that 
$\Ph_0(\Xi^{j'})\in{}_4\cm^{\preceq}(\tcb^4)$ for any $j'\in\ZZ$. This follows 
from the fact that $\Xi^{j'}\in{}_4\cm^{\preceq}(\tcb^2)$ (see 2.23(b)) and the
fact that $\Ph_0$ carries ${}_4\cd^{\preceq}(\tcb^2)$ to 
${}_4\cd^{\preceq}(\tcb^4)$. Thus (b) holds. A similar proof gives 
(c). We prove (d). We can assume
that $K$ is as in (a). Then the result follows from (a).

\subhead 3.17\endsubhead
In this subsection we study the functor ${}''d_!:\cd_m(\tcb^5)@>>>\cd_m(\tcb^4)$.
Let $\ww=(w_1,w_2,w_3,w_4)$, $\pmb\l=(\l_1,\l_2,\l_3,\l_4)$,
$\pmb\o=(\o_1,\o_2,\o_3,\o_4)$ (with $\o_i\in\k\i_q(w_i)$). Assume that
$w_1\cdo\l_1\preceq\boc$. Let $K=M_{\pmb\l}^{\pmb\o,[1,4]}\la|\ww|+5\r+\nu\ra\in\cd_m(\tcb^5)$.
We show:

(a) {\it If $h>a+\r$ then $({}''d_!K)^h\in{}'\cm^{\prec}(\tcb^4)$. Moreover,}
$$\align&\un{gr_{a+\r}(({}''d_!K)^{a+\r})}((a+\r)/2)=\op_{y'\in W;y'\cdo\l_2\in\boc}
\Hom_{\cc^\boc\tcb^2}(\LL_{\l_2}^{\dy'},\LL_{\l_1}^{\o_1}\un{\cir}
\LL_{\l_2}^{\o_2})\\&
\ot M_{\l_2,\l_3,\l_4}^{\dy',\o_3,\o_4,[1,3]}\la|w_3|+|w_4|+|y'|+4\r+\nu\ra.\endalign$$
We shall apply  \cite{\CONV, 1.12} with 
$\Ph:\cd_m(Y_1)@>>>\cd_m(Y_2)$ replaced by 
$\Ph_0:\cd_m(\tcb^2)@>>>\cd_m(\tcb^4)$,
$M\m p_{01}^*M\ot p_{12}^*L_{\l_3}^{\o_3\sha}\la|w_3|\ra\ot 
p_{23}^*L_{\l_4}^{\o_4\sha}\la|w_4|\ra\la\r-\nu\ra$
and with $\cd^{\preceq}(Y_1)$, $\cd^{\preceq}(Y_2)$ replaced by
${}_1\cd^{\preceq}(\tcb^2)$, ${}_1\cd^{\preceq}(\tcb^4)$, see 3.15. We shall 
take $\XX$ in {\it loc.cit.} equal to 
$\Xi:=L_{\l_1}^{\o_1\sha}\cir L_{\l_2}^{\o_2\sha}\la|w_1|+|w_2|+4\r+2\nu\ra$.
The conditions of {\it loc.cit.} are satisfied: those concerning $\XX$ are
satisfied with $c'=a-\nu$ (see 2.23); those concerning $\Ph_0$ are satisfied 
with $c=\r+\nu$ (using the definitions). Since ${}''d_!K=\Ph_0(\Xi)$ and 
$c+c'=a+\r$ we see that the first sentence in (a) holds; moreover we see that,
setting $K_1=\un{gr_{a+\r}(({}''d_!K)^{a+\r})}((a+\r)/2)$, we have
$$K_1=\un{gr_{\r+\nu}((p_{01}^*M\ot p_{12}^*L_{\l_3}^{\o_3\sha}\la|w_3|\ra\ot 
p_{23}^*L_{\l_4}^{\o_4\sha}\la|w_4|+\r-\nu\ra)^{\r+\nu})}((\r+\nu)/2)$$
where 
$$M=\LL_{\l_1}^{\o_1}\un{\cir}\LL_{\l_2}^{\o_2}=\op_{y'\cdo\eta\in\boc}
\Hom_{\cc^\boc\tcb^2}(\LL_\eta^{\dy'},
\LL_{\l_1}^{\o_1}\un{\cir}\LL_{\l_2}^{\o_2})\LL_\eta^{\dy'}.$$
From 2.13(c) we see that $\eta$ must satisfy $\eta=\l_2$. Thus we have
$$\align&K_1=\op_{y'\in W;y'\cdo\l_2\in\boc}\Hom_{\cc^\boc\tcb^2}(\LL_{\l_2}^{\dy'},  
\LL_{\l_1}^{\o_1}\un{\cir}\LL_{\l_2}^{\o_2})\\&\ot\un{gr_{\r+\nu}
(M_{\l_2,\l_3,\l_4}^{\dy',\o_3,\o_4,[1,3]}
\la|w_3|+|w_4|+|y'|+3\r\ra)^{\r+\nu})}((\r+\nu)/2).\endalign$$
It remains to use that
$$\align&\un{gr_{\r+\nu}(M_{\l_2,\l_3,\l_4}^{\dy',\o_3,\o_4,[1,3]}
\la|w_3|+|w_4|+|y'|+3\r\ra)^{\r+\nu})}((\r+\nu)/2)\\&=\un{gr_{\r+\nu}
(M_{\l_2,\l_3,\l_4}^{\dy',\o_3,\o_4,[1,3]}\la|w_3|+|w_4|+|y'|+4\r+\nu\ra)^0
(-(r+\nu)/2))}((\r+\nu)/2)\\&=
M_{\l_2,\l_3,\l_4}^{\dy',\o_3,\o_4,[1,3]}\la|w_3|+|w_4|+|y'|+4\r+\nu\ra.\endalign$$
We state the following properties of the functor
${}''d_!:\cd_m(\tcb^5)@>>>\cd_m(\tcb^4)$.

(b) {\it If $K\in{}_1\cd^{\preceq}(\tcb^5)$ then 
${}''d_!(K)\in{}_1\cd^{\preceq}(\tcb^4)$.}

(c) {\it If $K\in{}_1\cd^{\prec}(\tcb^5)$ then 
${}''d_!(K)\in{}_1\cd^{\prec}(\tcb^4)$.}

(d) {\it If $K\in{}_1\cm^{\preceq}(\tcb^5)$ and $h>a+\r$ then 
$({}''d_!(K))^h\in{}_1\cm^{\prec}(\tcb^4)$.}
\nl
The proof of (b),(c),(d) is completely similar to that of 3.16(b),(c),(d).

\subhead 3.18\endsubhead
Let $w\cdo\l\in\overset\smile\to{W\fs}$, $u\cdo\l'\in\boc$. We shall apply  
\cite{\CONV, 1.12} with $\Ph:\cd_m(Y_1)@>>>\cd_m(Y_2)$ replaced by 
${}'d_!:\cd_m(\tcb^5)@>>>\cd_m(\tcb^4)$ and with $\cd^{\preceq}(Y_1)$, 
$\cd^{\preceq}(Y_2)$ replaced by ${}_4\cd^{\preceq}(\tcb^5)$, 
${}_4\cd^{\preceq}(\tcb^4)$, see 3.15. We shall take $\XX$ in
{\it loc.cit.} equal to $\Xi={}'\vt_!({}'K)$ as in 3.15, $(w_2,w_4)=(w,u)$, 
$(\l_2,\l_4)=(\l,\l')$. The conditions of {\it loc.cit.} are satisfied: those 
concerning $\XX$ are satisfied with $c'=k=|w|+|u|+3\nu+5\r$ (see 3.8(c)); 
those concerning $\Ph$ are satisfied with $c=a+\r$ (see 3.16). We see that
$$\align&\un{gr_{a+\r+k}(({}'d_!{}'\vt_!({}'K))^{a+\r+k})}((a+\r+k)/2)\\&=
\un{gr_{a+\r}(({}'d_!\un{gr_k(({}'\vt_!({}'K))^k)}(k/2))^{a+\r})}((a+\r)/2).\endalign$$
Using 3.11(a) we have:
$$\align&gr_k({}'\vt_!({}'K))^k)(k/2)=\op_{y\in W}M_{\l,\l,y(\l),\l'}^{\dy,\dw,\dy\i,
\du,[1,4]}\la2|y|+|w|+|u|+5\r+\nu\ra\\&=\un{gr_k({}'\vt_!({}'K))^k}(k/2).\endalign$$
Hence using 3.16(a) we have
$$\align&\un{gr_{a+\r}(({}'d_!\un{gr_k(({}'\vt_!({}'K))^k)}(k/2))^{a+\r})}((a+\r)/2)\\&=
\op_{y\in W}\op_{y'\in W;y'{}\i\cdo\l'\in\boc}\Hom_{\cc^\boc\tcb^2}
(\LL_{\l'}^{\dy'{}\i},\LL_{y(\l)}^{\dy\i}\un{\cir}\LL_{\l'}^{\du})\\&
\ot M_{\l,\l,\l'}^{\dy,\dw,\dy'{}\i,[1,3]}\la|y|+|w|+|y'|+4\r+\nu\ra.\endalign$$
Since $y'{}\i\cdo\l'\in\boc$, $u\cdo\l'\in\boc$, for $y\in W$ we have 
$$\Hom_{\cc^\boc\tcb^2}(\LL_{\l'}^{\dy'{}\i},\LL_{y(\l)}^{\dy\i}\un{\cir}
\LL_{\l'}^{\du})=0$$ unless $y\i\cdo y(\l)\in\boc$ (see 2.26) or equivalently 
(see 1.9(Q10), 1.11), $y\cdo\l\in\boc$. Thus we have
$$\align&\un{gr_{a+\r+k}(({}'d_!{}'\vt_!({}'K))^{a+\r+k})}((a+\r+k)/2)\\&=
\op_{y\in W;y\cdo\l\in\boc}
\op_{y'\in W;y'{}\i\cdo\l'\in\boc}\Hom_{\cc^\boc\tcb^2}(\LL_{\l'}^{\dy'{}\i},
\LL_{y(\l)}^{\dy\i}\un{\cir}\LL_{\l'}^{\du})\\&
\ot M_{\l,\l,\l'}^{\dy,\dw,\dy'{}\i,[1,3]}\la|y|+|w|+|y'|+4\r+\nu\ra.\endalign$$
The last $\Hom$-space is zero unless $y'{}\i(\l')=\l$ hence
$$\align&\un{gr_{a+\r+k}(({}'d_!{}'\vt_!({}'K))^{a+\r+k})}((a+\r+k)/2)\\&=
\op_{y\in W;y\cdo\l\in\boc}
\op_{y'\in W;y'{}\i\cdo y'(\l)\in\boc}\Hom_{\cc^\boc\tcb^2}
(\LL_{y'(\l)}^{\dy'{}\i},\LL_{y(\l)}^{\dy\i}\un{\cir}\LL_{\l'}^{\du})\\&
\ot M_{\l,\l,\l'}^{\dy,\dw,\dy'{}\i,[1,3]}\la|y|+|w|+|y'|+4\r+\nu\ra.\tag a\endalign$$
     
\subhead 3.19\endsubhead
In the setup of 3.18 we shall apply 
\cite{\CONV, 1.12} with $\Ph:\cd_m(Y_1)@>>>\cd_m(Y_2)$ replaced by 
${}''d_!:\cd_m(\tcb^5)@>>>\cd_m(\tcb^4)$ and with $\cd^{\preceq}(Y_1)$, 
$\cd^{\preceq}(Y_2)$ replaced by ${}_1\cd^{\preceq}(\tcb^5)$, 
${}_1\cd^{\preceq}(\tcb^4)$, see 3.15. We shall take $\XX$ in {\it loc.cit.} 
equal to $\Xi={}''\vt_!({}''K)$ as in 3.15, $(w_1,w_3)=(u,w)$, 
$(\l_1,\l_3)=(\l',\l)$. The conditions of {\it loc.cit.} are satisfied: those 
concerning $\XX$ are satisfied with $c'=k=|w|+|u|+3\nu+5\r$ (see 3.8(c)); 
those concerning $\Ph$ are satisfied with $c=a+\r$ (see 3.17). We see that
$$\align&\un{gr_{a+\r+k}(({}''d_!{}''\vt_!({}''K))^{a+\r+k})}((a+\r+k)/2)\\&=
\un{gr_{a+\r}(({}''d_!\un{gr_k(({}''\vt_!({}''K))^k)}(k/2))^{a+\r})}((a+\r)/2).\endalign$$
Using 3.11(a) we have:
$$\align&gr_k({}''\vt_!({}''K))^k)(k/2)=\op_{y\in W}  
M_{\l',\l,\l,y(\l)}^{\du,\dy,\dw,\dy\i,[1,4]}\la2|y|+|w|+|u|+5\r+\nu\ra\\&
=\un{gr_k({}''\vt_!({}''K))^k}(k/2).\endalign$$
Hence using 3.17(a) we have
$$\align&\un{gr_{a+\r}(({}''d_!\un{gr_k(({}''\vt_!({}''K))^k)}(k/2))^{a+\r})}
((a+\r)/2)\\&=\op_{y\in W}\op_{y'\in W;y'\cdo\l\in\boc}
\Hom_{\cc^\boc\tcb^2}(\LL_\l^{\dy'},\LL_{\l'}^{\du}\un{\cir}\LL_\l^{\dy})\\&
\ot M_{\l,\l,y(\l)}^{\dy',\dw,\dy\i,[1,3]}\la|y|+|w|+|y'|+4\r+\nu\ra.\endalign$$
Since $y'\cdo\l\in\boc$, $u\cdo\l'\in\boc$, for $y\in W$ we have 
$\Hom_{\cc^\boc\tcb^2}(\LL_\l^{\dy'},\LL_{\l'}^{\du}\un{\cir}\LL_\l^{\dy})=0$
unless $y\cdo\l\in\boc$ (see 2.26). Thus we have
$$\align&\un{gr_{a+\r+k}(({}''d_!{}''\vt_!({}''K))^{a+\r+k})}((a+\r+k)/2)\\&=
\op_{y\in W;y\cdo\l\in\boc}\op_{y'\in W;y'\cdo\l\in\boc}
\Hom_{\cc^\boc\tcb^2}(\LL_\l^{\dy'},\LL_{\l'}^{\du}\un{\cir}\LL_\l^{\dy})\\&
\ot M_{\l,\l,y(\l)}^{\dy',\dw,\dy\i,[1,3]}\la|y|+|w|+|y'|+4\r+\nu\ra.\endalign$$
The last $\Hom$-space is zero unless $y(\l)=\l'$ hence (with the change of 
notation $(y,y')\m(y',y)$):
$$\align&\un{gr_{a+\r+k}(({}''d_!{}''\vt_!({}''K))^{a+\r+k})}((a+\r+k)/2)\\&=
\op_{y\in W;y\cdo\l\in\boc}\op_{y'\in W;y'\cdo\l\in\boc}
\Hom_{\cc^\boc\tcb^2}(\LL_\l^{\dy},\LL_{\l'}^{\du}\un{\cir}\LL_\l^{\dy'})\\&
\ot M_{\l,\l,\l'}^{\dy,\dw,\dy'{}\i,[1,3]}\la|y|+|w|+|y'|+4\r+\nu\ra.\tag a\endalign$$

\subhead 3.20\endsubhead
Let $y_1\cdo\l_1\in\boc$, $y_2\cdo\l_2\in\boc$, $y_3\cdo\l_3\in\boc$. We show:

(a) {\it We have canonically}
$$\Hom_{\cc^\boc\tcb^2}(\LL_{y_2(\l_2)}^{\dy_2\i},
\LL_{y_1(\l_1)}^{\dy_1\i}\un{\cir}\LL_{\l_3}^{\dy_3})
=\Hom_{\cc^\boc\tcb^2}(\LL_{\l_1}^{\dy_1},\LL_{\l_3}^{\dy_3}\un{\cir}
\LL_{\l_2}^{\dy_2}).$$
When $\l_1\ne\l_2$, both sides of the last equality (to be proved) are zero and
the result is clear. In the rest of the proof we assume $\l_1=\l_2=\l$. We set
$u\cdo\l'=y_3\cdo\l_3$. Choose $w\in W$ such that $w\cdo\l\in\overset\smile\to{W\fs}$. 

Applying 3.18(a), 3.19(a) to our $w\cdo\l$, $u\cdo\l'$ and using the equality
$$\align&\un{gr_{a+\r+k}(({}'d_!{}'\vt_!({}'K))^{a+\r+k})}((a+\r+k)/2)\\&=
\un{gr_{a+\r+k}(({}''d_!{}''\vt_!({}''K))^{a+\r+k})}((a+\r+k)/2)\endalign$$
which comes from ${}'d_!{}'\vt_!({}'K)={}''d_!{}''\vt_!({}''K)$, see 3.15(a),
we deduce
$$\align&\op_{y\in W;y\cdo\l\in\boc}\op_{y'\in W;y'\cdo\l\in\boc}\Hom_{\cc^\boc\tcb^2}
(\LL_{y'(\l)}^{\dy'{}\i},\LL_{y(\l)}^{\dy\i}\un{\cir}\LL_{\l'}^{\du})\\&
\ot M_{\l,\l,\l'}^{\dy,\dw,\dy'{}\i,[1,3]}\la|y|+|w|+|y'|+4\r+\nu\ra\\&=
\op_{y\in W;y\cdo\l\in\boc}\op_{y'\in W;y'\cdo\l\in\boc}
\Hom_{\cc^\boc\tcb^2}(\LL_\l^{\dy},\LL_{\l'}^{\du}\un{\cir}\LL_\l^{\dy'})\\&
\ot M_{\l,\l,\l'}^{\dy,\dw,\dy'{}\i,[1,3]}\la|y|+|w|+|y'|+4\r+\nu\ra.\tag b\endalign$$
Considering the coefficient of
$$M_{\l,\l,\l'}^{\dy_1,\dw,\dy_2\i,[1,3]}\la|y_1|+|w|+|y_2|+4\r+\nu\ra$$
in the two sides of (b) we obtain (a). From the proof one can see that the 
identification in (a) does not depend on the choice of $w$. 

\subhead 3.21\endsubhead
We assume that $w\cdo\l\in\overset\smile\to\boc$, $u\cdo\l'\in\boc$. We apply 
$p_{03!}$ and $\la N\ra$ for some $N$ to the two sides of 3.20(b). (Recall 
that $p_{03}:\tcb^4@>>>\tcb^2$.) We obtain
$$\align&\op_{y\in W;y\cdo\l\in\boc}\op_{y'\in W;y'\cdo\l\in\boc}\Hom_{\cc^\boc\tcb^2}
(\LL_{y'(\l)}^{\dy'{}\i},\LL_{y(\l)}^{\dy\i}\un{\cir}\LL_{\l'}^{\du})
\ot\LL_\l^{\dy}\cir\LL_\l^{\dw}\cir\LL_{y'(\l)}^{\dy'{}\i}\\&
=\op_{y\in W;y\cdo\l\in\boc}\op_{y'\in W;y'\cdo\l\in\boc}
\Hom_{\cc^\boc\tcb^2}(\LL_\l^{\dy},\LL_{\l'}^{\du}\un{\cir}\LL_\l^{\dy'})
\ot\LL_\l^{\dy}\cir\LL_\l^{\dw}\cir\LL_{y'(\l)}^{\dy'{}\i}.\endalign$$
(We have replaced $\LL_{\l'}^{\dy'{}\i}$ by $\LL_{y'(\l)}^{\dy'{}\i}$; in the 
last equality the terms with $\l'\ne y'(\l)$ contribute $0$.) Applying 
$\un{()^{\{2(a-\nu)\}}}$ to both sides and using 2.24(a) we obtain
$$\align&\op_{y\in W;y\cdo\l\in\boc}\op_{y'\in W;y'\cdo\l\in\boc}\Hom_{\cc^\boc\tcb^2}
(\LL_{y'(\l)}^{\dy'{}\i},\LL_{y(\l)}^{\dy\i}\un{\cir}\LL_{\l'}^{\du})
\ot\LL_\l^{\dy}\un{\cir}\LL_\l^{\dw}\un{\cir}\LL_{y'(\l)}^{\dy'{}\i}\\&
=\op_{y\in W;y\cdo\l\in\boc}\op_{y'\in W;y'\cdo\l\in\boc}
\Hom_{\cc^\boc\tcb^2}(\LL_\l^{\dy},\LL_{\l'}^{\du}\un{\cir}\LL_\l^{\dy'})
\ot\LL_\l^{\dy}\un{\cir}\LL_\l^{\dw}\un{\cir}\LL_{y'(\l)}^{\dy'{}\i}\endalign$$
or equivalently
$$\op_{y\in W;y\cdo\l\in\boc}\LL_\l^{\dy}\un{\cir}\LL_\l^{\dw}\un{\cir}
\LL_{y(\l)}^{\dy\i}\un{\cir}\LL_{\l'}^{\du}
=\op_{y'\in W;y'\cdo\l\in\boc}\LL_{\l'}^{\du}\un{\cir}\LL_\l^{\dy'}\un{\cir}
\LL_\l^{\dw}\un{\cir}\LL_{y'(\l)}^{\dy'{}\i}.$$
Using 3.13(d), this can be rewritten as follows:
$$\un{\fb}(\Bbb L_\l^{\dw})\un{\cir}\LL_{\l'}^{\du}
=\LL_{\l'}^{\du}\un{\cir}\un\fb(\Bbb L_\l^{\dw}).\tag a$$
Another identification of the two sides in (a) is given by 3.14(j) with $L=\Bbb L_\l^{\dw}$,
$L'=\LL_{\l'}^{\du}$ (note that $\un\fb(L)=\un\fb'(L)$ by 3.14(d)). 
In fact, the 
arguments in 3.13-3.21 show that

(b) {\it these two identifications of the two sides of (a) coincide.}

\subhead 3.22\endsubhead
Let 
$$\align&V=\{(B_0,B_1,B_2,gU_{B_0},g'U_{B_1});\\&(B_0,B_1,B_2)\in\cb^3,g\in G,g'\in G,
gB_0g\i=B_1,g'B_1g'{}\i=B_2\}.\endalign$$
Define $p_{01}:V@>>>Z$, $p_{12}:V@>>>Z$, $p_{02}:V@>>>Z$ by 
$$p_{01}:(B_0,B_1,B_2,gU_{B_0},g'U_{B_1})\m(B_0,B_1,gU_{B_0}),$$
$$p_{12}:(B_0,B_1,B_2,gU_{B_0},g'U_{B_1})\m(B_1,B_2,g'U_{B_1}),$$
$$p_{02}:(B_0,B_1,B_2,gU_{B_0},g'U_{B_1})\m(B_0,B_2,g'gU_{B_0}).$$
For $L,L'$ in $\cd(Z)$ we set 
$$L\bul L'=p_{02!}(p_{01}^*L\ot p_{12}^*L')\in\cd(Z).$$ 
This operation is associative. Hence if ${}^1L,{}^2L,\do,{}^rL$ are in $\cd(Z)$
then \lb ${}^1L\bul{}^2L\bul\do\bul{}^rL\in\cd(Z)$ is well defined. We show: 

(a) {\it For $L,L'$ in $\cd(Z)$ we have canonically 
$\e^*(L\bul L')=\e^*(L)\cir\e^*(L')$.}
\nl
Let 
$$Y=\{(x\UU,y\UU,gU_{x\BB x\i});x\UU\in\tcb,y\UU\in\tcb;g\in G\}.$$
Define $j:Y@>>>\tcb^2$, $j_1:Y@>>>Z$, $j_2:Y@>>>Z$ by
$$\align &j(x\UU,y\UU,gU_{x\BB x\i})=(x\UU,y\UU),\\&
j_1(x\UU,y\UU,gU_{x\BB x\i})=(x\BB x\i,gx\BB x\i g\i,gU_{x\BB x\i}),\\&
j_2(x\UU,y\UU,gU_{x\BB x\i})=(gx\BB x\i g\i,y\BB y\i,y\UU x\i g\i).\endalign$$
From the definitions we have 
$$\e^*(L\bul L')=j_!(j_1^*(L)\ot j_2^*(L'))=\e^*(L)\cir\e^*(L')$$
and (a) follows.

\subhead 3.23\endsubhead
Let $L,L'\in\cd^\spa Z$. We show:

(a) {\it If $L\in\cd^{\preceq}Z$ or $L'\in\cd^{\preceq}Z$ then 
$L\bul L'\in\cd^{\preceq}Z$. If $L\in\cd^{\prec}Z$ or $L'\in\cd^{\prec}Z$ then
$L\bul L'\in\cd^{\prec}Z$.}
\nl
For the first assertion of (a) we can assume that $L=\Bbb L_\l^{\dw}$,
$L'=\Bbb L_{\l'}^{\dw'}$ with $w\cdo\l,w'\cdo\l'$ in $\overset\smile\to{W\fs}$ and either 
$w\cdo\l\preceq\boc$ or $w'\cdo\l'\preceq\boc$. Assume that 
$w_1\cdo\l_1\in\overset\smile\to{W\fs}$ and $\Bbb L_{\l_1}^{\dw_1}$ is a composition factor 
of $(L\bul L')^j$. Then 
$\LL_{\l_1}^{\dw_1}=\ti\e\Bbb L_{\l_1}^{\dw_1}$ is a composition factor of
$$\align&\e^*(L\bul L')^j\la\r\ra=(\e^*(L\bul L'))^{j+\r}(\r/2)
=(\e^*L\cir\e^*L')^{j+\r}(\r/2)\\&
=(\e^*L\la\r\ra\cir\e^*L'\la\r\ra)^{j-\r}(-\r/2)
=(\LL_\l^{\dw}\cir\LL_{\l'}^{\dw'})^{j-\r}(\r/2).\endalign$$
From 2.23(b) we see that $w_1\cdo\l_1\preceq\boc$. This proves the first
assertion of (a). The second assertion of (a) can be reduced to the first 
assertion.

(b) {\it Assume that $L,L'\in\cm^\spa Z$ and that either $L$ or $L'$ is in 
$\cd^{\preceq}Z$. If $j>a+\r-\nu$ then $(L\bul L')^j\in\cm^{\prec}Z$.}
\nl
We can assume that $L=\Bbb L_\l^{\dw}$, $L'=\Bbb L_{\l'}^{\dw'}$ with 
$w\cdo\l,w'\cdo\l'$ in $\overset\smile\to{W\fs}$ and either $w\cdo\l\in\boc$ or 
$w'\cdo\l'\in\boc$. Assume that $w_1\cdo\l_1\in\overset\smile\to{W\fs}$ and 
$\Bbb L_{\l_1}^{\dw_1}$ is a composition factor of $(L\bul L')^j$. Then as in 
the proof of (a), $\LL_{\l_1}^{\dw_1}$ is a composition factor of 
$$\ti e(L\bul L')^j=(\LL_\l^{\dw}\cir\LL_{\l'}^{\dw'})^{j-\r}(-\r/2).$$
Since $j-\r>a-\nu$ we see from 2.23(a) that $w_1\cdo\l_1\prec\boc$. This 
proves (b).

\subhead 3.24\endsubhead
For $L,L'\in\cc^\boc_0Z$ we set
$$L\un{\bul}L'=\un{(L\bul L')^{\{a+\r-\nu\}}}\in\cc^\boc_0Z.$$
Using 3.23(a),(b) we see as in 2.24 that for $L,L',L''\in\cc^\boc_0Z$ we have
$$L\un{\bul}(L'\un{\bul}L'')
=(L\un{\bul}L')\un{\bul}L''=\un{(L\bul L'\bul L'')^{\{2a+2\r-2\nu\}}}.$$
We see that $L,L'\m L\un{\bul}L'$ defines a monoidal structure on 
$\cc^\boc_0Z$. Hence if 
$${}^1L,{}^2L,\do,{}^rL$$ 
are in $\cc^\boc_0Z$, then ${}^1L\un{\bul}{}^2L\un{\bul}\do\un{\bul}{}^rL\in\cc^\boc_0Z$
is well defined; we have
$${}^1L\un{\bul}{}^2L\un{\bul}\do\un{\bul}{}^rL=
\un{({}^1L\bul{}^2L\bul\do\bul{}^rL)^{\{(r-1)(a+\r-\nu)\}}}.\tag a$$
For $L,L'\in\cc^\boc_0Z$ we have 
$\ti\e L,\ti\e L'\in\cc^\boc_0\tcb^2$. We show:
$$\ti\e(L\un{\bul}L')=(\ti\e L)\un{\cir}(\ti\e L').\tag b$$
It is enough to show that
$$\align&\e^*(gr_0((L\bul L')^{a+\r-\nu})((a+\r-\nu)/2))[\r](\r/2)\\&
=gr_0((\e^*L[\r](\r/2)\cir\e^*L'[\r](\r/2))^{a-\nu})((a-\nu)/2))).\endalign$$
The left hand side is equal to
$$gr_0(\e^*((L\bul L')^{a+\r-\nu})((a+\r-\nu)/2))[\r](\r/2))$$
hence it is enough to show:
$$\align&\e^*((L\bul L')^{a+\r-\nu})((a+\r-\nu)/2))[\r](\r/2)\\&
=(\e^*L[\r](\r/2)\cir\e^*L'[\r](\r/2))^{a-\nu}((a-\nu)/2))\endalign$$
that is,
$$\e^*((L\bul L')^{a+\r-\nu})[\r]=(\e^*L[\r]\cir\e^*L'[\r])^{a-\nu},$$
or, after using 3.3(b):
$$(\e^*(L\bul L'))^{a+2\r-\nu}=(\e^*L\cir\e^*L')^{a+2\r-\nu}.$$
It remains to use that $\e^*(L\bul L')=\e^*L\cir\e^*L'$, see 3.22(a).

\subhead 3.25\endsubhead
In the setup of 3.14 let
$${}^\di\cy=\TT^2\bsl\{((x_0\UU,x_1\UU,x_2\UU,x_3\UU),g)\in\tcb^4\T G;
g\in x_3\UU x_0\i,g\in x_2\UU x_1\i\}$$
where $\TT^2$ acts freely by
$$(t_1,t_2):((x_0\UU,x_1\UU,x_2\UU,x_3\UU),g)\m
((x_0t_1\UU,x_1t_2\UU,x_2t_2\UU,x_3t_1\UU),g).$$
We define ${}^\di\eta:{}^\di\cy@>>>Z$ by
$$((x_0\UU,x_1\UU,x_2\UU,x_3\UU),g)\m\e(x_1\UU,x_2\UU).$$
We define $d:{}^\di\cy@>>>Z$ by
$$((x_0\UU,x_1\UU,x_2\UU,x_3\UU),g)\m\e(x_0\UU,x_3\UU).$$
We define $\fb'':\cd(Z)@>>>\cd(Z)$ and $\fb'':\cd_m(Z)@>>>\cd_m(Z)$ by
$$\fb''(L)=d_!({}^\di\eta)^*L.$$
From the definitions it is clear that
$$\fb'(L)=\e^*\fb''(L).\tag a$$
Using (a) we see that 3.14(a),(b),(c) imply the following statements.

(b) {\it If $L\in\cd^{\preceq}(Z)$ then $\fb''(L)\in\cd^{\preceq}Z$. If 
$L\in\cd^{\prec}(Z)$ then $\fb''(L)\in\cd^{\prec}Z$.}

(c) {\it If $L\in\cm^{\preceq}(Z)$ and $h>2\nu+2a$ then 
$(\fb''(L))^h\in\cm^{\prec}\tcb^2$.}
\nl
We define $\un{\fb''}:\cc^\boc_0(Z)@>>>\cc^\boc_0(Z)$ by
$$\un{\fb''}(L)=\un{gr_{2\nu+2a}((\fb''(L))^{2\nu+2a})}(\nu+a).$$ 
Using results in 3.3 we see that, if $L\in\cc^\boc_0Z$, then

(d) $\un\fb'(L)=\ti\e(\un{\fb''}(L))$.

\head 4. The monoidal category $\cc^\boc\tcb^2$ and its centre\endhead
\subhead 4.1\endsubhead
We consider the inclusion $\tco_1@>>>\tcb^2$ where $\tco_1=\{(x\UU,y\UU)\in\tcb^2,x\i y\in\BB\}$. 
Let $w\cdo\l\in W\fs$ be such that $w\cdo\l\preceq\boc$; let $i\in\ZZ$. From 2.12(a),(b) we deduce:
$$\sum_{i\in\ZZ}\rk(\ch^iL_\l^{\dw\sha}|_{\tco_1})v^i\text{ is }v^{|w|}p^\l_{1,w}\text{ if }w\in W_\l
\text{ and is $0$ if }w\n W_\l;\tag a$$
where $p^\l_{1,w}\in\ZZ[v\i]$ (as in 1.8) belongs to $v^{-a(w)}\ZZ[v\i]$ (see \cite{\HEC, \S14}) 
and $a(w)$ is the value at $w$ of the $a$-function of the Coxeter group $W_\l$ (so that $a(w)\ge a$); 
moreover,
$$\ch^iL_\l^{\dw\sha}|_{\tco_1}\text{ is a local system of pure weight }i.\tag b$$
From (a) it follows that
$$\sum_{i\in\ZZ}\rk(\ch^i(L_\l^{\dw\sha}[|w|])|_{\tco_1})v^i\text{ is in }v^{-a(w)}\ZZ[v\i]
\text{ if }w\in W_\l\text{ and is $0$ if }w\n W_\l$$
so that 
$$\ch^i(L_\l^{\dw\sha}[|w|])|_{\tco_1}) \text{ is $0$ if }i>-a.$$
Define $\d:\tcb@>>>\tcb^2$ by $x\UU\m(x\UU,x\UU)$. Then the image of $d$ is contained in $\tcb_1$
and we deduce

(c) $\ch^i(\d^*L_\l^{\dw\sha}[|w|])=0$ if $i>-a$.
\nl
We show:

(d) {\it If $L\in\cm^{\preceq}\tcb^2$ and $j>-a-\r$ then $(\d^*L)^j=0$.}
\nl
We can assume that $L=\LL_\l^{\dw}$ with $w\cdo\l$ as above. It is enough 
to show that for any $k$ we have $(\ch^k(\d^*L)[-k])^j=0$ that is
$$(\ch^k(\d^*(L_\l^{\dw\sha}[|w|+\nu+2\r]))[\nu+\r])^{j-k-\nu-\r}=0.$$
Now $\ch^k(\d^*(L_\l^{\dw\sha}[|w|+\nu+2\r]))$ is a local system on $\tcb$ hence
$$\ch^k(\d^*(L_\l^{\dw\sha}[|w|+\nu+2\r]))[\nu+\r]$$ 
is a perverse sheaf on $\tcb$ so that
we can assume that $j-k-\nu-\r=0$. Thus it is enough to show that
$$\ch^{j-\nu-\r}(\d^*(L_\l^{\dw\sha}[|w|+\nu+2\r]))=0$$ 
or that $\ch^{j+\r}(\d^*(L_\l^{\dw\sha}[|w|))=0$. This is indeed true by (c).

We show:

(e)  {\it If $L\in\cm^{\preceq}_m\tcb^2$ is pure of weight $0$ and $j\in\ZZ$ then 
$(\d^*L)^j$ is pure of weight $j$.}
\nl
We can assume that $L=\LL_\l^{\dw}$ with $\l,w$ as in (c). It is enough to 
prove that for any $k$, $(\ch^k(\d^*L)[-k])^j$ is pure of weight $j$ that is,
$$(\ch^k(\d^*(L_\l^{\dw\sha}[|w|+\nu+2\r]))[\nu+\r])^{j-k-\nu-\r}((|w|+\nu+2\r)/2)$$
is pure of weight $j$. As in the proof of (d) we can assume that $j-k-\nu-\r=0$. 
Thus it is enough to show that
$$\ch^{j-\nu-\r}(\d^*(L_\l^{\dw\sha}[|w|+\nu+2\r]))[\nu+\r]((|w|+\nu+2\r)/2)$$
is pure of weight $j$ or that $\ch^{j+|w|+\r}(\d^*(L_\l^{\dw\sha}))\la\nu+\r\ra((|w|+\r)/2)$ 
is pure of weight $j$. This follows from (b).

\mpb

For $w\cdo\l\in\boc$ we have as in the proof of (e):
$$(\d^*\LL_\l^{\dw})^{-a-\r}=\ch^{-a+|w|}(\d^*(L_\l^{\dw\sha}))\la\nu+\r\ra((|w|+\r)/2).$$
We set
$$\b_{w\cdo\l}=\ch^{-a+|w|}(\d^*(L_\l^{\dw\sha}))((-a+|w|)/2)$$
This is a $G$-equivariant local system on $\tcb$ hence can be identified with a $\bbq$-vector space
which by (b) is pure of weight $0$; we have

We show:

(f) {\it $\dim\b_{w\cdo\l}$ is $1$ if $w\cdo\l\in\DD_\boc$ and is $0$ if $w\cdo\l\n\DD_\boc$. We have}
$$(\d^*\LL_\l^{\dw})^{-a-\r}=\b_{w\cdo\l}\la\nu+\r\ra((a+\r)/2).$$ 
By (a), $\dim\b_{w\cdo\l}$ is $0$ if $w\n W_\l$ while if $w\in W_\l$ it is
equal to the coefficient of $v^{-a}$ in $p^\l_{1,w}$, which by \cite{\HEC, 14.2, P5} is $1$ 
if $w$ is a distinguished involution of $W_\l$ and is $0$ otherwise. This proves (f).

\subhead 4.2\endsubhead
Let $\p':\tcb@>>>\pp$ be the obvious map. We show:

(a) {\it Assume that $L\in\cm_m(\tcb)$ is $G$-equivariant so that $L=V\ot\bbq\la\nu+\r\ra$
where $V$ is a mixed vector space. If $j>\nu+\r$ then $(\p'_!L)^j=0$.}
\nl
We have $\ch^j(p'_!L)=V\ot H_c^{j+\nu+\r}(\tcb,\bbq)((\nu+\r)/2)$. This is zero if $j+\nu+\r>2\nu+2\r$
since $\tcb$ is irreducible of dimension $\nu+\r$. 

We show:

(b) {\it If $L\in\cm^{\preceq}_m\tcb^2$ and $j>\nu-a$ then $(\p'_!\d^*L)^j=0$ (with $\d$ as in 4.1).
Moreover we have canonically $(\p'_!\d^*L)^{\nu-a}=(\p'_!((\d^*L)^{-a-\r}))^{\nu+\r}$.}
\nl
The proof has much in common with that of \cite{\CONV, 8.2}.

Let $\XX=\d^*L$. For any $i$ we have a distinguished triangle $(\t_{<i}\XX,\t_{\le i}\XX,\XX^i[-i])$ 
where we
write $\t_{<i},\t_{\le i}$ for what in \cite{\BBD} is denoted by ${}^p\t_{<i},{}^p\t_{\le i}$.
We deduce a distinguished triangle $(\p'_!(\t_{<i}\XX),\p'_!(\t_{\le i}\XX),\p'_!(\XX^i)[-i])$
hence an exact sequence
$$\align&(\p'_!(\XX^i))^{h-1}@>>>(\p'_!(\t_{<i}\XX))^{i+h}@>>>(\p'_!(\t_{\le i}\XX))^{i+h}@>>>
(\p'_!(\XX^i))^h\\&@>>>(\p'_!(\t_{<i}\XX))^{i+h+1}.\tag c\endalign$$
We show by induction on $i$ that 

(d) $(\p'_!(\t_{\le i}\XX))^j=0$ for $j>\nu-a$.
\nl
For sufficiently negative $i$, (d) is obvious. Thus we can assume that (d) is known when $i$ is replaced by 
$i-1$. From (c) with $h=j-i$ we see using the induction hypothesis that we have an exact sequence 
$0@>>>(\p'_!(\t_{\le i}\XX))^j@>>>(\p'_!(\XX^i))^{j-i}$. It is then enough to show that 
$(\p'_!(\XX^i))^{j-i}=0$. If $i>-a-\r$ then $\XX^i=0$ by 4.1(d). Thus we can
assume that $i\le-a-\r$ so that $j-i>\nu+\r$ and the equality $(\p'_!(\XX^i))^{j-i}=0$
follows from (a) with $L=\XX^i$. This proves (d). 

In particular the first assertion of (b) holds.
We now take $h=\nu-a-i$ in (c). Assuming that $i>-a-\r$ we obtain (using 4.1(d))
$$(\p'_!(\t_{<i}\XX))^{\nu-a}@>\si>>(\p'_!(\t_{\le i}\XX))^{\nu-a}.$$ 
Hence
$$(\p'_!(\t_{\le-a-\r}\XX))^{\nu-a}@>\si>>(\p'_!\XX)^{\nu-a}.$$

We show by induction on $i$:

(e) $(\p'_!(\t_{\le i}\XX))^j=0$ if $i<-a-\r$, $j=\nu-a$.
\nl
For sufficiently negative $i$, (e) is obvious. Thus we can assume that (e) is known when $i$ is replaced by 
$i-1$. From (c) with $h=j-i$ we see using the induction hypothesis that we have an exact sequence 
$0@>>>(\p'_!(\t_{\le i}\XX))^j@>>>(\p'_!(\XX^i))^{j-i}$. It is then enough to show that 
$(\p'_!(\XX^i))^{j-i}=0$. We have $j-i>\nu+\r$ and the equality $(\p'_!(\XX^i))^{j-i}=0$
follows from (a) with $L=\XX^i$. This proves (e). 

From (e) we see that $(\p_!(\t_{\le-a-\r-1}\XX))^{\nu-a}=0$. 
From (c) with $i=-a-\r$, $h=\nu+\r$ we deduce an exact sequence
$$0@>>>(\p'_!\XX)^{\nu-a}@>>>(\p'_!(\XX^{-a-\r}))^{\nu+\r}@>>>0.$$
This completes the proof of (b).

\mpb

We show:

(f) {\it If $L\in\cc^\boc_0\tcb^2$ then $(\p'_!((\d^*L)^{-a-\r}))^{\nu+\r}$ and $(\p'_!\d^*L)^{\nu-a}$ are 
pure of weight $\nu-a$.}
\nl
We can assume that $L=\LL_\l^{\dw}$ where $w\cdo\l\in\boc$. 
Using 4.1(f) we have 
$$\align&(\p'_!((\d^*L)^{-a-\r}))^{\nu+\r}((-a+\nu)/2)=
(\p'_!(\b_{w\cdo\l}\la\nu+\r\ra)((a+\r)/2))^{\nu+\r}((-a+\nu)/2)\\&=
\b_{w\cdo\l}\ot(\p'_!\bbq\la\nu+\r\ra)^{\nu+\r}((\nu+\r)/2)=
\b_{w\cdo\l}\ot(\p'_!\bbq)^{2\nu+2\r}(\nu+\r)=\b_{w\cdo\l}.\endalign$$
Since $b_{w\cdo\l}$ is pure of weight $0$, we see that $(\p'_!((\d^*L)^{-a-\r}))^{\nu+\r}$ is pure of weight
$\nu-a$. Using (b) we deduce that $(\p'_!\d^*L)^{\nu-a}$ is pure of weight $\nu-a$. 

\subhead 4.3\endsubhead
We set
$$\bold1'=\op_{d\cdo\l\in\DD_\boc}\b^*_{d\cdo\l}\ot\LL_\l^{\dot d}\in\cc^\boc_0\tcb^2.$$
Here $\b^*_{d\cdo\l}$ is the vector space dual to $\b_{d\cdo\l}$, see 4.1(f).
For $L\in\cc^\boc\tcb^2$ we show
$$\Hom_{\cc^\boc\tcb^2}(\bold1',L)=(\p'_!((\d^*L)^{-a-\r}))^{\nu+\r}((-a+\nu)/2).\tag a$$ 
We can assume that $L=\LL_\l^{\dw}$ where $w\cdo\l\in\boc$.
Assume first that $w\cdo\l\in\DD_\boc$. Then 
$\Hom_{\cc^\boc\tcb^2}(\bold1',L)=\b_{w\cdo\l}$ and (as in the proof of 4.2(f)):
$$(\p'_!((\d^*L)^{-a-\r}))^{\nu+\r}((-a+\nu)/2)=\b_{w\cdo\l}.$$
Thus (a) holds in this case. Next assume that $w\cdo\l\n\DD_\boc$. Then
both sides of (a) are zero (we use 4.1(f)).  

\subhead 4.4\endsubhead
For $L\in\cd_m(\tcb^2)$ we set $L^\da=\ti\fh^*L$ where $\ti\fh:\tcb^2@>>>\tcb^2$ 
is as in 2.1. Let $\p'':\tcb^2@>>>\pp$ be the obvious map. For $L,L'$ in $\cd_m(\tcb^2)$
we have from the definitions
$$\p'_!\d^*(L\circ L')=\p''_!(L\ot L'{}^\da).\tag a$$
We show:

(b) {\it If $L\in\cc^\boc_0\tcb^2$ then $\fD(L^\da)\in\cc^\boc_0\tcb^2$.}
\nl
We can assume that $L=\LL_\l^{\dw}$ where $w\cdo\l\in\boc$.
Using 2.2(a) and the definitions  we have 
$$\fD(L^\da)=\fD(\LL_{w(\l\i)}^{\dw\i})=\LL_{w(\l)}^{\dw\i}.\tag c$$
It remains to use that $w\i\cdo w(\l)\in\boc$, by property Q10 in 1.9 for $\HH$.

\subhead 4.5\endsubhead
We show:

(a) {\it For $L,L'$ in $\cc^\boc\tcb^2$ we have canonically}
$$\Hom_{\cc^\boc\tcb^2}(\bold1',L\un\circ L')=\Hom_{\cc^\boc\tcb^2}(\fD(L'{}^\da),L).$$
We can assume that $L=\LL_\l^{\dw}$, $L'=\LL_{\l'}^{\dw'}$ where $w\cdo\l\in\boc$,
$w'\cdo\l'\in\boc$. Using 4.3(a) and 4.2(b), 4.2(f), we have
$$\align&\Hom_{\cc^\boc\tcb^2}(\bold1',L\un\circ L')
=(\p'_!((\d^*(L\un\circ L'))^{-a-\r}))^{\nu+\r}((-a+\nu)/2)\\&
=(\p'_!\d^*(L\un\circ L'))^{\nu-a}((-a+\nu)/2)
=(\p'_!\d^*(L\un\circ L'))^{\{\nu-a\}}.\tag b\endalign$$                     
Applying \cite{\CONV, 8.2} with $\Ph:\cd^{\preceq}_mZ@>>>\cd_\pp$, $\tL\m\p'_!\d^*\tL$,
$c=\nu-a$, $c'=-\nu+a$ (see 4.2(b)) we see that we have canonically
$$(\p'_!\d^*(L\un\circ L'))^{\{\nu-a\}}\sub(\p'_!\d^*(L\bul L'))^{\{0\}}.\tag c$$
From \cite{\CSII, 7.4} we see that we have canonically
$$(\p''_!(L\ot L'{}^\da))^0=(\p''_!(L\ot L'{}^\da))^{\{0\}} 
=\Hom_{\cd(\tcb^2)}(\fD(L'{}^\da),L).\tag d$$
By 4.4(a) we have $\p'_!\d^*(L\bul L')=\p''_!(L\ot L'{}^\da)$. Hence by combining 
(b),(c),(d) we have
$$\Hom_{\cc^\boc\tcb^2}(\bold1',L\un\circ L')\sub\Hom_{\cc^\boc\tcb^2}(\fD(L'{}^\da),L).\tag e$$
The dimension of the left hand side of (e) is the sum over $z\cdo\l_1\in\DD_\boc$ of the coefficient of
$t_{z\cdo\l_1}$ in $t_{w\cdo\l}t_{w'\cdo\l'}\in\HH^\iy$ and by the properties Q1,Q2,Q4 (in 1.9) of $\HH$,
this sum is equal to $1$ if $w'\cdo\l'=w\i\cdo w(\l)$ and is $0$ if $w'\cdo\l'\ne w\i\cdo w(\l)$; hence 
it is equal to the  dimension of the right hand side of (e). It follows that (e) is an equality and 
(a) follows.

\subhead 4.6\endsubhead
The bifunctor $\cc^\boc_0\tcb^2\T\cc^\boc_0\tcb^2@>>>\cc^\boc_0\tcb^2$, $L,L'\m L\un{\circ}L'$ in 2.24 
gives rise to a bifunctor $\cc^\boc\tcb^2\T\cc^\boc\tcb^2@>>>\cc^\boc\tcb^2$ denoted again by 
$L,L'\m L\un{\cir}L'$ as follows. Let $L\in\cc^\boc\tcb^2$, $L'\in\cc^\boc\tcb^2$; we choose mixed 
structures of pure weight $0$ on $L,L'$, we define $L\un{\cir}L'$ as in 2.24 in terms of these mixed 
structures and we then disregard the mixed structure on $L\un{\cir}L'$. The resulting object of 
$\cc^\boc\tcb^2$ is denoted again by $L\un{\cir}L'$; it is independent of the choices made.

Similarly the bifunctor $\cc^\boc_0Z\T\cc^\boc_0Z@>>>\cc^\boc_0Z$, $L,L'\m L\un{\bul}L'$ in 3.24 
gives rise to a bifunctor $\cc^\boc Z\T\cc^\boc Z@>>>\cc^\boc Z$ denoted again by 
$L,L'\m L\un{\bul}L'$.

The operation $L\un{\bul}L'$ (resp. $L\un{\circ}L'$) makes $\cc^\boc Z$ (resp. $\cc^\boc\tcb^2$) into
a monoidal abelian category (see 2.24, 3.24).

The following result can be deduced from 2.22(c).

(a) {\it Let $w_i\cdo\l_i\in\boc$, $i=1,2$. In $\cc^\boc\tcb^2$ we have
$$\LL_{\l_1}^{\dw_1}\un\cir\LL_{\l_2}^{\dw_2}\cong\op_{w\cdo\l\in\boc}(\LL_\l^{\dw})^{\op f(w\cdo\l)}
$$
where $f(w\cdo\l)\in\NN$ are given by}
$$t_{w_1\cdo\l_1}t_{w_2\cdo\l_2}=\sum_{w\cdo\l\in\boc}f(w\cdo\l)t_{w\cdo\l}\in\HH^\iy_\boc.$$
It follows that:

(b) {\it in (a) we have $f(w\cdo\l)=0$ unless $\l=\l_2$, $w(\l)=w_1(\l_1)$, $\l_1=w_2(\l_2)$.}
\nl
To see this, we use that, setting for any $\l\in\fo$:

$1_\l^\boc=\sum_dt_{d\cdo\l}\in\HH^\iy_\boc$ (sum over all distinguished involutions $d$ of $W_\l$),
\nl
we have $1_\l^\boc1_{\l'}^\boc=\d_{\l,\l'}1_\l^\boc$ for $\l,\l'$ in $\fo$ and

$t_{w\cdo\l}=t_{w\cdo\l}1_\l^\boc=1_{w(l)}^\boc t_{w\cdo\l}$ for any $w\cdo\l\in\boc$.

\mpb

For any $\l_1,\l_2$ in $\fo$ let $\cc^\boc_{\l_1,\l_2}\tcb^2$ be the subcategory of $\cc^\boc\tcb^2$ 
consisting of objects which are direct sums of objects of the form $\LL_{\l_2}^{\dw}$ for various 
$w\in W$ such that $w(\l_2)=\l_1$ and $w\cdo\l_2\in\boc$. Clearly, any object $L\in\cc^\boc\tcb^2$ is
canonically of the form 

(c) {\it $\op_{\l_1,\l_2\in\fo}L_{\l_1,\l_2}$ where $L_{\l_1,\l_2}\in\cc^\boc_{\l_1,\l_2}$.}
\nl
From (b) we see that 

(d) {\it if $\l_1,\l_2,\l_3,\l_4$ are elements of $\fo$ and $L\in\cc^\boc_{\l_1,\l_2}$, 
$L'\in \cc^\boc_{\l_3,\l_4}$ then $L\un\cir L'\in\cc^\boc_{\l_1,\l_4}$; moreover, $L\un\cir L'=0$ 
unless $\l_2=\l_3$.}

\subhead 4.7\endsubhead
We set 
$$\bold1=\op_{d\cdo\l\in\DD_\boc}\b_{d\cdo\l}\ot\LL_\l^{(\dot d)\i}\in\cc^\boc_0\tcb^2.\tag a$$
Here $\b_{d\cdo\l}$ is as in 4.1(f). 

Let $y_2\cdo\l_2\in\boc$, $y_3\cdo\l_3\in\boc$. From 3.20(a) we have for any
$d\cdo\l\in\DD_\boc$:
$$\Hom_{\cc^\boc\tcb^2}(\LL_{y_2(\l_2)}^{\dy_2\i},\LL_\l^{(\dot d)\i}\un{\cir}\LL_{\l_3}^{\dy_3})
=\Hom_{\cc^\boc\tcb^2}(\LL_\l^{\dot d},\LL_{\l_3}^{\dy_3}\un{\cir}\LL_{\l_2}^{\dy_2}).$$
It follows that
$$\Hom_{\cc^\boc\tcb^2}(\LL_{y_2(\l_2)}^{\dy_2\i},\bold1\un{\cir}\LL_{\l_3}^{\dy_3})
=\Hom_{\cc^\boc\tcb^2}(\bold1',\LL_{\l_3}^{\dy_3}\un{\cir}\LL_{\l_2}^{\dy_2}).\tag b$$
From 4.5(a) we have:
$$\Hom_{\cc^\boc\tcb^2}(\bold1',\LL_{\l_3}^{\dy_3}\un{\circ}\LL_{\l_2}^{\dy_2})
=\Hom_{\cc^\boc\tcb^2}(\LL_{y_2(\l_2)}^{\dy_2\i},\LL_{\l_3}^{\dy_3}).$$
(We have used 2.2(a) and the equality
$$\fD(\LL_\l^\o)=\LL_{\l\i}^\o\tag c$$ 
for any $w\cdo\l\in W\fs$ and $\o\in\k\i_q(w)$.)
Using this, (b) becomes
$$\Hom_{\cc^\boc\tcb^2}(\LL_{y_2(\l_2)}^{\dy_2\i},\bold1\un{\cir}\LL_{\l_3}^{\dy_3})
 =\Hom_{\cc^\boc\tcb^2}(\LL_{y_2(\l_2)}^{\dy_2\i},\LL_{\l_3}^{\dy_3}).$$
Since this holds for any $y_2\cdo\l_2\in\boc$, we see that we have canonically
$$\bold1\un{\cir}\LL_{\l_3}^{\dy_3}=\LL_{\l_3}^{\dy_3}.$$
Since this holds for any $y_3\cdo\l_3\in\boc$, we see that we have canonically
$$\bold1\un{\cir}L=L\text{ for any }L\in\cc^\boc_0\tcb^2\tag d$$
for any $L$ in $\cc^\boc\tcb^2$. Now $\cc^\boc\tcb^2@>>>\cc^{\ti\boc}\tcb^2$, 
$L\m L^\da$, satisfies
$$(L\un{\cir}L')^\da=L'{}^\da\un{\cir}L^\da\tag e$$
for any $L,L'$ in $\cc^\boc\tcb^2$. Applying $L\m \fD(L^\da)$ to (d)
and using (e) and 2.25(a)  we get
$$L\un{\cir}\fD(\bold1^\da)=L\text{ for any }L\in\cc^\boc_0\tcb^2.\tag f$$
From (d),(f) we deduce that we have canonically $\bold1=\bold1\un{\cir}\fD(\bold1^\da)=\fD(\bold1^\da)$.
Using 4.4(c) we see that $\fD(\bold1'{}^\da)=\bold1$ hence $\fD(\bold1^\da)=\bold1'$. We see that 

(g) {\it $\bold1=\bold1'=\fD(\bold1^\da)$ is a unit object of the monoidal category $\cc^\boc\tcb^2$.}

\subhead 4.8\endsubhead
For $L\in\cc^\boc\tcb^2$ let $L^*=\fD(L^\da)$. We say that $L^*$ is the dual of $L$.
Note that $L\m L^*$ is a contravariant functor $\cc^\boc\tcb^2@>>>\cc^\boc\tcb^2$ and that $L^{**}=L$. 
We show how $L\m L^*$ gives a rigid structure on $\cc^\boc\tcb^2$.

We have the following special case of 4.5(a) (we use that $\bold1=\bold1'$, see 4.7(g)):
$$\Hom_{\cc^\boc\tcb^2}(\bold1,L\un\cir\fD(L^\da))=\Hom_{\cc^\boc\tcb^2}(L,L)\tag a$$ 
for any $L$ in $\cc^\boc_0\tcb^2$. Let $\x_L\in\Hom_{\cc^\boc\tcb^2}(\bold1,L\un\cir\fD(L^\da))$
be the element corresponding under (a) to the identity homomorphism in $\Hom_{\cc^\boc\tcb^2}(L,L)$. 
Using (e) and 2.25(a) we have
$$\align&\Hom_{\cc^\boc\tcb^2}(L\un\cir\fD(L^\da),\bold1)
=\Hom_{\cc^\boc\tcb^2}(\fD(\bold1^\da), \fD((L\un\cir\fD(L^\da))^\da))\\&
=\Hom_{\cc^\boc\tcb^2}(\bold1,\fD(L^\da)\un\cir L).\tag b\endalign$$
Let $\x'_L\in\Hom_{\cc^\boc\tcb^2}(L\un\cir\fD(L^\da)\bold1)$ be the element corresponding under (b) to the 
element $\x_{\fD(L^\da)}\in\Hom_{\cc^\boc\tcb^2}(\bold1,\fD(L^\da)\un\cir L)$.
The elements $\x_L,\x'_L$ define the rigid structure on $\cc^\boc\tcb^2$.

\subhead 4.9\endsubhead
Let $\cz^\boc$ be the centre (in the sense of Joyal and Street \cite{\JS}, Majid \cite{\MA} and Drinfeld)
of the monoidal abelian category $\cc^\boc\tcb^2$. By a general result on semisimple rigid monoidal
categories in \cite{\ENO, Proposition 5.4}, for any $L\in\cc^\boc\tcb^2$ one can define directly a central
structure on the object 
$$I(L):=\op_{y\cdo\l\in\boc}\LL_\l^{\dy}\un\cir L\un\cir\LL_{y(\l)}^{\dy\i}$$
of $\cc^\boc\tcb^2$ such that, denoting by $\ov{I(L)}$ the corresponding object of $\cz^\boc$, we have
canonically
$$\Hom_{\cc^\boc\tcb^2}(L,L')=\Hom_{\cz^\boc}(\ov{I(L)},L')\tag a$$
for any $L'\in\cz^\boc$. (We use that for $y\cdo\l\in\boc$, the dual of the simple object
$\LL_\l^{\dy}$ is $\LL_{y(\l)}^{\dy\i}$.) The central structure on $I(L)$ can be described as
follows: for any $X\in\cc^\boc\tcb^2$ we have canonically
$$\align&X\un\cir I(L)=\op_{y\cdo\l\in\boc}X\un\cir\LL_\l^{\dy}\un\cir L\un\cir\LL_{y(\l)}^{\dy\i}\\&
=\op_{y\cdo\l\in\boc,z\cdo\l'\in\boc}\Hom_{\cc^\boc\tcb^2}(\LL_{\l'}^{\dz},X\un\cir\LL_\l^{\dy})\ot
\LL_{\l'}^{\dz}\un\cir L\un\cir\LL_{y(\l)}^{\dy\i}\\&
=\op_{y\cdo\l\in\boc,z\cdo\l'\in\boc}\Hom_{\cc^\boc\tcb^2}(\LL_{y(\l)}^{\dy\i},\LL_{z(\l')}^{\dz\i}\ot X)\ot
\LL_{\l'}^{\dz}\un\cir L\un\cir\LL_{y(\l)}^{\dy\i}\\&
=\op_{z\cdo\l'\in\boc}\LL_{\l'}^{\dz}\un\cir L\un\cir \LL_{z(\l')}^{\dz\i}\ot X=I(L)\un\cir X.
\endalign$$
(The third equality uses 3.20(a).)

We show:

(b) {\it If $z\cdo\l\in\boc$ and $I(\LL_\l^{\dz})\ne0$ then $z\cdo\l\in\overset\smile\to{W\fs}$.}
\nl
For some $y\cdo\l'\in\boc$ we have $\LL_{\l'}^{\dy}\un\cir\LL_\l^{\dz}\ne0$ (hence $\l'=z(l)$) and
$\LL_\l^{\dz}\un\cir\LL_{y(\l')}^{\dy\i}\ne0$ (hence $\l=\l'$). It follows that $z(\l)=\l$ and (b) is
proved.

\subhead 4.10\endsubhead
By 3.13(d), for $z\cdo\l\in\overset\smile\to\boc$ we have canonically 
$$\un{\fb}(\Bbb L_\l^{\dz})=I(\LL_\l^{\dz})\tag a$$ 
as objects of $\cc^\boc\tcb^2$.
Here $\un\fb:\cc^\boc_0Z@>>>\cc^\boc_0\tcb^2$ in 3.13 is viewed as a functor
$\un\fb:\cc^\boc Z@>>>\cc^\boc\tcb^2$ as in 4.6.
Now $I(\LL_\l^{\dz})$ has a natural central structure (by 4.9) and 
$\un\fb(\Bbb L_\l^{\dz})$ has a natural central structure (by 3.14(j)). By 3.21(b),

(b) {\it these two central structures are compatible with the identification (a).}
\nl
In view of (a),(b) we can reformulate 4.9(a) as follows.

\proclaim{Theorem 4.11}For any $z\cdo\l\in\overset\smile\to\boc$, $L'\in\cz^\boc$ we have canonically
$$\Hom_{\cc^\boc\tcb^2}(\LL_\l^{\dz},L')=\Hom_{\cz^\boc}(\ov{\un{\fb}(\Bbb L_\l^{\dz})},L')\tag a$$
where $\ov{\un{\fb}(\Bbb L_\l^{\dz})}$ is $\un{\fb}(\Bbb L_\l^{\dz})$ viewed as an object of $\cz^\boc$ 
with the central structure given by 3.14(j).
\endproclaim

\subhead 4.12\endsubhead
We will state some variants (a)-(j) of results in 4.1-4.5 which will be needed in 
Section 6. Let $\d_0:\cb@>>>Z$ be map $B\m(B,B,U_B)$. Let $\p'_0:\cb@>>>\pp$ be the obvious map. 

(a) {\it Let $\l\in\fs$ and let $w\in W'_\l$ be such that $w\cdo\l\preceq\boc$; let 
$i\in\ZZ$. If $i>a$ we have $\ch^i(\d_0^*\cl_\l^{\dw\sha}[|w|])=0$.}
\nl
This can be deduced from 4.1(c) using that $e^*\ch^i(\d_0^*\cl_\l^{\dw\sha}[|w|])=
\ch^i(\d^*L_\l^{\dw\sha}[|w|])$ where $e:\tcb@>>>\cb$ is the map $x\UU\m x\BB x\i$.

(b)  {\it If $L\in\cm^{\preceq}Z$ and $j>-a-\r$ then $(\d_0^*L)^j=0$.}
\nl
We argue as in the proof of 4.1(d). We can assume that $L=\Bbb L_\l^{\dw}$ with $w\cdo\l$ as in (a). 
It is enough to show that for any $k$ we have $(\ch^k(\d_0^*L)[-k])^j=0$ that is
$(\ch^k(\d_0^*(\cl_\l^{\dw\sha}[|w|+\nu+\r]))[\nu])^{j-k-\nu}=0$.
Now $\ch^k(\d_0^*(\cl_\l^{\dw\sha}[|w|+\nu+\r]))$ is a local system on $\cb$ hence
$\ch^k(\d_0^*(\cl_\l^{\dw\sha}[|w|+\nu+\r]))[\nu]$ is a perverse sheaf on $\cb$ so that
we can assume that $j-k-\nu=0$. Thus it is enough to show that
$\ch^{j-\nu}(\d_0^*(\cl_\l^{\dw\sha}[|w|+\nu+\r]))=0$ or that
$\ch^{j+\r}(\d^*(L_\l^{\dw\sha}[|w|))=0$. This is indeed true by (a).

(c)  {\it If $L\in\cm^{\preceq}_mZ$ is pure of weight $0$ and $j\in\ZZ$ then 
$(\d_0^*L)^j$ is pure of weight $j$.}
\nl
We argue as in the proof of 4.1(e).
We can assume that $L=\Bbb L_\l^{\dw}$ with $\l,w$ as in (a). It is enough to 
prove that for any $k$, $(\ch^k(\d_0^*L)[-k])^j$ is pure of weight $j$ that is,
$(\ch^k(\d_0^*(\cl_\l^{\dw\sha}[|w|+\nu+\r]))[\nu])^{j-k-\nu}((|w|+\nu+\r)/2)$ 
is pure of weight $j$. As in the proof of (b) we can assume that $j-k-\nu=0$. 
Thus it is enough to show that
$\ch^{j-\nu}(\d_0^*(\cl_\l^{\dw\sha}[|w|+\nu+\r]))[\nu]((|w|+\nu+\r)/2)$ 
is pure of weight $j$ or that $\ch^{j+|w|+\r}(\d_0^*(\cl_\l^{\dw\sha}))\la\nu\ra((|w|+\r)/2)$ 
is pure of weight $j$. This follows from 4.1(b).

(d) {\it Let $\l\in\fs$, $w\in W'_\l$ be such that $w\cdo\l\in\boc$. If $w\cdo\l\in\DD_\boc$
(see 1.12) then $(\d_0^*\Bbb L_\l^{\dw})^{-a-\r}=\b_{w\cdo\l}\la\nu\ra((a+\r)/2)$ where
$\b_{w\cdo\l}$ is as in 4.1(f). If $w\cdo\l\n\DD_\boc$  then $(\d_0^*\Bbb L_\l^{\dw})^{-a-\r}=0$.}
\nl
As in the proof of (c) we have
$$(\d_0^*\Bbb L_\l^{\dw})^{-a-\r}=\ch^{-a+|w|}(\d_0^*(\cl_\l^{\dw\sha}))\la\nu+\r\ra((|w|+\r)/2).$$
Setting $\b_{w\cdo\l;0}=\ch^{-a+|w|}(\d_0^*(\cl_\l^{\dw\sha}))\la\nu+\r\ra((-a+|w|)/2)$ we have
$(\d_0^*\Bbb L_\l^{\dw})^{-a-\r}=\b_{w\cdo\l;0}\la\nu\ra((a+\r)/2)$ where $\b_{w\cdo\l;0}$ is a mixed 
vector space. If $e:\tcb@>>>\cb$ is the obvious map, we have
$$\align&\b_{w\cdo\l;0}=e^*(\ch^{-a+|w|}(\d_0^*(\cl_\l^{\dw\sha}))\la\nu+\r\ra((-a+|w|)/2)\\&=
\ch^{-a+|w|}(\d^*(L_\l^{\dw\sha}))\la\nu+\r\ra((-a+|w|)/2)=\b_{w\cdo\l}\endalign$$ 
where $\b_{w\cdo\l}$ is as in 4.1(f). Hence the result follows from 4.1(f).

(e) {\it Assume that $L\in\cm_m(\cb)$ is $G$-equivariant so that $L=V\ot\bbq\la\nu\ra$
where $V$ is a mixed vector space. If $j>\nu$ then $(\p'_{0!}L)^j=0$.}
\nl
We argue as in the proof of 4.2(a). We have 
$\ch^j(p'_{0!}L)=V\ot H_c^{j+\nu}(\tcb,\bbq)(\nu/2)$. This is zero if $j+\nu>2\nu$
since $\cb$ is irreducible of dimension $\nu+\r$. 

(f) {\it If $L\in\cm^{\preceq}_mZ$ and $j>\nu-a-\r$ then $(\p'_{0!}\d_0^*L)^!=0$.
Moreover we have canonically $(\p'_{0!}\d_0^*L)^{\nu-a-\r}=(\p'_{0!}((\d_0^*L)^{-a-\r}))^\nu$.}
\nl
The proof is almost identical to that of 4.2(b), using (b),(e) instead of 4.1(d), 4.2(a).

(g) {\it Let $L\in\cc^\boc_0Z$. Then
$(\p'_{0!}\d_0^*L)^{\nu-a-\r}$ and $(\p'_{0!}((\d_0^*L)^{-a-\r}))^\nu$ are pure of weight $\nu-a-\r$.}
\nl
We argue as in the proof of 4.2(f). 
We can assume that $L=\Bbb L_\l^{\dw}$ where $w\cdo\l\in\overset\smile\to\boc$. Using (d) we have 
$$\align&(\p'_{0!}((\d_0^*L)^{-a-\r}))^\nu((-a-\r+\nu)/2)=
(\p'_{0!}(\b_{w\cdo\l}\la\nu\ra)((a+\r)/2))^\nu((-a-\r+\nu)/2)\\&=
\b_{w\cdo\l}\ot(\p'_{0!}\bbq\la\nu\ra)^\nu(\nu/2)=
\b_{w\cdo\l}\ot(\p'_{0!}\bbq)^{2\nu}(\nu)=\b_{w\cdo\l}.\endalign$$
Since $\b_{w\cdo\l}$ is pure of weight $0$, we see that $(\p'_{0!}((\d_0^*L)^{-a-\r}))^\nu$ is pure of weight
$\nu-a-\r$. Using (f) we deduce that $(\p'_{0!}\d_0^*L)^{\nu-a-\r}$ is pure of weight $\nu-a-\r$.

We set
$$\bold1'_0=\op_{d\cdo\l\in\DD_\boc}\b^*_{d\cdo\l}\ot\Bbb L_\l^{\dot d}\in\cc^\boc_0Z.$$

(h) {\it For $L\in\cc^\boc Z$ we have canonically}
$$\Hom_{\cc^\boc Z}(\bold1'_0,L)=(\p'_{0!}((\d_0^*L)^{-a-\r}))^\nu((-a-\r+\nu)/2).$$ 
\nl
The proof is similar to that of 4.3(a); it uses (d) (instead of 4.1(f)) and the proof of (g).

\mpb

For $L\in\cd_m(Z)$ we set $L^\da=\fh^*L$ where $\fh:Z@>>>Z$ is as in 3.2.

(i) {\it If $L\in\cc^\boc_0Z$ then $\fD(L^\da)\in\cc^\boc_0Z$.}
\nl
This can be deduced from 4.4(b).

(j) {\it For $L,L'$ in $\cc^\boc Z$ we have canonically}
$$\Hom_{\cc^\boc Z}(\bold1'_0,L\un\bul L')=\Hom_{\cc^\boc Z}(\fD(L'{}^\da),L).$$
This can be proved by the same method as 4.5(b) or it can be deduced from 4.5(b) using the fully 
faithfulness of $\ti\e:\cc^\boc_0Z@>>>\cc^\boc_0\tcb^2$, the equality $\ti\e\bold1'_0=\bold1$ and 3.22(a).

\subhead 4.13\endsubhead
Let $\l\in\fo$. Using the decomposition 4.6(c) of any object of $\cc^\boc\tcb^2$ we see that 
$\cc^\boc\tcb^2$ can be viewed as the category of ``matrices'' with entries in the abelian category 
$\cc^\boc_{\l,\l}\tcb^2$ (see 4.6). (This is a category version of the isomorphism 
$\Psi:\HH_{\fo}@>>>\EE$ in 1.11(v).) Using this and a result of M\"uger \cite{\MU} it follows that
$\cz^\boc_\l$ is equivalent to the categorical centre of the abelian category 
$\cc^\boc_{\l,\l}\tcb^2$ with the monoidal structure induced by $\un\cir$ (see 4.6(d)).

\head 5. Truncated induction, truncated restriction, truncated convolution on $G$\endhead
\subhead 5.1\endsubhead
Let $\dZ=\{(B,B',g)\in\cb\T\cb\T G;gBg\i=B'\}$. We have a diagram 
$$Z@<f<<\dZ@>\p>>G\tag a$$ 
where $f(B,B',g)=(B,B',gU_B)$, $\p(B,B',g)=g$. We define
$\c:\cd(Z)@>>>\cd(G)$ and $\c:\cd_m(Z)@>>>\cd_m(G)$ by 
$$\c(L)=\p_!f^*L.$$
For any $w\cdo\l\in W\fs$ we define $\fR_\l^{\dw}\in\cd_m(G)$, 
$R_\l^{\dw}\in\cd_m(G)$ by 
$$\fR_\l^{\dw}=\c(\cl_\l^{\dw}), R_\l^{\dw}=\c(\cl_\l^{\dw\sha})
\text{ if }w\cdo\l\in\overset\smile\to{W\fs},$$
$$\fR_\l^{\dw}=0, R_\l^{\dw}=0\text{ if }w\cdo\l\n\overset\smile\to{W\fs}.$$
Here $\overset\smile\to{W\fs}$ is as in 3.3.

We say that a simple perverse sheaf $A$ on $G$ is a character sheaf if the 
following equivalent conditions are satisfied:

-there exists $w\cdo\l\in W\fs$ and $j\in\ZZ$ such that $(A:(\fR_\l^{\dw})^j)\ne0$;

-there exists $w\cdo\l\in W\fs$ and $j\in\ZZ$ such that $(A:(R_\l^{\dw})^j)\ne0$.
\nl
(For the equivalence of these two conditions see \cite{\CSIII, 12.7}.) 
A character 
sheaf $A$ determines a $W$-orbit on $\fs$: the set of $\l\in\fs$ such that 
$(A:(\fR_\l^{\dw})^j)\ne0$ for some $w\in W$ and some $j$ (or equivalently
$(A:(R_\l^{\dw})^j)\ne0$ for some $w\in W$ and some $j$), see 
\cite{\CSIII, 11.2(a), 12.7}. We say that $A$ is an {\it $\fo$-character sheaf} if
the $W$-orbit on $\fs$ determined by $A$ is $\fo$ (as in 2.14). Let 
$CS_{\fo}$ be a set of representatives for the isomorphism classes of 
$\fo$-character sheaves on $G$.

By \cite{\CSIII, 14.11}, for any $\l\in\fo$ there exists a pairing 
$CS_{\fo}\T\Irr W'_\l@>>>\bbq$, $(A,e)\m b_{A,e}$ such that for any 
$A\in CS_{\fo}$, any $z\in W'_\l$ and any $j\in\ZZ$ we have
$$(A:(R_\l^{\dz})^j)
=(j-\D-|z|;(-1)^{j+\D}\sum_{e\in\Irr(W'_\l)}b_{A,e}\tr(c_{z\cdo\l},e^v)).$$
Here $e^v$ is as in 1.12. (When $z\cdo\l\n\overset\smile\to{W\fs}$, both sides are zero.)
By the results in 1.12 this can be reformulated as follows.

There exists a pairing $CS_{\fo}\T\Irr_{\fo}(W\TT_n)@>>>\bbq$, 
$(A,E)\m b_{A,E}$ such that for any $A\in CS_{\fo}$, any $\l\in\fo$, any
$z\in W$ and any $j\in\ZZ$ we have
$$(A:(R_\l^{\dz})^j)=(-1)^{j+\D}(j-\D-|z|;\sum_{E\in\Irr_{\fo}(W\TT_n)}
b_{A,E}\tr(c_{z\cdo\l},E^v))\tag a$$
where $E^v$ is as in 1.12; if $E'$ is an $\bbq[W\TT_n]$-module isomorphic to 
$\op_{E\in\Irr(W\TT_n)}E^{\op m_E}$ (with $m_E\in\NN$) we set
$$b_{A,E'}=\sum_{E\in\Irr_{\fo}(W\TT_n)}m_Eb_{A,E}.$$
In particular if $E'\in\Irr(W\TT_n)-\Irr_{\fo}(W\TT_n)$, we have $b_{A,E'}=0$.
Moreover, given $A\in CS_{\fo}$, there is a unique two-sided cell $\boc_A$ 
of $W\fs$ such that $b_{A,E}=0$ whenever $E\in\Irr_{\fo}(W\TT_n)$ satisfies 
$\boc_E\ne\boc_A$. (This follows from \cite{\CSIII, 16.7}.) We have 
necessarily $\boc_A\sub\{w\cdo\l\in W\fs;\l\in\fo\}$. As in 
\cite{\CDGIX, 41.8}, \cite{\CDGX, 44.18}, we see that:

(b) {\it $(A:(R_\l^{\dz})^j)\ne0$ for some $z\cdo\l\in\boc_A$, $j\in\ZZ$; 
conversely, if $(A:(R_\l^{\dz})^j)\ne0$ for $z\cdo\l\in W\fs$, $j\in\ZZ$, then
$\boc_A\preceq z\cdo\l$.}
\nl
Let $a_A$ be the value of the $a$-function on $\boc_A$. If $z\cdo\l\in W\fs$, 
$E\in\Irr_{\fo}(W\TT_n)$ satisfy $\tr(c_{z\cdo\l},E^v)\ne0$ then 
$\boc_E\preceq z\cdo\l$; if in addition we have $z\cdo\l\in\boc_E$ then from
the definitions we have
$$\tr(c_{z\cdo\l},E^v)=\sum_{h\ge0}\g_{z\cdo\l,E,h}v^{a_E-h}$$
where $\g_{z\cdo\l,E,h}\in\bbq$ is zero for large $h$, 
$\g_{z\cdo\l,E,0}=\tr(t_{z\cdo\l},E^\iy)$ and $a_E$ is as in 1.13.
Hence from (a) we see that for $A\in CS_{\fo}$ and
$\l\in\fo$, $z\in W$, $j\in\ZZ$, the following holds:

(c) {\it $(A:(R_\l^{\dz})^j)=0$ unless $\boc_A\preceq z\cdo\l$; if 
$z\cdo\l\in\boc_A$ then
$$(A:(R_\l^{\dz})^j)=(-1)^{j+\D}(j-\D-|z|;\sum_{E\in\Irr_{\fo}(W\TT_n);
\boc_E=\boc_A;h\ge0}b_{A,E}\g_{z\cdo\l,E,h}v^{a_A-h})$$
which is $0$ unless $j-\D-|z|\le a_A$.}
\nl
Recall that $\boc$, $a$ are as in 2.14.

Let $\cm^{\preceq}G$  (resp. $\cm^{\prec}G$) be the category of perverse 
sheaves on $G$ whose composition factors are all of the form 
$A\in CS_{\fo}$ with $\boc_A\preceq\boc$ (resp. $\boc_A\prec\boc$). Let 
$\cd^{\preceq}G$  (resp. $\cd^{\prec}G$) be the subcategory of $\cd(G)$ whose 
objects are complexes $K$ such that $K^j$ is in $\cm^{\preceq}G$ (resp. 
$\cm^{\prec}G$) for any $j$. Let $\cd^{\preceq}_mG$ (resp. $\cd^{\prec}_mG$) 
be the subcategory of $\cd_m(G)$ whose objects are also in $\cd^{\preceq}G$ 
(resp. $\cd^{\prec}G$). 

Let $\l\in\fo$, $z\in W$. From (c) we deduce:

(d) {\it If $z\cdo\l\preceq\boc$ then $(R_\l^{\dz})^j\in\cm^{\preceq}G$ for 
all $j\in\ZZ$.}

(e) {\it If $z\cdo\l\in\boc$ and $j>a+\D+|z|$ then 
$(R_\l^{\dz})^j\in\cm^{\prec}G$.}

(f) {\it If $z\cdo\l\prec\boc$ then $(R_\l^{\dz})^j\in\cm^{\prec}G$ for all 
$j\in\ZZ$.}

\subhead 5.2\endsubhead
Let $CS_\boc=\{A\in CS_{\fo};\boc_A=\boc\}$. For any $z\cdo\l\in W\fs$ we 
set 
$$n_z=a(z)+\D+|z|.$$
Let $A\in CS_\boc$ and let $z\cdo\l\in\boc$. We have
$$(A:(R_\l^{\dz})^{n_z})
=(-1)^{a+|z|}\sum_{E\in\Irr(W\TT_n);\boc_E=\boc}b_{A,E}\tr(t_{z\cdo\l},E^\iy).
\tag a$$
Indeed, from 5.1(a) we have
$$(A:(R_\l^{\dz})^{n_z})=(-1)^{a+|z|}\sum_{E\in\Irr(W\TT_n);\boc_E=\boc}b_{A,E}
        (a;\tr(c_{z\cdo\l},E^v))$$
and it remains to use that $(a;\tr(c_{z\cdo\l},E^v))=\tr(t_{z\cdo\l},E^\iy)$. 
                   
We show:

(b) {\it For any $A\in CS_\boc$ there exists $E\in\Irr(W\TT_n)$ such that 
$\boc_E=\boc$, $b_{A,E}\ne0$.}
\nl
Assume that this is not so. Then, using 5.1(a), for any $\l\in\fo$, any 
$z\in W$ and any $j\in\ZZ$ we have $(A:(R_\l^{\dz})^j)=0$. This contradicts 
the assumption that $A\in CS_{\fo}$.

(c) {\it For any $A\in CS_\boc$ there exists $z\cdo\l\in\boc$ such that 
$(A:(R_\l^{\dz})^{n_z})\ne0$.}
\nl
Assume that this is not so. Then, using (a), we see that
$$\sum_{E\in\Irr(W\TT_n);\boc_E=\boc}b_{A,E}\tr(t_{z\cdo\l},E^\iy)=0$$
for any $z\cdo\l\in\boc$. Using this and (b) we see that the linear functions 
$t_{z\cdo\l}\m\tr(t_{z\cdo\l},E^\iy)$ on $\JJ_\boc$ (for various 
$E\in\Irr(W\TT_n)$ such that $\boc_E=\boc$) are linearly dependent. This is a 
contradiction since
the $E^\iy$ form a complete set of simple modules for the semisimple algebra 
$\JJ_\boc$.

We show:

(d) {\it Let $z\cdo\l\in\boc$ be such that $(R_\l^{\dz})^{n_z}\ne0$. Then 
$z\cdo\l\underset\text{left}\to\si z\i\cdo z(\l)$. In particular we have $z\in W'_\l$ and
$z,z\i$ are in the same left cell of $W'_\l$.}
\nl
Using (a) we see that there exists $E\in\Irr(W\TT_n)$ such that $\tr(t_{z\cdo\l},E^\iy)\ne0$. We have 
$E^\iy=\op_{d\cdo\l_1\in\DD}t_{d\cdo\l_1}E^\iy$ and $t_{z\cdo\l}:E^\iy@>>>E^\iy$ maps the summand
$t_{d\cdo\l_1}E^\iy$ where $z\cdo\l\underset\text{left}\to\si d\cdo\l_1$ 
into $t_{d'\cdo\l'_1}E^\iy$ where 
$d'\cdo\l'_1\in\DD$, $d\cdo\l'_1\underset\text{left}\to\si z\i\cdo z(\l)$ 
and all other summands to $0$. Since
$\tr(t_{z\cdo\l},E^\iy)\ne0$, we must have $t_{d\cdo\l_1}E^\iy=t_{d'\cdo\l'_1}E^\iy\ne0$ hence
$d\cdo\l_1=d'\cdo\l'_1$ and $z\cdo\l\underset\text{left}\to\si z\i\cdo z(\l)$. This proves (d).

\subhead 5.3\endsubhead
We show:

(a) {\it If $L\in\cd^{\preceq}Z$ then $\c(L)\in\cd^{\preceq}G$. If 
$L\in\cd^{\prec}Z$ then $\c(L)\in\cd^{\prec}G$.}

(b)  {\it If $L\in\cm^{\preceq}Z$ and $j>a+\nu$ then 
$(\c(L))^j\in\cm^{\prec}G$.}
\nl
It is enough to prove (a),(b) assuming in addition that $L=\Bbb L_\l^{\dz}$ 
where $z\cdo\l\in\overset\smile\to{W\fs}$, $z\cdo\l\preceq\boc$. Then (a) follows from 
5.1(d),(f). In the setup of (b) we have
$$(\c(\Bbb L_\l^{\dz}))^j=(R_\l^{\dz})^{j+|z|+\nu+\r}((|z|+\nu+\r)/2)$$ and 
this is in $\cm^{\prec}G$ since $j+|z|+\nu+\r>a+\D+|z|$, see 5.1(e).

\subhead 5.4\endsubhead
Let $\cc^\spa G$ be the subcategory of $\cm(G)$ consisting of semisimple 
objects. Let $\cc_0^\spa G$ be the subcategory of $\cm_m(G)$ consisting of 
those $K$ such that $K$ is pure of weight zero. Let $\cc^\boc G$ be the 
subcategory of $\cm(G)$ consisting of objects which are direct sums of objects
of the form $A\in CS_{\boc}$. Let $\cc_0^\boc G$ be the subcategory of 
$\cc^\spa_0G$ consisting of those $K$ such that, as an object of $\cc^\spa G$,
$K$ belongs to $\cc^\boc G$. For $K\in\cc_0^\spa G$ let $\un{K}$ be the 
largest subobject of $K$ such that as an object of $\cc^\spa G$, we have 
$\un{K}\in\cc^\boc G$.

\subhead 5.5\endsubhead
For $L\in\cc^\boc_0Z$ we set 
$$\un\c(L)=\un{(\c(L))^{a+\nu}}((a+\nu)/2)
=\un{(\c(L))^{\{a+\nu\}}}\in\cc^\boc_0G.$$  
(The last equality uses that $\p$ in 5.1 is proper hence it preserves purity.)
The functor $\un\c:\cc^\boc_0Z@>>>\cc^\boc_0G$ is called {\it truncated 
induction}. For $z\cdo\l\in\overset\smile\to\boc$ we have
$$\un\c(\Bbb L_\l^{\dz})=\un{(R_\l^{\dz})^{n_z}}(n_z/2).\tag a$$
Indeed,
$$\align&\un\c(\Bbb L_\l^{\dz})
=\un{(\c(\Bbb L_\l^{\dz}))^{a+\nu}}((a+\nu)/2)
=\un{(\c(\cl^{\dz\sha}_\l\la|z|+\nu+\r\ra))^{a+\nu}}((a+\nu)/2)\\&
=\un{(\c(\cl^{\dz\sha}_\l))^{|z|+a+\D}}((|z|+a+\D)/2)
=\un{(\c(\cl^{\dz\sha}_\l))^{n_z}}(n_z/2)=\un{(R_\l^{\dz})^{n_z}}(n_z/2).\endalign$$
Using (a) and 5.2(d) we see that:

(d) {\it If $z\cdo\l\in\overset\smile\to\boc$ is such that $\un\c(\Bbb L_\l^{\dz})\ne0$ then
$z\cdo\l\underset\text{left}\to\si z\i\cdo z(\l)$. In particular we have $z\in W'_\l$ and
$z,z\i$ are in the same left cell of $W'_\l$.}

\subhead 5.6\endsubhead
As in 1.9 we shall denote by $\t:\HH^\iy@>>>\ZZ$ the group homomorphism such that 
$\t(t_{z\cdo\l})=1$ if $z\cdo\l\in\DD$ and $\t(t_{z\cdo\l})=0$ if 
$z\cdo\l\in W\fs-\DD$. For $z\cdo\l,z'\cdo\l'$ in $\overset\smile\to\boc$ we show:
$$\dim\Hom_{\cc^\boc G}
(\un\c(\Bbb L_\l^{\dz}),\un\c(\Bbb L_{\l'}^{\dz'}))
=\sum_{y\in W;y\cdo\l'\in\boc}
\t(t_{y\i\cdo y(\l')}t_{z\cdo\l}t_{y\cdo\l'}t_{z'{}\i\cdo\l'}).\tag a$$
Using 5.5(a) and the definitions we see that the left hand side of (a) equals
$$\sum_{A\in CS_\boc}(A:(R_\l^{\dz})^{n_z})(A:(R_{\l'}^{\dz'})^{n_{z'}}).$$
Using 5.2(a) and the analogous identity for $(A:(R_{\l'}^{\dz'})^{n_{z'}})$ in
which the field automorphism $()^\spa:\bbq@>>>\bbq$ (see 1.16) is applied to 
both sides (the left hand side is fixed by $()^\spa$), we see that the left 
hand side of (a) equals
$$\align&(-1)^{|z|+|z'|}\sum_{E,E'\in\Irr(W\TT_n);\boc_E=\boc_{E'}=\boc}\\&\sum_{A\in CS_\boc}b_{A,E}
b_{A,E'}^\spa\tr(t_{z\cdo\l},E^\iy)(\tr(t_{z'\cdo\l'},E'{}^\iy))^\spa.\endalign$$
Replacing in the last sum $\sum_{A\in CS_\boc}b_{A,E}b_{A,E'}^\spa$ by $1$ if 
$E'=E$ and by $0$ if $E'\ne E$ (see \cite{\CDGVII, 35.18(g)}) we see that the 
left hand side of (a) equals
$$(-1)^{|z|+|z'|}\sum_{E\in\Irr(W\TT_n);\boc_E=\boc}\tr(t_{z\cdo\l},E^\iy)
(\tr(t_{z'\cdo\l'},E^\iy))^\spa.$$
By \cite{\CDGVII, 34.17}, for $E\in\Irr(W\TT_n)$ and $h\in\HH$ we have 
$\tr(h^\flat,E^v)=\tr(h,E^v)^\spa$ where $()^\spa:\bbq(v)@>>>\bbq(v)$ is as in 
1.16; in particular, for $w\cdo\l\in W\fs$ we have 
$\tr(c_{w\cdo\l},E^v)=\tr(c_{w\i\cdo w(\l)},E^v)^\spa$. Taking the coefficient
of $v^a$ in the two sides of the last equality we deduce
$\tr(t_{w\cdo\l},E^\iy)=\tr(t_{w\i\cdo w(\l)},E^\iy)^\spa$. Thus the left hand
side of (a) equals
$$(-1)^{|z|+|z'|}\sum_{E\in\Irr(W\TT_n),\boc_E=\boc}\tr(t_{z\cdo\l},E^\iy)
\tr(t_{z'{}\i\cdo\l'},E^\iy).$$
(Recall that $z'(\l')=z'$.)
This is equal to $(-1)^{|z|+|z'|}$ times the trace of the linear map 
$\x\m t_{z\cdo\l}\x t_{z'{}\i\cdo\l'}$ on $\JJ_\boc$; hence it is
equal to the sum over $y\cdo\l_1\in\boc$ of the coefficient of 
$t_{y\cdo\l_1}$ in $t_{z\cdo\l}t_{y\cdo\l_1}t_{z'{}\i\cdo\l'}$; this 
coefficient is $0$ if $\l_1\ne\l'$ while if $\l_1=\l'$ it is equal to 
$$\t(t_{y\i\cdo y(\l')}t_{z\cdo\l}t_{y\cdo\l'}t_{z'{}\i\cdo\l'}).$$
(We use 1.9(a) for $\HH^\iy$.) Thus we have
$$\dim\Hom_{\cc^\boc G}(\un\c(\Bbb L_\l^{\dz}),
\un\c(\Bbb L_{\l'}^{\dz'})) =(-1)^{|z|+|z'|}\sum_{y\in W;y\cdo\l'\in\boc}
\t(t_{y\i\cdo y(\l')}t_{z\cdo\l}t_{y\cdo\l'}t_{z'{}\i\cdo\l'}).$$
Since $\dim\Hom_{\cc^\boc G}(\un\c(\Bbb L_\l^{\dz}),
\un\c(\Bbb L_{\l'}^{\dz'}))\in\NN$ and the last sum is in $\NN$,
it follows that (a) holds.

\mpb

The proof above shows also that 
$\dim\Hom_{\cc^\boc G}(\un\c(\Bbb L_\l^{\dz}),\un\c(\Bbb L_{\l'}^{\dz'}))=0$ 
whenever $(-1)^{|z|+|z'|}=-1$.

\subhead 5.7\endsubhead
Let $L\in\cc^\boc_0Z$. We show that $\fD(L)\in\cc^{\ti{\boc}}_0Z$. It is enough to note that
for $w\cdo\l\in\overset\smile\to\boc$ and $\o\in\k\i_q(w)$ we have

(a) $\fD(\Bbb L_\l^\o)=\Bbb L_{\l\i}^\o$. 
\nl
We show: 

(b) {\it We have canonically $\un\c(\fD(L))=\fD(\un\c(L))$ where the first $\un\c$ is
relative to $\ti{\boc}$ instead of $\boc$.}
\nl
By the relative hard Lefschetz theorem \cite{\BBD, 5.4.10} applied to the 
projective morphism $\p$ (see 5.1) and to $f^*L\la\nu\ra$ (a perverse sheaf of
pure weight $0$ on $\dZ$, see 5.1) we have canonically for any $i$:
$$(\p_!f^*L\la\nu\ra)^{-i}=(\p_!f^*L\la\nu\ra)^i(i).\tag c$$
We have used the fact that $f$ is smooth with fibres of dimension $\nu$. This 
also shows that
$$\fD(\c(\fD(L)))=\c(L)\la2\nu\ra.\tag d$$
Using (d) we have
$$\align&\fD(\un\c(\fD(L)))=\fD((\c(\fD(L)))^{a+\nu}((a+\nu)/2)))
=(\fD(\c(\fD(L))))^{-a-\nu}((-a-\nu)/2)\\&
=(\c(L)\la2\nu\ra)^{-a-\nu}((-a-\nu)/2)=(\c(L)\la\nu\ra)^{-a}(-a/2).\endalign$$
Hence using (c) we have
$$\fD(\un\c(\fD(L)))=(\c(L)\la\nu\ra)^a(a/2)
=(\c(L))^{a+\nu}((a+\nu)/2)=\un\c(L).$$
This proves (b).

\subhead 5.8\endsubhead
Let $z\cdo\l\in\DD_\boc$ and let $\L_{z\cdo\l}$ be the left cell of $W\fs$ 
containing $z\cdo\l$. We show:

(a) {\it $(A:\un\c(\Bbb L_\l^{\dz}))=b_{A,[\L_{z,\l}]}$ for any $A\in CS_\boc$.}
\nl
Using 1.12(a) we see that for any $E\in\Irr(W\TT_n)$, $\tr(t_{z\cdo\l},E^\iy)$
is equal to the multiplicity of $E$ in the $\bbq[W\TT_n]$-module 
$[\L_{z\cdo\l}]$. Hence, using 5.5(a) and 5.2(a), we have
$$\align&(A:\un\c(\Bbb L_\l^{\dz}))
=(-1)^{a+|z|}\sum_{E\in\Irr(W\TT_n);\boc_E=\boc}b_{A,E}(\text{ multiplicity of
$E$ in } [\L_\d])\\&=(-1)^{a+|z|}b_{A,[\L_{z,\l}]}.\endalign$$
It is enough to show that $a+|z|=0\mod2$. Since $z$ is a distinguished 
involution of $W_\l$, the coefficient of $v^{-a}$ in $p^\l_{1,z}$ (see 1.8) is
nonzero (see \cite{\HEC, 14.1}). Using now \cite{\HEC, 5.4(b)} we deduce that 
$|z|_\l=a\mod2$. It remains to note that $|z|_\l=|z|\mod2$. (Indeed, $W_\l$ is
generated by elements $u\in W_\l$ such that $|u|_\l=1$ and such $u$ are 
reflections in $W$ so that $|u|$ is odd.)

\subhead 5.9\endsubhead
We define $\z:\cd(G)@>>>\cd(Z)$ and $\z:\cd_m(G)@>>>\cd_m(Z)$ by
$\z(K)=f_!\p^*K$ where $Z@<f<<\dZ@>\p>>G$ is as in 5.1(a). We show:

(a) {\it For any $L\in\cd(Z)$ or $L\in\cd_m(Z)$ we have $\fb''(L)=\z(\c(L))$.}
\nl
We have $\z(\c(L))=f_!\p^*\p_!f^*(L)$. We have
$$\dZ\T_G\dZ=\{((B_0,B_1,B_2,B_3),g)\in\cb^4\T G;gB_0g\i=B_3,gB_1g\i=B_2\}.$$
We have a cartesian diagram
$$\CD
\dZ\T_G\dZ@>\ti\p_1>>\dZ\\
@V\ti\p_2VV          @V\p VV\\
\dZ    @>\p>>   G
\endCD$$
where $\ti\p_1((B_0,B_1,B_2,B_3),g)=(B_0,B_3,g)$,
$\ti\p_2((B_0,B_1,B_2,B_3),g)=(B_1,B_2,g)$.
It follows that $\p^*\p_!=\ti\p_{1!}\ti\p_2^*$. Thus
$$\z(\c(L))=f_!\ti\p_{1!}\ti\p_2^*f^*(L)=(f\ti\p_1)_!(f\ti\p_2)^*(L).$$
Define $\p'_1:\dZ\T_G\dZ@>>>Z$, $\p'_2:\dZ\T_G\dZ@>>>Z$ by 
$$\p'_1((B_0,B_1,B_2,B_3),g)=(B_0,B_3,gU_{B_0}),$$
$$\p'_2((B_0,B_1,B_2,B_3),g)=(B_1,B_2,gU_{B_1}).$$
Then $\p'_1=f\ti\p_1$, $\p'_2=f\ti\p_2$ and $\z(\c(L))=\p'_{1!}\p'_2{}^*(L)$. 
We have an isomorphism ${}^\di\cy@>>>\dZ\T_G\dZ$ induced by
$$((x_0\UU,x_1\UU,x_2\UU,x_3\UU),g)\m
((x_0\BB x_0\i,x_1\BB x_1\i,x_2\BB x_2\i,x_3\BB x_3\i),g).$$
We use this to identify ${}^\di\cy=\dZ\T_G\dZ$. Then $\p'_1,\p'_2$ become
$d,{}^\di\eta$ of 3.25. We see that (a) holds.

\subhead 5.10\endsubhead
Let $z\cdo\l\in\fo$. We set $\Si=\e^*\z(R_\l^{\dz})\la2\nu+|z|\ra\in
\cd_m(\tcb^2)$. Let $j\in\ZZ$. We show:

(a) {\it If $z\cdo\l\preceq\boc$ then $\Si^j\in\cm^{\preceq}\tcb^2$.}

(b) {\it If $z\cdo\l\prec\boc$ then $\Si^j\in\cm^{\prec}\tcb^2$.}

(c) {\it If $z\cdo\l\in\boc$ and $j>\nu+2\r+2a$ then $\Si^j\in\cm^{\prec}\tcb^2$.}
\nl
If $z\cdo\l\n\overset\smile\to{W\fs}$ then $\Si=0$ and there is nothing to prove.
Now assume that $z\cdo\l\in\overset\smile\to{W\fs}$. Then, using 5.9(a), we have
$$\Si=\e^*\z(\c(\cl_\l^{\dz\sha}))\la2\nu+|z|\ra=
\fb'(\cl_\l^{\dz\sha})\la2\nu+|z|\ra=\fb'(\Bbb L_\l^{\dz})\la\nu-\r\ra.$$
Now (a),(b) follow from 3.14(a),(b) and (c) follows from 3.14(c).
(If $j>\nu+2\r+2a$ then $j+\nu-r>2\nu+\r+2a$.)

\subhead 5.11\endsubhead 
We show:

(a) {\it If $K\in\cd^{\preceq}G$ then $\z(K)\in\cd^{\preceq}Z$.}

(b) {\it If $K\in\cd^{\prec}G$ then $\z(K)\in\cd^{\prec}Z$.}

(c) {\it If $K\in\cd^{\preceq}G$ and $j>\nu+a$ then $(\z(K))^j\in\cm^{\prec}Z$.}
\nl
We can assume in addition that $K=A\in CS_{\fo}$ where $A\in CS_{\boc'}$ 
for a two-sided cell $\boc'$ such that $\boc'\preceq\boc$. Assume first that 
$\boc'=\boc$. By 5.2(c) we can find $z\cdo\l\in\boc$ such that 
$(A:(R_\l^{\dz})^{n_z})\ne0$. Then $A[-n_z]$ (without mixed structure) is a 
direct summand of $R_\l^{\dz}$ (which is a semisimple complex). Hence 
$\e^*\z(A)[-n_z]$ is a direct summand of $\e^*\z(R_\l^{\dz})$ and
$\e^*\z(A)[-n_z+2\nu+|z|]$ is a direct summand of $\Si$ (in 5.10), that is,
$\e^*\z(A)[-a-\r]$ is a direct summand of $\Si$. By 5.10, if $j\in\ZZ$ (resp. 
$j>\nu+2\r+2a$) then $\Si^j\in\cm^{\preceq}\tcb^2$ (resp. 
$\Si^j\in\cm^{\prec}\tcb^2$) hence $(\e^*\z(A)[-a-\r])^j\in\cm^{\preceq}\tcb^2$
(resp. $(\e^*\z(A)[-a-\r])^j\in\cm^{\prec}\tcb^2$), that is 
$(\e^*\z(A))^{j-a-\r}\in\cm^{\preceq}\tcb^2$ (resp. 
$(\e^*\z(A))^{j-a-\r}\in\cm^{\prec}\tcb^2$). We see that if $j'\in\ZZ$ (resp. 
$j'>\nu+\r+a$) then $(\e^*\z(A))^{j'}\in\cm^{\preceq}\tcb^2$ (resp. 
$(\e^*\z(A))^{j'}\in\cm^{\prec}\tcb^2$) so that
$(\z(A))^{j'-\r}\in\cm^{\preceq}Z$ (resp. $(\z(A))^{j'-\r}\in\cm^{\prec}Z$); 
here we use 3.3(a). We see that if $j\in\ZZ$ (resp. $j>\nu+a$, so that 
$j+\r>\nu+\r+a$) then $(\z(A))^j\in\cm^{\preceq}Z$ (resp.
$(\z(A))^j\in\cm^{\prec}Z$). Thus the desired results hold when $\boc'=\boc$.

Assume now that $\boc'\prec\boc$. Applying the above argument with $\boc$ 
replaced by $\boc'$ we see that (a),(b) hold.

\subhead 5.12\endsubhead 
For $K\in\cc^\boc_0G$ we set
$$\un\z(K)=\un{(\z(K))^{\{\nu+a\}}}\in\cc^\boc_0Z.$$
We say that $\un\z(K)$ is the {\it truncated restriction} of $K$.

\subhead 5.13\endsubhead 
Let $L\in\cc^\boc_0Z$. We show:

(a) {\it We have canonically $\un\z(\un\c(L))=\un{\fb''}(L)$.}
\nl
We shall apply \cite{\CONV, 1.12} with 
$\Ph:\cd_m(Y_1)@>>>\cd_m(Y_2)$ replaced by 
$\z:\cd_m(G)@>>>\cd_m(Z)$ and with $\cd^{\preceq}(Y_1)$, 
$\cd^{\preceq}(Y_2)$ replaced by $\cd^{\preceq}G$, $\cd^{\preceq}Z$. We shall
take $\XX$ in {\it loc.cit.} equal to $\c(L)$. 
The conditions of {\it loc.cit.} are satisfied: those concerning $\XX$ are
satisfied with $c'=a+\nu$, see 5.3. The conditions concerning $\z$ are 
satisfied with $c=a+\nu$, see 5.11. We see that 
$$(\z(\c(L)))^j=0\text{ if }j>2a+2\nu\tag b$$
and
$$\un{gr_{2a+2\nu}((\z(\c(L)))^{2a+2\nu})}(a+\nu)=\un\z(\un\c(L)).\tag c$$
Since $\z(\c(L))=\fb''(L)$, we see that the left hand side of (c) equals
$\un{\fb''}(L)$. Thus (a) is proved.

\mpb

Combining (a) with 3.25(d) and 3.14(d) we obtain the following result.

(b) {\it We have canonically $\ti\e\un\z(\un\c(L))=\un{\fb}(L)$.}

\subhead 5.14\endsubhead
Let $K\in\cd_m(G)$ and let $L\in\cd^\spa_m\tcb^2$. We show that

(a) {\it there is a canonical isomorphism 
$L\cir\e^*\z(K)@>\si>>\e^*\z(K)\cir L$.}
\nl
Let $Y=\tcb^2\T G$. Define $j:Y@>>>G$ by $j(x_0\UU,x_1\UU,g)=g$. Define 
$j_1:Y@>>>\tcb^2$ by $j_1(x_0\UU,x_1\UU,g)=(x_0\UU,g\i x_1\UU)$. Define 
$j_2:Y@>>>\tcb^2$ by $j_2(x_0\UU,x_1\UU,g)=(gx_0\UU,x_1\UU)$. From the 
definitions we have $L\cir\e^*\z(K)=j_{2!}(j_1^*(L)\ot j^*(K))$, 
$\e^*\z(K)\cir L=j_{2!}(j_2^*(L)\ot j^*(K))$. By the $G$-equivariance of $L$ 
we have $j_1^*L=j_2^*L$; (a) follows.

\mpb

Now let $K\in\cc^\boc_0G$ and let $L\in\cc^\boc_0\tcb^2$. We show that

(b) {\it there is a canonical isomorphism 
$L\un\cir\ti\e\un\z(K)@>\si>>(\ti\e\un\z(K))\un\cir L$.}
\nl
We apply \cite{\CONV, 1.12} with $\Ph:\cd^{\preceq}_m \tcb^2@>>>\cd^{\preceq}_m \tcb^2$, $L'\m L'\cir L$,
$\XX=\ti\e\z(K)$ 
and with $(c,c')=(a-\nu,\nu+a)$, see 2.23(a), 5.11(c). We deduce that we have canonically
$$\un{\un{((\ti\e\z(K))^{\{a+\nu\}}}\cir L)^{\{a-\nu\}}}=\un{(\ti\e\z(K)\cir L)^{\{2a\}}}.\tag c$$
We apply \cite{\CONV, 1.12} with $\Ph:\cd^{\preceq}_m \tcb^2@>>>\cd^{\preceq}_m \tcb^2$, $L'\m L\cir L'$,
$\XX=\ti\e\z(K)$ and with $(c,c')=(a-\nu,\nu+a)$, see 2.23(a), 5.11(c). We deduce that we have 
canonically
$$\un{((L\cir\un{(\ti\e\z(K))^{\{a+\nu\}})^{\{a-\nu\}}}}=\un{(L\cir\ti\e\z(K))^{\{2a\}}}.\tag d$$
We now combine (c),(d) with (a); we obtain (b).

\subhead 5.15\endsubhead
Let $\mu:G\T G@>>>G$ be the multiplication map. For $K,K'$ in $\cd(G)$ (resp. 
in $\cd_m(G)$) we set $K*K'=\mu_!(K\bxt K')$; this is in $\cd(G)$ (resp. in
$\cd_m(G)$). For $K,K',K''$ in $\cd_m(G)$ we have canonically
$(K*K') *K''=K*(K'*K'')$ (and we denote this by $K*K'*K''$). Note that if 
$K\in\cd_m(G)$ and $K'\in\cm_m(G)$ is $G$-equivariant for the conjugation 
action of $G$ then as in \cite{\CONV, 4.1} we have a canonical isomorphism
$$K*K'@>\si>>K'*K.\tag a$$

\subhead 5.16\endsubhead
We show:

(a) {\it For $K\in\cd_m(G)$, $L\in\cd_m(Z)$ we have canonically 
$K*\c(L)=\c(L\bul\z(K))$.}
\nl
Let $Y=G\T G\T\cb$. Define $c:Y@>>>G\T Z$ by 
$c(g_1,g_2,B)=(g_1,(B,g_2Bg_2\i,g_2U_B))$; define $d:Y@>>>G$ by
$d(g_1,g_2,B)=g_1g_2$. From the definitions we see that both $K*\c(L)$, 
$\c(L\bul\z(K))$ can be identified with $d_!c^*(K\bxt L)$. This proves (a).

Now let $L,L'\in\cd_m(Z)$. Replacing in (a) $K,L$ by $\c(L),L'$ and using
5.9(a), we obtain
$$\c(L)*\c(L')=\c(L'\bul\fb''(L)).\tag b$$

\subhead 5.17\endsubhead
Let $L,L'\in\cd^\spa_m(Z)$, $j\in\ZZ$. We show:

(a) {\it If $L\in\cd^{\preceq}Z$ or $L'\in\cd^{\preceq}Z$ then 
$L'\bul\fb''(L)\in\cd^{\preceq}Z$.}

(b) {\it If $L\in\cd^{\prec}Z$ or $L'\in\cd^{\prec}Z$ then 
$L'\bul\fb''(L)\in\cd^{\prec}Z$.}

(c) {\it If $L\in\cm^{\preceq}Z$, $L'\in\cm^\spa Z$ and $j>3a+\r+\nu$ then
$(L'\bul\fb''(L))^j\in\cd^{\prec}Z$.}
\nl
Now (a),(b) follow from 3.25(b) and 3.23(a). To prove (c) we may assume that 
$L=\Bbb L_\l^{\dw}$, $L'=\Bbb L_{\l'}^{\dw'}$ with $w\cdo\l,w'\cdo\l'$ in 
$\overset\smile\to{W\fs}$ and $w\cdo\l\preceq\boc$. We apply \cite{\CONV, 1.12} with
$\Ph:\cd_m^{\preceq}Z@>>>\cd_m^{\preceq}Z$, $L_1\m L'\bul L_1$ and
$\XX=\fb''(L)$ and with $c'=2\nu+2a$ (see 3.25(c)), $c=a+\r-\nu$ (see 3.23(b)).
We have $c+c'=\nu+\r+3a$ hence (c) holds.

\subhead 5.18\endsubhead
Let $L,L'\in\cd^\spa_m(Z)$, $j\in\ZZ$. We show:

(a) {\it If $L\in\cd^{\preceq}Z$ or $L'\in\cd^{\preceq}Z$ then 
$\c(L'\bul\fb''(L))\in\cd^{\preceq}G$.}

(b) {\it If $L\in\cd^{\prec}Z$ or $L'\in\cd^{\prec}Z$ then 
$\c(L'\bul\fb''(L))\in\cd^{\prec}G$.}

(c) {\it If $L'\in\cm^\spa Z$, $L\in\cm^{\preceq}Z$ and $j>4a+2\nu+\r$ then
$(\c(L'\bul\fb''(L)))^j\in\cm^{\prec}G$.}
\nl
(a),(b) follow from 5.3(a) using 5.17(a),(b). To prove (c) we can assume that 
$L=\Bbb L_\l^{\dw}$, $L'=\Bbb L_{\l'}^{\dw'}$ with $w\cdo\l,w'\cdo\l'$ in 
$\overset\smile\to{W\fs}$ and $w\cdo\l\preceq\boc$.
 We apply \cite{\CONV, 1.12} with 
$\Ph:\cd_m^{\preceq}Z@>>>\cd_m^{\preceq}G$, $L_1\m\c(L_1)$,
$\XX=L'\bul\fb''(L)$ and with $c'=\nu+\r+3a$ (see 5.17(c)), $c=a+\nu$ (see 5.3(b)).
We have $c+c'=2\nu+\r+4a$ hence (c) holds. 

\subhead 5.19\endsubhead
Let $K,K'\in\cd_m^\spa(G)$. We show:

(a) {\it If $K\in\cd^{\preceq}G$ or $K'\in\cd^{\preceq}G$ then 
$K*K'\in\cd^{\preceq}G$.}

(b) {\it If $K\in\cd^{\prec}G$ or $K'\in\cd^{\prec}G$ then $K*K'\in\cd^{\prec}G$.}

(c) {\it If $K\in\cd^{\preceq}G$ or $K'\in\cd^{\preceq}G$ and $j>2a+\r$ then 
$(K*K')^j\in\cd^{\prec}G$.}
\nl
We can assume that $K=A\in CS_{\fo}$, $K'=A'\in CS_{\fo}$. 
Let $A''\in\cm(G)$ be a composition factor of $(A*A')^j$. By 5.2(c) we can 
find $w\cdo\l\in\boc_A$, $w'\cdo\l'\in\boc_{A'}$ such that 
$(A:(R_\l^{\dw})^{n_w})\ne0$, $(A':(R_{\l'}^{\dw'})^{n_{w'}})\ne0$. Then $A$ 
is a direct summand of $R_\l^{\dw}[n_w]$ and $A'$ is a direct summand of 
$R_{\l'}^{\dw'}[n_{w'}]$ (without mixed structures). Hence $A*A'$ is a direct 
summand of 
$$R_\l^{\dw}*R_{\l'}^{\dw'}[a(w\cdo\l)+a(w'\cdo\l')+|w|+|w'|+2\D]$$ 
and $(A*A')^j$ is a direct summand of 
$$\align&(R_\l^{\dw}*R_{\l'}^{\dw'}[|w|+|w'|+2\nu+2\r])^{j+a(w\cdo\l)+a(w'\cdo\l')
+2\nu}\\&=(\c(\Bbb L_\l^{\dw})*\c(\Bbb L_{\l'}^{\dw'}))^{j+a(w\cdo\l)+
a(w'\cdo\l')+2\nu}.\endalign$$
Using 5.16(b) we see that $(A*A')^j$ is a direct summand of 
$$(\c(\Bbb L_{\l'}^{\dw'}\bul\fb''(\Bbb L_\l^{\dw}))^{j+a(w\cdo\l)
+a(w'\cdo\l')+2\nu}.\tag d$$
Hence $A''$ is a composition factor of (d). Using 5.18(a) we see that 
$A''\in CS_{\fo}$ and that $\boc_{A''}\preceq w\cdo\l$ and 
$\boc_{A''}\preceq w'\cdo\l'$.
In the setup of (a) we have $w\cdo\l\preceq\boc$ or $w'\cdo\l'\preceq\boc$ 
hence $\boc_{A''}\le\boc$. Thus (a) holds. Similarly, (b) holds. In the setup 
of (c) we have $w\cdo\l\preceq\boc$ and $w'\cdo\l'\preceq\boc$. Hence 
$a(w\cdo\l)\ge a$, $a(w'\cdo\l')\ge a$. Assume that $\boc_{A''}=\boc$. Since 
$A''$ is a composition factor of (d), we see from 5.18(c) that 
$$j+a(w\cdo\l)+a(w'\cdo\l')+2\nu\le 4a+2\nu+\r$$
hence $j+2a+2\nu\le4a+2\nu+\r$ and $j\le2a+\r$. This proves (c).

\subhead 5.20\endsubhead
For $K,K'\in\cc^\boc_0G$ we set
$$K\un{*}K'=\un{(K*K')^{\{2a+\r\}}}\in\cc^\boc_0G.$$
We say that $K\un{*}K'$ is the {\it truncated convolution} of $K,K'$. Note 
that 5.15(a) induces for $K,K'\in\cc^\boc_0G$ a canonical isomorphism
$$K\un{*}K'@>\si>>K'\un{*}K.\tag a$$

Let $L\in\cc^\boc_0Z$, $K\in\cc^\boc_0G$. Using \cite{\CONV, 1.12} several 
times, we see that
$$K\un{*}\un\c(L)=\un{gr_k((K*\c(L))^k)}(k/2)$$
where $k=(a+\nu)+(2a+\r)=3a+\nu+\r$ and
$$\un\c(L\un{\bul}\un\z(K))=\un{gr_{k'}((\c(L\bul\z(K))^{k'})}(k'/2)$$
where $k'=(a+\nu)+(a+\nu)+(a+\r-\nu)=3a+\nu+\r$.
Using now 5.16(a) and the equality $k=k'$ we obtain
$$K\un{*}\un\c(L)=\un\c(L\un{\bul}\un\z(K)).\tag b$$

\mpb

Let $L,L'\in\cc^\boc_0Z$. Using \cite{\CONV, 1.12} several times, we see that
$$\un\c(L)\un{*}\un\c(L')=\un{gr_k((\c(L)*\c(L'))^k)}(k/2)$$
where $k=(a+\nu)+(a+\nu)+(2a+\r)=4a+2\nu+\r$ and
$$\un\c(L'\un{\bul}\un{\fb''}(L)=
\un{gr_{k'}((\c(L'\bul\fb''(L)))^{k'})}(k'/2)$$
where $k'=(2a+2\nu)+(a+\r-\nu)+(a+\nu)=4a+2\nu+\r$.
Using now 5.16(b) and the equality $k=k'$ we obtain
$$\un\c(L)\un{*}\un\c(L')=\un\c(L'\un{\bul}(\un{\fb''}(L))).\tag c$$

\mpb

We show

(d) {\it For $K,K',K''$ in $\cc^\boc_0G$ there is a canonical isomorphism
$(K\un{*}K')\un{*}K''@>\si>>K\un{*}(K'\un{*}K'')$.}
\nl
Indeed, just as in \cite{\CONV, 4.7} we can identify, using \cite{\CONV, 1.12}, both
$(K\un{*}K')\un{*}K''$ and $K\un{*}(K'\un{*}K'')$ with $\un{(K*K'*K'')^{\{4a+2\r\}}}$.

A similar argument shows that the associativity isomorphism provided by (d) satisfies the
pentagon property.

\subhead 5.21\endsubhead
For $K,K'$ in $\cd_m(G)$ we show:

(a) {\it We have canonically $\z(K*K')=\z(K')\bul\z(K)$.}
\nl
Let $Y=\{(B,gU_B,h_1,h_2);B\in\cb,g\in G,h_1\in G,h_2\in G;h_1h_2\in gU_B\}$. 
Define $j_e:Y@>>>G$ by $j_e(B,gU_B,h_1,h_2)=h_e$ ($e=1,2$). Define $j:Y@>>>Z$ 
by $j(B,gU_B,h_1,h_2)=(B,gBg\i,gU_B)$. From the definitions we have 
$\z(K*K')=j_!(j_1^*(K)\ot j_2^*(K'))=\z(K')\bul\z(K)$; (a) follows.

For $K,K'$ in $\cd^\boc_0(G)$ we show:

(b) {\it We have canonically $\un\z(K\un{*}K')=\un\z(K')\un{\bul}\un\z(K)$.}
\nl
Using \cite{\CONV, 1.12} we see that
$$\un\z(K\un{*}K')= \un{gr_k((\z(K*K'))^k)}(k/2)$$
where $k=(a+\nu)+(2a+\r)=3a+\nu+\r$ and that
$$\un\z(K')\un{\bul}\un\z(K)=\un{gr_{k'}((\z(K)\bul\z(K'))^{k'})}(k'/2)$$
where $k'=(a+\r-\nu)+(a+\nu)+(a+\nu)=3a+\nu+\r$. It remains to use (a) and the
equality $k=k'$.

\subhead 5.22\endsubhead
Define $h:G@>>>G$ $g\m g\i$. For $K\in\cd_m(G)$ we set $K^\da=h^*K$. We show:

(a) {\it For $L\in\cd_m(Z)$ we have $(\c(L))^\da=\c(L^\da)$ with $L^\da$ as in 3.2.}
\nl
This follows from the definition of $\c$ using the commutative diagram
$$\CD
Z@<f<<\dZ@>\p>>G\\
@V\fh VV @V\dot{\fh}VV @VhVV\\
Z@<f<<\dZ@>\p>>G
\endCD$$
where $f,\p$ are as in 5.1, $\fh$ is as in 3.2 and $\dot{\fh}:\dZ@>>>\dZ$ is $(B,B',g)\m(B',B,g\i)$.

\mpb

From (a) and 3.2(a) we see that if $\l\in\fs$, $w\in W'_\l$, then
$$(\c(\Bbb L_\l^{\dw}))^\da=\c(\Bbb L_{\l\i}^{\dw\i}).\tag b$$
We deduce that

(c) {\it if $A$ is a $\fo$-character sheaf with associated two-sided cell $\boc$
then $A^\da$ is a $\fo\i$-character sheaf  with associated two-sided cell $\ti{\boc}$.}
\nl
From (a),(c) we deduce:

(d) {\it For $L\in\cc^\boc_0Z$ we have $(\un\c(L))^\da=\un\c(L^\da)$ where the second $\un\c$ is
relative to $\ti{\boc}$ instead of $\boc$.}

\head 6. The main results\endhead
\subhead 6.1\endsubhead
Let $\YY=\{(B,g);B\in\cb,g\in B\}$. We define $\fk:\YY@>>>\TT$ by $(B,g)\m t$ 
where $t\in\TT$ is given by the conditions $x\i gx\in t\UU$, $x\in G$, 
$B=x\BB x\i$; note that $t$ is independent of the choice of $x$. Let 
$$\tYY=\{((B,g),t)\in\YY\T\TT;\fk(B,g)=t^n\}.$$
Let $G_{rs}$ be the variety of regular semisimple elements in $G$. Let 
$$\YY_{rs}=
\{(B,g)\in\YY;g\in G_{rs}\},\tYY_{rs}=\{((B,g),t)\in\tYY;g\in G_{rs}\}.$$
For any $w\in W$ we define $\t_w:\YY_{rs}@>>>\YY_{rs}$ by $(B,g)\m(B',g)$ 
where $B'\in\cb$ is uniquely defined by the condition that $g\in B'$, 
$(B',B)\in\co_w$; one verifies that $\ti\t_w:\tYY_{rs}@>>>\tYY_{rs}$,
$((B,g),t)\m(\t_w(B,g),w(t))$, is well defined. Now $w\m\t_w$ (resp. 
$w\m\ti\t_w$) is an action of $W$ on $\YY_{rs}$ (resp. on $\tYY_{rs}$). For 
any $t_1\in\TT_n$ we define $\ti\t^{t_1}:\tYY_{rs}@>>>\tYY_{rs}$ by
$((B,g),t)\m(\t_w(B,g),tt_1)$. The operators $\ti\t_w$, $\ti\t^{t_1}$ define
an action of the semidirect product $W\TT_n$ (see 1.12) on $\tYY_{rs}$. This 
action leaves stable each fibre of the map $\tvp_{rs}:\tYY_{rs}@>>>G_{rs}$, 
$((B,g),t)\m g$, hence it induces an action of
$W\TT_n$ on $\tvp_{rs!}\bbq$, a local system of rank $\sha(W)n^\r$ on $G_{rs}$
such that the induced $W\TT_n$-action on any of its stalks is isomorphic to
the regular representation of $W\TT_n$. We show:

(a) {\it The algebra homomorphism $h:\bbq[W\TT_n]@>>>\End(\tvp_{rs!}\bbq)$ 
defined by the action above is an isomorphism.} 
\nl
Note that $h$ is injective since the induced algebra homomorphism from
$\bbq[W\TT_n]$ to the space of linear endomorphisms of any stalk of
$\tvp_{rs!}\bbq$ is clearly injective. Since 
$$\End(\tvp_{rs!}\bbq)
\cong\Hom_{\cd(\tYY_{rs})}(\bbq,\tvp_{rs}^*\tvp_{rs!}\bbq),$$ 
it is enough to show that
$$\dim\Hom_{\cd(\tYY_{rs})}(\bbq,\tvp_{rs}^*\tvp_{rs!}\bbq)\le\sha(W)n^\r.$$
Since $\tvp_{rs}^*\tvp_{rs!}\bbq$ is a local system of rank $\sha(W)n^\r$ on 
$\tYY_{rs}$, it is enough to show that $\tYY_{rs}$ is connected. Since 
$\tYY_{rs}@>>>\cb,((B,g),t)\m B$ is a $G$-equivariant fibration with $G$ 
acting transitively on $\cb$, it is enough to show that its fibre over $\BB$ 
is connected or that $\{(g,t)\in G_{rs}\T\TT;g\in t^n\UU\}$ is connected, or 
that 
$$\{(t,u)\in\TT\T\UU;t^nu\in G_{rs}\}=\{t\in\TT;t^n\in G_{rs}\}\T\UU$$
is connected. It is enough to observe that $\{t\in\TT;t^n\in G_{rs}\}$ is
connected (it is a nonempty open subset of $\TT$). This proves (a).

\mpb

We define $\tvp:\tYY@>>>G$ by $((B,g),t)\m g$. We have $\tvp=\vp\p_1$ where 
$\p_1:\tYY@>>>\YY$ is $((B,g),t)\m(B,g)$ and $\vp:\YY@>>>G$ is $(B,g)\m g$.
From the cartesian diagram
$$\CD
\tYY@>>>\TT\\
@V\p_1VV   @V\io VV\\
\YY@>\fk>>\TT\endCD$$
with $\io:\TT@>>>\TT$ as in 1.4, we see that $\fk^*\io_!\bbq=\p_{1!}\bbq$, 
hence
$$\Xi:=\tvp_!\bbq=\vp_!(\p_{1!}\bbq)=\vp_!(\fk^*\io_!\bbq)
=\op_{\l\in\fs}\vp_!(\fk^*L_\l)=\op_{l\in\fs}\Xi_\l$$
where $\Xi_\l=\vp_!(\fk^*L_\l)$. Since $\Xi|_{G_{rs}}=\tvp_{rs!}\bbq$, we have 
$\tvp_{rs!}\bbq=\op_{\l\in\fs}\Xi_\l|_{G_{rs}}$. As observed in \cite{\GRE}, 
$\vp:\YY@>>>G$ is small and $\Xi_\l$ is the intersection cohomology complex of
$G$ with coefficients in $\Xi_\l|_{G_{rs}}$; hence $\Xi$ is the intersection 
cohomology complex of $G$ with coefficients in $\tvp_{rs!}\bbq$. It follows 
that $\End_{\cd(G)}\Xi=\End_{\cd(G_{rs})}(\tvp_{rs!}\bbq)$ hence, using (a), 
$$\End_{\cd(G)}\Xi=\bbq[W\TT_n].\tag b$$
For any $E\in\Irr(W\TT_n)$ we set
$$A_E=\Hom_{\bbq[W\TT_n]}(E,\Xi\la\D\ra).\tag c$$
We see that $A_E$ is a simple perverse sheaf on $G$ and that for
$E\ne E'$ in $\Irr(W\TT_n)$ we have $A_E\not\cong A_{E'}$. Moreover we have
$$\Xi\la\D\ra=\op_{E\in\Irr(W\TT_n)}E\ot A_E.\tag d$$
From the definitions, for $\l\in\fs$ we have $\Xi_\l=R_\l^{\dot1}$ (notation 
of 5.1). Using this and (d) we see that for any $E\in\Irr(W\TT_n)$, $A_E$ is a
character sheaf on $G$.

We state the following result.

\proclaim{Proposition 6.2} For any $E'\in\Irr(W\TT_n)$ we have $\boc_{A_{E'}}=\boc_{E'}$.
In particular, we have $a_{A_{E'}}=a_{E'}$. (Notation of 5.1, 1.13).
\endproclaim

The proof is given in 6.3, 6.4, assuming, to simplify the exposition, that $n=1$. 
It consists in a reduction to an analogous (known) statement in the representation theory 
of the finite group $G^F$ in \cite{\ORA}. (We denote by
$F:G@>>>G$, $F:\cb@>>>\cb$ the Frobenius maps corresponding to the $\FF_q$-structures on $G,\cb$.)

\subhead 6.3\endsubhead 
Until the end of 6.5 we assume that $n=1$. Then $W\TT_n=W$, $\fs=\{1\}$ hence we can identify 
$W\fs=W$; $\HH$ has a basis $\{T_w;w\in W\}$. 
As in \cite{\DL}, for each $w\in W$ we consider the
variety $X_w=\{B\in\cb;(B,F(B))\in\co_w\}$ on which $G^F$ acts by conjugation and
the resulting $G^F$-module $H^i_c(X_w,\bbq)$ for each $i\in\ZZ$.
Let $\Irr_u(G^F)$ be a set of representatives for the isomorphism classes
of irreducible representations $\rr$ of $G^F$ such that $\rr$ appears in the $G^F$-module
$H^i_c(X_w,\bbq)$ for some $w\in W$ and some $i\in\ZZ$, or equivalently (see \cite{\DL}), 
such that $\sum_i(-1)^i(\rr:H^i_c(X_w,\bbq))\ne0$ for some $w\in W$. (Here $(\rr:?)$ denotes
the multiplicity of $\rr$ in $?$.) In the terminology of \cite{\DL}, $\Irr_u(G^F)$ is the set of
unipotent representations of $G^F$.
For any $\rr\in\Irr_u(G^F)$ and any $E\in\Irr(W)$ we set
$$b_{\rr,E}=\sha(W)\i\sum_{w\in W}\tr(w,E)\sum_i(-1)^i(\rr:H^i_c(X_w,\bbq))\in\QQ.\tag a$$
By \cite{\ORA, 4.23}, given $\rr\in\Irr_u(G^F)$, there is a unique two-sided cell $\boc_\rr$ 
of $W$ such that $b_{\rr,E}=0$ whenever $E\in\Irr(W)$ satisfies $\boc_E\ne\boc_\rr$; let $a_\rr$
be the value of the $a$-function $a:W@>>>\NN$ on $\boc_\rr$.

For $w=1$, we have $X_1=\cb^F$ and $\cf:=H^0_c(X_1,\bbq)$ is the vector space of functions
$\cb^F@>>>\bbq$. This vector space is naturally an $G^F$-module and its space of
$G^F$-equivariants endomorphism can be naturally identified with the semisimple algebra
$\HH^{\sqrt{q}}:=\bbq\ot_\ca\HH$ where $\bbq$ is viewed as an $\ca$-algebra via $v\m\sqrt{q}$.
Hence for any simple $\HH^{\sqrt{q}}$-module $M$, the vector space $\rr_M=\Hom_{\HH^{\sqrt{q}}}(M,\cf)$
is either $0$ or an object of $\Irr_u(G^F)$; in fact, it is known that it is $\ne0$.
For any $E\in\Irr(W)$ let $E^{\sqrt{q}}$ be the simple $\HH^{\sqrt{q}}$-module corresponding to $E$ under the
algebra isomorphisms $\HH^{\sqrt{q}}@>\psi^{\sqrt{q}}>>\bbq\ot\HH^\iy@>(\psi^1)\i>>\HH^1=\bbq[W]$ 
obtained by extension of scalars from $\psi:\HH@>>>\ca\ot\HH^\iy$ (see 1.12); we write $\rr_E$ 
instead of $\rr_{E^{\sqrt{q}}}$. Thus we have an imbedding $\Irr(E)@>>>\Irr_u(G^F)$, $E\m\rr_E$.

\subhead 6.4\endsubhead 
We write $\fR^{\dz}$ instead of $\fR^{\dz}_\l$ (see 5.1) where $\l=1$, $z\in W$.
The following result can be deduced from \cite{\REFL, 2.1}.

(a) {\it Let $y\in W$, $E'\in \Irr(W)$, $i\in\ZZ$. Then $\ch^i\fR^{\dy}|_{G_{rs}}$ is a local system
and} 
$$\sum_i(-1)^i(A_{E'}[-\D]|_{G_{rs}}:\ch^i\fR^{\dy}|_{G_{rs}})
=\sum_i(-1)^i(\rr_{E'}:H^i_c(X_y,\bbq)).$$
(In the left hand side, $(:)$ denotes the multiplicity
of an irreducible local system on $G_{rs}$ in another local system on $G_{rs}$.)

Using 6.3(a) we deduce for any $E,E'$ in $\Irr(W)$:
$$b_{\rr_{E'},E}
=\sha(W)\i\sum_{y\in W}\tr(y,E)\sum_i(-1)^i(A_{E'}[-\D]|_{G_{rs}}:\ch^i\fR^{\dy}|_{G_{rs}}).\tag b$$
From \cite{\CSIII, (14.10.1)}, for any $E\in\Irr(W)$, $A\in CS_{\fo}$ we have
$$b_{A,E}=\sha(W)\i\sum_{y\in W}\tr(y,E)\sum_i(-1)^{i+\D}(A:(\fR^{\dy})^i).$$
($b_{A,E}$ is as in 5.1(a).) In particular for $E'\in\Irr(W)$ we have
$$\align&b_{A_{E'},E}=\sha(W)\i\sum_{y\in W}\tr(y,E)\sum_i(-1)^{i+\D}(A_{E'}:(\fR^{\dy})^i)\\&
=\sha(W)\i\sum_{y\in W}\tr(y\i,E)\sum_i(-1)^{i+\D}(A_{E'}|_{G_{rs}}:(\fR^{\dy})^i|_{G_{rs}}).\endalign$$
Since $\ch^i\fR^{\dy}_{G_{rs}}$ are local systems, we see that
$$\sum_i(-1)^{i+\D}(A_{E'}|_{G_{rs}}:(\fR^{\dy})^i|_{G_{rs}})=
\sum_i(-1)^i(A_{E'}[\D]|_{G_{rs}}:\ch^i\fR^{\dy}|_{G_{rs}})$$
(in the last sum $(:)$ refers to multiplicities of an irreducible local system in another
local system). Thus,
$$b_{A_{E'},E}=\sha(W)\i\sum_{y\in W}\tr(y,E)\sum_i(-1)^i(A_{E'}[\D]|_{G_{rs}}:\ch^i\fR^{\dy}|_{G_{rs}})$$
so that, using (b), we have
$$b_{A_{E'},E}=b_{\rr_{E'},E}.$$
Using the definitions we now see that
$$\boc_{A_{E'}}=\boc_{\rr_{E'}},\text{ hence }a_{A_{E'}}=a_{\rr_{E'}}$$
for any $E'\in\Irr(W)$. Thus we can restate 6.2 as follows:

(c) {\it For any $E'\in\Irr(W)$ we have $\boc_{\rr_{E'}}=\boc_{E'}$.}
\nl
We shall deduce (c) from the following result which is equivalent to \cite{\ORA, 12.2(i)}:

(d) {\it Let $\L$ be a left cell of $W$. Write $[\L]=\op_{E\in\Irr(W)}E^{\op f_E}$ where $f_E\in\NN$.
In the Grothendieck module of $G^F$-modules tensored by $\QQ$ we have}
$$\sum_{E\in\Irr(W)}f_E\rr_E=\sum_{E\in\Irr(W)}\sum_{\rr\in\Irr_uG^F}b_{\rr,E}\rr.$$
Let $E'\in\Irr(W)$.
We can find a left cell $\L$ of $W$ as in (d) such that $\L\sub\boc_{E'}$ and $E'$ appears in $\L]$, that is,
$f_{E'}>0$. Then $\rr_{E'}$ appears with nonzero coefficient in the left hand side of the identity in (d)
hence it appears with nonzero coefficient in the right hand side of the identity in (d). Thus there exists
$E\in\Irr(W)$ such that $b_{\rr_{E'},E}\ne0$. By definition this means that $\boc_{\rr_{E'}}=\boc_{E'}$,
proving (c), hence also Proposition 6.2 (assuming $n=1$). The proof for general $n$ goes along similar lines.

Note that the proof of 6.4(d) given in \cite{\ORA, 12.2(i)} is case by case. It is likely that a more
efficient proof can be obtained using the inductive description of $W$-modules carried by left cells
in terms of constructible representations given in \cite{\CLASSII}.

\subhead 6.5\endsubhead 
The following inequality is a special case of Proposition 6.2.

(a) {\it For any $E'\in\Irr(W)$ we have $a_{A_{E'}}\le a_{E'}$.}
\nl
We give an alternative proof of (a) which avoids the use of 6.4(d) hence of \cite{\ORA, 12.2(i)}.
We again assume for simplicity that $n=1$. As in 6.4 it is enough to prove:

(b) $a_{\rr_{E'}}\le a_{E'}$.
\nl
It is known \cite{\GP, 8.1.8} that
$$\dim(\rr_{E'})=\dim(E')\sum_{w\in W}q^{|w|}(\sum_{w\in W}q^{-|w|}\tr(T_w,E'{}^{\sqrt{q}})^2)\i.$$
From \cite{\MADIS, 3.14, 3.16, 3.17, 3.19} we have
$$\dim\rr_{E'}=\sha(W)\i\sum_{E\in\Irr(W)}b_{\rr_{E'},E}\sum_{i\ge0}\dim\Hom_W(E,\bS^iV)q^i.$$
Since $b_{\rr_{E'},E}=0$ unless $\boc_E=\boc_{\rr_{E'}}$, it follows that
$$\align&\dim(E')\sum_{w\in W}q^{|w|}(\sum_{w\in W}q^{-|w|}\tr(T_w,E'{}^{\sqrt{q}})^2)\i\\&
=\sha(W)\i\sum_{E\in\Irr(W);\boc_E=\boc_{\rr_{E'}}}b_{\rr_{E'},E}\sum_{i\ge0}\dim\Hom_W(E,\bS^iV)q^i.
\tag c\endalign$$
Since $b_{\rr_{E'},E}=b_{A_{E'},E}$ is independent of $q$, we may regard (c) as an equality of polynomials
with rational coefficients in an indeterminate $q^{1/2}$.
From 1.20(c${}'$) we see that the right hand side of (c) is in $q^c\QQ[q^{-1/2}]$ where
$$c=\max_{E;\boc_E=\boc_{\rr_{E'}}}(\nu-b_{E\ot\sg})=\nu-a_{\rr_{E'}}$$ (the last equality uses 1.19(a),(b)).
From \cite{\HEC, 20.11} we see that the left hand side of (c) is in
$q^{\nu-a_{E'}}(c_0+q^{-1/2}\QQ[q^{-1/2}])$ where $c_0\in\QQ-\{0\}$.
Hence from (c) we deduce that $\nu-a_{E'}\le\nu-a_{\rr_{E'}}$ so that $a_{E'}\ge a_{\rr_{E'}}$, as required.

\subhead 6.6\endsubhead
We now return to our general $n$.
Let $A$ be a character sheaf of $G$. By \cite{\CDGVI, 30.12}, there exists
a parabolic subgroup $P$ of $G$, a Levi subgroup $L$ of $P$ and a subset $S_1$
of $L$ which is a single conjugacy class of $L$ times the connected centre of 
$L$ such that the support of $A$ is the union of $G$-conjugates of elements in
the closure of $S_1$ times the unipotent radical of $P$; moreover, if $P\in\cb$
then $A=A_E$ for some $E\in\Irr(W\TT_n)$ while if $P\n\cb$ then $A|_{\{1\}}=0$.
(Here we use the cleanness of cuspidal character sheaves, see \cite{\CLEAN} and 
its references. Actually we only use a weak form of the cleanness property which is
more elementary than what appears in \cite{\CLEAN}.) 

\subhead 6.7\endsubhead
Let $\ph:\pp@>>>G$ be the map with image $1$. For any $K\in\cd_m(G)$ we have
$$(\ph^*K)^j=\ch^j_1K.$$
The identification 6.1(b) induces for any $i$ an algebra homomorphism
$$\bbq[W\TT_n]@>>>\End(\ch^i_1\Xi);$$
thus, 
$$\ch^i_1\Xi=H^i_c(\tvp\i(1),\bbq)=H^i_c(\cb\T\TT_n,\bbq)=
H^i_c(\cb,\bbq)\ot\bbq[\TT_n]$$
becomes naturally a $W\TT_n$-module; one verifies that the action of $wt$ 
(with $w\in W$, $t\in\TT_n$) is given by $wt:e\ot t_1\m w(e)\ot w(tt_1)$ (here 
$e\in H^i(\cb,\bbq)$, $t_1\in\TT_n$ and $e\m w(e)$ is the $W$-action on 
$H^i(\cb,\bbq)=H^i(G/\BB,\bbq)=H^i(G/\TT,\bbq)$ induced by the conjugation
action of $N\TT$. (In the case where $n=1$ this is proved as in
\cite{\SPALT, \S2}; the proof in the general case is along similar lines.)
Note also that $\ch^i_1\Xi=0$ if $i$ is odd. We show:

(a) {\it Let $E$ be an irreducible $W\TT_n$-module. We have 
$(\ch^{2i}_1\Xi)^E=0$ for $i>\nu-a_E$. Moreover, $\dim(\ch^{2\nu-2a_E}_1\Xi)^E$
is $1$ if $E$ is special and $0$ if $E$ is not special.}
\nl
Let $\bS V=\op_{i\ge0}\bS^iV$ be as in 1.20. It is well known that for $i\ge0$
we have canonically $\bS^iV=H^{2i}_c(\cb,\bbq)$ compatibly with the 
$W$-actions. This extends to an identification
$$\bS^iV\ot\bbq[W\TT_n]=H^{2i}_c(\cb,\bbq)\ot\bbq[W\TT_n]=\ch^{2i}_1\Xi$$
which is compatible with the $W\TT_n$-actions. Hence (a) follows from 1.21(a).

\subhead 6.8\endsubhead
We show: 

(a) {\it Let $A\in CS_\boc$. If $j>-2a-\r$ then $(\ph^*A)^j=\ch^j_1A=0$.}
\nl
We can assume that $A\cong A_E$ for some $E\in\Irr(W\TT_n)$. (If $A$ is not of this 
form, the result holds by 6.6.) Since $a\le a_E$ (see 6.5(a)) we have $j>-2a_E-\r$.
By definition we have $\ch^j_1A_E=(\ch^j_1(\Xi\la\D\ra))^E=(\ch^{j+\D}_1\Xi)^E(\D/2)$
hence by 6.7(a), $\ch^j_1A_E=0$ if $j+\D>2\nu-2a_E$ that is if $j>-2a_E-\r$.

\mpb

We show: 

(b) {\it Let $E=E_\boc$, see 1.19. Then $(\ph^*A_E)^{-2a-\r}=\ch^{-2a-\r}_1A_E$ is a 
$1$-dimensional mixed vector space of pure weight $-2a-\r$. }
\nl
As in the proof of (a) we have
$$\ch^{-2a-\r}_1A_E=(\ch^{-2a-\r}_1(\Xi\la\D\ra))^E=(\ch^{-2a+2\nu}_1\Xi)^E(\D/2).$$
By 6.7(a) we have $\dim(\ch^{2\nu-2a}_1\Xi)^E=1$. It remains to note that 
$(\ch^{2\nu-2a}1q\Xi)^E$ is pure of weight $2\nu-2a$ (indeed,
$\ch^{2\nu-2a}_1\Xi=H^{2\nu-2a}_c(\cb,\bbq)\ot\bbq[\TT_n]$ is pure of weight $2\nu-2a$).

We show:

(c) {\it Let $A\in CS_\boc$ be such that $A\not\cong A_{E_\boc}$. 
Then $(\ph^*A)^{-2a-\r}=\ch^{-2a-\r}_1A=0$.}
\nl
As in the proof of (a) we can assume that $A\cong A_E$ for some $E\in\Irr(W\TT_n)$.
By 6.2 we have $\boc_E=\boc$ and by assumption we have $E\not\cong E_\boc$. Hence 
$E$ is not special. By 6.7(a) we have $(\ch^{2\nu-2a_E}_1\Xi)^E=0$ hence 
$\ch^{-2a_E-\r}_1A_E=0$. We have $a_E=a_\boc=a$ hence $\ch^{-2a-\r}_1A_E=0$. 
This proves (c).

From (b),(c) we see that if $K\in\cc^\boc G$ then

(d) $\dim\Hom_{\cc^\boc G}(A_{E_\boc},K)=\dim(\ph^*K)^{-2a-\r}=\dim\ch^{-2a-\r}_1K$.

\subhead 6.9\endsubhead
Let $\d_0:\cb@>>>Z$, $\p'_0:\cb@>>>\pp$ be as in 4.12.
We show that for $L\in\cd_m(Z)$ we have

(a) $\ph^*\c(L)=\p'_{0!}\d_0^*L$.
\nl
Define $\ph':\cb@>>>\dZ$ (see 4.1) by $B\m(B,B,1)$. 
We have a commutative diagram in which the right square is cartesian:
$$\CD\cb@>=>> \cb@>\p'_0>>\pp \\
@V\d_0 VV @V\ph'VV  @V\ph VV\\
Z@<f<<\dZ@>\p>>G
\endCD$$
(Here $f,\p$ are as in 4.1.) It follows that for $L\in\cd_m(Z)$ we have 
$\ph^*\c(L)=\ph^*\p_!f^*L=\p'_{0!}\ph'{}^*f^*L=\p'_{0!}\d_0^*L$. This proves (a).

Let $L\in\cc^\boc_0Z$. Applying \cite{\CONV, 8.2} with
$\Ph:\cd^{\preceq}_mG@>>>\cd_m\pp$, $K_1\m\ph^*K_1$, $c=-2a-\r$ (see 6.8(a)),
$K$ replaced by $\c(L)$ and $c'=a+\nu$, we see that we have canonically
$$(\ph^*(\un\c(L)))^{\{-2a-\r\}}\sub(\ph^*\c(L))^{\{-a-\r+\nu\}}
=(\p'_{0!}\d_0^*(L))^{\{-a-\r+\nu\}}.\tag b$$
(The last equality follows from (a).)
By 4.12(c), $\d_0^*L$ is pure of weight $0$ hence $\p'_{0!}\d_0^*L$ is pure of weight $0$ hence
$(\p'_{0!}\d_0^*L)^{\nu-a-\r}$ is pure of weight $\nu-a-\r$ so that
$$(\p'_{0!}\d_0^*L)^{\{-a-\r+\nu\}}=(\p'_{0!}\d_0^*L)^{-a-\r+\nu}((-a-\r+\nu)/2).$$
From 4.12(f) we have $(\p'_{0!}\d_0^*L)^{\nu-a-\r}=(\p'_{0!}((\d_0^*L)^{-a-\r}))^\nu$. Hence
$$(\p'_{0!}\d_0^*L)^{\{-a-\r+\nu\}}=(\p'_{0!}((\d_0^*L)^{-a-\r}))^\nu((-a-\r+\nu)/2).$$
Thus (b) becomes
$$(\ph^*(\un\c(L)))^{\{-2a-\r\}}\sub(\p'_{0!}((\d_0^*L)^{-a-\r}))^\nu((-a-\r+\nu)/2)$$
and using 4.12(h):
$$(\ph^*(\un\c(L)))^{\{-2a-\r\}}\sub\Hom_{\cc^\boc Z}(\bold1'_0,L).\tag c$$
We show that (c) is an equality:
$$(\ph^*(\un\c(L)))^{\{-2a-\r\}}=\Hom_{\cc^\boc Z}(\bold1'_0,L).\tag d$$
To prove this we can assume that $L=\Bbb L_\l^{\dw}$ for some $w\cdo\l\in\overset\smile\to\boc$.
If $w\cdo\l\n\DD_\boc$ then the right hand side of (d) is zero, hence by (c), the left hand side
of (d) is zero and (d) holds. Assume now that $w\cdo\l\in\DD_\boc$. Then the right hand side of (c)
has dimension $1$. Hence the left hand side of (c) has dimension $0$ or $1$; it is enough to prove
that it has dimension $1$. By 6.8(d) with $K=\un\c(L)$ we see that the left hand side of (c) has
dimension equal to $(A_{E_\boc}:\un\c(L))$. (We have also used 6.8(b).) In particular we have
$$(A_{E_\boc}:\un\c(\Bbb L_\l^{\dw}))\text{ is $0$ or $1$}\tag e$$
and we must prove that
$$(A_{E_\boc}:\un\c(\Bbb L_\l^{\dw}))=1.\tag f$$
In the rest of the proof we set $A=A_{E_\boc}$.
Using 5.8(a) we can reformulate (e) as $b_{A,[\L_{w,\l}]}\in\{0,1\}$ for any $w\cdo\l\in\DD_\boc$; 
we must prove that $b_{A,[\L_{w,\l}]}=1$ for any $w\cdo\l\in\DD_\boc$. 
Since $\boc_A=\boc$ we have $b_{A,[\L]}=0$ for any left cell $\L$ not contained in 
$\boc$. Hence it is enough to show that $\sum_\L b_{A,[\L]}=\sha(\DD_\boc)$ where
$\L$ runs over the left cells in $W\fs$. We have $\sum_\L b_{A,[\L]}=b_{A,Reg}$ where
$Reg$ is the regular representation of $W\TT_n$. Hence it is enough to show that
$b_{A,Reg}=\sha(\DD_\boc)$. From 5.1(a) with $z=1,j=\D$ we have
$$\align&\sum_{\l\in\fs}(A:(R_\l^1)^\D)
=(0;\sum_{E\in\Irr(W\TT_n)}b_{A,E}\sum_{\l\in\fs}\tr(1_\l,E^v))\\&=
(0;\sum_{E\in\Irr(W\TT_n)}b_{A,E}\dim(E))=b_{A,Reg}.\endalign$$
Hence it is enough to show that $\sum_{\l\in\fs}(A:(R_\l^1)^\D)=\sha(\DD_\boc)$ or equivalently
(see 6.1) that $(A:\Xi^\D)=\sha(\DD_\boc)$. By 6.1(d) we have $(A:\Xi^\D)=\dim(E_\boc)$. It
remains to show that $\dim(E_\boc)=\sha(\DD_\boc)$. The left hand side is 
$\sum_{\l\in\fo}\dim(1_\l E_\boc)$ where $1_\l E_\boc$ is the special representation of
$W'_\l$ attached to the two-sided cell $\boc_\l$ of $W'_\l$ determined by $\boc$; the right hand side is
$\sum_{\l\in\fo}n_\l$ where $n_\l$ is the number of distinguished involutions of $W'_\l$ contained
in $\boc_\l$. It is then enough to show that $\dim(1_\l E_\boc)=n_\l$ for any $\l\in\cl o$. This can be 
deduced from the following known property of a two-sided cell $\boc_0$ of $W_\l$: the dimension of the 
special representation of $W_\l$ corresponding to $\boc_0$ is equal to the number of 
distinguished involutions of $W_\l$ contained in $\boc_0$.
This completes the proof of (f) hence that of (d).

\mpb

We now state the following complement to (f).

(g) {\it If $w\cdo\l\in\overset\smile\to\boc$ and $w\cdo\l\n\DD_\boc$ then
$(A_{E_\boc}:\un\c(\Bbb L_\l^{\dw}))=0$.}
\nl
Let $L=\Bbb L_\l^{\dw}$. By 6.8(d) it is enough to show that $(\ph^*(\un\c L))^{-2a-\r}=0$.
By 6.8(b),(c), $(\ph^*(\un\c L))^{-2a-\r}$ is pure of weight $-2a-\r$ hence it is enough to
show that $(\ph^*(\un\c(L)))^{\{-2a-\r\}}=0$. Using (c) it is enough to note that, by our
assumption we have $\Hom_{\cc^\boc Z}(\bold1'_0,L)=0$.

\subhead 6.10\endsubhead
Let $u:G@>>>\pp$ be the obvious map. From \cite{\CSII, 7.4} we see that for $K,K'$ in $\cm_mG$
we have canonically
$$(u_!(K\ot K'))^0=\Hom_{\cm(G)}(\fD(K),K'),\qua (u_!(K\ot K'))^j=0 \text{ if }j>0.$$
We deduce that if $K,K'$ are also pure of weight $0$ then $(u_!(K\ot K'))^0$  is pure of weight $0$
that is $(u_!(K\ot K'))^0=gr_0(u_!(K\ot K'))^0$. From the definitions we see that we have
$u_!(K\ot K')=\ph^*(K^\da*K')$ where $K^\da$ as in 5.22. Hence for
$K'$ in $\cc^\boc_0G$ and $K$ in $\cc^{\ti{\boc}}_0G$ (so that $K^\da\in\cc^\boc_0G$, see 5.22(c))
we have
$$\Hom_{\cm(G)}(\fD(K),K')=(\ph^*(K^\da*K'))^0=(\ph^*(K^\da*K'))^{\{0\}}.\tag a$$
Applying \cite{\CONV, 8.2} with $\Ph:\cd^{\preceq}_mG@>>>\cd_m\pp$, $K_1\m\ph^*K_1$, $c=-2a-\r$ (see 6.8(a)),
$K$ replaced by $K^\da*K'$ and $c'=2a+\r$ we see that we have canonically
$$(\ph^*(K^\da\un{*}K'))^{\{-2a-\r\}}\sub(\ph^*(K^\da*K'))^{\{0\}}.$$
In particular, if $L,L'$ are in $\cc^\boc_0Z$ then we have canonically
$$(\ph^*(\un\c(L')\un{*}\un\c(L)))^{\{-2a-\r\}}\sub(\ph^*(\un\c(L')*\un\c(L)))^{\{0\}}.$$
Using the equality
$$(\ph^*(\un\c(L')\un{*}\un\c(L)))^{\{-2a-\r\}}=\ph^*(\un\c(L\un\bul\un\z(\un\c(L')))))^{-2a-\r}$$
which comes from 5.20(b), we deduce that we have canonically
$$\ph^*(\un\c(L\un\bul\un\z(\un\c(L')))))^{-2a-\r}\sub(\ph^*(\un\c(L')*\un\c(L)))^{\{0\}},$$
or equivalently, using (a) with $K, K'$ replaced by $\un\c(L')^\da$, $\un\c(L)$,
$$\align&\ph^*(\un\c(L\un\bul\un\z(\un\c(L')))))^{-2a-\r}\sub
\Hom_{\cc^\boc G}(\fD(\un\c(L')^\da),\un\c(L))\\&=
\Hom_{\cc^\boc G}(\fD(\un\c(L)^\da),\un\c(L')).\endalign$$
Using now 6.9(d) with $L$ replaced by $L\un\bul\un\z(\un\c(L'))$ we deduce that we have canonically
$$\Hom_{\cc^\boc Z}(\bold1',L\un\bul\un\z(\un\c(L'))))\sub
\Hom_{\cc^\boc G}(\fD(\un\c(L)^\da),\un\c(L'))$$
or equivalently (using 4.12(j)):
$$\Hom_{\cc^\boc Z}(\fD(\un\z(\un\c(L'))^\da),L)\sub\Hom_{\cc^\boc G}(\fD(\un\c(L)^\da),\un\c(L')).$$
Now we have
$$\align&\Hom_{\cc^\boc Z}(\fD(\un\z(\un\c(L'))^\da),L)
=\Hom_{\cc^{\ti{\boc}}Z}(\fD(L),\un\z(\un\c(L'))^\da)\\&=
\Hom_{\cc^\boc Z}((\fD(L))^\da,\un\z(\un\c(L')))\endalign$$
hence 
$$\Hom_{\cc^\boc Z}((\fD(L))^\da,\un\z(\un\c(L')))\sub\Hom_{\cc^\boc G}(\fD(\un\c(L)^\da),\un\c(L')).$$
We set ${}^1L=\fD(L^\da)=(\fD(L))^\da$ and note that 
$$\fD(\un\c(L)^\da)=\fD(\un\c(L^\da))=\un\c(\fD(L^\da))=\un\c({}^1L),$$
see 5.22(d), 5.7(b). We obtain
$$\Hom_{\cc^\boc Z}({}^1L,\un\z(\un\c(L')))\sub\Hom_{\cc^\boc G}(\un\c({}^1L),\un\c(L'))\tag b$$
for any ${}^1L,L'$ in $\cc^\boc_0Z$.

We show that (b) is an equality:
$$\Hom_{\cc^\boc Z}({}^1L,\un\z(\un\c(L')))=\Hom_{\cc^\boc G}(\un\c({}^1L),\un\c(L')).\tag c$$
Let $N'$ (resp. $N''$) be the dimension of the left (resp. right) hand side of (b). It is enough to show 
that $N'=N''$.
We can assume that ${}^1L=\Bbb L_\l^{\dw}$, $L'=\Bbb L_{\l'}^{\dw'}$ where $w\cdo\l\in\overset\smile\to\boc$,
$w'\cdo\l'\in\overset\smile\to\boc$. 
By 5.13(a), $N'$ is the multiplicity of ${}^1L$ in $\un{\fb''}(L')$; by the fully faithfulness of
$\ti\e$ this is the same as the multiplicity of $\ti\e{}^1L$ in 
$\ti\e\un{\fb''}(L')=\un{\fb'}(L')=\un\fb(L')$ (the last two equalities use
3.25(d) and 3.14(d)). By 3.13(d) this is the same as the multiplicity of $\LL_\l^{\dw}$ in 
$$\op_{y\in W;y\cdo\l'\in\boc}\LL_{\l'}^{\dy}\un{\cir}\LL_{\l'}^{\dw'}\un{\cir}\LL_{y(\l')}^{\dy\i}.$$
Using now 2.22(c) we see that $N'$ is the coefficient of $t_{w\cdo\l}$ in
$$\sum_{y\in W;y\cdo\l'\in\boc}t_{y\cdo\l'}t_{w'\cdo\l'}t_{y\i\cdo y(\l')}\in\HH^\iy.$$
Hence if $\t:\HH^\iy@>>>\ZZ$ is as in 4.6 (see also 1.9) then
$$N'=\sum_{y\in W;y\cdo\l'\in\boc}\t(t_{y\cdo\l'}t_{w'\cdo\l'}t_{y\i\cdo y(\l')}t_{w\i\cdo\l}).$$
This can be rewritten as
$$N'=\sum_{y\cdo\l_1\in\boc}\t(t_{y\cdo\l_1}t_{w'\cdo\l'}t_{y\i\cdo y(\l_1)}t_{w\i\cdo\l}).$$
(In the last sum, the terms corresponding to $y\cdo\l_1$ with $\l_1\ne\l'$ are equal to zero.)
By 5.6(a) we have
$$N''=\sum_{y\cdo\l_1\in\boc}\t(t_{y\i\cdo y(\l_1)}t_{w\cdo\l}t_{y\cdo\l_1}t_{w'{}\i\cdo\l'}).$$
Since $\t(\x^\flat)=\t(\x)$ for all $\x\in\HH^\iy$ and $\x\m\x^\flat$ is the ring antiautomorphism in 1.9
we have also
$$N''=\sum_{y\cdo\l_1\in\boc}\t(t_{w'\cdo\l'}t_{y\i\cdo y(\l_1)}t_{w\i\cdo\l}t_{y\cdo\l_1}).$$
To show that $N'=N''$ it is enough to show that for any $y\cdo\l_1\in\boc$ we have
$$\t(t_{y\cdo\l_1}t_{w'\cdo\l'}t_{y\i\cdo y(\l_1)}t_{w\i\cdo\l})=
\t(t_{w'\cdo\l'}t_{y\i\cdo y(\l_1)}t_{w\i\cdo\l}t_{y\cdo\l_1}).$$
This follows by taking $\x=t_{w'\cdo\l'}t_{y\i\cdo y(\l_1)}t_{w\i\cdo\l}$, $\x'=t_{y\cdo\l_1}$ in
the identity $\t(\x\x')=\t(\x'\x)$ which (as we see from 1.9(a)) holds for any $\x,\x'$ in $\HH^\iy$.
This completes the proof of the equality $N'=N''$ and hence that of (c).

\subhead 6.11\endsubhead
In the reminder of this section we assume that the $\FF_q$-rational structure on $G$ in 2.8 is such that 

(a) {\it any $A\in CS_\boc$ admits a mixed structure of pure weight $0$.}
\nl
(This can be achieved by replacing if necessary $q$ by a power of $q$.)

The bifunctor $\cc^\boc_0G\T\cc^\boc_0G@>>>\cc^\boc_0G$, $K,K'\m K\un{*}K'$ in 5.20
defines a bifunctor $\cc^\boc G\T\cc^\boc G@>>>\cc^\boc G$ denoted again by $K,K'\m K\un{*}K'$
as follows. Let $K\in\cc^\boc G$, $K'\in\cc^\boc G$; we choose mixed structures of pure weight $0$
on $K,K'$ (this is possible by (a)), we define $K\un{*}K'$ as in 5.20 in terms of these mixed structures
and we then disregard the mixed structure on $K\un{*}K'$. The resulting object of $\cc^\boc G$ is denoted
again by $K\un{*}K'$; it is independent of the choices made.

In the same way the functor $\un\c:\cc^\boc_0Z@>>>\cc^\boc_0G$ gives rise to a functor
$\cc^\boc Z@>>>\cc^\boc G$ denoted again by $\un\c$; the functor $\un\z:\cc^\boc_0G@>>>\cc^\boc_0Z$ 
gives rise to a functor $\cc^\boc G@>>>\cc^\boc Z$ denoted again by $\un\z$.

The operation $K\un{*}K'$ is again called truncated convolution. It has a canonical associativity
isomorphism (deduced from that in 5.20(d)) which again satisfies the pentagon property. Thus $\cc^\boc G$
becomes a monoidal category; it has a braiding coming from 5.20(a).

\subhead 6.12\endsubhead
If $K\in\cc^\boc G$ then the isomorphisms 5.14(b) provide a central structure on
$\ti\e\un\z(K)\in\cc^\boc\tcb^2$ so that
$\ti\e\un\z(K)$ can be naturally viewed as an object of $\cz^\boc$ denoted by 
$\ov{\ti\e\un\z(K)}$. (Here $\ti\e$ is as in 3.3, $\un\z$ is as in 5.9, 
$\cz^\boc$ is as in 4.9.) Then $K\m\ov{\ti\e\un\z(K)}$ is a functor $\cc^\boc G@>>>\cz^\boc$.
We shall prove the following result.

\proclaim{Theorem 6.13} The functor $\cc^\boc G@>>>\cz^\boc$, $K\m\ov{\ti\e\un\z(K)}$ is an equivalence
of categories.
\endproclaim
From 5.13(a), 3.14(d), 3.25(d) we have canonically for any $z\cdo\l\in\overset\smile\to\boc$:
$$\ti\e\un\z(\un\c(\Bbb L_\l^{\dz}))=\un\fb(\Bbb L_\l^{\dz})\tag a$$
as objects of $\cc^\boc\tcb^2$.
From the definitions we see that the central structure on the left hand side of (a) provided by 6.12 is the
same as the central structure on the right hand side of (a) provided by 3.14(j). Hence we have
$$\ov{\ti\e\un\z(\un\c(\Bbb L_\l^{\dz}))}=\ov{\un\fb(\Bbb L_\l^{\dz})}\tag b$$
as objects of $\cz^\boc$. Using this and 4.11(a) with
$L'=\ti\e\un\z(\un\c(\Bbb L_{\l'}^{\dw}))$ 
(where $z\cdo\l,w\cdo\l'$ are in $\overset\smile\to\boc$), we have
$$\Hom_{\cc^\boc\tcb^2}(\LL_\l^{\dz},\ti\e\un\z(\un\c(\Bbb L_{\l'}^{\dw})))=       
\Hom_{\cz^\boc}(\ov{\ti\e\un\z(\un\c(\Bbb L_\l^{\dz}))},\ov{\ti\e\un\z(\un\c(\Bbb L_{\l'}^{\dw}))}).$$
Combining this with the equalities
$$\Hom_{\cc^\boc G}(\un\c(\Bbb L_\l^{\dz}),\un\c(\Bbb L_{\l'}^{\dw}))=
\Hom_{\cc^\boc Z}(\Bbb L_l^{\dz},\un\z(\un\c(\Bbb L_{\l'}^{\dw})))=
\Hom_{\cc^\boc\tcb^2}(\LL_l^{\dz},\ti\e\un\z(\un\c(\Bbb L_{\l'}^{\dw}))),$$
of which the first comes from 6.10(c) and the second comes from the fully faithfulness of $\ti\e$,
we obtain
$$\Hom_{\cc^\boc G}(\un\c(\Bbb L_\l^{\dz}),\un\c(\Bbb L_{\l'}^{\dw}))=
\Hom_{\cz^\boc}(\ov{\ti\e\un\z(\un\c(\Bbb L_\l^{\dz}))},
\ov{\ti\e\un\z(\un\c(\Bbb L_{\l'}^{\dw}))}).$$
In other words, setting 
$$\AA_{z\cdo\l,w\cdo\l'}=\Hom_{\cc^\boc G}(\un\c(\Bbb L_\l^{\dz}),\un\c(\Bbb L_{\l'}^{\dw})),$$
$$\AA'_{z\cdo\l,w\cdo\l'}=\Hom_{\cz^\boc}(\ov{\ti\e\un\z(\un\c(\Bbb L_\l^{\dz}))},
\ov{\ti\e\un\z(\un\c(\Bbb L_{\l'}^{\dw}))}),$$
we have 
$$\AA_{z\cdo\l,w\cdo\l'}=\AA'_{z\cdo\l,w\cdo\l'}.\tag c$$
Note that the identification (c) is induced by the functor $K\m\ov{\ti\e\un\z(K)}$.
Let $\AA=\op\AA_{z\cdo\l,w\cdo\l'}$, $\AA'=\op\AA_{z\cdo\l,w\cdo\l'}$ (both direct sums are taken over
all $z\cdo\l,w\cdo\l'$ in $\overset\smile\to\boc$). Then from (c) we have $\AA=\AA'$. Note that this
identification is compatible with the obvious algebra structures of $\AA,\AA'$.

For any $A\in CS_\boc$ we denote by $\AA_A$ the set of all $f\in\AA$ such that for any
$z\cdo\l,w\cdo\l'$, the $(z\cdo\l,w\cdo\l')$-component of $f$ maps the $A$-isotypic component of 
$\un\c(\Bbb L_\l^{\dz})$ to the $A$-isotypic component of $\un\c(\Bbb L_{\l'}^{\dw})$ and any other
isotypic component of $\un\c(\Bbb L_\l^{\dz})$ to $0$. Thus, $\AA=\op_{A\in CS_\boc}\AA_A$ is the
decomposition of $\AA$ into a sum of simple algebras. 
(Each $\AA_A$ is nonzero since, by 5.2(c) and 5.5(a), any
$A$ is a summand of some $\un\c(\Bbb L_\l^{\dz})$.)

From \cite{\MUG}, \cite{\ENO} we see that $\cz^\boc$ is a semisimple abelian category with finitely
many simple objects up to isomorphism. Let $\fS$ be a set of representatives for the isomorphism
classes of simple objects of $\cz^\boc$. For any $\s\in\fS$ we denote by $\AA'_\s$ the set of all
$f'\in\AA'$ such that for any
$z\cdo\l,w\cdo\l'$, the $(z\cdo\l,w\cdo\l')$-component of $f'$ maps the $\s$-isotypic component of 
$\ov{\ti\e\un\z(\un\c(\Bbb L_\l^{\dz}))}$ to the $\s$-isotypic component of 
$\ov{\ti\e\un\z(\un\c(\Bbb L_{\l'}^{\dw}))})$ and all other isotypic components of 
$\ov{\ti\e\un\z(\un\c(\Bbb L_\l^{\dz}))}$ to zero. Then $\AA'=\op_\s\AA'_\s$ is the decomposition
of $\AA'$ into a sum of simple algebras. (Each $\AA'_\s$ is nonzero since any $\s$ is a summand of some
$\ov{\ti\e\un\z(\un\c(\Bbb L_\l^{\dx}))}$ with $x\cdo\l\in\boc$.
Indeed, we can find $x\cdo\l\in\boc$ such that 
$\LL_\l^{\dx}$ is a summand of $\s$, viewed as an object of $\cc^\boc\tcb^2$; 
then by 4.9(a), $\s$ is a summand of $\ov{I(\LL_\l^{\dx})}$. If in addition, 
$x\cdo\l\in\overset\smile\to{W\fs}$ then,
by 4.10(a),(b), we have $\ov{I(\LL_\l^{\dx})}=\ov{\un\fb(\Bbb L_\l^{\dx})}$ 
hence  $\s$ is a summand of $\ov{\un\fb(\Bbb L_\l^{\dx})}$ hence, 
by (a), $\s$ is a summand of $\ov{\ti\e\un\z(\un\c(\Bbb L_\l^{\dx}))}$, as required.
If $x\cdo\l\n\overset\smile\to{W\fs}$ 
then, by 4.9(b) we have $I(\LL_\l^{\dx})\ne0$ which is a contradiction.)

Since $\AA=\AA'$, from the uniqueness of decomposition of a semisimple algebra as a direct sum of
simple algebras, we see that there is a unique bijection $CS_\boc\lra\fS$, $A\lra\s_A$ such that
$\AA_A=\AA'_{\s_A}$ for any $A\in CS_\boc$. From the definitions we now see that for any $A\in CS_\boc$
we have $\ov{\ti\e\un\z(K)}\cong\s_A$. Therefore, Theorem 6.13 holds.

\proclaim{Theorem 6.14} Let $L\in\cc^\boc Z$, $K\in\cc^\boc G$. We have canonically
$$\Hom_{\cc^\boc Z}(L,\un\z(K))=\Hom_{\cc^\boc G}(\un\c(L),K).\tag a$$
\endproclaim
We can assume that $L=\Bbb L_\l^{\dz}$ where $z\cdo\l\in\overset\smile\to\boc$.
From 6.13 and its proof we see that
$$\Hom_{\cc^\boc G}(\un\c(L),K)=\Hom_{\cz^\boc}(\ov{\ti\e\un\z(\un\c(L))},\ov{\ti\e\un\z(K)})=
\Hom_{\cz^\boc}(\ov{I(\LL_\l^{\dz})},\ov{\ti\e\un\z(K)}).$$
Using 4.9(a) we see that 
$$\Hom_{\cz^\boc}(\ov{I(\LL_\l^{\dz})},\ov{\ti\e\un\z(K)})=
\Hom_{\cc^\boc\tcb^2}(\LL_\l^{\dz},\ti\e\un\z(K))=\Hom_{\cc^\boc Z}(L,\un\z(K)).$$
This proves the theorem.

\subhead 6.15\endsubhead
We show that for $K\in\cc^\boc G$ we have canonically
$$\fD(\un\z(\fD(K))))=\un\z(K).\tag a$$
Here the first $\un\z$ is relative to $\ti\boc$.
It is enough to show that for any $L\in\cc^\boc Z$ we have canonically
$$\Hom_{\cc^\boc Z}(L,\fD(\un\z(\fD(K)))))=\Hom_{\cc^\boc Z}(L,\un\z(K)).$$
Here the left side equals
$$\align&\Hom_{\cc^{\ti{\boc}}Z}(\un\z(\fD(K)),\fD(L))=\Hom_{\cc^\boc G}(\fD(K),\un\c(\fD(L)))\\&=
\Hom_{\cc^\boc G}(\fD(K),\fD(\un\c(L))).\endalign$$
(We have used 6.14(a) for $\ti\boc$ and 5.7(b).) The right hand side equals
$$\Hom_{\cc^\boc G}(\un\c(L),K)=\Hom_{\cc^\boc G}(\fD(K),\fD(\un\c(L))).$$
(We have again used 6.14(a).) This proves (a).

\subhead 6.16\endsubhead
The monoidal structure on $\cc^\boc\tcb^2$ induces a monoidal structure on $\cz^\boc$.
Using 5.21(b) and 3.24(b) we see that the equivalence of categories in 6.13
is compatible with the monoidal structures. Since $\cz^\boc$ has a unit object, it follows
that the monoidal category $\cc^\boc G$ also has a unit object, say $A$. We show:
$$A\cong A_{E_\boc}\tag a$$
with $A_{E_\boc}$ as in 6.8(c). From 6.9(f),(g) we see that for $w\cdo\l\in\overset\smile\to\boc$,
$(A_{E_\boc}:\un\c(\Bbb L_\l^{\dw}))$ is $1$ if $w\cdo\l\in\DD_\boc$ and is $0$ if
$w\cdo\l\n\DD_\boc$. Using 6.13 we deduce that 
$$\dim\Hom_{\cc^\boc\tcb^2}(\LL_\l^{\dw},\ti\e\un\z(A_{E_\boc}))$$
is $1$ if $w\cdo\l\in\DD_\boc$ and is $0$ if
$w\cdo\l\n\DD_\boc$. Thus $\ti\e\un\z(A_{E_\boc})$ is isomorphic in $\cc^\boc\tcb^2$ to the unit
object $\bold1$ of the monoidal category $\cc^\boc\tcb^2$. Then 
$\ti\e\un\z(A_{E_\boc})$ viewed as an object of $\cz^\boc$ is also the unit object of $\cz^\boc$ hence
is isomorphic in $\cz^\boc$ to $\ti\e\un\z(A)$. Using 6.13 we deduce that (a) holds.

\mpb

{\it Acknowledgement.}
This research was supported in part by National Science Foundation grant DMS-1303060 and by a 
Simons Fellowship. A part of the writing was done during a visit (April, May 2015) to the 
Mittag-Leffler Institute, Djursholm, Sweden, whose hospitality and support is hereby acknowledged.

\enddocument